\documentclass[12 pt]{amsart}

\usepackage{hyperref}
\usepackage{etex}
\usepackage[shortlabels]{enumitem}
\usepackage{amsmath}
\usepackage{amsxtra}
\usepackage{amssymb}
\usepackage{enumitem}
\usepackage{array}

\usepackage{amscd}
\usepackage{amsthm}
\usepackage{adjustbox}
\usepackage{amsfonts}
\usepackage{amssymb}
\usepackage{eucal}
\usepackage[all]{xy} 
\usepackage{graphicx}
\usepackage{tikz-cd}
\usepackage{mathrsfs}
\usepackage{subfiles}
\usepackage{mathpazo}
\usepackage[colorinlistoftodos, textsize=tiny]{todonotes}
\usepackage{morefloats}
\usepackage{pdfpages}
\usepackage{thm-restate}
\usepackage[percent]{overpic}
\usepackage[utf8]{inputenc}
\usepackage{epigraph}
\usepackage{csquotes}
\usepackage[margin=1in]{geometry}
\usepackage{adjustbox}
\usepackage{microtype}
\usepackage{stmaryrd}

\usepackage{bm}
\usepackage{verbatim}
\usepackage{stmaryrd}
\usepackage{scalerel}
\usepackage{stackengine}
\stackMath
\newcommand\reallywidehat[1]{%
\savestack{\tmpbox}{\stretchto{%
  \scaleto{%
    \scalerel*[\widthof{\ensuremath{#1}}]{\kern-.6pt\bigwedge\kern-.6pt}%
    {\rule[-\textheight/2]{1ex}{\textheight}}%WIDTH-LIMITED BIG WEDGE
  }{\textheight}% 
}{0.5ex}}%
\stackon[1pt]{#1}{\tmpbox}%
}
\usepackage{mathtools}

\graphicspath{ {images/} }

\RequirePackage{color}
\definecolor{myred}{rgb}{0.75,0,0}
\definecolor{mygreen}{rgb}{0,0.5,0}
\definecolor{myblue}{rgb}{0,0,0.65}

\usepackage{color}

\usepackage{hyperref}
\hypersetup{citecolor=blue}
\usepackage{tikz}
\usetikzlibrary{matrix,arrows,decorations.pathmorphing}

\theoremstyle{plain}
\newtheorem{theorem}[subsubsection]{Theorem}

\newtheorem{proposition}[subsubsection]{Proposition}
\newtheorem{prop/constr}[subsubsection]{Proposition/Construction}
\newtheorem{lemma}[subsubsection]{Lemma}
\newtheorem{corollary}[subsubsection]{Corollary}

\theoremstyle{definition}
\newtheorem{definition}[subsubsection]{Definition}
\newtheorem{remark}[subsubsection]{Remark}
\newtheorem{example}[subsubsection]{Example}

\newtheorem{assumption}[subsubsection]{Assumption}

\newtheorem{conjecture}[subsubsection]{Conjecture}

\newtheorem{setup}[subsubsection]{Setup}

\newtheorem*{claim*}{Claim}
\theoremstyle{remark}
\newtheorem{notation}[subsubsection]{Notation}
\newtheorem{construction}[subsubsection]{Construction}
\numberwithin{equation}{section}
\newcommand\nc{\newcommand}
\nc\on{\operatorname}
\nc\renc{\renewcommand}

\DeclareMathOperator\id{id}
\DeclareMathOperator\Hom{Hom}

\DeclareMathOperator\End{End}

\DeclareMathOperator\gr{gr}

\DeclareMathOperator\At{At}
\DeclareMathOperator\PD{PD}

\DeclareMathOperator\Sch{Sch}
\DeclareMathOperator\ad{ad}
\DeclareMathOperator\KS{KS}
\DeclareMathOperator\Def{Def}

\nc\mf\mathfrak
\nc\mc\mathcal
\nc\mb\mathbb
\nc\msf\mathsf
\nc\mscr\mathscr

\def\Spec{{\rm Spec}}
\newcommand{\defeq}{\vcentcolon=}

\newcommand{\Qbar}{\overline{\mathbb{Q}}}

\newcommand{\AR}{\mathsf{AR}}
\newcommand{\X}{\mscr{X}}

\newcommand{\fracr}{\mscr{K}}
\newcommand{\Xbar}{\overline{X}}

\newcommand{\Rint}{R^{int}}

\newcommand{\fbar}{\overline{f}}

\newcommand{\Binthat}{\hat{B} ^{int}}
\newcommand{\Ainthat}{\hat{A} ^{int}}

\newcommand{\vecu}[1]{\underline{#1}}

\makeatletter
\renewcommand\part{%
   \if@noskipsec \leavevmode \fi
   \par
   \addvspace{4ex}%
   \@afterindentfalse
   \secdef\@part\@spart}

\def\@part[#1]#2{%
    \ifnum \c@secnumdepth >\m@ne
      \refstepcounter{part}%
      \addcontentsline{toc}{part}{\thepart\hspace{1em}#1}%
    \else
      \addcontentsline{toc}{part}{#1}%
    \fi
    {\parindent \z@ \raggedright
     \interlinepenalty \@M
     \normalfont
     \ifnum \c@secnumdepth >\m@ne
       \Large\bfseries \partname\nobreakspace\thepart
       \par\nobreak
     \fi
     \huge \bfseries #2%
     \par}%
    \nobreak
    \vskip 3ex
    \@afterheading}
\def\@spart#1{%
    {\parindent \z@ \raggedright
     \interlinepenalty \@M
     \normalfont
     \huge \bfseries #1\par}%
     \nobreak
     \vskip 3ex
     \@afterheading}
\makeatother

\renewcommand{\thepart}{\Roman{part}}

\title{Algebraicity and integrality of solutions to differential equations}
\author{Yeuk Hay Joshua Lam and Daniel Litt}
\date{\today}

\begin{document}

\begin{abstract}
We formulate a conjecture classifying algebraic solutions to (possibly non-linear) algebraic differential equations, in terms of the primes appearing in the denominators of the coefficients of their Taylor expansion at a non-singular point. For linear differential equations, this conjecture is a strengthening of the Grothendieck-Katz $p$-curvature conjecture. We prove the conjecture for many differential equations and initial conditions of algebro-geometric interest. For linear differential equations, we prove it for Picard-Fuchs equations at initial conditions corresponding to cycle classes, among other cases. For non-linear differential equations, we prove it for isomonodromy differential equations, such as the Painlev\'e VI equation and Schlesinger system, at initial conditions corresponding to Picard-Fuchs equations. We draw a number of algebro-geometric consequences from the proofs.
\end{abstract}

\maketitle

\setcounter{tocdepth}{1}
\tableofcontents

%\pagebreak

\section{Introduction}\label{section:introduction}
%\subsubsection{Notation}\label{notation}

When is a solution to an algebraic differential equation itself an algebraic function? This question goes back at least to the latter half of the 19th century, when Schwarz answered it for hypergeometric functions \cite{schwarz1873ueber}, answering a question of Kummer \cite{kummer1887}, and Fuchs asked it explicitly for general linear ODE \cite{fuchs1875uber}. For linear ODE, there is a classical conjectural arithmetic answer, due to Grothendieck and Katz \cite{katz-p-curvature, katz1982conjecture}: conjecturally, the holomorphic solutions to the ODE $$\left(\frac{\partial}{\partial z}-A\right)\vec{f}(z)=0, ~A\in \on{Mat}_{r\times r}(\mathbb{Q}(z))$$ are algebraic if and only if the algebraic differential operator $\left(\frac{\partial}{\partial z}-A\right)$ satisfies $$\left(\frac{\partial}{\partial z}-A\right)^p\equiv 0 \bmod p$$ for almost all primes $p$. The purpose of this paper is to formulate an analogous conjecture for possibly \emph{non-linear} ODE, and to provide some evidence for this conjecture by proving it, in \autoref{thm:NA-main-intro}, for a large class of non-linear ODE---those coming from isomonodromy, including, for example, the Painlev\'e VI equation and Schlesinger systems---at initial conditions arising from algebraic geometry, which we refer to as \emph{Picard-Fuchs initial conditions}. We view this result as an analogue, for non-linear ODE, of Katz's proof of the classical $p$-curvature conjecture for Picard-Fuchs equations \cite{katz-p-curvature}. And in fact, proving our conjecture for only  certain isomonodromy ODE (and all initial conditions) would suffice to prove the classical Grothendieck-Katz $p$-curvature conjecture in full generality. 

It turns out that in the case of isomonodromy ODE and Picard-Fuchs initial conditions our conjecture is closely related to a variational form of the relative Fontaine-Mazur conjecture, and as a byproduct of our proof, we verify some predictions of that conjecture.

When specialized to linear ODE, our conjecture is apparently strictly stronger than the $p$-curvature conjecture. In this case a form of it has been attributed to Andr\'e and Christol (see e.g.~\cite[Remarque 5.3.2]{andre2004conjecture} and \cite[Conjecture 1.12]{bostan2023algebraic}). It seems that even for linear ODE little was known; in addition to our main result on non-linear ODE, we verify the conjecture for many classical linear ODE and initial conditions of algebro-geometric interest.

\subsection{The conjecture}
The simplest setting in which our conjecture applies is the following. Let $g\in \mathbb{Q}(z,y_0, \cdots, y_{n-1})$ be a rational function in $n+1$ variables, and consider the differential equation \begin{equation}\label{eqn:non-linear-ODE-1} f^{(n)}(z)=g(z,f(z), f'(z), \cdots, f^{(n-1)}(z)).\end{equation} Let $t_0,\cdots, t_{n-1} \in \mathbb{Q}$ be such that $g(0, t_0,\cdots, t_{n-1})$ is defined\footnote{that is, $g=p/q$ with $p, q$ polynomials and $q(0, t_0, \cdots, t_{n-1})\neq 0$}; then there is a unique formal power series $$f(z)=\sum_{n\geq 0} a_nz^n\in \mathbb{Q}[[z]]$$ with $f^{(i)}(0)=t_i$ for $0\leq i\leq n-1$, satisfying \eqref{eqn:non-linear-ODE-1}.  A weak form of our conjecture is:

\begin{conjecture}\label{conjecture:fg-power-series}
	The following are equivalent:
	\begin{enumerate}
	\item (algebraicity) the power series $f(z)$ is algebraic over $\mathbb{Q}(z)$;
	\item (integrality) there exists $N$ such that $a_i\in \mathbb{Z}[\frac{1}{N}]$ for all $i$;
	\item ($\omega(p)$-integrality) There exists a function $\omega(p): \text{Primes}\to \mathbb{Z}$ with $$\lim_{p\to \infty} \frac{\omega(p)}{p}=\infty$$ such that, for each prime $p$, the rational numbers $a_0, a_1, \cdots, a_{\omega(p)}$ are in $\mathbb{Z}_{(p)}$, the ring of rational numbers whose denominators are prime to $p$. 
	\end{enumerate}
\end{conjecture}
Here the integrality condition (2) says that there are only finitely many primes dividing the denominators of the Taylor coefficients $a_i$ of $f$; condition (3), which is a priori weaker, only asserts that for all $m>0$ and sufficiently large $p$ (in terms of $m$), the first $m p$ coefficients of the Taylor expansion of $f$ have denominators which are prime to $p$. For example, one could take $\omega(p)=\lfloor p\log\log\log\log p\rfloor$. That (1) implies (2) is a theorem of Eisenstein from 1852 \cite{eisenstein1852allgemeine}; (2) trivially implies (3). 

\begin{remark}
Eisenstein's paper \cite{eisenstein1852allgemeine} was the final paper he published before he died of tuberculosis at the age of 29. In it he writes: 
\begin{quotation}
	Die wichtigsten Anwendungen der so erhaltenen Sätze habe ich auf Fälle gemacht, in denen die algebraischen Funktionen als Integrale von Differential-Gleichungen definirt werden, und diese Differential-Gleichungen für einfache Reihen-Entwicklung geeignet sind, während die vielleicht sehr complicirte Darstellung in endlicher Form ganz unbekannt bleibt und für diesen Zweck auch wirklich ganz aus dem Spiele gelassen werden kann. Das Einzelne der hieranf bezüglichen Untersuchungen mag für eine künftige Mittheilung vorbehalten bleiben.\footnote{Or, in English (translation by the second author, with some help from Franz Lemmermeyer \cite{MOQ}):
\begin{quotation}
	The most important applications of the theorems thus obtained I have made to cases in which the algebraic function is defined as the integral of a differential equation suitable for simple series expansion, while a perhaps-very-complicated representation in finite form remains completely unknown and can actually be completely ignored. The details of the investigations relating to this may be reserved for a future communication.
\end{quotation}
}

\end{quotation}It is apparently unknown \cite{MOQ} what applications to differential equations Eisenstein had in mind due to his untimely passing; the authors found some pleasure in using his ideas for an application towards classifying their algebraic solutions.
\end{remark}

\begin{remark}\label{remark:hypergeometric-globally-bounded}
Note that if the ODE \eqref{eqn:non-linear-ODE-1} is linear, $0$ is by assumption a non-singular point of the ODE. This assumption is necessary---there exist linear ODE with non-algebraic solutions whose Taylor series expansions about a \emph{singular} point are integral in the sense of \autoref{conjecture:fg-power-series}(2). For example, the non-algebraic hypergeometric function  $$F(t)=\frac{1}{\pi}\int_0^1 \frac{dx}{\sqrt{x(1-x)(1-tx)}}=\sum_{n=0}^\infty \frac{1}{16^n}{2n \choose n}^2t^n$$ has only powers of $2$ in the denominators of the coefficients of its Taylor expansion about $0$. This does not contradict \autoref{conjecture:fg-power-series} as the linear differential equation satisfied by $F(t)$ is singular at $0$.\footnote{This example is taken from \cite{roberts2022hypergeometric}.} See \autoref{remark:taylor-expansion-hypergeometric} for a brief discussion of the Taylor expansion of $F$ about other points.

 See ~\autoref{section:elliptic-curves} and \autoref{section:hypergeometric} for further related discussion, and see \cite{christol1986fonctions} for a systematic study of hypergeometric functions with integral Taylor expansions.
\end{remark}
\begin{remark}
	In the case \eqref{eqn:non-linear-ODE-1} is linear, some form of the equivalence between (1) and (2) in \autoref{conjecture:fg-power-series} is sometimes referred to as the Andr\'e-Christol conjecture; see e.g.~\cite[Remark 5.3.2]{andre2004conjecture} and \cite[Conjecture 1.12]{bostan2023algebraic}. This is a variant of Christol's \cite[Conjecture 4]{gilles2006globally}, which is primarily concerned with the case of \emph{singular} points of linear ODE. For linear ODE, \autoref{conjecture:fg-power-series} is a priori much stronger than the $p$-curvature conjecture, and in fact it makes concrete predictions about linear ODE that we do not know how to deduce from the $p$-curvature conjecture; for example, see \autoref{conj:singular-points}. In the non-linear case, our conjecture is closely related to \cite[Conjecture 3.2]{movasati2024leaf}.
\end{remark}
\begin{remark}
Condition (3) of \autoref{conjecture:fg-power-series} is sharp in the sense that, for any integer $m>0$, there exist non-algebraic functions satisfying a (linear) algebraic differential equation whose first $m p$ Taylor coefficients have denominators prime to $p$, for all $p$. For example, this holds true for the functions $$e^{z^m}=\sum_n \frac{1}{n!}z^{nm}.$$
\end{remark}

In fact we formulate a generalization of \autoref{conjecture:fg-power-series} for arbitrary foliations. Putting off the precise definitions until \autoref{section:foliations}, we conjecture:
\begin{conjecture}\label{conjecture:general-foliation-conjecture-intro}
	Let $R\subset \mathbb{C}$ be a finitely-generated $\mathbb{Z}$-algebra and $X$ a smooth $R$-scheme equipped with a foliation $\mathscr{F}\subset T_X$. Let $x\in X(\mathbb{C})$ be a point. Then the following are equivalent:
	\begin{enumerate}
		\item (algebraicity) The leaf of $\mathscr{F}_{\mathbb{C}}$ through $x$ is algebraic;
		\item (integrality) The formal leaf of $\mathscr{F}_{\mathbb{C}}$ through $x$ descends to a finitely-generated $R$-algebra;
		\item ($\omega(p)$-integrality) The formal leaf  $\mathscr{F}_{\mathbb{C}}$ through $x$ is $\omega(p)$-integral.
	\end{enumerate}
\end{conjecture}
Here the integrality condition (2) means, loosely speaking, that the coefficients of the formal power series used to define the formal leaf lie in a finitely-generated $R$-algebra. Condition (3) is an analogue of condition (3) of \autoref{conjecture:fg-power-series}; we defer it to \autoref{section:foliations}. We explain why \autoref{conjecture:general-foliation-conjecture-intro} implies \autoref{conjecture:fg-power-series} in \autoref{subsection:algebraic-sections}.

\begin{remark}
There already exists a well-known analogue of the $p$-curvature conjecture for foliations---the Ekedahl--Shepherd-Barron--Taylor F-conjecture \cite{esbt}, which predicts that a foliation closed under $p$-th powers mod $p$ for almost all $p$ has \emph{all} of its leaves algebraic. But a foliation may have sporadic algebraic leaves (respectively, a non-linear ODE may have sporadic algebraic solutions) and these are often of great interest; see for example the quest to classify algebraic solutions to the Painlev\'e VI and Schlesinger equations over much of the 20th and early 21st century.\footnote{The classification of algebraic solutions to the Painlev\'e VI equation was completed in \cite{lisovyy-tykhyy}, who proved that the list of algebraic solutions found over the previous century (by Fuchs, Hitchin, Boalch, Dubrovin-Mazzocco, and others) was complete; see \cite{P6survey, boalch2007towards} for the first complete list of algebraic solutions. Algebraic solutions to the rank 2 Schlesinger system were recently completely classified by combining \cite{lam2023finite} and \cite{bronstein2024tykhyy}.}    One motivation for this paper is to try to formulate (and prove cases of) a conjecture that characterizes such sporadic algebraic leaves. That said, we will return to the Ekedahl--Shepherd-Barron--Taylor conjecture in the context of isomonodromy foliations in upcoming work.
\end{remark}

%We summarize schematically the analogies between linear and non-linear differential equations in the following table, highlighting the concepts most relevant for this paper.
%\begin{table}[h!]
%    \centering
%    \renewcommand{\arraystretch}{1.0} % Adjust row height
%    \begin{tabular}{|c | c |}
%        \hline
%        \textbf{Linear} & \textbf{Non-linear} \\ \hline
%        flat connections & foliations \\ \hline
%        $\pi_1$-representations  & $\pi_1$-actions \\ \hline
%        algebraic solution  & algebraic leaf \\ \hline
%        $p$-curvature conjecture & Ekedahl--Shepherd-Barron--Taylor conjecture\\ \hline
%        Andr\'e--Christol conjecture & \autoref{conj:main-conjecture}\\ \hline
%        Picard-Fuchs equations & isomonodromy equations \\ \hline
%        $(p,p)$-classes in cohomology & variations of Hodge structure \\ \hline
%    \end{tabular}
%    \caption{Analogies}
%\end{table}

%For the rest of the introduction we ignore the somewhat technical \autoref{conjecture:fg-power-series}(3); we will return to it in \autoref{section:foliations}.

\subsection{Main results on non-linear differential equations}
We now discuss the main cases in which we can verify \autoref{conjecture:fg-power-series} and its generalizations, beginning with the class of \emph{non-linear} differential equations and initial conditions where our results apply.

\subsubsection{Picard-Fuchs equations and the non-abelian Gauss-Manin connection} \label{subsubsection:NAGM}
Let $k$ be a field of characteristic zero, $\overline{X}/k$ a smooth projective variety, and $D\subset \overline{X}$ a simple normal crossings divisor. Let $X=\overline{X}\setminus D$. A flat bundle $(\mathscr{E}, \nabla: \mathscr{E}\to \mathscr{E}\otimes \Omega^1_{\overline{X}}(\log D))$ with regular singularities along $D$ is a \emph{Picard-Fuchs equation} if there exists $\pi: Y\to X$ smooth projective so that $$(\mathscr{E}, \nabla)|_X\simeq R^i\pi_*\Omega^\bullet_{dR, Y/X}$$ with its Gauss-Manin connection $\nabla_{GM}$, and if the eigenvalues of the residues of $\nabla$ about the components of $D$ are rational numbers in $[0,1)$. (Note that summands of Picard-Fuchs equations are not necessarily Picard-Fuchs equations.) Katz \cite{katz-p-curvature} famously proved the Grothendieck-Katz $p$-curvature conjecture for such flat bundles.

The non-abelian counterpart to this notion is the ``non-abelian Gauss-Manin connection," or isomonodromy foliation, as we now explain.

\begin{notation}\label{notation:X/S}
	Let $R\subset \mb{C}$ be a ring which is finitely generated over $\mb{Z}$, and denote by $\fracr$ the fraction field of $R$. Suppose  $S\rightarrow \Spec(R)$ is a smooth scheme. Let $\overline{f}: \overline{{X}}\to S$ be a smooth projective morphism and ${D}\subset {\overline{X}}$ a simple normal crossings divisor relative to $S$. Let ${X}=\overline{{X}}\setminus {D}$, and let $f: {X}\to S$ be the structure morphism. Fix $s\in S(R)$.
\end{notation}

 For $r\geq 1$,  $\mathscr{M}_{dR}({X}/S, r)$ is the stack of flat vector bundles of rank $r$ on $\overline{{X}}/S$ with regular singularities along ${D}$. Informally, this stack carries a foliation, whose complex-analytic leaves parametrize families of flat bundles $({X}_t, \mathscr{E}_t, \nabla_t)$, such that the conjugacy class of the monodromy representation $$\pi_1(\overline{{X}}_t\setminus {D}_t)\to GL_r(\mathbb{C})$$ associated to $(\mathscr{E}_t, \nabla_t)$ is independent of $t$.\footnote{This description is informal in a number of ways---most notably that $\mathscr{M}_{dR}({X}/S)$ is not smooth, nor even a scheme. The correct analogue of a foliation in this setting is a crystal (of functors or stacks). See the appendices for more detail.}  Here we use that the homotopy type of $\overline{{X}}_t\setminus {D}_t$ is independent of $t$ to identify $\pi_1(\overline{{X}}_t\setminus {D}_t)$ for different $t$.

The leaves of this foliation admit a simple algebraic description in characteristic zero. Namely, let $k$ be a field of characteristic zero and $T$ a smooth (formal) $k$-scheme equipped with a map to $S$. A logarithmic flat bundle $$(\mathscr{E}, \nabla: \mathscr{E} \to \mathscr{E}\otimes \Omega^1_{\overline{{X}}_T/T}(\log {D}_T))$$ on $(\overline{{X}}_T, {D}_T)/T$ with regular singularities along ${D}_T$ is \emph{isomonodromic} if $\nabla$ extends to an absolute flat connection $$\widetilde\nabla: \mathscr{E} \to \mathscr{E}\otimes \Omega^1_{\overline{{X}}_T/k}(\log {D}_T).$$ Given $0\in T$ and a flat bundle $(\mathscr{E}_0, \nabla_0)$ on $(\overline{{X}}_0, {D}_0)$, we say that an isomonodromic bundle $(\mathscr{F}, \nabla)$ on  $(\overline{{X}}_T, {D}_T)/T$ equipped with an isomorphism $(\mathscr{F}, \nabla)_0\overset{\sim}{\to} (\mathscr{E}_0, \nabla_0)$ is an \emph{isomonodromic deformation} of $(\mathscr{E}_0, \nabla_0)$; such isomonodromic deformations, if they exist, are unique up to isomorphism. These definitions work equally well if $T$ is a smooth complex-analytic space, and it is not hard to see that in this case isomonodromic deformations always exist if $T$ is contractible. Likewise, if $T=\on{Spf} k[[t_1, \cdots, t_n]]$ is the formal spectrum of a power series ring over a field of characteristic zero, and $0$ is the the unique closed point of $T$, isomonodromic deformations always exist (see \autoref{prop:isomonodromy-existence-uniqueness}). 

Many classical non-linear differential equations come from isomonodromy:
\begin{example}[Schlesinger equations and Painlev\'e VI]\label{example:schlesinger}
Let $\Delta\subset (\mathbb{P}_{\mathbb{C}}^1)^n$ be the big diagonal, i.e.~the subscheme $$\Delta:=\{(a_1, \cdots, a_n)\mid a_i= a_j \text{ for some } i\neq j\}.$$ Let $S=(\mathbb{P}_{\mathbb{C}}^1)^n\setminus \Delta$ and let $X=\mathbb{P}^1_S$. Let $D\subset X$ be the divisor (\'etale over $S$) whose fiber over $(a_1, \cdots, a_n)\in S$ 
is $a_1+\cdots+a_n$. Fix $s=(a_1, \cdots, a_n)\in S$, and suppose for simplicity that no $a_i=\infty\subset\mathbb{P}^1_{\mathbb{C}}$. Let $(\mathscr{E}, \nabla)$ be a Fuchsian flat connection on $(\mathbb{P}^1_{\mathbb{C}}, a_1+\cdots+a_n),$ 
i.e.~$\mathscr{E}=\mathscr{O}^{\oplus r}$ and $$\nabla=d+\sum_{i=1}^n A_i\cdot \frac{dz}{z-a_i},$$ with $A_i\in \on{Mat}_{r \times r}(\mathbb{C}), \sum A_i=0.$ As one varies the $a_i$, how should one vary $A_i$ so that  the monodromy of $(\mathscr{E}, \nabla)$ does not change, up to conjugacy? Such a variation is \emph{isomonodromic} by definition; the answer is given by the Schlesinger system \cite{schlesinger1912klasse, malgrange1983deformations}:
\begin{equation}
    \begin{dcases}
      \frac{\partial A_i}{\partial a_j}=\frac{[A_i, A_j]}{a_i-a_j} & i \neq j\\
   	  \sum_i \frac{\partial A_i}{\partial a_j}=0 & \forall j
    \end{dcases}\,.
\end{equation}
Here we view the $A_i$ as $\on{Mat}_{r\times r}(\mathbb{C})$-valued functions of $a_1, \cdots, a_n$. If all the $A_i$ lie in $\mathfrak{sl}_2(\mathbb{C})$ and $n=4$, this is, up to a change of coordinates, the Painlev\'e VI equation 	\cite{fuchs1905quelque, jimbo1981monodromy}.
\end{example}

Our main theorem in the non-linear setting verifies \autoref{conjecture:general-foliation-conjecture-intro} for isomonodromy foliations, at initial conditions which are Picard-Fuchs equations.
\begin{theorem}\label{thm:NA-main-intro}
With notation as in \autoref{notation:X/S}, suppose that $[(\mathscr{E}, \nabla)]\in \mathscr{M}_{dR}({X}_s, r)$ is a Picard-Fuchs equation. Then the following are equivalent:
\begin{enumerate}
	\item The leaf of the isomonodromy foliation on $\mathscr{M}_{dR}({X}/S, r)$ through $[(\mathscr{E}, \nabla)]$ is algebraic;
	\item The formal leaf of the isomonodromy foliation through $[(\mathscr{E}, \nabla)]$ is integral;
	\item The formal leaf of the isomonodromy foliation through $[(\mathscr{E}, \nabla)]$  is $\omega(p)$-integral.
\end{enumerate}
\end{theorem}
See \autoref{thm:NA-main} for a precise statement; we complete the proof in \autoref{section:proof-of-main-theorem}.

We view this result as an analogue, for \autoref{conjecture:fg-power-series}, of Katz's verification of the $p$-curvature conjecture for Picard-Fuchs equations. While we've written the result in some generality, note that it can be specialized to many concrete and classical non-linear ODE, e.g.~those in \autoref{example:schlesinger}.
\begin{remark}
The results on isomonodromy discussed here, including \autoref{thm:NA-main-intro}, were announced	in \cite[\S4]{litt2024motives}, written for the 2024 Current Developments in Mathematics conference at Harvard.
\end{remark}

\subsubsection{Relationship to the $p$-curvature conjecture}\label{subsubsection:intro-relationship-to-groth-katz}
It is not hard to see that \autoref{conjecture:fg-power-series} or \autoref{conjecture:general-foliation-conjecture-intro} imply the classical Grothendieck-Katz $p$-curvature conjecture; indeed, if $(\mathscr{E},\nabla)$ has vanishing $p$ curvature mod $p$ for almost all $p$, then its flat sections are $\omega(p)$-integral---see e.g.~\cite[Proposition 3.9]{bost2001algebraic}. What may be more surprising is that verifying \autoref{conjecture:general-foliation-conjecture-intro} for certain isomonodromy foliations---and arbitrary initial conditions---would imply the $p$-curvature conjecture for all flat bundles on all varieties. Indeed, we have:
\begin{proposition}\label{prop:isomonodromy-and-groth-katz}
	Let $r$ be a positive integer.  Suppose \autoref{conjecture:general-foliation-conjecture-intro} holds true for the isomonodromy foliation on $\mathscr{M}_{dR}(\mathscr{C}_g/\mathscr{M}_g, r)$ for all $g\gg r$, where $\mathscr{M}_g$ is the Deligne-Mumford moduli stack of smooth proper curves of genus $g$ and $\mathscr{C}_g$ is the universal curve. Then the $p$-curvature conjecture holds for all rank $r$ flat bundles on all varieties.
\end{proposition}
We give a precise statement and proof of \autoref{prop:isomonodromy-and-groth-katz} in \autoref{section:isomonodromy-and-groth-katz}.
\subsection{Main results for linear ODE}
We were surprised to discover that even for linear ODE, not much is known about \autoref{conjecture:fg-power-series}. We have included a number of results in this direction; in the cases they apply, they refine the main result of \cite{katz-p-curvature}.

For the reader's convenience we restate our conjecture in this case.
\begin{conjecture}\label{conjecture:flat-bundles}
	Let $R\subset\mathbb{C}$ be a finitely-generated $\mathbb{Z}$-algebra,  $X$  a smooth $R$-scheme, and $(\mathscr{E}, \nabla)$ a flat vector bundle on $X/R$. Let $x\in X(\mathbb{C})$ be a point, and fix $v\in \mathscr{E}_x$. Then the formal flat section to $\mathscr{E}$ through $v$ is algebraic if and only if it is integral, if and only if it is $\omega(p)$-integral.
\end{conjecture}

If $X$ is an open subset of $\mathbb{A}^1_{R}$ containing $0$, integrality of a formal flat section $s$ boils down to the statement that the coefficients of the power series defining $s$ generate a finitely-generated $\mathbb{Z}$-algebra; $\omega(p)$-integrality boils down to a version of \autoref{conjecture:fg-power-series}(3). See \autoref{section:foliations} for precise definitions. 

Given $x\in X$, we refer to the fiber $\mathscr{E}_x$ of $\mathscr{E}$ at $x$ as the space of \emph{initial conditions} for $(\mathscr{E}, \nabla)$. A formal flat section $s$ at $x$ has initial condition $s(x)\in \mathscr{E}_x$. 
\subsubsection{Picard-Fuchs equations}
Our main result towards \autoref{conjecture:flat-bundles} in the linear case is that we can verify the conjecture for Picard-Fuchs equations at initial conditions in the image of the cycle class map. 
\begin{theorem}\label{thm:intro-cycle-class-initial-conditions}
Let $R\subset\mathbb{C}$ be a finitely-generated $\mathbb{Z}$-algebra, $X$ a smooth variety over $R$, and $f: Y\to X$ a smooth projective morphism. Suppose $(\mathscr{E},\nabla)=(R^{2i}f_*\Omega_{dR, Y/X}^\bullet, \nabla_{GM})$ is a Picard-Fuchs equation. Let $x\in X(\mathbb{C})$ be a point at and $s$ an integral (resp.~$\omega(p)$-integral) formal flat section to $(\mathscr{E},\nabla)$ at $x$. If the initial condition $s(x)\in H^{2i}_{dR}(Y_x)$ is in the image of the cycle class map $$Z^i(Y_x)_{\mathbb{C}}\to H^{2i}_{dR}(Y_x),$$ then $s$ is algebraic.
\end{theorem}

See \autoref{thm:cycle-class-initial-conditions} for a precise statement. Arguably \autoref{thm:NA-main-intro} is a non-abelian analogue of this result.

We also verify \autoref{conjecture:flat-bundles} for a number of differential equations of algebro-geometric interest, with no assumptions on the initial conditions---for example, for Picard-Fuchs equations arising from families of elliptic curves:
\begin{theorem}\label{thm:ell-curves-intro}
	Let $X$ be a smooth variety over $\mathbb{C}$ and $f: Y\to X$ a family of elliptic curves. Let $(\mathscr{E},\nabla)=(R^{1}f_*\Omega_{dR, Y/X}^\bullet, \nabla_{GM})$. Let $x\in X(\mathbb{C})$ be a point at and $s$ an integral formal flat section to $(\mathscr{E},\nabla)$. Then $s$ is algebraic.
\end{theorem}
See \autoref{theorem:elliptic-curves} for a more general statement, allowing for symmetric powers of bundles as in \autoref{thm:ell-curves-intro}. In particular, this latter theorem applies to the solutions to Picard-Fuchs equations associated to modular forms; see e.g.~\autoref{remark:taylor-expansion-hypergeometric} for an example.

One pleasant aspect of \autoref{thm:ell-curves-intro} is that it can be made very concrete---indeed, one may explicitly write linear recurrences controlling the coefficients of Taylor expansions of solutions to Picard-Fuchs equations as in \autoref{thm:ell-curves-intro}. Then the theorem says that there are infinitely many primes dividing the denominators of the terms of these recurrences. See e.g.~\autoref{prop:ell-curve-corollary} for an explicit example.
\subsubsection{Singular points}
The following is a local prediction of our main conjecture.
\begin{conjecture}\label{conj:singular-points}
	Let $X/\mathbb{C}$ be a smooth projective curve, $D\subset X$ a reduced effective divisor, and $(\mathscr{E}, \nabla: \mathscr{E}\to \mathscr{E}\otimes \Omega^1_X(\log D))$ a logarithmic flat bundle on $(X, D)$. Let $x\in D$ be a point and let $s$ be an integral formal flat section to $(\mathscr{E}, \nabla)$ at $x$. Then either $s$ is algebraic, or $(\mathscr{E}, \nabla)$ has infinite monodromy at $x$.
\end{conjecture}
Indeed, suppose we are in the situation of \autoref{conj:singular-points}, and $(\mathscr{E}, \nabla)$ has finite monodromy at $x$. Then passing to a ramified cover $X'$, we may assume $(\mathscr{E}, \nabla)$ has trivial monodromy; after replacing $\mathscr{E}$ with a modification we may assume the residue of $\nabla$ at $x$ is zero, whence we may remove it from $D$. Then our conjecture predicts the pullback of $s$ to $X'$ is algebraic, and hence the same is true for $s$.

By combining work of Christol \cite{christol1986fonctions} and Beukers-Heckman \cite{beukersheckman}, we are able to check that \autoref{conj:singular-points} holds true for hypergeometric functions, in \autoref{section:hypergeometric}; see \autoref{prop: hyp-result}.
\subsection{Proofs}\label{subsection:intro-outline}
We now briefly outline the proofs of \autoref{thm:NA-main-intro} and \autoref{thm:intro-cycle-class-initial-conditions}. See \autoref{subsection:outline-of-proof} for a much more detailed outline of the proof of \autoref{thm:NA-main-intro}.

Both proofs have an algebraic step and an analytic step. In the algebraic step, we show that some property of our initial conditions---underlying a $\mathbb{Z}$-variation of Hodge structure in the case of \autoref{thm:NA-main-intro}; being a $\mathbb{C}$-linear combination of Hodge classes in the case of \autoref{thm:intro-cycle-class-initial-conditions}---extends to the formal leaf of the foliation we are studying. In the analytic step, we show that this implies some boundedness, which we may use to conclude that the leaf of our foliation is finitely branching, and hence (ultimately by some application of the theory of regular singularities) algebraic. As the analytic steps here are largely contained in previous work, we focus on the algebraic steps in the sketch below.
\subsubsection{Proof of \autoref{thm:intro-cycle-class-initial-conditions}}
Here we start with a $\mathbb{C}$-linear combination of cycle classes  $s(x)\in H^{2i}_{dR}(Y_x)$; in particular, $s(x)$ is a $\mathbb{C}$-linear combination of Hodge classes. Under the assumption the formal solution $s$ to the Gauss-Manin connection is ($\omega(p)$-)integral, we would like to show that this property extends over the formal neighborhood of $x$ in $X$. Ultimately the proof of this fact relies on some (relative) Fontaine-Laffaille theory. We show that we may write $s(x)$ as a $\mathbb{C}$-linear combination of $\mathbb{Q}$-cycle classes $e_i$ such that each $e_i$ remains a Hodge class in the formal neighborhood of $x$. Letting $s_i$ be the formal solution to the Gauss-Manin connection through $e_i$, we wish to show that $s_i$ remains in the appropriate piece of the Hodge filtration, $F^i_{\text{Hodge}}R^{2i}f_*\Omega^\bullet_{dR, Y/X}$. We do so after reduction mod $p$ for infinitely many $p$ by an analysis of the action of the crystalline Frobenius on $e_i$. See the proof of \autoref{prop:crystalline-step} for details.

Now, the analytic step---finite branching of the corresponding holomorphic solution to the Gauss-Manin connection---follows from the existence of a polarization on $R^{2i}f_*\Omega^\bullet_{dR, Y/X}$.

\subsubsection{Proof of \autoref{thm:NA-main-intro}} We now give a brief sketch of the proof of our main theorem; see \autoref{subsection:outline-of-proof} for a substantially more detailed sketch. With notation as in \autoref{notation:X/S}, we start with a Picard-Fuchs equation $(\mathscr{E},\nabla)$ on $(\overline{X}_s, D_s)$; in particular, after base change to $\mathbb{C}$, $(\mathscr{E},\nabla)$ underlies a  polarizable $\mathbb{Z}$-variation of Hodge structure. The algebraic step, performed in \autoref{theorem:extending-the-Hodge-filtration}, amounts to showing that the Hodge filtration on $(\mathscr{E},\nabla)$ extends Griffiths-transversally to the formal isomonodromic deformation of $(\mathscr{E},\nabla)$.

We achieve this by comparing the obstruction to doing so (modulo $p$) to the $p$-curvature of the isomonodromy foliation, which we compute in \autoref{appendix:NAGM}. The comparison relies on the positive characteristic non-Abelian Hodge theory of Ogus-Vologodsky \cite{ogus-vologodsky} and Schepler \cite{schepler2005logarithmic}; the proof of vanishing when the leaf of the isomonodromy foliation is ($\omega(p)$-)integral ultimately relies on the Higgs-de Rham flow of Lan-Sheng-Zuo \cite{lan2019semistable}, as exposited by Esnault-Groechenig \cite{eg_revisit}.

\begin{remark}
	The computation of the $p$-curvature of the isomonodromy foliation  in \autoref{appendix:NAGM} may be thought of as a non-abelian analogue of Katz's formula for the $p$-curvature of a Gauss-Manin connection \cite[Theorem 3.2]{katz-p-curvature}. We will make this precise in upcoming work.
\end{remark}
Having extended the Hodge filtration as desired, we find that the leaf of the isomonodromy foliation through $(\mathscr{E}, \nabla)$ is finitely-branching via the closedness of so-called non-abelian Noether-Lefschetz (or non-abelian Hodge) loci; see \autoref{section:hodge-preliminaries} for details. Ultimately finiteness comes from a generalization of Deligne's finiteness theorem for $\mathbb{C}$-local systems underlying a $\mathbb{Z}$-variation of Hodge structure \cite{deligne1987theoreme}, due to Simpson \cite{simpson62hodge}, Esnault-Kerz \cite{esnault2024non}, and others; see also \cite[Corollary 4.4.3]{litt2024motives}.
\subsection{Algebro-geometric consequences}
In \autoref{section-variational-conjectures}, we discuss some implications of our results for various conjectures on local systems in algebraic geometry. In particular we study consequences of \autoref{conjecture:variational-motivic}, itself a consequence of the relative Fontaine-Mazur conjecture \cite[Conjecture 1 bis]{petrov2023geometrically}. \autoref{conjecture:variational-motivic} predicts that if $f: X\to S$ is a smooth proper morphism of (say) smooth $\mathbb{C}$-schemes, $s\in S(\mathbb{C})$, and $\mathbb{V}$ is a local system on $X$ with $\mathbb{V}|_{X_s}$ of geometric origin, then $\mathbb{V}|_{X_t}$ is of geometric origin for all $t\in S(\mathbb{C})$. (See \autoref{defn:geometric-origin} for the definition of ``geometric origin," which loosely means that a local system appears in the cohomology of a family of varieties over the base.)

Combined with \autoref{conjecture:general-foliation-conjecture-intro}, this predicts that if $(\mathscr{E}, \nabla)$ is of geometric origin and admits an ($\omega(p)$-)integral (and hence algebraic) formal isomonodromic deformation, then the isomonodromic deformation is also of geometric origin. In particular, any properties of flat bundles that hold for bundles of geometric origin---for example, the property of underling a $\mathbb{Z}$-variation of Hodge structure, or having nilpotent $p$-curvature mod $p$ for almost all $p$---should persist to the isomonodromic deformation. 

Under the assumption that $(\mathscr{E}, \nabla)$ admit an algebraic isomonodromic deformation, \cite{katzarkov1999non, jost2001harmonic, esnault2024non} show that the isomonodromic deformation underlies a variation of Hodge structure. In \autoref{theorem:extending-the-Hodge-filtration} we prove that (assuming $(\mathscr{E}, \nabla)$ is a Picard-Fuchs equation) ($\omega(p)$-)integrality is enough to guarantee that the Hodge filtration extends to the formal isomonodromic deformation. Finally, we show that the conjugate filtration mod $p$ extends to the isomonodromic deformation, and in particular the $p$-curvature of the isomonodromic deformation is nilpotent mod $p$ for almost all $p$.

\subsection{Related work}
The primary antecedent to this work is Katz's paper \cite{katz-p-curvature}, in which he proves the $p$-curvature conjecture for Gauss-Manin connections---the bulk of this paper is devoted to a ``non-abelian" analogue of his result. We were inspired by the PhD theses \cite{papaioannou2013algebraic, menzies2019p}, which study the $p$-curvature of the isomonodromy foliation. The papers \cite{shankarsimpleloops, patel2021rank} are also arguably in this spirit, and in fact the former is what got the second author interested in this subject. To our knowledge Kisin was the first person to consider the $p$-curvature conjecture in the context of isomonodromy, and this paper owes him a substantial debt.

Arguably this work belongs to a mathematical tradition that goes back to Siegel \cite{siegel2014einige}---the theory of $G$-functions. That theory is too large to summarize here, but see \cite{andre1989g} for a beautiful overview. The work of Andr\'e \cite{andre2004conjecture}, Bost \cite{bost2001algebraic}, and Chudnovsky-Chudnovsky \cite{chudnovsky2006applications} on the solvable case of the $p$-curvature conjecture, of which we make use in \autoref{section:elliptic-curves}, is part of this tradition. We will return to the theory of $G$-functions in the context of isomonodromy in future work. Our main conjecture is closely related to some work in this area, particularly that of Christol \cite{christol1986fonctions}.

One goal of this work was to introduce arithmetic to the study of algebraic solutions to non-linear differential equations. That study is a field in its own right; in the context of isomonodromy, classifications of algebraic solutions are due in different cases to Dubrovin-Mazzocco \cite{dubrovin-mazzocco}, Boalch \cite{P6survey, boalch2007towards}, Lisovyy-Tykhyy \cite{lisovyy-tykhyy} in the case of Painlev\'e VI, and \cite{lam2023finite, bronstein2024tykhyy} in the case of general rank 2 Schlesinger systems.

While our methods are largely unrelated to his, Bost's beautiful paper \cite{bost2001algebraic} proves algebraicity of leaves of certain foliations under arithmetic conditions related to ours, assuming some additional Archimedean hypotheses. See also \cite{movasati2024leaf} for  conjectures on algebraicity of leaves of foliations closely related to our main conjectures here.  
\subsection{Acknowledgements}
We are extremely grateful to Greg Baldi, Yohan Brunebarbe, Fran\c{c}ois Charles, Vesselin Dimitrov, H\'el\`ene Esnault, Michael Groechenig, Nick Katz, Moritz Kerz, Bruno Klingler, Mark Kisin, Raju Krishnamoorthy, Aaron Landesman, Ngaiming Mok, Hossein Movasati, Sasha Petrov, Ananth Shankar, Mao Sheng, Yunqing Tang, Vadim Vologodsky, and Kang Zuo for many useful conversations over the course of this project. In particular, discussions with Moritz Kerz and H\'el\`ene Esnault were crucial to \autoref{section:hodge-preliminaries}. Lam was supported by a BMS Dirichlet Fellowship and a DFG Walter Benjamin grant. Litt was supported by an NSERC Discovery Grant, a Sloan Research Fellowship, and an Ontario Early Researcher Award.

\section{Conjectures and basic properties of foliations}\label{section:foliations}
We now prepare to precisely state our main conjecture in the language of foliations, and to develop some basic notions around integrality of leaves, $p$-curvature of foliations, and so on.
\subsection{Foliations}
\begin{definition}
	\label{definition:foliation} Let $R$ be a ring and let $X$ be a smooth $R$-scheme. A subbundle $\mathscr{F}\subset T_X$ is a \emph{foliation} if it is closed under the Lie bracket. If $f: X\to S$ is a smooth morphism of smooth $R$-schemes, an integrable  foliation $\mathscr{F}\subset T_X$ is a \emph{horizontal foliation} (with respect to $f$) if the composition $$\mathscr{F}\hookrightarrow T_X\overset{df}{\longrightarrow} f^*T_S$$ is an isomorphism. Put another way, a horizontal foliation is a splitting of the natural exact sequence $$0\to T_f\to T_X\to f^*T_S\to 0$$ compatible with the Lie algebra structure.
\end{definition}
\begin{definition} A (formal) leaf of $\mathscr{F}$ is a smooth (formal) $R$-scheme $Z$ equipped with a map $l: Z\to X$ such that the induced map $$dl: T_Z\to l^*T_X$$ factors through $l^*\mathscr{F}$ and induces an isomorphism $T_Z\overset{\sim}{\to} l^*\mathscr{F}$.	 If $x$ is a $T$-point of $X$ for some $R$-algebra $T$, we say that $Z$ is a \emph{leaf through} $x$ if $x$ factors through $Z$.
\end{definition}
\begin{example}[Analytic foliations]
	The above definitions make equal sense for complex-analytic spaces, where there is a simple construction of many horizontal foliations. Let $(S,s)$ be a smooth pointed connected complex-analytic space with universal cover $\widetilde S$ and let $V$ be a smooth complex-analytic space equipped with an action of $\pi_1(S,s)$. Then setting $X=(\widetilde S \times V)/\pi_1(S,s)$, where the action on $\widetilde S$ is via deck transformations, $X$ is  equipped with a foliation horizontal with respect to the natural map $X\to S$ inherited from projection onto the first coordinate. Indeed, the subsheaf $p_1^*T_{\widetilde S}\subset T_{\widetilde S \times V}$, where $p_1$ is projection onto $\widetilde s$, descends to a horizontal foliation on $X$. The leaves of this foliation are precisely the images of $\widetilde S\times \{v\}\to X$ for $v\in V$. 
	
	Up to mild issues of stackiness and non-smoothness, the complex-analytifications of all foliations considered in this paper will have this form.
\end{example}

\subsection{Integrality of leaves}\label{subsection:integrality-of-leaves}
We now say what it means for a (formal) leaf of a foliation to be integral. 

Let $R$ be an integral domain of characteristic zero. Let $X$ be a smooth $R$-scheme, and let $\mathscr{F}\subset T_X$ be a foliation of rank $d$.
\begin{proposition}  \label{proposition:formal-leaves-exist}
Suppose $\mathbb{Q}\subset R$, and $x$ is an $R$-point of $X$. Then there always exists a unique (up to unique isomorphism) formal leaf $\on{Spf}(R[[x_1, \cdots, x_d]])\to X$ through $x$.	
\end{proposition}
\begin{proof}
This follows from Taylor's formula, see e.g.~\cite[\S3.4.1]{bost2001algebraic}.	
\end{proof}
\begin{definition}
	Suppose that $R\subset \mathbb{C}$ is a finitely-generated $\mathbb{Z}$-algebra. We say that the leaf through $x$ is \emph{integral} if the formal leaf $\on{Spf}(R_\mathbb{Q}[[x_1, \cdots, x_d]])\to X$ factors through $\on{Spf}(R[1/N][[x_1,\cdots, x_d]])$ for some integer $N$. 
\end{definition}

This is the analogue of condition (2) of \autoref{conjecture:fg-power-series} for foliations. Loosely speaking it means that the formal power series defining the leaf through $x$ have coefficients lying in some fixed finitely-generated $\mathbb{Z}$-algebra $R'$. 

\begin{remark}
Note that in our definition of integrality of a leaf, there is no assumption that the integral model is itself a leaf, i.e.~there need not be an isomorphism between $T_{\on{Spf}(R[1/N][[x_1,\cdots, x_d]])}$ and the pullback of $\mathscr{F}$. But as the following proposition shows, this is automatic after some further localization.
\end{remark}

\begin{proposition}[Integral models of leaves]\label{proposition:integral-models-of-leaves}
	Let $R\subset \mathbb{C}$ be a finitely-generated $\mathbb{Z}$-algebra, and set $\mathscr{K}=\on{Frac}(R)$. Let $X$ be a smooth $R$-scheme and $\mathscr{F}\subset T_X$ a foliation. Let $\iota: \on{Spf} \mathscr{K}[[x_1, \cdots, x_{\dim(S)}]]\to X$ be a leaf of $\mathscr{F}_{\mathscr{K}}$. Suppose there exists an integral model of $\iota$, i.e.~a map $\iota^\natural: \on{Spf} R[[x_1, \cdots, x_{\dim{S}}]]\to X$ such that $\iota^\natural \widehat\otimes \mathscr{K}=\iota$. Then there exists $f\in R$ such that the induced map $\iota^\natural[\frac{1}{f}]: \on{Spf} R[1/f][[x_1, \cdots, x_{\dim{S}}]]\to X$ is a leaf of $\mathscr{F}$.
\end{proposition}
	\begin{proof}
	  Set $\mathscr{L}=\on{Spf} R[[x_1, \cdots, x_{\dim{S}}]]$. First, observe that $d\iota: T_{\mathscr{L}}\to \iota^{\natural *}T_X$ factors through $\iota^{\natural*}\mathscr{F}$. Indeed, this is equivalent to showing that the composite map $$T_{\mathscr{L}}\to \iota^{\natural *}T_X\to \iota^{\natural *}T_X/\iota^{\natural *}\mathscr{F}$$ is zero. This composition is zero after taking the completed tensor product $-\widehat{\otimes}\mathscr{K}$, whence it follows that it was already zero as the sheaves in question are torsion-free.
	  
	  Given $f\in R$, set $\mathscr{L}[\frac{1}{f}]=\on{Spf} R[\frac{1}{f}][[x_1, \cdots, x_{\dim{S}}]]$. We wish to find nonzero
	   $f\in R$ such that $d\iota^\natural[\frac{1}{f}]: T_{\mathscr{L}[\frac{1}{f}]}\to \iota^{\natural}[\frac{1}{f}]^*\mathscr{F}$ is an isomorphism. Now consider $\det d\iota^\natural\in \det T_{\mathscr{L}}^\vee\otimes \det \iota^\natural\mathscr{F}.$ It suffices to take $f$ so that this element vanishes nowhere on $\mathscr{L}[\frac{1}{f}]$; such a non-zero $f$ exists as $\iota^\natural$ is invertible on $\mathscr{L}_{\mathscr{K}}$.
	\end{proof}

We now formulate the analogue of \autoref{conjecture:fg-power-series}(3).

Let $\widehat{\mathscr{I}}_d\subset R[[x_1, \cdots, x_d]]$ be the ideal $(x_1, \cdots, x_d)$. For each function $$g: \text{Primes}\to \mathbb{Z}_{\geq 0}$$ we define the subring of
 $g$-integral elements $$R[[x_1, \cdots, x_d]]^g\subset R[[x_1, \cdots, x_d]]\widehat{\otimes} \mathbb{Q}$$ to be the completion of the subring generated by $R[[x_1, \cdots, x_d]]$ and $$\{\widehat{\mathscr{I}}_d^{g(p)+1}\widehat{\otimes}\mathbb{Z}[1/p]\}_{p\in \text{Primes}}$$ at the ideal $$\widehat{\mathscr{I}}_d+\sum_{p\in \text{Primes}}\widehat{\mathscr{I}}_d^{g(p)+1}\widehat{\otimes}\mathbb{Z}[1/p].$$ 
 
 %Note that $R[[x_1, \cdots, x_d]]^g$ is equipped with a natural map to $R$; we denote the kernel of this map by $\widehat{I}_{d,g}$.

%\begin{example}
%Suppose $S=\on{Spf} R[[x_1, \cdots, x_n]],$ and $s$ is defined by the ideal $(x_1, \cdots, x_n)$, so that $\widehat{P}_s=R[[x_1, \cdots, x_n]], \widehat{\mathscr{I}}_s=(x_1, \cdots, x_n)$. Then $\widehat{P}_s^g$ is the subring of $R_\mathbb{Q}[[x_1, \cdots, x_n]]$ consisting of those power series $$f(x_1, \cdots, x_n)=\sum a_Ix^I$$ such that $$a_I\in R\otimes \bigotimes_{g(p_i)\leq |I|}\mathbb{Z}[1/p_i]$$ for all $I$.
%In other words, $\widehat{P}_x^g$ consists of those power series which are $p$-integral to order $g(p)$, for all primes $p$. In fact this is the only example we will need to consider.
%\end{example}

\begin{definition}
	Let $g: \text{Primes}\to \mathbb{Z}_{\geq 0}$ be a function. Let $R\subset \mathbb{C}$ be a finitely generated $\mathbb{Z}$-algebra, $X$ a smooth  $R$-scheme, and $\mathscr{F}\subset T_X$ a foliation. Let $x\in X(R)$ be a point. We say that the leaf through $x$ is \emph{$g$-integral} if the formal leaf $\on{Spf}(R_{\mathbb{Q}}[[x_1, \cdots, x_d]])\to X$ through $x$ factors through $\on{Spf}(R[[x_1, \cdots, x_d]]^g)$, the formal spectrum of the subring of $R_{\mathbb{Q}}[[x_1, \cdots, x_d]]$ of $g$-integral elements.
	
	If $f: \on{Primes}\to \mathbb{Z}_{\geq 0}$ is a function, we say that the leaf through $x$ is $\omega(f)$-integral if there exists some $g$ with $$\lim_{p\to \infty} \frac{g(p)}{f(p)}=\infty$$ such that the leaf through $x$ is $g$-integral. 
\end{definition}
The evident analogue of \autoref{proposition:integral-models-of-leaves} holds for $\omega(f)$-integrality in place of integrality, with an identical proof.

\begin{remark}
Sometimes (as in the introduction) we will allow $x$ to be a $\mathbb{C}$-point of $X$, or begin with a foliated variety $(X, \mathscr{F})$ over $\mathbb{C}$, and still (by abuse of notation) reference integrality or $\omega(f)$-integrality of the leaf through $x$. By this we mean that that there exists a finitely-generated $\mathbb{Z}$-algebra $R$ to which $(X, \mathscr{F}, x)$ descends as a smooth pointed foliated scheme $(X_R, \mathscr{F}_R, x_R)$, and the leaf through $x_R$ is integral or $\omega(f)$-integral.
\end{remark}

Suppose $S$ is an open subscheme of $\mathbb{A}^1_R$ containing zero, $X\to S$ is a smooth morphism and $x\in X(R)$ lies over $0\in S(R)$. Let $\mathscr{F}\subset T_X$ be a foliation horizontal over $S$. That the leaf through $x$ is $g$-integral means, loosely speaking, that the first $g(p)$ coefficients of the Taylor expansion of the functions defining the leaf have denominators prime to $p$.

Suppose $X$ is a scheme over a finitely-generated $\mathbb{Z}$-algebra $R$ and $s$ is an $R$-point of $X$. Suppose we are given sheaves $\mathfrak{E}_{\mathfrak{p}}$ on $X_{\mathfrak{p}}$ for each closed point $\mathfrak{p}$ of $R$, and sections $f_{\mathfrak{p}}\in \Gamma(X_{\mathfrak{p}}, \mathfrak{E}_{\mathfrak{p}})$ for each $\mathfrak{p}$. We will say that $f_{\mathfrak{p}}$ vanishes to order $\omega(f)$ at $s_{\mathfrak{p}}$ if there exists $g: \text{Primes}\to \mathbb{Z}_{\geq 0}$ with $$\lim_{p\to \infty} \frac{g(p)}{f(p)}=\infty$$ such that $f_{\mathfrak{p}}\in \mathscr{I}_{s_\mathfrak{p}}^{g(p)}\mathscr{E}_{\mathfrak{p}}$ for almost all $p$. In this case we will also say $f_\mathfrak{p}$ lies in $\mathscr{I}_{s_\mathfrak{p}}^{\omega(p)}\mathscr{E}_{\mathfrak{p}}$.

\subsection{Algebraic sections to horizontal foliations}\label{subsection:algebraic-sections}
We can now formally state our main conjecture; this is a slightly more refined statement of \autoref{conjecture:general-foliation-conjecture-intro}.
\begin{conjecture}\label{conjecture:general-algebraic-leaves}
	Let $R\subset \mb{C}$ be an integral domain  which is finitely generated over $\mb{Z}$. Let $(X, \mathscr{F})$ be a smooth scheme with a foliation over $R$. Let $x\in X(R)$ be a point. Then the following are equivalent:
	\begin{enumerate}
		\item (algebraicity) the leaf of $\mathscr{F}$ through $x_{\mb{C}}$ is algebraic,
		\item (integrality) the leaf through $x$ is integral, and
		\item ($\omega(p)$-integrality) the leaf through $x$ is $\omega(p)$-integral.
	\end{enumerate}
\end{conjecture} 
We briefly explain why this conjecture implies \autoref{conjecture:fg-power-series}. Suppose we are given an ODE as in \eqref{eqn:non-linear-ODE-1}. Set $S=\mathbb{A}^1_{\mathbb{Z}}$, with coordinate $z$, and $X'=\Spec~\mathbb{Z}[z, y_0, \cdots, y_{n-1}],$ and let $\pi: X'\to S$ be the evident map induced by the coordinate $z$. Let $X\subset X'$ be the open subscheme on which $g$ is defined, and set $$\vec v=\frac{\partial}{\partial z}+\sum_{i=0}^{n-2}y_{i+1}\frac{\partial}{\partial y_i}+g(z, y_0,\cdots, y_{n-1})\frac{\partial}{\partial y_{n-1}}.$$ Set $\mathscr{F}\subset T_X$ to be the span of $\vec v$; this defines a foliation on $X$ horizontal over $S$. (In the language of a first calculus course, the vector field $\vec v$ is the ``phase portrait " of the ODE \eqref{eqn:non-linear-ODE-1}.) 

Then the leaves of this foliation correspond to the solutions to \eqref{eqn:non-linear-ODE-1}: a solution $f$ yields a leaf $z\mapsto (z, f(z), f'(z), \cdots, f^{n-1}(z))$, and the $y_0$ coordinate of a leaf is a solution to \eqref{eqn:non-linear-ODE-1}. Under this correspondence, our notions of algebraicity, integrality, and $\omega(p)$-integrality for leaves of foliations correspond precisely to conditions (1-3) of \autoref{conjecture:fg-power-series}.

\begin{remark}
While we state conjecture \autoref{conjecture:general-foliation-conjecture-intro} for arbitrary foliations, all of our results are for foliations horizontal over some base. The reader might be tempted to believe that the conjecture should only hold in this horizontal setting. However, there is in fact no difference between considering all foliations and foliations horizontal over some base; given a smooth variety $X$ equipped with a foliation $\mathscr{F}\subset T_X$, one can always locally choose a map $X\to S$ such that $\mathscr{F}$ becomes horizontal over $S$. See e.g.~\cite[\S3.4.1]{bost2001algebraic}.
\end{remark}

\subsection{$p$-curvature of a foliation}
If $X$ is a smooth variety over a field $k$ of characteristic $p>0$, and $v$ is a local section to $T_X$ (i.e.~a local derivation), then $v^p$ is also a derivation. 
\begin{definition}\label{defn:p-curv-foliation}
	Let $X$ be a smooth variety over a field $k$ of characteristic $p>0$ and let $\mathscr{F}\subset T_X$ be a foliation. The \emph{$p$-curvature} of $\mathscr{F}$ is the $\mathscr{O}_X$-linear map $\psi_p$ given as the composition 
    $$F_{\text{abs}}^*\mathscr{F}\to T_X \to T_X/\mathscr{F},$$
    where the first map is induced by taking $p$-th powers on $T_X$, and the second is the natural projection. 
\end{definition}
See e.g.~\cite[\S2.1.1]{bost2001algebraic} for further discussion of the $p$-curvature of a foliation. Loosely speaking, the $p$-curvature measures the failure of $\mathscr{F}$ to be preserved under the $p$-th power operation. Foliations with vanishing $p$-curvature are called \emph{$p$-closed}.

\begin{definition}\label{definition:p-power-leaves}
	Let $k$ be a field of characteristic $p>0$, $S$ a smooth $k$-scheme, and $X\to S$ a smooth morphism with $\mathscr{F}\subset T_X$ a foliation, horizontal over $S$. Let $Z\subset S$ be a closed subscheme, and set $Z^{[p]}\subset S$ to be $Z^{[p]}=F_{\text{abs}}^{-1}(Z)$. A $p$-power leaf of $\mathscr{F}$ about $Z$  is a section $\iota: Z^{[p]}\to X$ which satisfies the definition of a leaf of $\mathscr{F}$ apart from smoothness, i.e.~such that the induced map $T_{Z^{[p]}}\to \iota^*T_X$ factors through $\iota^*\mathscr{F}$ and induces an isomorphism $T_{Z^{[p]}}\to \iota^*\mathscr{F}$.
\end{definition}
	\begin{remark}
		Note that $\Omega^1_{Z^{[p]}/k}$ is locally free, despite the fact that $Z^{[p]}$ is almost never smooth. Indeed, the local equations defining $Z^{[p]}$ are, by definition, $p$-th powers, and hence have vanishing differential. Hence $T_{Z^{[p]}}$ is a locally free sheaf of rank equal to $\dim S$.
	\end{remark}
	The following criterion will be useful for checking the vanishing of $p$-curvature of a foliation on certain subschemes of $X$. See also \autoref{prop:crystal-p-power-leaves} for a variant we will use in the proof of \autoref{thm:NA-main-intro}.
\begin{proposition}[Vanishing of $p$-curvature along $p$-power leaves]\label{prop:p-power-leaves}
	With notation as in \autoref{definition:p-power-leaves}, let $\iota: Z^{[p]}\to X$ be a $p$-power leaf of $\mathscr{F}$. Then the $p$-curvature of $\mathscr{F}$, $$\psi_p: F_{\text{abs}}^*\mathscr{F}\to T_X/\mathscr{F}$$ vanishes when restricted to $\iota(Z^{[p]})\subset X$.
\end{proposition}
\begin{proof}
We wish to compute the map $\iota^*F^*_{\text{abs}}\mathscr{F}\to \iota^*T_X/\mathscr{F}$ induced by taking $p$-th powers of derivations. As $\iota$ is a $p$-power leaf, this is the same as the map $T_{Z^{[p]}}\to \iota^*T_X/T_Z^{[p]}$ induced by taking $p$-th powers; but this latter map is zero as $T_{Z^{[p]}}$ is closed under $p$-th powers.
\end{proof}

%\daniel{Do we want to state this as a conjecture?}
%\subsection{Invariant subvarieties}\daniel{Remove this?}
%\begin{conjecture}\label{conjecture:invariant-subvarieties}
%	A subvariety is invariant under a horizontal foliation iff $p$-curvature at all points is contained in the tangent space.
%\end{conjecture}
%\begin{proposition}
%	If $\mathscr{D}=\emptyset$ (i.e.~if $\mathscr{X}/S$ is proper) then \autoref{conjecture:invariant-subvarieties} follows for $\mathscr{X}/S$ from \autoref{conjecture:general-algebraic-leaves} for all proper varieties.
%\end{proposition}
%\begin{proof}
%Hilbert scheme argument, need to worry about singularities on the Hilbert scheme.	
%\end{proof}
\pagebreak

\part{Linear differential equations}
We now  prepare to state and prove our main results on linear differential equations.  This is not meant to be an exhaustive list of the linear differential equations and initial conditions for which \autoref{conjecture:general-algebraic-leaves} can be verified, but rather a sampling of representative results meant to persuade the reader of its plausibility and introduce some of the objects that will be relevant for our results on non-linear (isomonodromy) differential equations in \autoref{part:non-linear} of this paper. We expect that we will study \autoref{conjecture:general-algebraic-leaves} for linear differential equations further in future work.
\section{Picard-Fuchs equations and cycle class initial conditions}
\subsection{The main theorem}
The main goal of this section is to prove \autoref{conjecture:general-algebraic-leaves} for Picard-Fuchs equations, with cycle class initial conditions.
\begin{theorem}\label{thm:cycle-class-initial-conditions}
Let $R\subset \mathbb{C}$ be a finitely-generated $\mathbb{Z}$-algebra with fraction field $K$ and let $f: X\to S$ be a smooth projective morphism over $R$, with $S$ smooth. Fix $s\in S(K)$ and $x\in H^{2i}_{dR}(X_s)_s$ in the image of the cycle class map $$\on{cl}: Z^i(X_s)_{K}\to H^{2i}_{dR}(X_s)=(R^{2i}f_*\Omega^\bullet_{dR, X/S})_s.$$ Then the formal flat section $v$ to $R^{2i}f_*\Omega^\bullet_{dR, X/S}$ through $x$ is algebraic if and only if it is integral (resp.~if and only if it is $\omega(p)$-integral).
\end{theorem}
\subsection{Preliminaries}\label{subsection:berthelot-ogus-preliminaries}
In this section we briefly summarize work of Berthelot, Faltings, Mazur, and Ogus, which explains how to read off the Hodge filtration on the cohomology of a family of algebraic varieties (mod $p$) from the associated $F$-crystal arising from the formalism of crystalline cohomology, at least in good situations. This theory is explained nicely, in slightly different language, in \cite[8.26]{berthelot2015notes}, \cite{faltings1989crystalline}, and \cite{ogus-griffiths}.

Let $k$ be a perfect field of characteristic $p>0$, $S$ a smooth $p$-adically complete formal $W(k)$-scheme, and let $f: X\to S$ be a smooth projective morphism. Let $f_0: X_0\to S_0$ be the reduction of $f$ mod $p$. Let $F: S\to S$ be a lift of absolute Frobenius, and set $X'=X\times_{S, F} S$, so that we have a diagram
$$\xymatrix{
X' \ar[r]^\sigma \ar[d]^{f'} & X \ar[d]^f \\
S \ar[r]^F & S
}$$ lifting the analogous diagram mod $p$ (which we denote by appending zeroes to all the objects above, so for example $X_0'$ is the Frobenius twist of $X_0$). The formalism of crystalline cohomology gives $$\mathscr{H}_{dR}^i(X/S):=R^if_*\Omega^\bullet_{dR, X/S}$$ the structure of an $F$-crystal, i.e.~there is a natural ``crystalline Frobenius" $$\Phi: \mathscr{H}_{dR}^i(X'/S)\simeq F^*\mathscr{H}_{dR}^i(X/S)\to \mathscr{H}_{dR}^i(X/S),$$ horizontal for the Gauss-Manin connection, arising by functoriality from the relative Frobenius morphism $F_{X_0/S}: X_0\to X_0'$. The map $\Phi$ is an isogeny, i.e.~its kernel and cokernel are $p$-power torsion. We also have a natural map $$\sigma^*: \mathscr{H}^i_{dR}(X/S)\to \mathscr{H}^i_{dR}(X'/S)$$ via pullback along $\sigma$ in the Cartesian square above.

We now explain how to recover the Hodge filtration on $\mathscr{H}^i_{dR}(X_0/S_0)$ from $\Phi$. 

Let $$M^k\mathscr{H}^i_{dR}(X'/S)=\Phi^{-1}(p^k\mathscr{H}^i_{dR}(X/S)).$$
\begin{theorem}[Berthelot-Mazur-Ogus, {\cite[8.26]{berthelot2015notes}}, Faltings {\cite[Theorem 6.2]{faltings1989crystalline}}]\label{thm:berthelot-mazur-ogus}
	Suppose that 
	\begin{enumerate}
		\item for all $s\in S_0$, the Hodge-de Rham spectral sequence of $X_{0,s}$ degenerates at $E_1$, and
		\item $\mathscr{H}^i_{dR}(X/S)$ is locally free for all $i$.	
	\end{enumerate}
Then the reduction mod $p$ of $M^k\mathscr{H}^i_{dR}(X'/S)$ is precisely $F^k_{\text{Hodge}}\mathscr{H}^i_{dR}(X'_0/S_0).$ This filtration is flat for the Gauss-Manin connection on $\mathscr{H}^i_{dR}(X'_0/S_0)$ and descends to the Hodge filtration on $\mathscr{H}^i_{dR}(X_0/S_0)$ through the natural base-change isomorphism
$$\sigma_0^*: F_{\text{abs}}^*\mathscr{H}^i_{dR}(X_0/S_0)\overset{\sim}{\to} \mathscr{H}^i_{dR}(X'_0/S_0).$$ That is, $F^k_{\text{Hodge}}\mathscr{H}^i_{dR}(X_0/S_0)$ consists of exactly those $\omega$ such that $\sigma_0^*\omega$ lies in $F^k_{\text{Hodge}}\mathscr{H}^i_{dR}(X'_0/S_0)$.

If in addition $p\geq i+2$, then $$\mathscr{H}^i_{dR}(X/S)=\sum_k \frac{\Phi(\sigma^*F^k_{\text{Hodge}}\mathscr{H}^i_{dR}(X/S))}{p^k}.$$
\end{theorem}  
\begin{remark}
Note that if $S$ is a finite-type smooth $\mathbb{Z}$-scheme and $f: X\to S$ is a smooth projective morphism, the $p$-adic completion of $f$ satisfies the hypotheses of \autoref{thm:berthelot-mazur-ogus} for all sufficiently large primes $p$, see~\cite[Corollaire 4.1.5 and Remarque 4.1.6]{deligne-illusie}.

The last condition says essentially that $\mathscr{H}^i_{dR}(X/S)$ is a Fontaine-Laffaille module (as generalized by Faltings to the relative setting) \cite{faltings1989crystalline}. See also \cite{liu2024relative} for a recent reference on relative Fontaine-Laffaille theory.
\end{remark}
%\begin{corollary}\label{cor:fr-mod-p}
%	Suppose $X/S$ satisfies the hypotheses of \autoref{thm:berthelot-mazur-ogus}, with $S$ irreducible. Suppose further that $F: S\to S$ fixes a $W(k)$-point $s\in S(W(k))$, and that $v\in \Gamma(S, \mathscr{H}_{dR}^{2r}(X/S))$ satisfies $\Phi(\sigma^*v)|_s=p^rv|_s$ and is flat for the Gauss-Manin connection. Then the reduction of $v$ mod $p$ lies in $\on{F}^r_{\text{Hodge}}\mathscr{H}_{dR}^{2r}(X_0/S_0)$.
%\end{corollary}
%\begin{proof}
%It suffices, by \autoref{thm:berthelot-mazur-ogus}, to show that $\sigma_0^*v$ lies in $\on{F}^r_{\text{Hodge}}\mathscr{H}_{dR}^{2r}(X_0'/S_0)$, and hence, by the definition of $M^r\mathscr{H}^i_{dR}(X'/S)$, that $\Phi(\sigma^*(v))$ lies in $p^r\mathscr{H}^i_{dR}(X/S)$. We claim that in fact $\Phi(\sigma^*(v))=p^rv$. This is true by assumption at $s$. By the horizontality of $\Phi$ and $\sigma^*$, $\Phi(\sigma^*(v))$ is flat for the Gauss-Manin connection; thus it is determined by its value at $s$. But $p^rv$ is another flat section to the Gauss-Manin connection with the same value at $s$; thus we must have $\Phi(\sigma^*(v))=p^rv$, as desired.
%\end{proof}
\subsection{Proof of \autoref{thm:cycle-class-initial-conditions}}
We now prepare for the proof of \autoref{thm:cycle-class-initial-conditions}. 
\begin{proposition}[Crystalline step]\label{prop:crystalline-step}
	Notation as in \autoref{thm:cycle-class-initial-conditions}, and assume in addition that $K$ is a number field. Let $V_x\subset H^{2i}_{dR}(X_s)$ be the intersection $$\bigcap_{W\subset Z^i(X_s)_{\mathbb{Q}}\mid x\in \on{cl}(W)_K} \on{cl}(W),$$ i.e.~the intersection  of the images under the cycle class map of all $\mathbb{Q}$-subspaces of $Z^i(X_s)_{\mathbb{Q}}$ whose $K$-span contains $x$. Suppose the formal flat section $v$ to $R^{2i}f_*\Omega^\bullet_{dR, X/S}$ through $x$ is integral (resp.~$\omega(p)$-integral). Then for any $x'\in (V_x)_K$, the formal flat section through $x'$ factors through $$F^i_{\text{Hodge}}R^{2i}f_*\Omega^\bullet_{dR, X_K/S_K}.$$
\end{proposition}
\begin{proof}
	We give the proof assuming $v$ is integral; the proof with $\omega(p)$-integrality is only notationally more complicated, due to the necessity of working over Artin neighborhoods of $s$.

	After perhaps replacing $R$ with a localization and a finite extension, we may assume that $s$ extends to an $R$-point of $S$, $v$ is a formal flat section to $R^{2r}f_*\Omega_{dR, X/S}^\bullet$, and that $K$ is Galois over $\mathbb{Q}$.
	
%	We first argue that $v$ itself factors through $F^i_{\on{Hodge}}R^{2i}f_*\Omega^\bullet_{dR, X_K/S_K}$. 
Fix $e_1, \cdots, e_r\in V_x$ a $\mathbb{Q}$-basis, and write $x=\sum_j a_j e_j$ with $a_j\in K$. (Note that the $e_j$ are necessarily $K$-linearly independent, not just $\mathbb{Q}$-linearly independent, so the $a_j$ are uniquely determined.)
%	
%	First, consider the set of primes $p$ which split completely in $K$ and  such that the $p$-adic completion of $f$ satisfies the hypotheses of \autoref{thm:berthelot-mazur-ogus}; there are infinitely many such by the Chebotarev density theorem. For a place $\lambda$ above such a prime $p$ we obtain by $\lambda$-adic completion an embedding $R\hookrightarrow \mathbb{Z}_p$; let $S_p$ be the (smooth) formal scheme obtained by base-changing to $\mathbb{Z}_p$ and completing $S$ at the $k$-point $\bar s$ corresponding to $s$; let $X_p$ be the base change of $X$ to $S_p$. By the Cohen structure theorem, we may write $S_p\simeq \on{Spf}(\mathbb{Z}_p[[t_1, \cdots, t_n]])$ where $(t_1, \cdots, t_n)$ is the ideal of $s$. Let $F: S_p\to S_p$ be a Frobenius lift such that $F(t_i)=t_i^p$ for all $i$, i.e.~$F(s)=s$.
%	
%	By the integrality of $v$, we are in the setting of \autoref{cor:fr-mod-p}, i.e.~$v$ gives rise to an element of $\Gamma(S, \mathscr{H}^{2i}(X_p/S_p))$ flat for the Gauss-Manin connection. But $v|_s=x$ is a $K$-linear combination of cycle classes, and $p$ splits in $K$, so we have $$\Phi(\sigma^*v)|_s=\Phi(\sigma^*v|_s)=\Phi(\sigma^*(\sum_j a_je_j))p^iv=\sum_j a_j \Phi(\sigma^*(e_j))=\sum_j a_j p^ie_j=p^iv|_s$$. Thus \autoref{cor:fr-mod-p} applies and the reduction of $v$ mod $p$ lies in $F^i_{\text{Hodge}}$.
	
%	But this was the case for infinitely many $p$, so the claim that that $v$ factors through $F^i_{\on{Hodge}}R^{2i}f_*\Omega^\bullet_{dR, X_K/S_K}$ follows by coherence of $F^i_{\text{Hodge}}R^{2i}f_*\Omega^\bullet_{dR, X/S}$.
	
	 We now argue that the formal flat section through any $x'\in (V_x)_K$ factors through $F^i_{\text{Hodge}}R^{2i}f_*\Omega^\bullet_{dR, X_K/S_K}$. By the definition of $V_x$, it suffices to show that for each $\gamma\in \on{Gal}(K/\mathbb{Q})$, the formal flat section $v^\gamma$ through $$x^\gamma:=\sum \gamma(a_j)e_j$$ factors through $F^i_{\text{Hodge}}R^{2i}f_*\Omega^\bullet_{dR, X_K/S_K}$.
	
	By Chebotarev, there exist infinitely many places $\lambda$ of $K$ such that the Frobenius at $\lambda$ is $\gamma$, and such that the completion of $X/S$ at $v$ satisfies the hypotheses of \autoref{thm:berthelot-mazur-ogus}. Let $k$ be the residue field of such a $\lambda$, with residue characteristic $p$, and $W(k)$ the completion of $R$ at $\lambda$. Let $S_{\lambda}$ be the smooth formal $W(k)$ scheme obtained by base-changing $S$ to $W(k)$ and completing at $s$. By the Cohen structure theorem, we may write $S_\lambda\simeq \on{Spf}(W(k)[[t_1, \cdots, t_n]])$ where $(t_1, \cdots, t_n)$ is the ideal of $s$. Let $F: S_\lambda\to S_\lambda$ be a Frobenius lift such that $F(t_i)=t_i^p$ for all $i$, i.e.~$F(s)=s$.  Let $X_\lambda$ be the base-change of $X$ to $S_\lambda$.
	
	We now compute $$\Phi(\sigma^*v^{\gamma^{-1}})|_s=\Phi(\sigma^*x^{\gamma^{-1}})=\Phi(\sum_j \gamma\gamma^{-1}(a_j)\sigma^*e_j)=\sum_j a_j p^ie_j=p^ix=p^iv|_s.$$
	
	Thus $\Phi(\sigma^*v^{\gamma^{-1}})$ and $p^iv$ are both flat sections to the Gauss-Manin connection with the same initial condition at $s$, so they are equal. But $v$ is integral, and hence $\Phi(\sigma^*v^{\gamma^{-1}})$ lies in $p^i\mathscr{H}^{2i}(X_\lambda/S_\lambda)$, whence $v^{\gamma^{-1}}$ lies in $$\sum_k \frac{F^k_{\text{Hodge}}\mathscr{H}^{2i} _{dR}(X_\lambda/S_\lambda)}{p^{k-i}}$$ by \autoref{thm:berthelot-mazur-ogus}. Thus we see that modulo $F^i_{\text{Hodge}}$, $v^{\gamma^{-1}}$ lies in $p\mathscr{H}^{2i}_{dR}(X_\lambda/S_\lambda).$ As this holds true for infinitely many $\lambda$, we must have that $v^{\gamma^{-1}}$ lies in $F^i_{\text{Hodge}}$, as desired.
	
	As $\gamma$ was arbitrary, this completes the proof.
\end{proof}
\begin{proposition}[Analytic step]\label{prop:analytic-step}
	Let $f: X\to S$ be a smooth projective morphism of smooth varieties over $\mathbb{C}$. Fix $s\in S$ and $x_{\text{Betti}}\in H^{2i}(X_s(\mathbb{C})^{\text{an}}, \mathbb{Z})$. Let $x_{dR}$ be the image of $x_{\text{Betti}}$ under the Betti-de Rham comparison morphism $$H^{2i}(X_s(\mathbb{C})^{\text{an}}, \mathbb{Z}){\to} H^{2i}_{dR}(X_s),$$ and suppose that the formal flat section to $R^if_*\Omega^\bullet_{X/S}$ with respect to  the Gauss-Manin connection  through  $x_{dR}$ factors through $F^i_{\text{Hodge}}$. Then this flat section is in fact algebraic.
\end{proposition}
\begin{proof}
We first claim that the orbit $\pi_1(S,s)\cdot x_{\text{Betti}}$ inside $H^{2i}(X_s(\mathbb{C})^{\text{an}}, \mathbb{C})$ is finite. As $x_{\text{Betti}}$ is a $\mathbb{Z}$-class, this orbit is is discrete in $H^{2i}(X_s(\mathbb{C})^{\text{an}}, \mathbb{C})$, so it suffices to show the orbit is compact. We will do so using the existence of a polarization on $H^{2i}(X_s(\mathbb{C})^{\text{an}}, \mathbb{C})$.

As $x_{\text{Betti}}$ is a $\mathbb{Z}$-class, it is its own complex conjugate; hence $v$ factors through $\overline{F^i_{\text{Hodge}}}$ as well as $F^i_{\text{Hodge}}$. Thus the orbit $\pi_1(S, s)\cdot x_{\text{Betti}}$ lies in $H^{i,i}(X_s)\subset H^{2i}(X_s, \mathbb{C})$. This orbit is compact (as the polarization on $H^{i,i}(X_s)$ is preserved by $\pi_1(S,s)$ and definite on each piece of the Lefschetz decomposition), hence is, by the paragraph above, finite.

Now let $\widetilde{S}$ be the finite cover of $S$ corresponding to the stabilizer of $x_{\text{Betti}}$ in $\pi_1(S,s)$. The pullback of $v$ to $\widetilde{S}$ extends to an (analytic) global section on $\widetilde{S}$, flat for the Gauss-Manin connection, by  construction. But it is in fact algebraic as the Gauss-Manin connection has regular singularities at infinity.
\end{proof}
\begin{proof}[Proof of \autoref{thm:cycle-class-initial-conditions}]
Again we only give the proof when $v$ is integral, as the case where it is $\omega(p)$-integral is only notationally more complicated. 

	The theorem is more or less immediate from \autoref{prop:crystalline-step} and \autoref{prop:analytic-step}. We first prove it in the case that $K$ is a number field. Defining $V_x$ as in \autoref{prop:crystalline-step}, choose a $\mathbb{Q}$-basis $e_1, \cdots, e_r$ of $V_x$ consisting of honest $\mathbb{Z}$-linear combinations of cycle classes. By \autoref{prop:crystalline-step}, the formal flat section to the Gauss-Manin connection through each $e_i$ factors through $F^i_{\text{Hodge}}$; thus by \autoref{prop:analytic-step}, each such formal flat section is algebraic. As $v$ is a $K$-linear combination of these, $v$ is algebraic as well.
	
	We now consider the case where $R\subset \mathbb{C}$ is an arbitrary finitely-generated $\mathbb{Z}$-algebra. By extending $R$ we may assume $s$ extends to an $R$-point of $S$ and $v$ is a formal flat section to $R^{2r}f_*\Omega_{dR, X/S}^\bullet$. Let $B=\on{Spec}(R)$. For each $L$-point $b$ of B, with $L$ a number field, we have that $v_b$ is algebraic by the previous paragraph.
	
	Choose an embedding $\iota: R\hookrightarrow \mathbb{C}$. It suffices (again by the fact that the Gauss-Manin connection has regular singularities at infinity) to show that the $\pi_1(S_{\mathbb{C}}^{\text{an}}, s)$-orbit of $x$ is finite. But this is true for any specialization $\iota': R\to \mathbb{C}$ of $\iota$ whose image is a number field; the corresponding monodromy representations are the same (as $f$ is defined over $R$), and so the proof is complete.
\end{proof}
\begin{remark}
If $x$ is in the image of $Z^i(X_s)_\mathbb{Q}$, as opposed to $Z^i(X_s)_K$, the proof becomes much simpler. Specializing the proof to this case, one does not use that $x$ is a cycle class; rather, one uses only that:
\begin{enumerate}
	\item $x$ is a Hodge cycle, i.e.~it lies in the intersection $H^{2i}(X, \mathbb{Q}(i))\cap F^i_{\text{Hodge}}H^{2i}_{dR}(X)$, and
	\item $x$ is an absolute Tate cycle, i.e.~viewed as an element of crystalline cohomology, the crystalline Frobenius at $p$ multiplies it by $p^i$ for almost all $p$.  (See \cite[VI, 4.1.3]{hodge-cycles-motives-etc} for a discussion of this notion.)
\end{enumerate}
Of course one expects these conditions to be equivalent to one another, and to be, in turn, equivalent to being the class of an algebraic cycle. In this case one does not need the full force of (relative) Fontaine-Laffaille theory. Instead, one may replace \autoref{prop:crystalline-step} by the observation that $\Phi(\sigma^*(v))=p^iv$, whence by \autoref{thm:berthelot-mazur-ogus}, $v$ lies in $F^i_{\text{Hodge}}$ mod $p$. As this holds for infinitely many $p$, $v$ is in fact in $F^i_{\text{Hodge}}$.

In this case one only needs the integrality conditions of the statement to hold at infinitely many primes; by contrast in our proof above one needs integrality at a density one set of primes due to the use of the Chebotarev density theorem.
\end{remark}

\begin{remark}[Relationship to Katz's proof]
	 \autoref{thm:cycle-class-initial-conditions} is related, but somewhat orthogonal to, Katz's proof of the $p$-curvature conjecture for Picard-Fuchs equations \cite{katz-p-curvature}, and in particular we do not know a proof of \autoref{thm:cycle-class-initial-conditions} which avoids the machinery of crystalline cohomology. The analytic step of Katz's proof is more or less identical to our \autoref{prop:analytic-step}. Arguably there is some relationship between the algebraic steps as well; the crystalline Frobenius operator we use ``lifts" the inverse Cartier operator used by Katz; see \cite[Theorem 8.26]{berthelot2015notes} for details. We will review some aspects of Katz's argument in the next section as we adapt them for our purposes.
\end{remark}

\begin{remark}
We expect that proving \autoref{conjecture:general-algebraic-leaves} for Picard-Fuchs equations and arbitrary initial conditions will require new ideas. Indeed, doing so would resolve the Grothendieck-Katz $p$-curvature conjecture for summands of Picard-Fuchs equations, which is open; see \cite[\S16]{andre2004conjecture} for some results towards this case of the $p$-curvature.
\end{remark}

\section{Families of elliptic curves; modular forms}\label{section:elliptic-curves}
We now prove our main conjecture, \autoref{conjecture:general-algebraic-leaves}, for differential equations arising from families of elliptic curves; for example, it is true for Taylor series expansions of solutions to the differential equation \begin{equation}\label{eqn:beukers} t(t^2-11t-1)\frac{\partial^2{f}}{\partial t^2}+(3t^2-22t-1)\frac{\partial f}{\partial t}+(t-3)f=0\end{equation} away from its singular points (namely $0$ and the roots of $t^2-11t-1=0$). Our primary reason for including this example is that such differential equations often admit a (non-algebraic) solution holomorphic at a \emph{singular} point of the ODE, whose Taylor expansion is integral. For example, \eqref{eqn:beukers} has a solution whose Taylor expansion is given by $$\sum u_nt^n,\ u_n=(-1)^n\sum_k {n \choose k}^2{n+k\choose k}.$$ See \cite[\S10]{beukers-stienstra} for further discussion and examples. Of course this does not contradict \autoref{conjecture:general-algebraic-leaves}, which only allows Taylor expansions at non-singular points.

Explicitly writing down the recurrence relation for the coefficients of the Taylor coefficients for a solution to \eqref{eqn:beukers} at $t=a$, we find:
\begin{proposition}[Corollary to \autoref{theorem:elliptic-curves} below] \label{prop:ell-curve-corollary} Fix  $a, b_0, b_1\in \mathbb{C}$ with $a(a^2-11a-1)\neq 0$. Let $b_i, i>1$ be the solution to the recurrence relation

 \begin{multline} \label{eqn:apery-recurrence} a(a^2-11a-1)(n+2)(n+1) b_{n+2} + (3a^2-22a-1)(n+1)^2 b_{n+1}  +  \\ [(3a-11)(n+1)n +a-3] b_n + n^2 b_{n-1} = 0.	
	\end{multline}

	Then the ring $\mathbb{Z}[b_0, b_1, b_2, \cdots] \subset \mathbb{C}$ generated by the $b_i$ is not a finitely-generated $\mathbb{Z}$-algebra.
\end{proposition}
For example, if $a, b_0, b_1\in \mathbb{Q}$, this says that there are infinitely many primes dividing the denominators of the $b_i$. 

Note that while \eqref{eqn:apery-recurrence} is an order $3$ linear recurrence, the value of $b_2, b_3, \cdots$ is determined uniquely by the value of $b_0, b_1$.
\subsection{Statement of the theorem for families of elliptic curves}
\begin{theorem}\label{theorem:elliptic-curves}
Let $R\subset \mathbb{C}$ be a finitely-generated $\mathbb{Z}$-algebra with fraction field $K$. Let $f: \mathscr{X}\to S/R$ be a non-isotrivial relative genus one curve (i.e.~a smooth projective morphism of relative dimension $1$ with geometric fibers curves of genus one, whose Jacobian has non-constant $j$-invariant), with $S$ smooth and connected. Let $(\mathscr{E},\nabla)=(R^1f_*\Omega^\bullet_{dR, \mathscr{X}/S},\nabla_{GM})$ be the relative de Rham cohomology of $f$, equipped with its Gauss-Manin connection.

Fix $n>0$. Let $s\in S(R)$ be a point and $v\in \on{Sym}^n H^1_{dR}(\mathscr{X}_s)$ be a non-zero element. Then the formal flat section to $\on{Sym}^n(\mathscr{E},\nabla)_{K}$ through $v_K$ is not integral (resp.~not $\omega(p)$-integral).
\end{theorem}
For example, this result applies to the solutions to Picard-Fuchs equations associated to modular forms, see e.g.~\cite[\S5.4]{bruinier2008elliptic}. 
\begin{remark}\label{remark:taylor-expansion-hypergeometric}
The function $F(t)$ of \autoref{remark:hypergeometric-globally-bounded} is an example of a solution to a Picard-Fuchs equation associated to a modular form; it is the Taylor expansion of the square of a certain theta function in the modular function $\lambda(z)$; see \cite[equation (73)]{bruinier2008elliptic} (though the expression there has a small typo; the binomial coefficient should be squared). Thus \autoref{theorem:elliptic-curves} implies that the coefficients of the Taylor expansion of $F(t)$ at any point $t_0\not\in\{0,1\}$ do not lie in a finitely-generated $\mathbb{Z}$-algebra.
\end{remark}

Note that $f$ is not isotrivial by assumption, and so $(\mathscr{E}, \nabla)$ has geometric monodromy group $\on{SL}_2$, and hence $\on{Sym}^n(\mathscr{E},\nabla)$ does not admit any nonzero algebraic flat sections. So to prove \autoref{conjecture:general-algebraic-leaves} in this case it is enough to show that no nonzero formal flat section is integral (resp.~$\omega(p)$-integral).
As in \autoref{prop:ell-curve-corollary} this is, when $S$ is an open subset of $\mathbb{A}^1$, really a very concrete statement about sequences $\{b_n\}$ defined by certain concrete recurrence relations, saying that infinitely many prime factors appear in the denominators of the $b_i$.
\subsection{Preliminaries on the conjugate filtration}\label{subsection:conjugate-preliminaries}
We briefly recall some of the basics on the conjugate filtration, $p$-curvature, and the Gauss-Manin connection that we will use in the proof.

Let $f: X\to S$ be a smooth projective morphism of varieties over a field $k$ of positive characteristic. Let $\tau_{\leq p}\Omega^\bullet_{dR, X/S}$ be the subcomplex which consists of $\Omega^i_{dR, X/S}$ in degree $i$ for $i<p$, $\ker d$ in degree $i=p$, and $0$ in degree $i>p$. This filtration induces the ``second hypercohomology spectral sequence" $$E_2^{ij}=R^if_*\mathscr{H}^j(\Omega^\bullet_{dR, X/S})\implies R^{i+j}f_*\Omega^\bullet_{dR, X/S}.$$ The induced (decreasing) filtration on  $R^{i+j}f_*\Omega^\bullet_{dR, X/S}$ is called the conjugate filtration, i.e.~$$F^r_{\text{conj}}R^{s+r}f_*\Omega^\bullet_{dR,X/S}=\on{Im}(R^{s+r}f_*\tau_{\leq s}\Omega^\bullet_{dR,X/S}\to R^{s+r}f_*\Omega^\bullet_{dR,X/S}).$$
In particular the graded pieces of the conjugate filtration on $R^sf_*\Omega^\bullet_{dR, X/S}$ are supported in degrees $[0, s]$. This is nicely explained in e.g.~\cite[\S2.3]{katz-p-curvature}.

 If the conjugate spectral sequence degenerates at $E_2$ (which will often be the case for us---see e.g.~ \autoref{theorem:katz-big-diagram} below), the graded pieces of the conjugate filtration may be identified explicitly via Cartier theory; we have the Cartier isomorphism $$C^{-1}: R^af'_*\Omega^b_{X'/S} \to \on{gr}^{a}_{F_{\on{conj}}}R^{a+b}f_*\Omega^\bullet_{dR, X/S}$$ where $f': X'\to S$ is the base change of $f$ along the absolute Frobenius map $F_{\text{abs}}: S\to S$.

We now recall the main computation of \cite{katz-p-curvature}, variants of which will be used many times in the rest of this paper.
\begin{theorem}[{\cite[Theorem 3.2]{katz-p-curvature}}]\label{theorem:katz-big-diagram} 
Let $f: X\to S$ be a smooth projective morphism of $k$-schemes, where the characteristic of $k$ is $p>0$. Let $D$ be a relative simple normal crossings divisor on $X/S$. Suppose the Hodge cohomology sheaves $R^af_*\Omega^b_{X/S}(\log D)$ are locally free and that their formation commutes with base change, and that the Hodge-de Rham spectral sequence degenerates at $E_1$. Then the conjugate spectral sequence for $X/S$ degenerates at $E_2$, and the following diagram commutes:
$$\xymatrix{
F_{\text{abs}}^*\on{gr}^{b}_{F_{\on{Hodge}}}R^{a+b}f_*\Omega^\bullet_{dR, X/S} (\log D)\ar[rr]^-{F_{\text{abs}}^*\on{gr}_{F_{\text{Hodge}}}\nabla_{GM}} \ar[d]^\sim && F_{\text{abs}}^*\on{gr}^{b-1}_{F_{\on{Hodge}}}R^{a+b}f_*\Omega^\bullet_{dR, X/S}(\log D)\otimes F_{\text{abs}}^*\Omega^1_{S/k} \ar[d]^\sim \\
F_{\text{abs}}^*R^af_*\Omega^b_{X/S}(\log D) \ar[rr] \ar[d]^\sim &&  F_{\text{abs}}^*R^{a+1}f_*\Omega^{b-1}_{X/S}(\log D) \otimes F_{\text{abs}}^*\Omega^1_{S/k} \ar[d]^\sim\\
R^af'_*\Omega^b_{X'/S}(\log D') \ar[rr] \ar[d]^{C^{-1}} && R^{a+1}f'_*\Omega^{b-1}_{X'/S}(\log D')\otimes F_{\text{abs}}^*\Omega^1_{S/k} \ar[d]^{C^{-1}}\\
\on{gr}^{a}_{F_{\on{conj}}}R^{a+b}f_*\Omega^\bullet_{dR, X/S}(\log D) \ar[rr]^-{-\on{gr}_{F_{\on{conj}}}\psi_p} && \on{gr}^{a+1}_{F_{\on{conj}}}R^{a+b}f_*\Omega^\bullet_{dR, X/S}(\log D) \otimes  F_{\text{abs}}^*\Omega^1_{S/k}
}$$
Here the top two vertical arrows record the degeneration of the Hodge-de Rham spectral sequence, the middle two are base change isomorphisms, and $C^{-1}$ is the inverse Cartier operator.  The map $\on{gr}_{F_{\text{Hodge}}}\nabla_{GM}$ is the associated graded of the Gauss-Manin connection with respect to the Hodge filtration, which has degree $-1$ by Griffiths transversality. The map $\on{gr}_{F_{\on{conj}}}\psi_p$ is the associated graded of the $p$-curvature with respect to the conjugate filtration, which has degree $-1$ as $\on{gr}_{F_{\on{conj}}}  R^if_*\Omega^\bullet_{dR, X/S}$ has zero $p$-curvature.
\end{theorem}

The following is a special case of a question of Ogus \cite[Problem 2.1]{hodge-cycles-motives-etc}:
\begin{lemma}\label{lemma:bost-ogus-conjugate}
Let $R\subset \mathbb{C}$ be a finitely-generated integral $\mathbb{Z}$-algebra, and let $X/R$ be a smooth projective variety. Let $v\in H^1_{dR}(X/R)$ be any non-zero element. Then there exists a Zariski-dense set of closed points $\{\mathfrak{p}_i\}\subset \on{Spec}(R)$ such that the reduction of $v\bmod \mathfrak{p}_i$ does not lie in $$F^1_{\text{conj}}H^1_{dR}(X_{\mathfrak{p}_i}/k(\mathfrak{p}_i)).$$
\end{lemma}
\begin{proof}
In \cite[VI, Remark 2.5]{hodge-cycles-motives-etc}, it is explained how to reduce the desired statement to the Grothendieck-Katz $p$-curvature conjecture for flat bundles which are extensions of $(\mathscr{O}_X, d)$ by itself.  This latter statement has been proven by Andr\'e \cite[\S6]{andre2004conjecture}, Bost \cite[\S2.4.2]{bost2001algebraic}, and Chudnovsky-Chudnovsky \cite{chudnovsky2006applications}. For the reader's convenience, we include a self-contained sketch of the reduction here.

Without loss of generality, by replacing $R$ with a localization we may assume that for $s$ a closed point of $\on{Spec}(R)$, the local freeness, base change, and degeneration hypotheses of \autoref{theorem:katz-big-diagram} are satisfied for $X_s/s$, as they are satisfied in characteristic zero. Thus we may assume that the conjugate spectral sequence for $X_s/s$ degenerates at $E_2$ for all $s$.

Let $v\in H^1_{dR}(X/R)$ be a class, corresponding to an extension $$0\to (\mathscr{O}_X, d)\to (\mathscr{E}_v, \nabla)\to (\mathscr{O}_X, d)\to 0.$$ That the reduction of $v\bmod \mathfrak{p}_i$ lies in $F^1_{\text{conj}}H^1_{dR}(X_{\mathfrak{p}_i}/k(\mathfrak{p}_i))$ is equivalent to the statement that the image of $v\bmod \mathfrak{p}_i$ in $H^0(\mathscr{H}^1_{dR}(X_{\mathfrak{p}_i}/k(\mathfrak{p}_i)))$ vanishes. In other words, if $U\subset X_{\mathfrak{p}_i}$ is an affine open, so that $v$ is represented by a closed $1$-form $\omega$ on $U$, then $\omega$ is locally exact. Thus $(\mathscr{E}_v, \nabla)$ is (mod $\mathfrak{p}_i$) locally trivial as a bundle with flat connection. If this occurs for all $\mathfrak{p}_i$ in the complement of some closed subset of $\on{Spec}(R)$, we have by Andr\'e, Bost, Chudnovsky-Chudnovsky that $(\mathscr{E}_v, \nabla)$ is trivial, i.e.~$v=0$ as desired.
\end{proof}

\begin{lemma}\label{lemma:veronese-conjugate}
	Fix $n>0$ Let $R\subset \mathbb{C}$ be a finitely-generated integral $\mathbb{Z}$-algebra in which $n$ is invertible, and let $X/R$ be an elliptic curve. Let $v\in \on{Sym}^n H^1_{dR}(X/R)$ be any non-zero element. Then there exists a Zariski-dense set of closed points $\{\mathfrak{p}_i\}\subset \on{Spec}(R)$ such that the reduction of $v\bmod \mathfrak{p}_i$ does not lie in $$F^n_{\text{conj}} \on{Sym}^n H^1_{dR}(X_{\mathfrak{p}_i}/k(\mathfrak{p}_i)).$$ Here $F^\bullet_{\text{conj}}$ is the filtration on $\on{Sym}^n H^1_{dR}(X_{\mathfrak{p}_i}/k(\mathfrak{p}_i))$ induced by the conjugate filtration on $H^1_{dR}(X_{\mathfrak{p}_i}/k(\mathfrak{p}_i))$.
\end{lemma}
\begin{proof}
For $n=1$, this is immediate from \autoref{lemma:bost-ogus-conjugate}. For $n>1$, we may without loss of generality assume that there exists a Zariski-dense set of closed points $\{\mathfrak{q}_j\}\subset \on{Spec}(R)$ such that the reduction of $v\bmod \mathfrak{q}_j$ \emph{lies} in $$F^n_{\text{conj}} \on{Sym}^n H^1_{dR}(X_{\mathfrak{q}_j}/k(\mathfrak{q}_j)),$$ as otherwise we're done. 

For such $\mathfrak{q}_j$, $v\bmod \mathfrak{q}_j$ lies in the image of the Veronese embedding $$\mathbb{P}H^1_{dR}(X_{\mathfrak{q}_j}/k(\mathfrak{q}_j))\hookrightarrow \mathbb{P}\on{Sym}^n H^1_{dR}(X_{\mathfrak{q}_j}/k(\mathfrak{q}_j)),$$ since $$F^n_{\text{conj}} \on{Sym}^n H^1_{dR}(X_{\mathfrak{q}_j}/k(\mathfrak{q}_j))=(F^1_{\text{conj}}H^1_{dR}(X_{\mathfrak{q}_j}/k(\mathfrak{q}_j)))^{\otimes n}.$$ As the $\mathfrak{q}_j$ are Zariski-dense in $\on{Spec}(R)$, the same is true (after possibly replacing $R$ with a finitely generated extension) over $R$, i.e.~$v$ is in the span of $w^{\otimes n}$ for some $w\in H^1_{dR}(X_{\mathfrak{q}_j}/k(\mathfrak{q}_j))$. But now applying \autoref{lemma:bost-ogus-conjugate}, we find some Zariski-dense set of closed points $\{\mathfrak{p}_i\}\subset \on{Spec}(R)$ such that the reduction of $w\bmod \mathfrak{p}_i$ does not lie in $F^1_{\text{conj}}H^1_{dR}(X_{\mathfrak{p}_i}/k(\mathfrak{p}_i))$, which suffices.
\end{proof}

\subsection{The proof of \autoref{theorem:elliptic-curves}}
We may now prove \autoref{theorem:elliptic-curves}.
\begin{proof}[Proof of \autoref{theorem:elliptic-curves}]
	Let $\widetilde v_K$ be the formal flat section to $\on{Sym}^n(\mathscr{E},\nabla)_K$ through $v_K$. It suffices to show that $\widetilde v_K$ is not $\omega(p)$-integral. Suppose to the contrary that it is; then we will show that $\mathscr{X}/S$ is isotrivial, a contradiction. For a closed point $\mathfrak{p}$ of $R$, let $\widetilde v_\mathfrak{p}$ be the reduction of $\widetilde{v}_K$ modulo $\mathscr{I}_s^{\omega(p)}, \mathfrak{p}$. (See \autoref{subsection:integrality-of-leaves} for this notation.)
	
	As $v$ is non-zero, \autoref{lemma:veronese-conjugate} gives us a Zariski-dense set of closed points $\{\mathfrak{p}_i\}_{i\in I}\subset \on{Spec}(R)$ such that the reduction of $v$ mod $\mathfrak{p}_i$ is not contained in $F^n_{\text{conj}}\on{Sym}^nH^1_{dR}(X_{\mathfrak{p}_i}/k(\mathfrak{p}_i))$, for all $i\in I$. Thus there exists some $r<n$ and some Zariski-dense set of closed points $\{\mathfrak{p}'_i\}_{i\in I}\subset \on{Spec}(R)$ such that the reduction of $v$ mod $\mathfrak{p}'_i$ is contained in $F^r_{\text{conj}}$ but not in $F^{r+1}_{\text{conj}}$. Hence for $\mathfrak{p}_i'$ of residue characteristic $p$, $\on{gr}^r_{F_{\text{conj}}}\psi_p$ vanishes to order $\omega(p)$ on the image of $\widetilde{v}_{\mathfrak{p}'_i}$ in $\on{gr}^r_{F_{\text{conj}}}\on{Sym}^n(\mathscr{E}, \nabla)|_{V(\mathscr{I}_s^{\omega(p)}, \mathfrak{p}_i')}$ (for example, by \autoref{prop:p-power-leaves}).  As, for $\mathfrak{p}'_i$ of sufficiently large characteristic $p$, this class generates $\on{gr}^r_{F_{\text{conj}}}\on{Sym}^n(\mathscr{E}, \nabla)|_{V(\mathscr{I}_s^{\omega(p)}, \mathfrak{p}_i')}$, we thus have that $\on{gr}_{F_{\text{conj}}}\psi_p$ is zero when restricted to $V(\mathscr{I}_s^{\omega(p)}+ \mathfrak{p}'_i)$, i.e.~$\on{gr}^r_{F_{\text{conj}}}\psi_p$ vanishes to order $\omega(p)$ at $s \bmod \mathfrak{p}'_i$. 
	
	It follows from \autoref{theorem:katz-big-diagram} that $\on{gr}^{n-r}_{F_{\text{Hodge}}}\nabla_{GM}$ vanishes to order $\omega(p)/p$ at $s \bmod \mathfrak{p}'_i$. As $\omega(p)/p\to \infty$ and the $\mathfrak{p}'_i$ are Zariski-dense, this implies that $\on{gr}^{n-r}_{F_{\text{Hodge}}}\nabla_{GM}$ is identically zero, i.e.~$\nabla$ preserves $F^{n-r}_{\text{Hodge}}\on{Sym}^n(\mathscr{E},\nabla)$, and so the monodromy representation associated to $\on{Sym}^n(\mathscr{E},\nabla)$ is reducible. This contradicts the assumption that $f$ was non-isotrivial. Indeed, for any non-isotrivial family of elliptic curves the geometric monodromy group of $(\mathscr{E}, \nabla)$ is $SL_2$; hence $\on{Sym}^n(\mathscr{E},\nabla)$ is irreducible. 
\end{proof}
\section{Hypergeometric functions}\label{section:hypergeometric}
\subsection{The local conjecture}
In this section we will state an  elementary corollary of \autoref{conjecture:general-algebraic-leaves} about solutions of algebraic  linear differential equations, and verify it in the case of hypergeometric functions.

\begin{definition}
For $p_1(z), \cdots, p_n(z)\in \mb{C}(z)$, we say that 
    \[
    y^{(n)}+p_1(z)y^{(n-1)}+\cdots +p_n(z)y=0
    \]
    has a regular singularity at $0$ if for  $i=1, \cdots, n$, $p_i$ has a pole of order at most $i$. We say that the differential equation is Fuchsian if it has  regular singularities at all $t\in \mb{P}^1(\mb{C})$: i.e. when written   in  a local parameter at $t$, the differential equation  has a regular singularity at $0$.
\end{definition}
The following is a translation of \autoref{conj:singular-points} to this setting:
\begin{conjecture}\label{conjecture:local-conjecture-fuchsian}
 Suppose the power series $f(z)\in \Qbar[[z]]$ is a solution to the Fuchsian differential equation 
\begin{equation}\label{eqn: fuchsian-de}
   y^{(n)}+p_1(z)y^{(n-1)}+\cdots +p_n(z)y=0,
\end{equation}
with $p_i(z)\in \Qbar(z)$ for $i=1, \cdots, n$.
If $f(z)\in \overline{\mathbb{Z}}[\frac{1}{N}][[z]]$ for some integer $N\geq 1$, then either $f(z)$ is algebraic, or the monodromy of \eqref{eqn: fuchsian-de} has infinite order  about $0$.
\end{conjecture}

\subsection{Proof for hypergeometric functions}
While we do not know how to say much about \autoref{conjecture:local-conjecture-fuchsian} in general, we check it for Taylor expansions of hypergeometric functions about zero. Our conjecture seems to be closely related to \cite[Conjecture 4]{gilles2006globally}, though that conjecture remains open for hypergeometric functions.
\begin{definition}
\begin{enumerate}
    \item The hypergeometric functions are defined explicitly as  power series
\begin{equation}\label{eqn:explicithyp}
 _{p+1}F_{p}[a_1, \cdots , a_{p+1}; b_1, \cdots , b_p; z] \defeq \sum_{n=0}^\infty \frac{(a_1)_n\cdots (a_{p+1})_n}{(1)_n(b_1)_n\cdots (b_p)_n} z^n,
\end{equation}
for parameters $a_i \in \mb{Q}$ and $ b_i \in \mb{Q} \setminus \mb{Z}_{\leq 0}$, where $(x)_n$ is the Pochhammer symbol, or rising factorial, given by $(x)_0=1$, and $(x)_n= x(x+1)\cdots (x+n-1)$ for $n\geq 1$. When the context is clear, we sometimes simplify the notation to $_{p+1}F_{p}[a_i; b_i; z]$.
\item 
For $a=(a_1, \cdots, a_{n}), b=(b_1, \cdots, b_n)$ tuples of complex numbers, we define also the hypergeometric differential operator
\[
    \mathrm{Hyp}(a;b):=\prod_{i=1}^{n} (z\partial_z+b_i-1)-  z \prod_{j=1}^{n} (z\partial_z+a_j). 
\]
Note that this operator defines a Fuchsian differential equation of order $n$, with singularities at $0,1,\infty$.  Moreover,  $_{p+1}F_{p}[a_1, \cdots , a_{p+1}; b_1, \cdots , b_p; z]$ is annihilated by $\mathrm{Hyp}(a_1, \cdots, a_{p+1}; b_1, \cdots , b_p, 1)$; see \cite[\S2]{beukersheckman}. 
\end{enumerate}
  
\end{definition}
 \autoref{prop: hyp-result} below  verifies \autoref{conjecture:local-conjecture-fuchsian} for hypergeometric functions:
\begin{proposition}\label{prop: hyp-result}
    For $a_1, \cdots, a_{p+1}\in \mb{Q}$ and $b_1, \cdots, b_{p}\in \mb{Q}\setminus \mb{Z}_{\leq 0}$, write  $F= _{p+1}F_{p}[a_i; b_i; z]$. Suppose that
    \begin{itemize}
        \item $F$ has coefficients in $\mathbb{Z}[1/M]$ for some positive integer $M$, and 
        \item $F$ is not algebraic.
    \end{itemize}
    Then the monodromy of $\mathrm{Hyp}(a; b)$ has infinite order around $0$.
\end{proposition}
The proof is essentially a combination of some of the main results of \cite{beukersheckman} and \cite{christol1986fonctions}, which we now recall.

\begin{proposition}[{\cite[Prop.~3.2, Thm. 3.5]{beukersheckman}}]\label{prop:hyp local mono}
    We draw based loops around $0,1,\infty$ as in \cite[p. 328]{beukersheckman}, and write $g_{0}, g_1, g_{\infty}$ for the monodromy matrices for $\on{Hyp}(a;b)$ around these loops. Write $\alpha_i=\exp(2\pi\sqrt{-1}a_i)$ and $\beta_i=\exp(2\pi\sqrt{-1}b_i)$. Then the elements  $g_0, g_1, g_{\infty}$ satisfy
    \begin{align*}
        det(t-g_{\infty}) &= \prod_{i=1}^{p+1} (t-\alpha_i),\\
         det(t-g_{0}^{-1}) &= \prod_{i=1}^{p+1} (t-\beta_i),\\
         \mathrm{rk}(g_1-\mathrm{Id})&=1, \ \det(g_1)=\exp(2\pi \sqrt{-1}\sum_{i=1}^{p+1} (b_i-a_i)).
    \end{align*}
Moreover, for each generalized eigenvalue of $g_{\infty}$, there is a unique Jordan block with this generalized eigenvalue, and the same is true of $g_0$.       
\end{proposition}

For any  $r\in \mb{R}$, write $\{r\}$ for the unique element of $(0,1]$ such that $r-\{r\} \in \mb{Z}$. 
\begin{lemma}\label{lemma: hyp-criterion}
    For $a_1,\cdots , a_{p+1}\in \mb{Q}, b_1, \cdots, b_p \in \mb{Q}\setminus \mb{Z}_{\leq 0}$, write  $F= _{p+1}F_{p}[a_1, \cdots , a_{p+1}; b_1, \cdots , b_p; z]$  and let $N$ be the least common multiple of the denominators of the $a_i, b_i$'s. Set $b_{p+1}:=1$. 
     For any integer $\Delta$ with $(\Delta, N)=1$, set 
    \[
    a_{\Delta, i} = \{\Delta a_i\}, b_{\Delta, i} = \{\Delta b_i\},
    \]
    
    \begin{enumerate}
       \item\label{criterion-christol} \cite[\S IV]{christol1986fonctions} $F$ has coefficients in $\mathbb{Z}[\frac{1}{M}]$ for some positive integer $M$ if and only if the following holds:
   
    for any $\Delta$ with $(\Delta, N)=1$, and   any $\lambda \in (0,1]$, 
    \[
    \#\{a_{\Delta, i}| a_{\Delta, i}\leq \lambda\} \geq \# \{b_{\Delta, i}| b_{\Delta, i}\leq \lambda\}.
    \]
    Colloquially, this condition says that, as we go from $0$ to $1$, the number of $a_{\Delta, i}$'s we cross is always at least that of $b_{\Delta, i}$'s. 
    \item\label{criterion-bh} \cite[Theorem 4.8]{beukersheckman} $F$ is algebraic if and only if the following holds: for any $\Delta$ with $(\Delta, N)=1$, the $a_{\Delta, i}$'s and $b_{\Delta, i}$'s interlace. More precisely, we may reindex the   $a_{\Delta, i}$'s and $b_{\Delta, i}$'s such that 
    \[
    0< a_{\Delta, 1}<b_{\Delta, 1}< a_{\Delta, 2}< b_{\Delta, 2} < \cdots < a_{\Delta, p+1}< b_{\Delta, p+1}=1.
    \]
    \end{enumerate} 
\end{lemma}

\begin{proof}[Proof of \autoref{prop: hyp-result}]
    We assume that $F$ has coefficients in $\mathbb{Z}[\frac{1}{M}]$ for some $M>0$ and that $F$ is not algebraic. By the last part of  \autoref{prop:hyp local mono}, it suffices to show that there is $i\neq p+1$ such that $\beta_i=1$, or in other words, $b_i\in \mb{Z}$. Let us suppose that there is no such $i$ and aim to derive a contradiction; note that this  implies that there is no $i<p+1$ such that $b_{\Delta, i}=1$.

Let us reorder the $a_{1, i}$'s and $b_{1, i}$'s so that
\[
a_{1, 1}\leq a_{1, 2} \leq \cdots \leq a_{1, p+1}, \ \ b_{1, 1}\leq b_{1, 2} \leq \cdots \leq b_{1, p+1}.
\]
By our assumption,  $b_{1, p}\neq 1$. Taking $\Delta=-1$ in \autoref{lemma: hyp-criterion} (\ref{criterion-christol}) we deduce that $b_{1, p}< a_{1, p+1}$. On the other hand, taking $\Delta=1$ in \autoref{lemma: hyp-criterion} (\ref{criterion-christol}), we deduce that $a_{1, p}<b_{1, p}$; continuing in this way, we conclude that the $a_{1,i}$'s and $b_{1,i}$'s interlace. 

The same argument shows that the $a_{\Delta, i}$'s and $a_{\Delta, i}$'s interlace for any $(\Delta, N)=1$. Therefore $F$ is algebraic by \autoref{lemma: hyp-criterion} (\ref{criterion-bh}), which is a contradiction. Therefore $b_{1,p}=1$, as required.
\end{proof}
\begin{remark}
Landau \cite{landau1904anwendung, landau1911zahlentheoretischen} reproduced Schwarz's list of algebraic hypergeometric functions ${}_2F_1[a_1, a_2; b; z]$ using Eisenstein's theorem on integrality of Taylor expansions of algebraic functions \cite{eisenstein1852allgemeine}; Katz rediscovered this approach from a more conceptual point of view in \cite[\S6]{katz-p-curvature}. The main result of this section, \autoref{prop: hyp-result}, amounts to the observation that, a posteriori, Eisenstein's theorem \cite{eisenstein1852allgemeine} combined with a local monodromy computation suffices to classify algebraic ${}_{n+1}F_n$ hypergeometric functions in general.
\end{remark}

\pagebreak

\part{Non-linear differential equations}\label{part:non-linear}
We now prepare to recall and prove our main results about non-linear (isomonodromy) differential equations.
\section{Main theorem and outline of proof}
\subsection{Statement of the conjecture for isomonodromy foliations, and the main theorem}
We begin by specializing \autoref{conjecture:general-algebraic-leaves} to isomonodromic deformations, using the definitions laid out in \autoref{subsubsection:NAGM}. We defer the precise definition of the stack $\mathscr{M}_{dR}(X/S,r)$, the isomonodromy foliation  etc.~until \autoref{subsection:M_dR-and-isomonodromy}, and for now concretely describe the statement of our main theorem and the implications of \autoref{conjecture:general-algebraic-leaves}  in the setting of isomonodromy in classical language.  

We use notation as in \autoref{notation:X/S}, i.e. $R\subset \mathbb{C}$ is a finitely-generated $\mathbb{Z}$-algebra and $f: X\to S$ is a smooth morphism of smooth $R$-schemes equipped with a relative normal crossings compactification $(\overline{X}, D)$ over $S$. Let $\mathscr{K}$ be the fraction field of $R$. Fix $s\in S(R)$, and let $$(\mathscr{E}, \nabla: \mathscr{E}\to \mathscr{E}\otimes \Omega^1_{\overline{X}_s}(\log D_s))$$ be a flat vector bundle on $\overline{X}_s$ with regular singularities along $D_s$. Our goal is to understand when $(\mathscr{E}, \nabla)$ has an algebraic isomonodromic deformation over $S$. That is, when does there exist a smooth connected pointed $\mathbb{C}$-scheme $(T,t)$ equipped with a dominant map to to $(S_{\mathbb{C}}, s_\mathbb{C})$ such that $(\mathscr{E}, \nabla)|_{\overline{X}_t}$ extends to an isomonodromic flat bundle on $X_T$?

In this setting, \autoref{conjecture:general-algebraic-leaves} predicts the following:

\begin{conjecture}\label{conj:main-conjecture} 
The following are equivalent:
\begin{enumerate}
	\item (algebraicity) $(\mathscr{E}, \nabla)$ admits an algebraic isomonodromic deformation over $S$.
	\item (integrality) $(\mathscr{E},\nabla)$ has an integral formal isomonodromic deformation.
	\item ($\omega(p)$-integrality) $(\mathscr{E},\nabla)$ has an $\omega(p)$-integral formal isomonodromic deformation.
\end{enumerate}
\end{conjecture}
Let us spell out precisely what conditions (2) and (3) of \autoref{conj:main-conjecture} mean. Let $\widehat S$, $\widehat X, \widehat{\overline X}, \widehat D$ be the formal $R$-schemes obtained by completing $S, X, \overline{X}, D$, respectively, at $s$, and $\widehat{S}_\mathbb{Q}, \widehat{X}_\mathbb{Q}, \widehat{\overline X}_\mathbb{Q}, \widehat{D}_\mathbb{Q}$ the formal schemes obtained by completing $S_\mathbb{Q}, X_\mathbb{Q}, \overline{X}_\mathbb{Q}, D_\mathbb{Q}$ at $s_\mathbb{Q}$; for an integer $N>0$ we analogously define $\widehat S[\frac{1}{N}], \widehat {X}[\frac{1}{N}]$, etc. After replacing $R$ with a localization we may assume that $\widehat{S}\simeq \on{Spf} R[[x_1, \cdots, x_{\dim(S)}]]$ with the ideal of $s$ given by $(x_1, \cdots, x_{\dim S})$; fix such an isomorphism.

Similarly, if $$g: \text{Primes}\to \mathbb{Z}_{\geq 0}$$ is a function, we define $\widehat S^g=\on{Spf}(R[[x_1, \cdots, x_{\dim S}]]^g)$ as in \autoref{subsection:integrality-of-leaves}; note that there are natural maps $\widehat S_{\mathbb{Q}}\to \widehat S^g\to \widehat S$. We define $\widehat  {\overline X}^g, \widehat D^g$ as the base change of $\widehat  {\overline X}, \widehat D$ to $\widehat S^g$.

By \autoref{prop:isomonodromy-existence-uniqueness} below, $(\mathscr{E},\nabla)$ always admits a (unique up to isomorphism) formal isomonodromic deformation $(\widetilde{\mathscr{E}}, \widetilde\nabla)$ to $(\overline{X}_\mathbb{Q}, D_\mathbb{Q})/\widehat S_\mathbb{Q}$. 
\begin{definition}
	We say that $(\mathscr{E}, \nabla)$ has an \emph{integral formal isomonodromic deformation} if there exists an integer $N>0$ and a flat bundle $(\mathscr{E}', \nabla')$ on $(\widehat  {\overline X}[\frac{1}{N}], \widehat D[\frac{1}{N})/\widehat S[\frac{1}{N}]$ equipped with an isomorphism $(\mathscr{E}', \nabla')\widehat\otimes\mathbb{Q}\simeq (\widetilde{\mathscr{E}}, \widetilde\nabla).$ Given $$g: \text{Primes}\to \mathbb{Z}_{\geq 0},$$ we say that $(\mathscr{E}, \nabla)$ has $g$-integral formal isomonodromic deformation if $(\widetilde{\mathscr{E}}, \widetilde\nabla)$  has a descent to $(\widehat  {\overline X}^g, \widehat D^g)/\widehat S^g$, and we say that it is $\omega(p)$-integral if it is $g$-integral for some $g$ with $$\lim_{p\to \infty} \frac{g(p)}{p}=\infty.$$
	
	We will often drop the word ``formal" and refer to ($\omega(p)$-)integral isomonodromic deformations if the meaning is clear.
\end{definition}
Loosely speaking integrality means that only finitely many primes appear in the denominators of the coefficients of the formal power series appearing as the structure constants of the formal isomonodromic deformation $(\widetilde{\mathscr{E}}, \widetilde\nabla)$; $\omega(p)$-integrality means that there exists $g$ as above so that $p$ does not divide the denominators of the first $g(p)$ coefficients of these power series. Note that we do not require that $\nabla'$ extend to an absolute connection integrally. We will return to and study these notions in \autoref{subection:integral-isomonodromic}.

Now recall the definition of a \emph{Picard-Fuchs equation} from \autoref{subsubsection:NAGM}. Namely, a flat bundle $(\mathscr{E}, \nabla: \mathscr{E}\to \mathscr{E}\otimes \Omega^1_{\overline{X}_s}(\log D_s))$ with regular singularities along $D_s$ is a \emph{Picard-Fuchs equation} if there exists $\pi: Y\to X_s$ smooth projective, and $i\geq 0$, so that $$(\mathscr{E}, \nabla)|_{X_s}\simeq R^i\pi_*\Omega^\bullet_{dR, Y/X}$$ with its Gauss-Manin connection $\nabla_{GM}$, and if the eigenvalues of the residues of $\nabla$ about the components of $D$ are rational numbers in $[0,1)$.\footnote{That the eigenvalues are rational is automatic from the quasi-unipotent local monodromy theorem.} Loosely speaking, a Picard-Fuchs equation is a flat bundle arising as the cohomology of a family of algebraic varieties.

The main result of this paper is:
\begin{theorem}\label{thm:NA-main}
	\autoref{conj:main-conjecture} is true if $(\mathscr{E},\nabla)$ is a Picard-Fuchs equation.
\end{theorem}
As the structure constants of isomonodromic deformations are in many cases controlled by classical non-linear ODE---see for example \autoref{example:schlesinger}---this implies \autoref{conjecture:fg-power-series} for those ODE at initial conditions corresponding to Picard-Fuchs equations. 
\subsection{Outline of proof}\label{subsection:outline-of-proof}
We now give a sketch of the proof of \autoref{thm:NA-main}, in somewhat more detail than we did in \autoref{subsection:intro-outline}. We complete the proof in \autoref{section:proof-of-main-theorem}; most of the hard work happens in \autoref{section:deformation-of-hodge} and \autoref{section:smoothness}.

We fix $\pi: Y\to X_s$ such that $(\mathscr{E}, \nabla)|_{X_s}=(R^i\pi_*\Omega^\bullet_{dR, Y/X_s}, \nabla_{GM})$. By replacing $\overline{X}, D, S$ by covers (possibly ramified along $D$), we may assume that $\mathbb{V}:=R^i\pi_*\mathbb{C}$ has unipotent monodromy around $D$ and hence that $(\mathscr{E}, \nabla)$ has nilpotent residues along $D$. By \autoref{lemma:passing-to-cover-trick} and \autoref{prop:finite-orbit-criterion}, it will suffice to show that the isomorphism class of $\mathbb{V}$ has finite $\pi_1(S,s)$-orbit, under the natural action of $\pi_1(S,s)$ on isomorphism classes of local systems on $X_s$.

Let $(\widetilde{\mscr{E}}, \widetilde\nabla)$ denote the formal isomonodromic deformation of $(\mscr{E}, \nabla)$ to $(\widehat X, \widehat D)_{\mathscr{K}}/\widehat S_{\mathscr{K}}$, where $\mathscr{K}$ is the fraction field of $R$. The main point will be to show that the Hodge filtration extends Griffiths-transversally to $(\widetilde{\mathscr{E}}, \widetilde{\nabla})$. We do so in \autoref{theorem:extending-the-Hodge-filtration}: 
\begin{theorem}[see \autoref{theorem:extending-the-Hodge-filtration}]\label{thm:extendhodge}
    There exists a Griffiths-transverse filtration on $(\widetilde{\mscr{E}}, \widetilde\nabla)$ (relative to $\widehat{S}_{\mathscr{K}}$) deforming the Hodge filtration on $(\mathscr{E},\nabla)_{\mathscr{K}}$.   
\end{theorem}

This will suffice by \autoref{theorem:esnault-kerz}, essentially due to Esnault-Kerz \cite[Theorem 1.1]{esnault2024non} (see also \cite[Corollary 4.4.3]{litt2024motives}), which explains that such a Griffiths-transverse deformation of the Hodge filtration guarantees finiteness of the $\pi_1(S,s)$-orbit of $\mathbb{V}$. See \autoref{section:hodge-preliminaries} for an explanation.

We now explain the proof of \autoref{thm:extendhodge} (see \autoref{theorem:extending-the-Hodge-filtration} for the precise statement), eliding some details when $(\mathscr{E},\nabla)$ is reducible.
%\begin{lemma}
%    Suppose that $(\mscr{E}_{\hat{S}_s}, \nabla_{\hat{S}_s})$ admits a Griffiths transverse filtration deforming $F^{\bullet}$.   Then the orbit of $\rho$ under $\pi_1(S, s)$ is finite. Therefore, \autoref{thm:extendhodge} implies \autoref{thm:NA-main}. 
%\end{lemma}
%\begin{proof}
%    Take a contractible analytic neighborhood $s\in W\subset S$. Let $\Flag\rightarrow W$ be the moduli space of Griffiths transverse filtrations on $\mscr{E}_{\hat{S}_s}|_W$, where we work in the category of analytic spaces. Our assumption that $F^{\bullet}$ deforms in the formal neighborhood of $s$ implies that $\Flag\rightarrow W$ is dominant, and therefore  $F^{\bullet}$ deforms in an analytic neighborhood of $s\in S$. 
%
%    Let $\tilde{S}$ be the universal cover of $S$, and $(\mscr{E}^{univ}, \nabla^{univ})$ the universal isomonodromic deformation of $(\mscr{E}, \nabla)$ to $\X_{\tilde{S}}$. The locus in $\tilde{S}$ where $(\mscr{E}^{univ}, \nabla^{univ})$ admits a Griffiths transverse filtration is a closed analytic subset; on the other hand, by the above, this locus contains an open neighborhood of $s$, and therefore must be $\tilde{S}$ itself. 
%
%    We deduce that every $\rho'$ in the $\pi_1(S, s)$-orbit of $\rho$ underlies a $\mb{Z}$-VHS; by \autoref{thm:delignefinite}\josh{we should perhaps move Deligne's theorem}, this is a finite set, as required.
%\end{proof}
Let $S_n\subset \widehat{S}_\mathscr{K}$ be the local Artin $\mathscr{K}$-scheme defined by the $n$-th power of the ideal sheaf of $s_{\mathscr{K}}$, and suppose for the sake of induction that the Hodge filtration on $(\mathscr{E},\nabla)$ extends Griffiths-transversally to $(\widetilde{\mathscr E}, \widetilde{\nabla})|_{(\overline{X}_{S_n}, D_{S_n})}/S_n$. By a version of the $T^1$-lifting theorem (\autoref{cor:t1-lifting-rank-1}), it will suffice for the induction step to show that for any map $$q: S_n[\epsilon]/\epsilon^2\to S$$ reducing to the natural inclusion $S_n \hookrightarrow S$ mod $\epsilon$, the Hodge filtration further extends Griffiths-transversally over the isomonodromic deformation to $$q^*(\overline{X}, D)=(\overline{X}, D)\times_{S, q} S_n[\epsilon]/\epsilon^2.$$

We explicitly compute the obstruction to doing so in \autoref{prop: gt-transverse-obstruction}. It is the image of the Kodaira-Spencer class $v$ of the deformation $q^*(\overline{X}, D)$ of $(\overline{X}_{S_n}, D_{S_n})$ under the map $$H^1(T_{\overline{X}_{S_n}/{S_n}}(-\log D_{S_n}))\to H^1(\on{End}(\widetilde{\mathscr E}, \widetilde{\nabla})_{dR}|_{\overline{X}_{S_n}}/F_{\text{Hodge}}^0)$$ induced by the connection $\widetilde \nabla$; see \autoref{prop: gt-transverse-obstruction} for the definition of the complex $F_{\text{Hodge}}^0$. 

The main idea is to compare this obstruction to the $p$-curvature of the isomonodromy foliation (in characteristic $p$); we may reduce mod $p$ for sufficiently large $p$ due to the assumption of the existence of an ($\omega(p)$-)integral isomonodromic deformation. We will make this comparison using Ogus-Vologodsky's \cite{ogus-vologodsky} and Schepler's \cite{schepler2005logarithmic} (logarithmic) non-abelian Hodge theory in positive characteristic,\footnote{This latter is mostly worked out only under the assumption of nilpotent residues at infinity, which is why we reduce to this case. We expect that most of the arguments of this paper would go through without this reduction given the existence of a good parabolic theory.} as well as the Higgs-de Rham flow of Lan-Sheng-Zuo \cite{lan2019semistable}, as explained by Esnault-Groechenig \cite{eg_revisit}.

We explain this comparison now. Under the assumption that the Hodge-de Rham spectral sequence for  $\on{End}(\widetilde{\mathscr E}, \widetilde{\nabla})_{dR}|_{\overline{X}_{S_n}}$ degenerates at $E_1$, which we prove inductively in the course of proving \autoref{theorem:smoothness}, it will suffice to show that the image of $v$ under the natural map 
\begin{equation}\label{eqn:sketch-ob-hodge}
	H^1(T_{\overline{X}_{S_n}/{S_n}}(-\log D_{S_n}))\to H^1(\on{gr}^{-1}_{F_{\text{Hodge}}}\on{End}(\widetilde{\mathscr E}, \widetilde{\nabla})_{dR}|_{\overline{X}_{S_n}})
\end{equation}
 induced by $\widetilde{\nabla}$ (i.e.~by the Higgs field on the associated graded bundle) vanishes. We now pass to characteristic $p$ for sufficiently large $p$ and show that this element vanishes there. Let $A$ be a spreading-out of $S_n$ and $A_\mf{p}$ the fiber of $A$ over some closed point $\mf{p}$ of $\on{Spec}(R)$. Let $B_\mf{p}$ be the preimage of $A_\mf{p}$ under the absolute Frobenius on $S_{\mf p}$.

Using non-abelian Hodge theory in characteristic $p$, we construct (in \autoref{defn:lifting-and-conjugate-filtration}) a flat bundle $(\mathscr{E}_B, \nabla)$ on $(\overline{X}, D)_{B_{\mf p}}$ equipped with a conjugate filtration $F_{\text{conj}}$. In \autoref{appendix:NAGM} (summarized in \autoref{subsection:appendix-summary}) we give an explicit formula for the $p$-curvature of the isomonodromy foliation on $\mathscr{M}_{dR}(X/S)$ (suitably understood and restricted to our given $B_{\mf p}$-point of $\mathscr{M}_{dR}(X/S)$) as the map 
\begin{equation}\label{eqn:sketch-p-curvature} F_{\text{abs}}^*T_{S, B_{\mf p}}\to H^1(T_{\overline{X}'_{B_{\mf p}}/B_{\mf p}}(-\log D'_{B_{\mf p}})) \to  \mb{H}^1(F^*_{\text{abs}}(T_{\overline{X}'_{B_{\mf p}}/B_{\mf p}}(-\log D_{B_{\mf p}}))_{dR}) \overset{\psi_p}{\to} 
        \mb{H}^1(\End(\mscr{E}_B)_{dR})
        \end{equation}
       where the superscript $'$ denotes the Frobenius twist, the first map is the pullback of the Kodaira-Spencer map for $(\overline{X}, D)/S$, and $\psi_p$ is the $p$-curvature of $(\mathscr{E}_B, \nabla)$. Characteristic $p$ non-abelian Hodge theory provides a comparison between the $\on{gr}^1_{F_{\text{conj}}}$-component of \eqref{eqn:sketch-p-curvature} and the reduction of \eqref{eqn:sketch-ob-hodge} mod $\mf p$.
       
  Finally, by the degeneration of the conjugate spectral sequence (\autoref{prop:conjugate-ss-degenerates-e2}), it will suffice to show that the $p$-curvature of the isomonodromy foliation restricted to the $B_{\mf p}$-point of $\mathscr{M}_{dR}(X/S)$ corresponding to $(\mathscr{E}_B, \nabla)$---that is, \eqref{eqn:sketch-p-curvature}---vanishes. This will follow as in \autoref{prop:p-power-leaves} once we show that $(\mathscr{E}_B, \nabla)$ gives rise (loosely speaking) to a $p$-power leaf of the isomonodromy foliation; see \autoref{prop:crystal-p-power-leaves} for the precise statement.
  
  We check this by constructing an isomonodromic lift of $(\mathscr{E}_B, \nabla)$ to characteristic zero, essentially using the fact that one may, in good situations, check if a power series (in positive characteristic) is a solution to a differential equation by lifting to characteristic zero. This is an analogue of \autoref{proposition:integral-models-of-leaves}; see \autoref{cor:vanishing-p-curvature-for-p-power-leaves-isomonodromy} for a precise statement in the setting of isomonodromy. This lift is constructed in \autoref{subsection:proof-of-key-lemma} using the Higgs-de Rham flow of \cite{lan2019semistable}, closely following the presentation of \cite{eg_revisit}.

The key point in the construction of this lift of $(\mathscr{E}_B,\nabla)$ to characteristic zero is as follows. Set ${A}_{\widehat R}$ to be the $\mf{p}$-adic completion of $A$. We choose an appropriate Frobenius lift $F_S$ on the $\mf p$-adic completion of $\widehat S$, and set ${B}_{\widehat R}$ to be the preimage of $A_{\widehat R}$ under $F_S$. By induction we have a Griffiths-transverse filtration on an integral model of $(\widetilde{\mathscr{E}}, \widetilde{\nabla})$ over $(\overline{X}, D)_{A_{\widehat R}}$, deforming the Hodge filtration on $(\mathscr{E},\nabla)$: we denote these by $(\mscr{E}_{A_{\widehat R}}, \nabla_{A_{\widehat R}}, F_{A_{\widehat R}}^{\bullet})$. The Higgs-de Rham flow takes this filtered flat bundle as input and constructs a new  flat bundle $\Phi(\mscr{E}_{A_{\widehat R}}, \nabla_{A_{\widehat R}}, F_{A_{\widehat R}}^{\bullet})$ on $(\overline{X}, D)_{B_{\widehat R}}/B_{\widehat R}$, lifting $(\mathscr{E}_B, \nabla)$. We now observe that (essentially by Fontaine-Laffaille theory), as $(\mathscr{E},\nabla)$ was a Picard-Fuchs equation, we have $$\Phi(\mscr{E}_{A_{\widehat R}}, \nabla_{A_{\widehat R}}, F_{A_{\widehat R}}^{\bullet})|_{(\overline{X}_s, D_s)}\simeq (\mathscr{E},\nabla).$$ That is, $(\mathscr{E},\nabla)$ is a ``Higgs-de Rham fixed point." 

Combined with the fact (checked in \autoref{prop:special-point-and-isomonodromic} and \autoref{lemma:frobenius-pullback-isomonodromic}) that the base change of $\Phi(\mscr{E}_{A_{\widehat R}}, \nabla_{A_{\widehat R}}, F_{A_{\widehat R}}^{\bullet})$ to characteristic zero is isomonodromic, we see that $\Phi(\mscr{E}_{A_{\widehat R}}, \nabla_{A_{\widehat R}}, F_{A_{\widehat R}}^{\bullet})$ is itself an integral model of our isomonodromic deformation; as these are essentially unique, by \autoref{cor:unique-integral-isomonodromic}, it must agree with the base change of our given ($\omega(p)$-)integral model to $B_{\widehat R}$. As a result the hypotheses of \autoref{cor:vanishing-p-curvature-for-p-power-leaves-isomonodromy} are satisfied, and the $p$-curvature map \eqref{eqn:sketch-p-curvature} vanishes as desired.
  \begin{remark}
  Crucial to the proof is the notion of the $p$-curvature of a crystal of functors, developed in \autoref{appendix:crystals-and-functors} and studied in the case of the functor $\mathscr{M}_{dR}(X/S)^\natural$ (the functor of isomorphism classes of objects of $\mathscr{M}_{dR}(X/S)$) in \autoref{appendix:NAGM}. That said, we have tried to write the paper so it can be read without dealing with crystals if one is willing to accept certain results as black boxes; all crystalline computations have been exiled to the appendices, which are summarized in \autoref{subsection:appendix-summary}.
  \end{remark}
\section{Preliminaries on isomonodromy}
\subsection{Basic properties in characteristic zero}\label{subsection:isomonodromy-basic-properties}
In this section we briefly prove the existence and uniqueness of (formal) isomonodromic deformations in characteristic zero, and discuss some of their basic properties. Everything in this section is surely well-known to experts, but we were unable to find an appropriate reference.

Let $R$ be a ring containing $\mathbb{Q}$, and set $S=\on{Spf}(R[[t_1, \cdots, t_n]])$. Let $\overline X$ be a $(t_1, \cdots, t_n)$-adically complete formally smooth (formal) $S$-scheme, and let $D\subset \overline X$ be a  (relative) simple normal crossings divisor over $S$. Let $0: \on{Spec}(R)\to S$ be the $R$-point defined by $(t_1, \cdots, t_n)$. In the following,  for any object over $S$ we append the subscript $0$ to denote the basechange along this map; for example, we have the smooth simple normal crossings pair $(\overline X_0, D_0)$ over $R$.
\begin{definition}\label{defn:naive-isomondromy}
	A flat bundle $(\mathscr{E}, \nabla: \mathscr{E}\to \mathscr{E}\otimes \Omega^1_{\overline X/S}(\log D))$ on $(\overline X, D)/S$ is \emph{isomonodromic} if the connection $\nabla$ extends to an absolute connection $\widetilde{\nabla}: \mathscr{E}\to \mathscr{E}\otimes \Omega^1_{\overline X}(\log D)$. 
	
	Let $(\mathscr{E}_0, \nabla_0)$ be a flat bundle with regular singularities on $(\overline X_0, D_0)$. 
	%Recall that a flat bundle $$(\mathscr{E}, \nabla: \mathscr{E}\to \mathscr{E}\otimes \Omega^1_{X/S}(\log D))$$ on $(X,D)/S$ is said to be \emph{isomonodromic} if the connection $\nabla$ extends to a flat connection $$\widetilde{\nabla}: \mathscr{E}\to \mathscr{E}\otimes \Omega^1_{X/R}(\log D).$$  
	An isomonodromic flat bundle $(\mathscr{E}, \nabla) $ on $(\overline X,D)/S$ equipped with an isomorphism $\varphi: (\mathscr{E}, \nabla)|_{\overline X_0}\overset{\sim}{\to} (\mathscr{E}_0, \nabla_0)$ is said to be an \emph{isomonodromic deformation} of $(\mathscr{E}_0, \nabla_0)$.
\end{definition}
\begin{proposition}\label{prop:isomonodromy-existence-uniqueness}
	Let $(\mathscr{E}_0, \nabla_0)$ be a flat bundle with regular singularities on $(\overline X_0, D_0)$. There exists an isomonodromic deformation of $(\mathscr{E}_0, \nabla_0)$; any two isomonodromic deformations are (non-canonically) isomorphic to one another.
\end{proposition}
Before proving the proposition, we need a lemma on Taylor series.
\begin{lemma}\label{lemma:char-0-taylor-series}
Let $U$ be a $(t_1, \cdots, t_n)$-adically complete formally smooth (formal) $S$-scheme, and let $D\subset U$ be a relative simple normal crossings divisor over $S$. Let $(\mathscr{E}, \nabla: \mathscr{E}\to \mathscr{E}\otimes \Omega^1_U(\log D))$ be a logarithmic flat bundle on $U$. Let $s\in \Gamma(U_0, \mathscr{E}|_{U_0})$ be a flat section to $(\mathscr{E}, \nabla)|_{U_0}$. Then there exists a unique flat section $\tilde s\in \Gamma(U, \mathscr{E})$ extending $s$.
\end{lemma}
\begin{proof}
Let $d=\dim(U/S)$. We first assume $U$ is equipped with an \'etale morphism $U\to \mathbb{A}^d_S$, where we give $\mathbb{A}^d_S$ coordinates $x_1, \cdots, x_d$, and that $D$ is the preimage of $x_1\cdots x_e=0$. Note that the operators $\nabla_{\frac{\partial}{\partial t_i}}, i=1, \cdots,  n,  \nabla_{x_j\frac{\partial}{\partial x_j}}, j=1, \cdots, e, \nabla_{\frac{\partial}{\partial x_j}}, j=e+1, \cdots, d$ all commute.

Let $u$ be any section to $\mathscr{E}$ such that $u|_{U_0}=s$. Now set $$\tilde s=\sum_{i_1, \cdots, i_n} (-1)^{\sum_j i_j}\frac{t_1^{i_1}\cdots t_n^{i_n}}{I!}\left(\prod_{j=1}^n\nabla_{\frac{\partial}{\partial t_j}}^{i_j}\right)(u),$$ where $I!=i_1!i_2!\cdots i_n!$.

Then $\tilde s$ extends $s$ by direct computation. We claim it is flat. Indeed, $\nabla_{\frac{\partial}{\partial t_i}}\tilde s=0$ by direct computation, as the sum telescopes. To see that $\nabla_{\frac{\partial}{\partial x_j}}\tilde s=0$, for $j>e$, first observe that, as $s$ is flat, we have that $\nabla_{\frac{\partial}{\partial x_j}}\tilde s\in (t_1, \cdots, t_n)\mathscr{E}$. Moreover, for all $i$, 
$$\nabla_{\frac{\partial}{\partial t_i}}\nabla_{\frac{\partial}{\partial x_j}}\tilde s=\nabla_{\frac{\partial}{\partial x_j}}\nabla_{\frac{\partial}{\partial t_i}}\tilde s=0.$$ 
But if $\nabla_{\frac{\partial}{\partial x_j}}\tilde s$ were non-zero there would, by the Leibniz rule, be some $i_1, \cdots, i_n$, not all zero, such that $$0=\frac{1}{I!}\left(\prod_{j=1}^n\nabla_{\frac{\partial}{\partial t_j}}^{i_j}\right)(\nabla_{\frac{\partial}{\partial x_j}}\tilde s)$$ was non-zero modulo $(t_1, \cdots, t_n)$, a contradiction. An identical argument shows that $\tilde s$ is annihilated by $\nabla_{x_j\frac{\partial}{\partial x_j}}$ for $j=1, \cdots, e$.

Uniqueness follows as flat sections are determined by their value at any point.

Now for general $U$ we may find an \'etale cover by formal subschemes $U_i$ equipped with \'etale maps to $\mathbb{A}^d_S$ as above. On each $U_i$ we obtain sections $s_i$ extending $s$. By uniqueness, these sections glue to a flat section on all of $U$ as desired.
\end{proof}
\begin{proof}[Proof of \autoref{prop:isomonodromy-existence-uniqueness}]
Let $\{U_i\}_{i\in I}$ be a cover of $X$ by affine formal schemes, and for each $i$ pick a splitting $s_i$ of the inclusion $(U_i)_0\hookrightarrow U_i$ such that $s_i^{-1}(D\cap (U_i)_0)=D\cap U_i$ (these exist by formal smoothness). Given indices $i,j, k$, let $U_{ij}=U_i\cap U_j$ and $U_{ijk}=U_i\cap U_j\cap U_k$. Let $(\mathscr{E}_i, \nabla_i)=s_i^*(\mathscr{E}_0, \nabla_0)$. We claim that for all $i,j$ there is a canonical isomorphism $\varphi_{ij}$ between $(\mathscr{E}_i, \nabla_i)|_{U_{ij}}$ and $(\mathscr{E}_j, \nabla_j)|_{U_{ij}}$ restricting to the identity on $(U_{ij})_0$; this follows immediately by applying \autoref{lemma:char-0-taylor-series} to the bundle $$\on{Hom}((\mathscr{E}_i, \nabla_i)|_{U_{ij}}, (\mathscr{E}_j, \nabla_j)|_{U_{ij}})$$
and the identity section in the special fiber. That these isomorphisms satisfy the cocycle condition, i.e.~$$\varphi_{jk}|_{U_{ijk}}\circ \varphi_{ij}|_{U_{ijk}}=\varphi_{ik}|_{U_{ijk}},$$ follows from the same statement over $(U_{ijk})_0$ and the uniqueness in \autoref{lemma:char-0-taylor-series}. Thus we have constructed gluing data for the $(\mathscr{E}_i, \nabla_i)$ to a bundle $(\mathscr{E}, \widetilde{\nabla}: \mathscr{E}\to \mathscr{E}\otimes \Omega^1_X)$ with (absolute) flat connection on $X$; remembering only the relative connection over $S$ yields an isomonodromic deformation.

Given two isomonodromic deformations $(\mathscr{E}_1, \nabla_1)$ and $(\mathscr{E}_2, \nabla_2)$ of $(\mathscr{E}_0, \nabla_0)$, one may construct an isomorphism between them by choosing absolute connections $\widetilde{\nabla}_1, \widetilde{\nabla}_2$ lifting $\nabla_1, \nabla_2$, and applying \autoref{lemma:char-0-taylor-series} to $\on{Hom}((\mathscr{E}_1, \widetilde{\nabla}_1), (\mathscr{E}_2, \widetilde{\nabla}_2)).$ The isomorphism thus constructed depends in general on the choice of $\widetilde{\nabla}_1, \widetilde{\nabla}_2$.
\end{proof}

\begin{remark}\label{remark:general-isomonodromic-in-char-0}
From the construction of \autoref{prop:isomonodromy-existence-uniqueness}, we immediately deduce another criterion by which we can check if a family of flat bundles is isomonodromic; this latter criterion makes sense even over non-smooth bases and will be important for our inductive arguments, which will take place over a non-smooth base $S$. Namely, let $R$ be a Noetherian ring containing $\mathbb{Q}$. Suppose $\overline{X}\to S$ is a smooth and quasiprojective  morphism of finite type $R$-schemes, and $D\subset \overline{X}$ is a relative simple normal crossings divisor over $S$, with $S_{\text{red}}=\on{Spec}(R)$. Suppose we are given a flat bundle $(\mathscr{E},\nabla)$ on $(\overline{X}, D)/S$. Let $0$ be the $R$-point of $S$ coming from $S_{\text{red}}$. Choose an affine cover $\{U_i\}$ of $\overline{X}$ such that each $(U_i, U_i\cap D)\simeq ((U_i)_0, (U_i)_0\cap D)\times S$; such a cover exists by smoothness. Let $s_i: U_i\to \overline{X}_0$ be the corresponding map splitting the inclusion $(U_i)_0\hookrightarrow U_i$. We say that $(\mathscr{E},\nabla)$ is isomonodromic if for each $i$, $(\mathscr{E},\nabla)|_{U_i}\simeq s_i^*(\mathscr{E},\nabla)|_{\overline{X}_0}$.

In other words, $(\mathscr{E}, \nabla)$ is isomonodromic if it is locally pulled back from $(\overline{X}_0, D_0)$.

Even over such $S$ isomonodromic deformations  exist and are unique, as follows by embedding $S$ in a smooth formal $R$-scheme $S'\simeq \on{Spf}(R[[t_1, \cdots, t_n]])$, embedding $(\overline{X}, D)$ in some smooth simple normal crossings pair $(\overline{X'}, D')/S'$, replacing $(X,D)$ with the completion of $\widehat{(\overline{X}', D')}$ of $(\overline{X}', D')$ along $\overline{X}_0$ and applying the recipes above as in the proof of \autoref{prop:isomonodromy-existence-uniqueness}. That is, locally on $\widehat{(\overline{X}', D')}$ we may, by the smoothness of $(\overline{X}_0, D_0)$, write it as $(\overline{X}_0, D_0)\times \on{Spf}(R[[u_1, \cdots, u_m]])$ and pull back $(\mathscr{E}, \nabla)$ to get a local isomonodromic deformation; these deformations glue over $\widehat{(\overline{X}', D')}$ by \autoref{lemma:char-0-taylor-series} as above. This yields an absolute connection on $\widehat{(\overline{X}', D')}$. We may then restrict to $(\overline{X}, D)/S$ to get a global isomonodromic deformation.
\end{remark}
\begin{remark}
While the construction in \autoref{remark:general-isomonodromic-in-char-0} (embedding $(\overline{X}, D)/S$, with $S$ possibly non-smooth, in a smooth $(\overline{X}', D')/S'$ with $S'$ smooth, and then completing) may seem involved, in all of our situations it will be particularly simple. Namely, $S$ will start life as a closed subscheme of some smooth $S'$ (typically $\on{Spf}(R[[t_1, \cdots, t_n]])$) and $(\overline{X}, D)$ will be pulled back from a family $(\overline{X}', D')$ over $S'$. In this case it is easy to say what it means for a flat bundle $(\mathscr{E},\nabla)$ on $(\overline{X},D)/S$  to be isomonodromic. It simply means that $(\mathscr{E},\nabla)$ is the pullback of an isomonodromic flat bundle on $(\overline{X}', D')$, defined as in \autoref{defn:naive-isomondromy}.

In fact this is always how we will check that flat bundles are isomonodromic, in characteristic zero. In mixed or positive characteristic (where we study several different notions of isomonodromy, all of which are equivalent in good situations, to which we will always reduce) we will almost always compare to this notion and check isomonodromy in characteristic zero, typically with the aid of \autoref{cor:unique-integral-isomonodromic} below. The only exception is in \autoref{appendix:NAGM}, where we perform certain formal calculations explicitly with the crystalline notion of isomonodromy developed there.
\end{remark}

\begin{remark}
See \autoref{construction:analytic-isomonodromy} for a sketch of an analytic construction of isomonodromic deformations over Artin $\mathbb{C}$-schemes.
\end{remark}

\begin{lemma}\label{lemma:aut-functor-smooth}
In the situation of \autoref{remark:general-isomonodromic-in-char-0}, let $(\mathscr{E}, \nabla)$ be an isomonodromic flat bundle on $(\overline{X}, D)/S$. Then $\on{Aut}(\mathscr{E}, \nabla)$ is formally smooth over $S$.
\end{lemma}
\begin{proof}
It suffices to show that $H^0(\on{End}(\mathscr{E})_{dR})$ is locally free on $S$ and that its formation commutes with arbitrary base change, as the functor $\on{Aut}(\mathscr{E}, \nabla)$ is given by the invertible elements of $H^0(\on{End}(\mathscr{E})_{dR})$. When $S$ is (formally) smooth, this is immediate from the fact that $H^*(\on{End}(\mathscr{E})_{dR})$ carries a (Gauss-Manin) connection (as $\on{End}(\mathscr{E}, \nabla)$ carries an absolute logarithmic connection on $(\overline{X}, D)$); see \cite[(3.0)]{katz1970nilpotent} for the existence of this connection and \cite[Proposition (8.8)]{katz1970nilpotent} for local freeness. In general the claim follows by embedding into the (formally) smooth setting as in \autoref{remark:general-isomonodromic-in-char-0}, and then restricting to $S$.
\end{proof}
\begin{remark}
See \autoref{remark:aut-functor-smooth} for a sketch of an analytic argument of a form of \autoref{lemma:aut-functor-smooth} for $S$ an Artin $\mathbb{C}$-scheme.
\end{remark}

\subsection{The moduli of flat bundles, and the isomonodromy foliation}\label{subsection:M_dR-and-isomonodromy}
Let $\overline{f}:\overline{X}\to S$ be a smooth morphism, and $D\subset \overline X$ a simple normal crossings divisor over $S$; set $X=\overline{X}\setminus D$ and $f=\overline{f}|_X$. The moduli stack $\mathscr{M}_{dR}(X/S, r)$ is the stack whose fiber over an $S$-scheme $T$ is the groupoid of flat vector bundles $(\mathscr{E}, \nabla)$ on $\overline{X}_T/T$ of rank $r$ with logarithmic singularities along $D_T$. The stack $\mathscr{M}_{dR}^0(X/S, r)$ is the stack whose fiber over an $S$-scheme $T$ is the groupoid of triples $(\mathscr{E}, \nabla, \phi)$ where $(\mathscr{E}, \nabla)$ is as above and $\phi: \on{det}(\mathscr{E}, \nabla)\overset{\sim}{\to}(\mathscr{O}_{\overline X_T}, d)$ is an isomorphism of flat bundles. We denote by $\mathscr{M}_{dR}(X/S)$ the disjoint union $$\mathscr{M}_{dR}(X/S)=\bigcup_r \mathscr{M}_{dR}(X/S, r),$$ and similarly with $\mathscr{M}_{dR}^0(X/S, r)$. Set $\mathscr{M}_{dR}(X/S, r)^\natural$ to be the functor $\on{Sch}/S\to \on{Sets}$ sending an $S$-scheme $T$ to the set of isomorphism classes of objects of $\mathscr{M}_{dR}(X/S, r)(T)$, and similarly with $\mathscr{M}_{dR}^0(X/S, r)^\natural$, etc.

In \autoref{appendix:NAGM} we explain that $\mathscr{M}_{dR}(X/S)^\natural$ is a crystal over $S$ and, loosely speaking, compute the $p$-curvature of the corresponding foliation. Note that this construction gives rise to a notion of isomonodromic deformations over arbitrary (not necessarily smooth) PD-thickenings of our base; see \autoref{subsection:appendix-nagm}. 

In characteristic zero this can be understood via \autoref{remark:general-isomonodromic-in-char-0}.  That such isomonodromic deformations exist follows from that remark or from the recipe of \autoref{subsection:appendix-nagm}; we will construct them directly in the cases of interest in this paper.

We now define the complexes relevant for the computation of the $p$-curvature of the isomonodromy foliation and show they control the deformation theory of $\mathscr{M}_{dR}$.

\subsection{Deformation of flat bundles}

Suppose $Y\rightarrow T$ is a smooth morphism, and $(\mscr{E}, \nabla)$ is a flat bundle on $Y/T$. 
The Atiyah bundle is a vector bundle $\At_{Y/T}(\mscr{E})$ fitting in an exact sequence 
\begin{equation} \label{eqn:atiyah-exact-sequence} 0\rightarrow \End(\mscr{E})\rightarrow \At_{Y/T}(\mscr{E}) \rightarrow T_{Y/T}\rightarrow 0.
\end{equation}
known as the Atiyah exact sequence. 

More precisely, let $\mscr{D}^1(\mscr{E})$ be the sheaf of differential operators on $\mscr{E}$ of order $\leq 1$. This fits into a short exact sequence 
\[
0\rightarrow \End(\mscr{E})\rightarrow \mscr{D}^1(\mscr{E})\xrightarrow{\sigma} \End(\mscr{E})\otimes T_{Y/T}\rightarrow 0,
\]
where $\sigma$ is the map sending a differential operator $\partial$ to its \emph{symbol} $\sigma(\partial)$. Then 
\[
\At_{Y/T}(\mscr{E}):= \{\partial \in \mscr{D}^1_S(\mscr{E})| \sigma(\partial)\in \id \otimes T_{Y/T}\subset \End(\mscr{E})\otimes T_{Y/T}\}.
\]
A splitting $q^\nabla$ of the Atiyah exact sequence \eqref{eqn:atiyah-exact-sequence} compatible with Lie algebra structures  is the same as a flat connection $\nabla$ on $\mathscr{E}$. If $D\subset Y$ is a relative simple normal crossings divisor on $Y$ over $T$, we define $\on{At}_{(Y,D)/T}(\mathscr{E})\subset \At_{Y/T}(\mscr{E})$ to be the preimage of $T_{Y/T}(-\log D)$ under the symbol map; again a splitting of the natural surjection $$\on{At}_{(Y,D)/T}(\mathscr{E})\to T_{Y/T}(-\log D)$$ compatible with Lie algebra structures is the same as a logarithmic flat connection on $\mathscr{E}$.

See e.g.~\cite[\S3.1]{landesman2024geometric} for more background on Atiyah bundles.

Fix an isomorphism $\on{det}(\mathscr{E})\overset{\sim}{\to}\mathscr{O}_Y$ (assuming one exists). Then there is a natural induced map $$\on{tr}: \on{At}_{(Y, D)/T}(\mathscr{E})\to \on{At}_{Y/T}(\mathscr{O}_Y)$$ given by sending a differential operator on $\mathscr{E}$ to its induced action on $\det(\mathscr{E})$. As $\mathscr{O}_Y$ is equipped with the tautological connection $d$, there is a canonical isomorphism $$\on{At}_{Y/T}(\mathscr{O}_Y)\simeq \mathscr{O}_Y\oplus T_{Y/T}$$$$\delta\mapsto (\delta(1), \sigma(\delta)).$$ We set $$\on{At}^0_{(Y,D)/T}(\mathscr{E})=\ker(\on{At}_{(Y,D)/T}(\mathscr{E})\overset{\on{tr}}{\longrightarrow} \on{At}_{Y/T}(\mathscr{O}_Y)\to \mathscr{O}_Y),$$ where the second map sends a differential operator $\delta$ to $\delta(1)$.
\begin{definition}
    The Atiyah-de Rham complex is a complex associated to flat bundle $(\mathscr{E}, \nabla)$, denoted $\At_{(Y,D)/T, dR}(\mscr{E})$
    \[
    0\rightarrow \At_{(Y,D)/T}(\mscr{E})\xrightarrow{\ad_{\nabla}} \End(\mscr{E})\otimes \Omega^1_{Y/T}(\log D)\rightarrow \cdots 
    \]
    where each term and differential in degree greater than zero is same as that of the de Rham complex of $(\End(\mscr{E}), \nabla)$. The first differential $\ad_{\nabla}$ is given by the following recipe. We first define a map $\ell: \At_{(Y,D)/T}(\mscr{E}) \rightarrow \At_{(Y,D)/T}(\mscr{E}\otimes \Omega^1_{Y/T}(\log D))$ given by 
    \[
    \partial \mapsto \partial \otimes \id +\id \otimes L_{\sigma(\partial)},
    \]
    where, for a vector field $v$ tangent to $D$,  $L_v: \Omega^1_{Y/T}(\log D)\rightarrow \Omega^1_{Y/T}(\log D)$ is the Lie derivative with respect to $v$, defined as 
    \[
    L_v:= d\circ i_v+ i_v\circ d
    \]
    where $i_v$ denotes the contraction with $v$ map. Then we define $\ad_{\nabla}$ as  
    \[
    \ad_{\nabla}(\partial):= \nabla \circ \partial - \ell(\partial)\circ \nabla.
    \]
    
    There is a variant complex $\on{At}^0_{(Y,D)/T,dR}(\mathscr{E})$ where we replace $\on{At}_{(Y,D)/T}(\mathscr{E})$ with $\on{At}^0_{(Y,D)/T}(\mathscr{E})$ and $\on{End}(\mathscr{E})$ with the bundle $\on{End}^0(\mathscr{E})$ of trace-zero endomorphisms of $\mathscr{E}$.
    
    We denote by $\mathscr{E}_{dR}$ the complex $$\mathscr{E}\overset{\nabla}{\to} \mathscr{E}\otimes \Omega^1_{Y/T}(\log D)\to \overset{\nabla}{\to} \mathscr{E}\otimes \Omega^2_{Y/T}(\log D)\to \cdots.$$

\end{definition} 

Let us  recall the following piece of algebra needed for the deformation theory to come.
\begin{definition}\label{defn: split-thickening}
    For $T$ an arbitrary scheme, and $M$ an $\mscr{O}_T$-module, let $\mscr{O}_T\oplus M$ be the $\mscr{O}_T$-algebra with multiplication defined as 
    \[
    (t, m)\cdot (t', m') = (tt', tm'+t'm),
    \]
    for local sections $t, t'$ of $\mscr{O}_T$ and $m, m'$ of $M$. 
We write   $T\oplus M$ for the scheme defined by the $\mscr{O}_T$-algebra $\mscr{O}_T\oplus M$. This is a split square zero thickening of $T$, in the sense that there is a diagram
\[
T\xhookrightarrow[]{} T\oplus M \rightarrow T
\]
where the first arrow is a square-zero thickening, and the composition is the identity on $T$.
\end{definition}
As we now recall from \cite{katzarkov1999non}, the Atiyah-de Rham complexes control the deformation theory of the triple $(Y, \mathscr{E}, \nabla)$. 
\begin{proposition}\label{prop:atiyah-de-rham-deformation}
Let $f: \overline{X}\rightarrow S$ be a smooth morphism, and $D\subset \overline X$ a relative simple normal crossings divisor over $S$. Let $(\mathscr{E},\nabla)$ be a flat bundle on $\overline X$ with regular singularities along $D$. Then there is a natural  identification between the  deformation-obstruction theory  of $(X, D, \mscr{E}, \nabla)$, and the  cohomology of the Atiyah-de Rham complex.

That is, suppose $M$ is a $\mscr{O}_S$-module, and denote by $S\oplus M$ the corresponding square-zero thickening of $S$. Then the deformations of $(\overline{X}, D, \mathscr{E}, \nabla)$ to $S\oplus M$ are naturally in bijection with $\mathbb{H}^1(\on{At}_{(\overline{X}, D)/T, dR}(\mathscr{E})\otimes M)$; infinitesimal automorphisms of a given deformation are given by  $\mathbb{H}^0(\on{At}_{(\overline{X}, D)/T, dR}(\mathscr{E})\otimes M)$.

Let $(\X, \mathscr{D})$ be a fixed deformation of $(\overline{X}, D)$ over $S\oplus M$. Suppose the flat bundle $(\mscr{E}', \nabla')$ on $(\X, \mathscr{D})$ is an isomonodromic deformation of $(\mscr{E}, \nabla)$. Let $v\in H^1(T_{\overline{X}/S}(-\log D)\otimes M)$ be the Kodaira-Spencer class corresponding to $(\X, \mathscr{D})$. 

Then $(\mscr{E}', \nabla')$ corresponds to the class in $\mb{H}^1(\At_{(\overline{X}, D)/S, dR}(\mscr{E})\otimes M)$ given by the image of $v$ under 
\begin{equation}\label{eqn: induced-h1-atiyah-splitting}
    H^1(T_{\overline{X}/S}\otimes M) \rightarrow \mb{H}^1(\At_{(\overline{X}, D)/S, dR}(\mscr{E})\otimes M),
\end{equation}

where \eqref{eqn: induced-h1-atiyah-splitting} is given by the splitting of the Atiyah exact sequence \eqref{eqn:atiyah-exact-sequence} $$T_{\overline{X}/S}(-\log D)\rightarrow \At_{(\overline{X}, D)/S}(\mscr{E})$$ induced by the connection on $\mathscr{E}$.
\end{proposition}
\begin{proof}
	The first statement (without the presence of a simple normal crossings divisor) is \cite[Lemma 4.4]{katzarkov1999non}; the proof in the presence of an snc divisor is identical. The second statement is \cite[Lemma 4.5]{katzarkov1999non}.
\end{proof}

\begin{corollary}\label{cor:de-rham-complex-deformations}
	Notation as in \autoref{prop:atiyah-de-rham-deformation}. The deformation theory of a fixed flat bundle $(\mathscr{E}, \nabla)$ on $(\overline{X}, D)/S$ is controlled by the cohomology of the de Rham complex $\on{End}(\mathscr{E})_{dR}$of $\on{End}(\mathscr{E})$. That is, suppose $M$ is a $\mscr{O}_S$-module, and denote by $S\oplus M$ the corresponding square-zero thickening of $S$. Let  Then the deformations of $(\mathscr{E}, \nabla)$ to $(\overline{X}, D)_{S\oplus M}$ are naturally in bijection with $\mathbb{H}^1(\on{End}(\mathscr{E})_{dR}\otimes M)$; infinitesimal automorphisms of a given deformation are given by  $\mathbb{H}^0(\on{End}(\mathscr{E})_{dR}\otimes M)$.
 \end{corollary}

\subsection{Deformation of filtered flat bundles}
The flat bundles $(\mathscr{E},\nabla)$ we will be deforming will be equipped with a Griffiths-transverse Hodge filtration, i.e.~a decreasing filtration $F^\bullet$ such that $$\nabla(F^i)\subset F^{i-1}\otimes \Omega^1(\log D).$$ We now explain how to understand the obstruction to deforming such a filtration (Griffiths-transversally) to an isomonodromic deformation of $(\mathscr{E},\nabla)$.

\begin{proposition}\label{prop: gt-transverse-obstruction}
Let $f: \overline{X}\rightarrow S$ be a smooth morphism, and $D\subset \overline X$ a relative simple normal crossings divisor over $S$. Let $(\mathscr{E},\nabla)$ be a flat bundle on $\overline X$ with regular singularities along $D$. Suppose $F^i$ is a decreasing Griffiths-transverse filtration on $\mscr{E}$. This induces a filtration on $\End(\mscr{E})$, which we denote by $F^i\End(\mscr{E})$.

Let $(\X, \mathscr{D})$ be a fixed deformation of $(\overline{X}, D)$ over $S\oplus M$, corresponding to a class in $v\in H^1(T_{\overline{X}/S}(-\log D)\otimes M)$.  Let $(\mathscr{E}', \nabla')$ be an isomonodromic deformation of $(\mathscr{E},\nabla)$ to $(\X, \mathscr{D})$.

% Suppose that the Hodge-de Rham spectral sequence for $(\mscr{E}, \nabla, F)$ degenerates.  
 Then the obstruction to deforming the $F^i$'s in a  Griffiths-transverse manner to $\mathscr{E}'$ is   the image of $v$ under 
 \begin{equation}\label{eqn:obstruction-map-gt} H^1(T_{\overline X/S}(-\log D)\otimes M) \rightarrow \mb{H}^1(\on{End}{(\mscr{E})_{dR}}/F^0\otimes M),\end{equation}
where the filtration on $\on{End}{(\mscr{E})_{dR}}\otimes M$ has $F^i$ given by 
\[
F^i\End(\mscr{E})\otimes M\rightarrow F^{i-1}\End(\mscr{E})\otimes M \otimes \Omega^1_{\overline{X}/S}(\log D) \rightarrow \cdots ,
\]
 and the map \eqref{eqn:obstruction-map-gt} is induced by the $\mathscr{O}_{\overline{X}}$-linear(!) map $$T_{\overline X/S}(-\log D)\otimes M \rightarrow \on{End}({\mscr{E}}_{dR}/F^0\otimes M$$ 
 arising from the connection on $\mathscr{E}$. If this obstruction vanishes, the set of such deformations as acted on transitively by 
 $$\mb{H}^0(\on{End}{(\mscr{E})_{dR}}/F^0\otimes M).$$
\end{proposition}
\begin{proof}
    The proof is the same as that of \cite[Lemma 4.4]{katzarkov1999non} with the addition  of one step. We use freely the notation of loc.cit., and we give only the proof in the case when $M=\mscr{O}_T$; the general case is identical. For a complex $C^{\bullet}$, we write $Z^i(C^{\bullet})$ for the group of $i$-th hypercocycles of $C^{\bullet}$.
    
Suppose the isomonodromic deformation $(\mathscr{E}', \nabla')$ of $(\mathscr{E}, \nabla)$ to $\mathscr{X}$  is represented by   a hypercocycle $(\mf{x}_{UV}, \mf{m}_{U})  \in Z^1(\At_{(\overline{X}, D)/S, dR}(\mscr{E}))$, i.e. for a covering of $\overline{X}$ by opens $U_i$'s
    \[
    \mf{x}_{UV} \in \Gamma(U\cap V, \At_{(\overline{X}, D)/S}(\mscr{E})), \ \mf{m}_U\in \Gamma(U, \End(\mscr{E})\otimes \Omega^1_{\overline{X}/S}(\log D))
    \]
    satisfying the hypercocycle conditions. Roughly, $\mf{x}_{UV}$ gives the deformation of the bundle $\mscr{E}$, and $\mf{m}$ gives the deformation of $\nabla$; we refer the reader to loc.cit. for details. 
    
    Now $(\mf{x}_{UV})$ preserve $F^{\bullet}$ if and only if 
    $\mf{x}_{UV} \in F^0 \At_{(\overline{X}, D)/S}(\mscr{E})$, the subsheaf consisting of differential operators preserving the filtration $F^\bullet$. Similarly, $(\mf{m}_{U})$ satisfy Griffiths transversality if and only if $\mf{m}_{U} \in F^{-1}\End(\mscr{E})\otimes \Omega^1_{Y/T}$. Therefore the filtration $F^{\bullet}$ on  $(\mscr{E}, \nabla)$ extends to  $(\mscr{E}', \nabla')$ if and only if we may take   
    $(\mf{x}_{UV}, \mf{m}_U)\in Z^1(F^0\At_{(\overline X, D)/S, dR}(\mscr{E}))$, where $F^0\At_{(\overline X, D)/S, dR}(\mathscr{E})$ is defined identically
     to $F^0(\on{End}(\mathscr{E})_{dR})$, except with $F^0 \At_{(\overline{X}, D)/S}(\mscr{E})$ in degree zero. Finally, we have a short exact sequence of complexes 
    \[
    0\rightarrow F^0\At_{(\overline X, D)/S, dR}(\mscr{E})\rightarrow \At_{(\overline X, D)/S, dR}(\mscr{E}) \rightarrow \on{End}{(\mscr{E})_{dR}}/F^0\rightarrow 0,
    \]
    and hence the obstruction  to extending $F^{\bullet}$ is the image of $[(\mf{x}_{UV}, \mf{m}_U)]$ under 
    \[
    \mb{H}^1(\At_{(\overline X, D)/S, dR}(\mscr{E}))
    \rightarrow 
    \mb{H}^1(\on{End}{(\mscr{E})_{dR}}/F^0),
    \]
    as claimed. A completely analogous computation gives the final claim: we have shown that a Griffiths transverse deformation gives rise to a class in $\mb{H}^1(F^0\on{At}_{(\overline X, D)/S, dR}(\mathscr{E}))$; the set of such Griffiths transverse deformations on a fixed deformation of $(\mathscr{E}, \nabla)$ is exactly the preimage of a fixed class in $\mb{H}^1(\on{At}_{(\overline X, D)/S, dR}(\mathscr{E}))$, on which there is a transitive action of $\mb{H}^0(\on{End}{(\mscr{E})_{dR}}/F^0)$ by the long exact sequence in (hyper)cohomology.
\end{proof}
\begin{corollary}\label{cor:griffiths-transverse-def-unique}
	With notation as in \autoref{prop: gt-transverse-obstruction}, suppose in addition that $H^0(\on{End}(\mathscr{E})_{dR})$ is locally free of rank one and that its formation commutes with arbitrary base change. Suppose moreover that the Hodge-de Rham spectral sequence for $(\mathscr{E}, \nabla, F^\bullet)$ degenerates at $E_1$. Finally, assume $M=\mathscr{O}_S$. Then a Griffiths-transverse deformation of $F^\bullet$ to $(\mathscr{E}', \nabla')$ is unique, if it exists. 
\end{corollary}
\begin{proof}
The assumption implies that 	$H^0(\on{End}(\mathscr{E})_{dR})$ is spanned by the identity map, which is in $F^0$. Thus by the degeneration of the Hodge-de Rham spectral sequence, $\mb{H}^0(\on{End}{(\mscr{E})_{dR}}/F^0)=0$, completing the proof.
\end{proof}

\subsection{Integral isomonodromic deformations}\label{subection:integral-isomonodromic}

We now study some basic properties of integral isomonodromic deformations. Let $R$ be a Noetherian domain containing $\mathbb{Z}$, $B$ a finitely generated $R$-algebra with $B_{\text{red}}=R$, $I$ the nilradical of $B$, and $S=\on{Spf}(B)$. Let $0: \on{Spec}(R)\hookrightarrow S$ be the inclusion induced by $I$. Let $\overline X/S$ be a smooth $S$-scheme, and $D$ a relative simple normal crossings divisor on $\overline X/S$. Let $I=I_1\supset I_2 \supset \cdots \supset I_N=0$ be a decreasing filtration of $B$ by ideals such that $I\cdot I_n\subset I_{n+1}$ for all $n$. Finally, assume:
\begin{assumption}\label{assumption:torsion-free-ideals}
	The $R$-modules $I_m/I_{m+1}$ are torsion-free.
\end{assumption}

Set $S_{\mathbb{Q}}=\on{Spec}(B_{\mathbb{Q}})$, and let $\overline{X}_{\mathbb{Q}}, D_{\mathbb{Q}}$ denote the base-change of $\overline{X}, D$ to $S_{\mathbb{Q}}$.

Let $(\mathscr{E}_0, \nabla_0)$ be a (logarithmic) flat bundle  on $(\overline{X}_0, D_0)/R$. Let $(\widetilde{\mathscr{E}}, \widetilde{\nabla})$ be an isomonodromic deformation of $(\mathscr{E}_0, \nabla_0)|_{X_{0,\mathbb{Q}}}$ to $\overline{X}_{\mathbb{Q}}$. Say that this flat bundle \emph{admits an integral isomonodromic deformation over $S$} if (possibly after replacing $R$ with a localization), there exists a flat bundle on $(\overline{X}, D)/S$ equipped with an isomorphism to $(\widetilde{\mathscr{E}}, \widetilde{\nabla})$ on base change to $\overline{X}_\mathbb{Q}$, i.e.~if $(\widetilde{\mathscr{E}}, \widetilde{\nabla})$ descends to  $(\overline{X}, D)/S$.

Note that in the definition above we do not require that this descent satisfy any additional condition; in particular the connection $\nabla$ need not extend to an absolute connection integrally. This is because integral models of isomonodromic deformations (and flat bundles more generally) are essentially unique.
\begin{proposition}\label{prop:special-fiber-char-0-comparison}
	Let $(\mathscr{E}_1,\nabla_1)$, $(\mathscr{E}_2,\nabla_2)$ be logarithmic flat bundles on $(\overline{X},D)/S$ equipped with isomorphisms $$\varphi: (\mathscr{E}_1,\nabla_1)|_{X_0}\overset{\sim}{\to}(\mathscr{E}_2,\nabla_2)|_{X_0}$$ and $$\psi: (\mathscr{E}_1,\nabla_1)_{\overline{X}_{\mathbb{Q}}}\overset{\sim}{\to} (\mathscr{E}_2,\nabla_2)_{\overline{X}_{\mathbb{Q}}}$$ which are compatible in the evident sense, i.e. $$\psi \bmod I=\varphi_{X_{0, \mathbb{Q}}}.$$ Suppose further that $\on{Aut}(\mathscr{E}_1, \nabla_1)_{S_{\mathbb{Q}}}$ is formally smooth over $S_\mathbb{Q}$, and that $H^1_{dR}(\on{End}(\mathscr{E}_1, \nabla_1)|_{\overline{X}_0})$ is locally free, and that its formation  commutes with arbitrary base change. Finally, suppose \autoref{assumption:torsion-free-ideals} holds. Then $(\mathscr{E}_1,\nabla_1)$, $(\mathscr{E}_2,\nabla_2)$ are isomorphic to one another.
\end{proposition} 
\begin{proof}
By induction it suffices to show that, given an isomorphism  $\varphi_n$ between $(\mathscr{E}_1,\nabla_1)$, $(\mathscr{E}_2,\nabla_2)$ over $B/I_n$, reducing to $\varphi\bmod I$, it may be lifted to an isomorphism over $B/I_{n+1}$. 

Let $\psi_n=\psi|_{X_{B_{\mathbb{Q}}/I_n}}$. By assumption the automorphism functor of $(\mathscr{E}_1, \nabla_1)_{S_{\mathbb{Q}}}$ is smooth over $S_{\mathbb{Q}}$; hence we may lift $\psi_n^{-1} \circ \varphi_{n,\mathbb{Q}}$ to an automorphism $\alpha_n$ of $(\mathscr{E}_1, \nabla_1)|_{\overline{X}_{\mathbb{Q}}}.$ Replacing $\psi$ with $\psi\circ \alpha_n$, we may assume  that $\psi$ and $\varphi_n$ agree over $B_{\mathbb{Q}}/I_n$.

Now $(\mathscr{E}_1, \nabla_1)|_{B/I_{n+1}}$, $(\mathscr{E}_2, \nabla_2)|_{B/I_{n+1}}$, are two a priori distinct deformations of $(\mathscr{E}_1, \nabla_1)|_{B/I_n}\overset{\sim}{\to}(\mathscr{E}_2, \nabla_2)|_{B/I_n}$, where the identification here is given by $\varphi_n$. Hence by \autoref{cor:de-rham-complex-deformations}, their difference is controlled by a class in $H^1_{dR}(\on{End}(\mathscr{E}_1, \nabla_1)|_{\overline{X}_0})\otimes I_n/I_{n+1}$. To complete the induction it suffices to show that this class is zero. But the existence of $\psi$ implies that the class vanishes on base change to $R_{\mathbb{Q}}$; the result follows from the torsion-free-ness of $H^1_{dR}(\on{End}(\mathscr{E}_1, \nabla_1)|_{\overline{X}_0})$ and \autoref{assumption:torsion-free-ideals}.
\end{proof}
Note that in general there exists an element $f$ of $R$ such that $(\mathscr{E}_0, \nabla_0)_{R[1/f]}$ satisfies the cohomological hypotheses of \autoref{prop:special-fiber-char-0-comparison}. As isomonodromic deformations are unique in characteristic zero and have formally smooth deformations in characteristic zero by \autoref{lemma:aut-functor-smooth}, we have the following corollary; compare to \autoref{proposition:integral-models-of-leaves}.
\begin{corollary}\label{cor:unique-integral-isomonodromic}
	Suppose $(\mathscr{E}_0, \nabla_0)$ admits two integral isomonodromic deformations $(\mathscr{E}_1, \nabla_1), (\mathscr{E}_2, \nabla_2)$ over $S$. Assume $S$ satisfies \autoref{assumption:torsion-free-ideals}. Then there exists an element $f$ of $R$ (not depending on $B$) such that after passing to $B[1/f]$, $(\mathscr{E}_1, \nabla_1), (\mathscr{E}_2, \nabla_2)$ are isomorphic. If $(\mathscr{E}_0, \nabla_0)$ satisfies the cohomological conditions of \autoref{prop:special-fiber-char-0-comparison}---i.e.~if $H^1_{dR}(\on{End}(\mathscr{E}_0, \nabla_0)$ is locally free and  its formation  commutes with arbitrary base change---the localization is unnecessary.
\end{corollary} 
\begin{proof}
After replacing $R$ by a localization, we may assume the cohomological conditions of \autoref{prop:special-fiber-char-0-comparison} are satisfied; the condition on the smoothness of the automorphism functor is satisfied by \autoref{lemma:aut-functor-smooth}. It thus suffices to produce a map $\psi$ in characteristic zero as in the statement of that proposition; see \autoref{remark:general-isomonodromic-in-char-0} for an explanation of how to construct this map.
\end{proof}

Finally, we record the behavior of existence of integral isomonodromic deformations under passing to a summand; a completely identical argument works for $\omega(p)$-integral isomonodromic deformations.
\begin{lemma}\label{lemma:integral-isomonodromic-summand}
	Let $(\mathscr{E}_0, \nabla_0)=\bigoplus_{i=1}^n (\mathscr{F}_i, \nabla_i)$ be a logarithmic flat bundle on $(X_0, D_0)/R$. Suppose $(\mathscr{E}_0, \nabla_0)$ admits an integral isomonodromic deformation. Then the same is true for each $(\mathscr{F}_i, \nabla_i)$.
\end{lemma}
\begin{proof}
	Let $j: \overline{X}_{\mathbb{Q}}\hookrightarrow \overline{X}$ be the natural inclusion. Let $(\widetilde{\mathscr{E}},\widetilde\nabla)$ be an integral isomonodromic deformation of $(\mathscr{E}_0, \nabla_0)$ over $\overline{X}$, and let $(\widetilde{\mathscr{F}}_i, \widetilde{\nabla}_i)$ be an isomonodromic deformation of $ (\mathscr{F}_i, \nabla_i)$ over $\overline{X}_\mathbb{Q}$. Consider the composition \begin{equation}\label{eqn: composition-sum-decomp}(\widetilde{\mathscr{E}},\widetilde\nabla)\to j_*j^*(\widetilde{\mathscr{E}},\widetilde\nabla)\to j_*(\widetilde{\mathscr{F}}_i, \widetilde{\nabla}_i)
    \end{equation}
    where the latter map is obtained by applying \autoref{lemma:char-0-taylor-series} to a projection $(\mathscr{E}_0, \nabla_0)|_{\overline{X}_\mathbb{Q}}\to (\mathscr{F}_i, \nabla_i)_{\overline{X}_\mathbb{Q}}.$ The image of \eqref{eqn: composition-sum-decomp} is the desired integral isomonodromic deformation.
\end{proof}

\subsection{The $p$-curvature of the isomonodromy foliation}\label{subsection:appendix-summary}
In this section we briefly summarize the content of  \autoref{appendix:NAGM} (which specializes the formalism of \autoref{appendix:crystals-and-functors}, on crystals of functors, to $\mathscr{M}_{dR}(X/S)^\natural$). For this section we will pretend that $\mathscr{M}_{dR}(X/S)^\natural$ is a smooth $S$-scheme equipped with a horizontal foliation, and let $q:\mathscr{M}_{dR}(X/S)^\natural\to S$ be the structure map; in \autoref{appendix:crystals-and-functors} and \autoref{appendix:NAGM} we make everything rigorous and explain how to understand this in terms of the notion of a crystal of functors. 

Our goal is to understand the $p$-curvature of the isomonodromy foliation and give a criterion for its vanishing in terms of the existence of $(\omega(p)$-)integral isomonodromic deformations.

Let $R$ be a Noetherian ring, $S$ a smooth $R$-scheme, and $(\overline{X}, D)/S$ a smooth projective $S$-scheme equipped with a simple normal crossings divisor over $S$.

Let $T$ be an affine scheme of characteristic $p$ equipped with a map $h: T\to S$. The first main result of the appendices (\autoref{lemma:cocycles-agree}) computes the $p$-curvature of the isomonodromy foliation at a characteristic $p$ point $[(\mathscr{E},\nabla)]\in \mathscr{M}_{dR}(X/S)^\natural(T)$ in terms of the $p$-curvature $\psi_p$ of $(\mathscr{E},\nabla)$, where here $(\mathscr{E},\nabla)$ is a flat bundle on $(\overline{X}, D)_T/T$. We let $t: T\to \mathscr{M}_{dR}(X/S)^\natural$ be the map corresponding to $(\mathscr{E},\nabla)$ and $\overline{f_T}: \overline{X}\to T$ be the structure map.

This is a map $$\Psi: h^*F_{\text{abs}}^*T_S\to t^*T_{\mathscr{M}_{dR}(X/S)^\natural/S}.$$ By \autoref{cor:de-rham-complex-deformations}, we may think of $t^*T_{\mathscr{M}_{dR}(X/S)^\natural/S}$ as $R^1\overline{f_T}_{*}\on{End}(\mathscr{E})_{dR}$.

\begin{remark}
Strictly speaking this is not quite right, and 	$$t^*T_{\mathscr{M}_{dR}(X/S)^\natural/S}=R^1\overline{f_T}_{*}\on{End}(\mathscr{E})_{dR}/\on{Aut}(\mathscr{E}, \nabla),$$ as we are working with $\mathscr{M}_{dR}(X/S)^\natural$ rather than the stack $\mathscr{M}_{dR}(X/S)$. But we will elide this issue for now and address it in the appendices. Since we are primarily interested in the vanishing of $p$-curvature, the quotient will not cause any serious issues.
\end{remark}

Thus the $p$-curvature of the isomonodromy foliation may be thought of as a map 
$$\Psi: h^*F_{\text{abs}}^*T_S\to R^1\overline{f_T}_{*}\on{End}(\mathscr{E})_{dR}/\on{Aut}(\mathscr{E}, \nabla).$$ In \autoref{lemma:cocycles-agree} we describe this map explicitly as follows. It is the composite map $$h^*F_{\text{abs}}^*T_S\to h^*F_{\text{abs}}^*R^1\overline{f_T}_{*}T_{\overline{X}_T/T}(-\log D_T)\xrightarrow{\psi_p} R^1\overline{f_{T}}_*\on{End}(\mathscr{E})_{dR}$$
where the first map above is the pullback of the Kodaira-Spencer map of $(\overline{X}, D)/S$, and the second is the $p$-curvature of $(\mathscr{E}, \nabla)$.
 
 Moreover, we give a criterion for the vanishing of this map, completely analogous to \autoref{prop:p-power-leaves} (namely \autoref{cor:vanishing-p-curvature-for-p-power-leaves-isomonodromy}). In \autoref{prop:p-power-leaves}, we argued that the $p$-curvature vanishes along a $p$-power leaf of a foliation, defined as in \autoref{definition:p-power-leaves}. \autoref{cor:vanishing-p-curvature-for-p-power-leaves-isomonodromy} gives a criterion for checking flatness of sections to $\mathscr{M}_{dR}(X/S)^\natural$ over $p$-power subschemes of $S$ by lifting to characteristic zero, as we now explain in the situation of a horizontal foliation on $Y/S$, where $Y\to S$ is a smooth map of smooth schemes.
 
 Suppose one has a horizontal foliation on some $Y\to S/R$, where $S$ is equipped with local coordinates $t_1, \cdots, t_n$. Let $\mf p$ be a maximal ideal of $R$ of residue characteristic $p$. In good (i.e.~$p$-torsion-free) situations, one may check if a section to $Y$  through some point $y$ over some $p$-power monomial subscheme $Z^{[p]}=V(p, t^{pA_1}, \cdots, t^{pA_m})$ is flat, where the $A_i$ are multi-indices, by lifting to characteristic zero. Namely, it suffices that the leaf through $y$ is integral over the subscheme cut out by $(t_1, \cdots, t_n)(t^{pA_1}, \cdots, t^{pA_m})$. That is, if the leaf through $y$ is integral to order greater than the order of $Z^{[p]}$, then one obtains a $p$-power leaf over $Z^{[p]}$ and hence the $p$-curvature, restricted to $Z^{[p]}$, vanishes. Some version of this is formalized in \autoref{cor:vanishing-p-curvature-for-p-power-leaves-isomonodromy}.
\subsection{Analytic aspects of the isomonodromy foliation}
In this section we give a brief discussion of some analytic aspects of isomonodromy, and an analytic criterion for a flat bundle to have algebraic isomonodromic deformation. 

We first briefly sketch an analytic construction of isomonodromic deformations in the case of primary interest to us.
\begin{construction}\label{construction:analytic-isomonodromy}
	Let $S$ be the spectrum of an Artin local $\mathbb{C}$-algebra, and $\overline{X}$ a smooth projective $S$-scheme equipped with a simple normal crossing divisor $D$ over $S$. Let $X=\overline{X}\setminus D$, and let $0$ be the unique $\on{Spec}(\mathbb{C})$-point of $S$. Let $(\mathscr{E},\nabla)$ be a flat bundle on $\overline{X}_0$ with regular singularities along $D_0$, which we suppose for simplicity has residue matrices along $D_0$ whose eigenvalues have real parts in $[0,1)$. Let $\mathbb{V}=\ker(\nabla^{\text{an}})$ be the associated complex local system on the analytification $X_0^{\text{an}}$ of $X_0$. 
	
	One may give a complex-analytic construction of the isomonodromic deformation of $(\mathscr{E},\nabla)$ as follows. On $X^{\text{an}}$, consider the flat bundle $\mathbb{V}\otimes \mathscr{O}_{X^{\text{an}}}$, equipped with the connection $$\on{id}\otimes d: \mathbb{V}\otimes \mathscr{O}_{X^{\text{an}}}\to \mathbb{V}\otimes \Omega^1_{X^{\text{an}}/S^{\text{an}}}.$$ By e.g.~\cite[Theorem 5.7]{hai2023regular}, this bundle extends to a flat bundle on $\overline{X}$ with regular singularities along $D$, whose residue matrices along $D$ whose eigenvalues have real parts in $[0,1)$. As $\overline{X}$ is projective, this flat bundle admits a unique algebraization. This gives a construction of isomonodromic deformations over singular (Artin) bases $S$.
\end{construction}
\begin{remark}\label{remark:aut-functor-smooth}
This construction gives another proof of \autoref{lemma:aut-functor-smooth}, at least conditional on an appropriate logarithmic Riemann-Hilbert correspondence; such exists at least when $(\mathscr{E}, \nabla)$ has quasi-unipotent monodromy along $D_0$, which will always be the case in our applications; see e.g.~\cite{illusie20051}. (In fact in our case we will always reduce to the case of nilpotent residues, whence more classical results apply.)

Indeed, in the proof of that lemma we needed to show that $H^0(\on{End}(\widetilde{\mathscr{E}})_{dR})$ is locally free on $S$ and that its formation commutes with arbitrary base change, where $(\widetilde{\mathscr{E}}_{dR}, \widetilde \nabla)$ is an isomonodromic deformation of $(\mathscr{E}, \nabla)$. But this follows immediately from the compatibility of logarithmic de Rham cohomology and log Betti cohomology, see e.g.~\cite[Theorem 6.2]{illusie20051}. 
\end{remark}

Now $\overline{f}: \overline{X}\to S$ be a smooth projective map of smooth complex varieties, and let $D\subset \overline{X}$ be a simple normal crossings divisor. Set $X=\overline{X}\setminus D$ and $f=\overline{f}|_X$. Fix $s\in S$ and $x\in X_s$. Then there is a natural outer action of $\pi_1(S,s)$ on $\pi_1(X_s, x)$, inducing an action of $\pi_1(S,s)$ on the set of conjugacy classes of representations of $\pi_1(X_s, x)$.
\begin{proposition}\label{prop:finite-orbit-criterion}
	Let $(\mathscr{E}, \nabla)$ be a logarithmic flat bundle on $(\overline{X}_s, D_s)$, with associated monodromy representation $$\rho_{(\mathscr{E}, \nabla)}: \pi_1(X_s)\to GL_r(\mathbb{C}).$$ Suppose $\rho$ is semisimple, and its irreducible summands have finite determinant. Assume the real parts of the eigenvalues of the residue matrices of $(\mathscr{E}, \nabla)$ along $D_s$ lie in $[0,1)$.  Then $(\mathscr{E}, \nabla)$ admits an algebraic isomonodromic deformation over $S$ if and only if the orbit of the conjugacy class of $\rho$ under $\pi_1(S)$ is finite.
\end{proposition}
\begin{proof}
	By replacing $S$ with a cover we may assume $f$ admits a section, so that $\pi_1(X)=\pi_1(X_s, x)\rtimes \pi_1(S, s)$. We may also assume that the  orbit of the conjugacy class of $\rho$ under $\pi_1(S,s)$ is a singleton. Finally, as $\rho$ is semisimple, it suffices without loss of generality to treat the case where it is irreducible with finite determinant.
	
	For each $\gamma\in \pi_1(S,s)$, we have that $\rho^\gamma$ is conjugate to $\rho$ via some matrix $A_\gamma$. By Schur's lemma, $A$ is well-defined up to scaling. Thus the map $$\pi_1(X)=\pi_1(X_s, x)\rtimes \pi_1(S, s)\to PGL_r(\mathbb{C})$$ $$(\xi, \gamma)\mapsto \rho(\xi)A_\gamma$$ is a well-defined homomorphism and extends $\mathbb{P}\rho$. An argument identical to the proof of \cite[Proposition 2.3.4]{landesman2022canonical} shows that after passing to a further cover $T$ of $S$, we may lift this representation to a $GL_r(\mathbb{C})$-representation $\widetilde\rho$ extending $\rho$. Fix $t\in T$ lying over $s$.
	
	Now consider the flat bundle $(\mathscr{E}', \nabla')$ on $X_T$ with monodromy $\widetilde\rho$, and 
	let $(\widetilde{\mathscr{E}}, \widetilde\nabla)$ be the Deligne canonical extension to $\overline{X}_T$, i.e. the unique-up-to-isomorphism extension whose residue matrices have eigenvalues 
	with real part in $[0,1)$. By uniqueness of the Deligne canonical extension, this bundle is isomorphic to $(\mathscr{E},\nabla)$ when restricted to $\overline{X}_t$. But by \cite[Th\'eor\`eme 5.9]{deligne2006equations} we see that $(\widetilde{\mathscr{E}}, \widetilde\nabla)$ is algebraic, completing the proof. (Explicitly, choosing an (absolute) snc compactification of $\overline{X}_T$, say $\overline{X}'$, and an extension of $(\widetilde{\mathscr{E}}, \widetilde\nabla)$ to $\overline{X}'$ with regular singularities at infinity, by GAGA we have that $(\widetilde{\mathscr{E}}, \widetilde\nabla)$ is algebraic.)
\end{proof}
\begin{remark}
See also \cite{cousin2021algebraic},  \cite[Remark 2.3.6]{landesman2022canonical}, and \cite[Theorem 1.1]{esnault2024non} for some  statements closely related to \autoref{prop:finite-orbit-criterion}.
\end{remark}

The following will be useful to check the criterion of \autoref{prop:finite-orbit-criterion}.
\begin{lemma}\label{lemma:passing-to-cover-trick}
Suppose we have a commutative square of smooth complex varieties
$$\xymatrix{
X' \ar[r] \ar[d] & X \ar[d] \\
S' \ar[r] & S
}$$
where $S$ is smooth, $S'\to S$ is dominant \'etale, and $X$, (resp.~$X'$) is the complement in a smooth projective $S$-scheme (resp.~$S'$-scheme) of a relative simple normal crossing divisor. Fix $s'\in S'$ lying over $s$ and suppose the map $\pi_1(X'_{s'})\to \pi_1(X_s)$ has image with finite index. Let $\rho$ be a semisimple complex representation of $\pi_1(X_s)$. Then the conjugacy class $\rho$ has finite orbit under $\pi_1(S)$ if and only if the pullback $\rho'$ of $\rho$ to $\pi_1(X'_{s'})$ has finite orbit under $\pi_1(S')$.
\end{lemma}
\begin{proof}
That $\rho$ has finite orbit under $\pi_1(S)$ immediately implies the same about the orbit of $\rho'$ under $\pi_1(X'_{s'})$. 

We prove the opposite direction now. Without loss of generality $\rho$ is irreducible, and by replacing $X$ with its base change to $S'$ (using that the image of $\pi_1(S',s)\to \pi_1(S,s)$ has finite index, as $S'\to S$ is dominant, see \cite[Lemma 4.19]{debarre2001higher}), we may assume $S=S', s=s'$. Let $\Gamma$ be the image of $\pi_1(X'_{s})\to \pi_1(X_s)$, so we may view $\rho'$ as a representation of $\Gamma$. But $\rho$ is a summand of the induction $\tau$ of $\rho'$ from $\Gamma$ to $\pi_1(X_s)$. By assumption $\tau$ has finite $\pi_1(S)$-orbit; as it has finitely-many irreducible summands, the same is true for each summand, completing the proof.
\end{proof}

\section{Preliminaries on complex Hodge theory}\label{section:hodge-preliminaries}
We now briefly recall some Hodge theoretic-terminology that will appear in the proof of our main theorem, \autoref{thm:NA-main}. After introducing this terminology, the main point of this section will be to give a Hodge-theoretic condition which guarantees that the criterion of \autoref{prop:finite-orbit-criterion} holds, and hence that a flat bundle (underlying a polarizable $\mathbb{Z}$-variation of Hodge structure) admits an algebraic isomonodromic deformation. In the case of a smooth projective map $\overline{X}\to S$, this is the main result of \cite{esnault2024non}; see also \cite[Corollary 4.4.3]{litt2024motives}. We will briefly indicate how to extend the argument to the case where $\overline{X}$ is equipped with a simple normal crossing divisor $D$ over $S$. We are very grateful to H\'el\`ene Esnault and Moritz Kerz for several useful discussions on this topic.

Throughout this section $S$ will be a finite type $\mathbb{C}$-scheme, $\overline{X}$ will be a smooth projective $S$-scheme, and $D$ will be a simple normal crossings divisor on $\overline{X}$ over $S$. Let $X=\overline{X}\setminus D$. All flat bundles will be assumed to have nilpotent residues along $D$.
\subsection{Variations of Hodge structure}
Let $(\mathscr{E}, \nabla)$ be a flat bundle on $\overline{X}/S$ with regular singularities along $D$. Recall that a decreasing filtration $F$ on $\mathscr{E}$ is {Griffiths-transverse} if $$\nabla(F^i)\subset F^{i-1}\otimes \Omega^1_{\overline{X}/S}(\log D)$$ for all $i$. The main source of flat bundles equipped with Griffiths-transverse filtrations is Hodge theory: if $g: Y\to X$ is a smooth projective morphism then $R^ig_*\mathscr{O}_{Y, dR}$ is equipped with a natural Griffiths-transverse filtration, namely the Hodge filtration. 

Suppose $S=\on{Spec}(\mathbb{C})$. A polarizable complex variation $(V, V^{p, q}, D)$ of Hodge structures on $X^{\text{an}}$ is (by definition) equipped with a Griffiths-transverse Hodge filtration. In this case, the Hodge filtration extends canonically to a Griffiths-transverse filtration on the Deligne canonical extension of the corresponding flat bundle to $\overline{X}$, i.e.~the unique extension whose residue matrices have eigenvalues with real parts in $[0,1)$. See e.g.~\cite[Proposition 4.1.4]{landesman2024geometric} and the references therein.
\subsection{Higgs bundles}\label{subsection:higgs-bundles}
Given a flat bundle $(\mathscr{E},\nabla)$ equipped with a Griffiths transverse filtration $F$, one obtains a (logarithmic) Higgs bundle on $(\overline{X}, D)/S$ as follows. Namely, we consider $\on{gr}_F\mathscr{E}$, equipped with the Higgs field $$\theta=\on{gr}_F\nabla: \bigoplus F^{i}/F^{i+1}\to \bigoplus F^{i-1}/F^i\otimes \Omega^1_{\overline{X}/S}(\log D).$$

By the Higgs complex of a Higgs bundle $(\mathscr{F}, \theta)$ on $(\overline{X}, D)/S$, we mean the complex of $\mathscr{O}_{\overline{X}}$-modules $$\mathscr{F}\to \mathscr{F}\otimes \Omega^1_{\overline{X}/S}(\log D)\to  \mathscr{F}\otimes \Omega^2_{\overline{X}/S}(\log D)\to \cdots.$$
The Higgs cohomology $H^*_{\text{Higgs}}$ of $(\mathscr{F}, \theta)$ will be the hypercohomology of this complex.

\subsection{Deligne's finiteness theorem and deformations of the Hodge filtration}

We first explain our motivation for the Hodge-theoretic criterion of this subsection. In \autoref{prop:finite-orbit-criterion}, we explained that one may check if a flat bundle $(\mathscr{E}, \nabla)$ on the fiber $(\overline{X}_s, D_s)$ of a map $(\overline{X}, D)\to S$ has integral isomonodromic deformation by analyzing the orbit of the monodromy representation $\rho$ of $(\mathscr{E}, \nabla)$ under $\pi_1(S,s)$. There is a simple Hodge-theoretic criterion for this orbit to be finite when $(\mathscr{E}, \nabla)$ underlies a $\mathbb{Z}$-variation of Hodge structure.

We will use the following beautiful theorem of Deligne.
\begin{theorem}[See {\cite[Th\'eor\`eme 0.5]{deligne1987theoreme}}]\label{thm:delignefinite}
Let $X$ be a smooth complex algebraic variety, and $r$ a positive integer. The set of isomorphism classes of $\mathbb{C}$-local systems of rank $r$ on $X$ underlying a polarizable $\mathbb{Z}$-variation of Hodge structure is finite.
\end{theorem}
The following corollary is immediate:
\begin{corollary}
	Let $\Gamma\subset \on{Out}(\pi_1(X))$ be a subgroup, and $\mathbb{V}$ a $\mathbb{C}$-local system on $X$. If for all $\gamma\in \Gamma$, the local system $\mathbb{V}^\gamma$ underlies a polarizable $\mathbb{Z}$-variation of Hodge structure, then $\Gamma\cdot[\mathbb{V}]$ is finite.
\end{corollary}
Thus, to verify the criterion of \autoref{prop:finite-orbit-criterion}, it will suffice to show that, if $\mathbb{V}$ is the local system on $\overline{X}_s/D_s$ corresponding to $(\mathscr{E}, \nabla)$, then $\mathbb{V}^\gamma$ underlies a $\mathbb{Z}$-variation of Hodge structure for every $\gamma\in \pi_1(S)$. The key observation of this section, essentially due to Esnault-Kerz \cite{esnault2024non}, is that there is a local criterion which guarantees that this property holds. The idea is that it is enough to check that the Hodge filtration extends Griffiths-transversally to the formal isomonodromic deformation of $(\mathscr{E}, \nabla)$. By the closedness of the non-abelian Hodge loci (see \cite[\S12]{simpson62hodge} or \cite[Appendix A]{esnault2024non}) this implies that every point on the leaf of the isomonodromy foliation through $(\mathscr{E}, \nabla)$ corresponds to a flat bundle underlying a $\mathbb{Z}$-variation of Hodge structure.
\begin{theorem}\label{theorem:esnault-kerz}
	Let $\overline{X}\to S$ be a smooth projective morphism of smooth complex varieties, and $D\subset \overline{X}$ a simple normal crossings divisor over $S$. Set $X=\overline{X}\setminus D$ and suppose that $X\to S$ is ``good" in the sense of \cite[\S1]{esnault2024non}. Fix a point $s\in S$, and let $\mathbb{V}$ be a local system on $X_s$ underlying a $\mathbb{Z}$-variation of Hodge structure, with unipotent monodromy along $D$. Let $(\mathscr{E}, \nabla)$ be the Deligne canonical extension of $(\mathbb{V}\otimes\mathscr{O}_X, \on{id}\otimes d)$ to $(\overline{X}_s, D_s)$, equipped with its natural Hodge filtration $F$. Then the following are equivalent:
	\begin{enumerate}
		\item The $\pi_1(S,s)$-orbit of the conjugacy class of $\mathbb{V}$ is finite, and
		\item The Hodge filtration on $(\mathscr{E},\nabla)$ extends Griffiths-transversally to the formal isomonodromic deformation of $(\mathscr{E},\nabla)$ to $(\widehat{\overline{X}}, \widehat D)/\widehat S$, where $(\widehat{\overline{X}}, \widehat D)/\widehat S$ denotes the formal completion at $s$.
	\end{enumerate}
\end{theorem}
\begin{proof}[Proof sketch]
	In the case $D$ is empty, this is precisely the equivalence $(1)\iff (5)$ of \cite[Theorem 1.1]{esnault2024non}. The only case where properness is used in the proof of that theorem is in the last paragraph of \cite[\S6.3]{esnault2024non}, where it is used to show the following: if $\Delta$ is a contractible open neighborhood of $s\in S$, and $(\mathscr{E},\nabla)$ is an irreducible isomonodromic flat bundle  on $\overline{X}_\Delta/\Delta$ equipped with a relatively Griffiths-transverse filtration $F$ and polarization $Q$, then the set of points $t\in \Delta$, such that $(\mathscr{E},\nabla, F, Q)$ gives rise to a polarized $\mathbb{C}$-VHS, is open. But the case where $D$ is non-empty is immediate from \cite[Proposition 7.1.1]{kato2008classifying}.\footnote{We are very grateful to Moritz Kerz for pointing us to this reference; we briefly explain why it implies the desired statement. We must check that the set of $t$ where $Q$ remains a polarization is open. In a neighborhood of $x\in \overline{X}\setminus D$, it is clear that the set of relatively Griffiths-transverse filtrations where $Q$ satisfies the definition of a polarization is open in the set of all relatively Griffiths-transverse filtrations. \cite[Proposition 7.1.1]{kato2008classifying} gives precisely the analogous statement for $x\in D$, and the result follows by compactness of $\overline{X}$.}
\end{proof}

\section{Preliminaries on non-abelian Hodge theory in positive characteristic}
In this section we briefly summarize the properties of the inverse Cartier transform of Ogus-Vologodsky and Schepler that we will use, as well as some material on spectral sequences. 

\subsection{The inverse Cartier operator}\label{subsection:inverse-Cartier}

Let $k$ be a perfect field of characteristic $p>0$, $S$ a scheme over $k$, $\Xbar \rightarrow S$ a smooth morphism, and $D\subset \Xbar$ a simple normal crossings divisor over $S$. 

Let $(\Xbar', D')\rightarrow S$ denote the Frobenius twist of $(\Xbar, D)\rightarrow S$. 
\begin{definition}
   \begin{enumerate}
       \item Let $\on{HIG}_{p-1}(\Xbar', D')$ denote the category of Higgs bundles on $\Xbar'/S$ with logarithmic poles along $D'$, and such that the Higgs field is nilpotent of order $\leq p-1$. 
   
\item 
Let $\on{MIC}_{p-1}(\Xbar, D)$ denote the category of flat bundles on $\Xbar/S$ with logarithmic poles along $D$, whose $p$-curvature is nilpotent of order $\leq p-1$, and whose residue along each component of $D$ is nilpotent of order $\leq p-1$.
\end{enumerate} 
\end{definition}

\begin{theorem}[{\cite[Theorem 2.8]{ogus-vologodsky}, \cite[Corollary 4.11]{schepler2005logarithmic}}]\label{theorem:inverse-cartier-transform}
Suppose $(\Xbar_{W_2}', D_{W_2}')\rightarrow S_{W_2}$ be a lift of $(\Xbar', D')\rightarrow S$ to $W_2(k)$.  Then there is an equivalence of categories 
    \[
    C^{-1}_{(\Xbar_{W_2}', D_{W_2}')}: \on{HIG}_{p-1}(\Xbar', D') \rightarrow  \on{MIC}_{p-1}(\Xbar, D).
    \] 
\end{theorem}

When the choice of $(\Xbar_{W_2}', D_{W_2}')$ is clear from context, we often omit it from the notation and simply denote the functor by $C^{-1}$, and its quasi-inverse by $C$.

If $(\mscr E', \theta)\in \on{HIG}_{p-1}(\Xbar',D')$ is a Higgs bundle with $\theta=0$, then this functor is particularly simple---we have $C^{-1}(\mscr E',0)= F_{\Xbar/S}^*\mscr E'$ equipped with its canonical (Frobenius-pullback) connection. Likewise, if $(\mathscr{E}, \nabla)\in \on{MIC}_{p-1}(\Xbar, D)$ has vanishing $p$-curvature, we have $C(\mscr E,\nabla)=(\ker \nabla, 0).$ That is, the (inverse) Cartier transform restricts to usual Cartier theory in the $p$-curvature zero case. 

In this paper all Higgs bundles will be graded, and equipped with a Higgs field of degree $-1$ (they will be ``systems of Hodge bundles" in the terminology of \cite{simpson1992higgs}). We briefly explain the compatibility of the inverse Cartier transform with the classical Cartier theory in this case. Let $$\mathscr{E}'=\bigoplus_{i=0}^{p-1} \mathscr{E}_i', \theta_i: \mathscr{E}_i'\to \mathscr{E}_{i-1}'\otimes \Omega^1_{\overline{X}'/S}(\log D')$$ be such a graded Higgs bundle on $(\overline{X}', D')/S$. $(\mathscr{E}', \theta)$ comes with a natural decreasing filtration by sub-Higgs bundles induced by the grading, namely $F^m=\bigoplus_{i\leq -m} \mathscr{E}_i'$. As $C^{-1}$ is exact, this induces a decreasing filtration on $(\mathscr{E},\nabla)=C^{-1}(\mathscr{E}', \theta)$, which we denote by $F_{\text{conj}}$ and refer to as the \emph{conjugate filtration}. Note that as $C^{-1}$ is exact, the $p$-curvature of $\on{gr}_{F_\text{conj}} (\mathscr{E}, \nabla)$ vanishes, and so the $p$-curvature of $(\mathscr{E}, \nabla)$ has degree $1$.

The following gives a more precise description of the $p$-curvature; compare to \autoref{theorem:katz-big-diagram}:

\begin{proposition}\label{prop:higgs-field-vs-p-curvature}
	Fix notation as above. For all $i$ there is a natural isomorphism $$\phi_i: \on{gr}^i_{F_{\text{conj}}} (\mathscr{E},\nabla)\xrightarrow{\sim} F_{\overline{X}/S}^*\mathscr{E}'_{-i},$$ where $F_{\overline{X}/S}^*\mathscr{E}'_{-i}$ is equipped with its canonical (Frobenius-pullback) connection. The following diagram commutes:
	$$\xymatrix{
	\on{gr}^i_{F_{\text{conj}}} (\mathscr{E},\nabla) \ar[r]^-{\phi_i} \ar[d]_{-\on{gr}^i_{F_{\text{conj}}}\psi_p} & F_{\Xbar/S}^*\mathscr{E}'_{-i}\ar[d]^{F_{\Xbar/S}^*\theta_{-i}} \\
	\on{gr}^{i+1}_{F_{\text{conj}}} (\mathscr{E},\nabla)\otimes F_{\text{abs}}^*\Omega^1_{\overline{X}/S}(\log D) \ar[r]^-{\phi_{i+1}} & F_{\Xbar/S}^*\mathscr{E}'_{-i-1}\otimes F_{\Xbar/S}^*\Omega^1_{\Xbar'/S}(\log D')
	}$$
	where $\on{gr}^i_{F_{\text{conj}}}\psi_p$ is the associated graded of the $p$-curvature with respect to the conjugate filtration.
\end{proposition}
\begin{proof}
	In the case $D=\emptyset$, this is a special case of \cite[Corollary 2.24]{ogus-vologodsky}. In general, the existence of the $\phi_i$ follows immediately from the exactness of $C^{-1}$ and its compatibility with classical Cartier theory. That the diagram commutes follows by restricting to $\overline{X}\setminus D$, where it commutes by the case where $D=\emptyset$, whence it commutes on all of $\Xbar$ as $\Xbar\setminus D$ is Zariski-dense in $\Xbar$ and the sheaves in the diagram above are locally free.
\end{proof}

%\begin{theorem}\label{thm:ovcohomologies}
%    Suppose  $h: Y\rightarrow T$ is a smooth morphism of schemes over $\mb{F}_q$, and that 
%    $\widetilde{Y}'\rightarrow \widetilde{T}$ is a lifting to $W_2$ of 
%    $Y'\rightarrow T$. Then for an object 
%    $(\mscr{E}, \theta) \in HIG_{p-1}(Y'/T)$, there is a canonical isomorphism
%    \[
%        R^qh_{*}'^{HIG}(\mscr{E}, \theta) \cong
%         R^qh_{*}^{dR}(C^{-1}(\mscr{E}, \theta)).
%        \]
%\end{theorem}
%\begin{proof}
%    This is a special case of \cite[Theorem~3.8]{ogus-vologodsky}.
%\end{proof}

\subsection{Recollections on spectral sequences}\label{section:recallss}
In the situation above, the conjugate filtration on $(\mathscr{E},\nabla)$ gives rise to a spectral sequence converging to the de Rham cohomology of $(\mathscr{E},\nabla)$, as we now recall.

We first recall some hypercohomology spectral sequences associated with filtered complexes, following \cite{paranjape}. Given an abelian category $\mscr{C}$ and a filtered complex $(K^{\bullet}, F)$ of objects in $\mscr{C}$, we have the  spectral sequence
\[\mscr{E}_{1}^{p,q} = \mc{H}^{p+q}(\gr^pK^{\bullet}) \Rightarrow \mc{H}^{p+q}(K^{\bullet}). \]

On the other hand, given a left exact functor $T: \mscr{C}\rightarrow \mscr{C}'$ to another abelian category $\mscr{C}'$, we can take  \emph{hypercohomologies} of the complexes $\mscr{E}^{\bullet ,q}_1$, giving rise to a hypercohomology spectral sequence:

\begin{definition}\label{defn:conjss}
There is a spectral sequence
\[E_{2}^{p,q}\defeq \mb{T}^{p+q}(\mscr{E}_{1, p,q}) \Rightarrow \mb{T}^{p+q}(K^{\bullet}),\]
where \begin{itemize}
    \item $\mb{T}^{i}$ denotes the $i$-th hyper-derived functor associated to $T$, and 
    \item $\mscr{E}_{1, p,q}$ denotes the complex 
\[\cdots \rightarrow \mscr{E}_1^{p-1,q} \rightarrow \mscr{E}_1^{p,q} \rightarrow \mscr{E}_1^{p+1,q}\rightarrow \cdots\]
with the middle term placed in degree $p+q$. 
\end{itemize}

    We  refer to this as the \emph{conjugate spectral sequence} attached to $(K^{\bullet}, F)$ and $T$.
\end{definition}
This will be the spectral sequence of primary interest to us.

\begin{remark}[Comparison with Ogus--Vologodsky's conjugate spectral sequence]\label{remark:conj-ss-paranjape}
    We will be particularly interested in the following special case, where the spectral sequence is the \emph{conjugate spectral sequence} of Ogus--Vologodsky; we have presented this material in the above manner since we will also analyze the intermediate spectral sequence $\mscr{E}_1^{p,q}$.  
    
    With notation as in \autoref{subsection:inverse-Cartier}, let $(\mscr{E}, \nabla)$ a flat bundle on $(\Xbar,D)/S$, and $N^i\subset \mscr{E}$ a decreasing filtration of $\nabla$-stable bundles such that the associated graded bundle has vanishing $p$-curvature. Let $h: \Xbar\to S$ be the structure morphism.
    
    We take $K^{\bullet}$ to be  the  de Rham complex 
    \[
    K^{\bullet}:= \mscr{E}_{dR},
    \]
    with filtration 
    \[
    F^i=N^i_{dR}, 
    \]
    and $T$ being the functor $h_*$.
    Then the spectral sequence $E_{2}^{p,q}$ in \autoref{defn:conjss} agrees with the conjugate spectral sequence of Ogus--Vologodsky, as defined in \cite[text immediately after (3.18.2)]{ogus-vologodsky}, after a change of indices. \footnote{Indeed, the conjugate spectral sequence of loc.cit. is an $E_1$-spectral sequence, so one must re-index to make it an $E_2$ one so that it has a chance of agreeing with that given by \autoref{defn:conjss}.}
    
    We spell out this identification for the reader's convenience. Indeed,  the conjugate spectral sequence of loc.cit. is obtained by taking the filtered complex $Rh_*(K^{\bullet}, F)$, and then taking the spectral sequence of this filtered complex. By \cite[Equation (2)]{paranjape}, we see that Ogus--Vologodsky's conjugate spectral sequence agrees with that defined in \autoref{defn:conjss}.
    
   Suppose $(G,\nabla)=C^{-1}(G', \theta)$, where $(G', \theta)$ is a graded Higgs bundle, and $F_{conj}$ the filtration on $G$ defined  as in \autoref{subsection:inverse-Cartier}. In this case \autoref{prop:higgs-field-vs-p-curvature} (or \cite[Corollary 5.1.1]{ogus2004higgs} when $D=\emptyset$, and \cite[Lemma~5.5]{schepler2005logarithmic} in general), tells us that there is a commutative diagram whose horizontal arrows are isomorphisms:
 \begin{equation} 
 \adjustbox{scale=.9,center}{%
    \begin{tikzcd}
{\mscr{E}}_{1}^{p,q} \arrow[r, "\cong"] \arrow[d] & G'_{-p}\otimes \Omega^{p+q}_{\Xbar'}(\log D') \arrow[r] \arrow[d, "\theta_{-p}"]
& \mc{H}^{p+q}(\on{gr}^p_{F_{conj}}(G)_{dR}) \arrow[d]  & \on{gr}^{p}_{F_{conj}}(G)^{\nabla}\otimes \Omega^{p+q}_{\Xbar'}(\log D') \arrow{l}[swap]{C^{-1}_X} \arrow{d}{-\bar{\psi}} \\
{\mscr{E}_{1}}^{p+1,q} \arrow[r, "\cong"] & G'_{-p-1}\otimes \Omega^{p+q+1}_{\Xbar'}(\log D') \arrow[r]
&  \mc{H}^{p+q+1}_{dR}(\on{gr}^{p+1}_{F_{conj}}(G)_{dR})  &\on{gr}^{p+1}_{F_{conj}}(G)^{\nabla}\otimes \Omega_{\Xbar'}^{p+q+1}(\log D') \arrow{l}[swap]{C^{-1}_X}.
\end{tikzcd}
}
\end{equation}
In the above, $\bar{\psi}$ denotes the associated graded of the map induced by the $p$-curvature.

\end{remark}

\begin{proposition}\label{prop:conjugate-ss-degenerates-e2}
	Suppose $(\mathscr{E},\nabla)=C^{-1}(\mathscr{E}', \theta)$, with $(\mathscr{E}', \theta)$ a graded Higgs bundle, and equip $(\mathscr{E},\nabla)$ with the conjugate filtration as in \autoref{subsection:inverse-Cartier}. Suppose the length of $F_{\text{conj}}$ is less than $p-\dim_S(\Xbar)$, and that $h: \Xbar\to S$ is proper. Then the conjugate spectral sequence for $(\mathscr{E},\nabla)$ degenerates at $E_2$.
\end{proposition}
\begin{proof}
	In the case $D=\emptyset$ this is a special case of \cite[Corollary 3.25]{ogus-vologodsky} (without the assumption $h$ is proper, and after reindexing as above). In general (again without this assumption), it should follow identically from a filtered version of \cite[Corollary 5.7]{schepler2005logarithmic} or a logarithmic version of \cite[Theorem 3.22]{ogus-vologodsky}; as neither of these have been written down, we give a different argument in the case $h$ is proper, using that the modules on the $E_\infty$ page are coherent.
	
	Without loss of generality $S$ is an Artin $k$-scheme. By \cite[Corollary 5.7]{schepler2005logarithmic}, there is a quasi-isomorphism of complexes of $\mathscr{O}_{\Xbar'}$-modules: 
	\begin{equation}\label{eqn:schepler-qi}
	 (\mathscr{E}',\theta)_\text{Higgs}\simeq (\mathscr{E}, \nabla)_{dR},\end{equation}
	 where $(\mathscr{E}',\theta)_\text{Higgs}$ is the Higgs complex of \autoref{subsection:higgs-bundles}. The filtration on $(\mathscr{E}', \theta)_{\text{Higgs}}$ of \autoref{subsection:inverse-Cartier} yields a spectral sequence $$Rh_*^{p+q}(\on{gr}^p(\mathscr{E}',\theta)_\text{Higgs})\implies Rh_*^{p+q}(\mathscr{E}',\theta)_\text{Higgs}.$$ This spectral sequence degenerates at $E_2$ as the filtration on $(\mathscr{E}',\theta)_\text{Higgs}$ splits by assumption.
	
	On the other hand, the conjugate spectral sequence for $(\mathscr{E}', \theta)$ has, by \autoref{remark:conj-ss-paranjape} and \eqref{eqn:schepler-qi} the same terms on the $E_2$-page and converges to the same thing. Thus by comparing lengths of the modules on the $E_\infty$-pages, all differentials after $E_2$ must vanish.
\end{proof}

\section{Deformations of Hodge filtrations}\label{section:deformation-of-hodge}

\subsection{The main technical result}\label{subsection:main-technical-result}
Our goal in this section and the next is to show that if a Picard-Fuchs equation admits an integral (resp.~$\omega(p)$-integral) isomonodromic deformation, then the  obstruction to deforming the Hodge filtration to the formal isomonodromic deformation in characteristic zero, in a Griffiths-transversal manner, vanishes. We perform the main computations for the argument in this section, and finish it in  \autoref{section:smoothness}.

\begin{setup}\label{setup: induction-base-case}
Throughout this section we fix a finitely-generated $\mathbb{Z}$-algebra $R\subset\mathbb{C}$, a smooth $R$-scheme $S$, a smooth projective map $\overline{f}: \overline{X}\to S$, $D\subset \overline{X}$ a simple normal crossings divisor over $S$, and $s\in S(R)$. We let $(\mathscr{E}, \nabla)$ be a Picard-Fuchs equation on $(\overline{X}, D)_s$ (defined as in \autoref{subsubsection:NAGM}), and assumed to have \emph{nilpotent residues} along $D$. 
\end{setup}

Let $\widehat S, \widehat{\overline{X}}, \widehat D$ be the formal schemes obtained by completing $S, \overline{X}, D$ at $s$, and assume $(\mathscr{E}, \nabla)$ admits an integral (resp.~$\omega(p)$-integral) isomonodromic deformation to $(\widehat{\overline{X}}, \widehat D)$. After replacing $R$ with a finitely-generated localization we can and do assume $\widehat S\simeq R[[t_1, \cdots, t_{\dim S}]]$. 

We append a subscript $\mathbb{Q}$ to indicate the formal scheme obtained by base change to $R_{\mathbb{Q}}$ and completion at $s_{\mathbb{Q}}$. We denote by $S_n, \overline{X}_n$ etc.~ the $n$-th order neighborhood of $s$, i.e.~the subscheme cut out by the $n$-th power of the ideal sheaf of $s$, and similarly with $S_{n,\mathbb{Q}}$ etc. Let $\mathscr{K}$ be the fraction field of $R$.

\begin{remark}
The assumption that $(\mathscr{E},\nabla)$ has nilpotent residues is probably unnecessary for many of our arguments, though it makes things notationally simpler; the literature on non-abelian Hodge theory in positive characteristic is more complete under this assumption, which is why we make it. For the proof of our main theorem, \autoref{thm:NA-main}, we will reduce to the case of nilpotent residues.
\end{remark}

\begin{theorem}\label{theorem:extending-the-Hodge-filtration}We use notation as in \autoref{setup: induction-base-case}. 
The Hodge filtration on 	$(\mathscr{E}, \nabla)|_{\overline{X}_{s, \mathscr{K}}}$ extends Griffiths transversally to an isomonodromic deformation of $(\mathscr{E}, \nabla)|_{\overline{X}_{s, \mathscr{K}}}$ to $(\widehat{\overline{X}}, \widehat D)_{\mathscr{K}}/\widehat{S}_{\mathscr{K}}$.
\end{theorem}

The main point will be to inductively extend the Hodge filtration to $n$-th order neighborhoods of $s$, by comparing the (Frobenius twist of the) mod $p$ obstruction to the existence of a Griffiths-transverse deformation of the Hodge filtration to the $p$-curvature of the isomonodromy foliation in characteristic $p$. The latter will vanish to order $\omega(p)$ by the assumption of \mbox{($\omega(p)$-)}integrality, so the former will vanish modulo $p$ to order $\omega(p)/p$. As this latter quantity tends to infinity with $p$ by definition of $\omega(p)$-integrality, we will see that the obstruction to a Griffiths-transverse deformation of the Hodge filtration vanishes in characteristic zero. 

We complete the proof of \autoref{theorem:extending-the-Hodge-filtration} in \autoref{subsection:deforming-last-step}.

At this point we suggest that the reader interested in the precise details of the proof read the Appendices; for the reader more interested in the general idea, we suggest referring to the summary of \autoref{appendix:NAGM}, namely \autoref{subsection:appendix-summary}, as necessary.

\subsection{Base case of induction}\label{subsection:base-case}
We let $F_{\text{Hodge}}^\bullet$ denote the Hodge filtration on $\mathscr{E}$, and $\on{gr}_{F^\bullet_{\text{Hodge}}}(\mathscr{E})$ the associated graded bundle. This bundle is equipped with a Higgs field $\theta$ given as the associated graded of $\nabla$, which has degree $-1$ by Griffiths transversality.  Our goal in this section is to prove:

\begin{lemma}\label{lemma:ksvanishing}
We keep the notation of \autoref{setup: induction-base-case}. After possibly replacing $R$ by $R[\frac{1}{f}]$ for some $f\in R$, the map 
    \[KS: H^1(T_{\overline{X}_s}(-\log D_{s}))\overset{H^1(\theta)}{\longrightarrow} H^1_{\text{Higgs}}(\overline{X}_{s}, \on{End}(\on{gr}_{F_{\text{Hodge}}^\bullet}({\mathscr{E}}),\theta)/F^0)\]
    vanishes on the image of the Kodaira-Spencer map $T_{S,s}\rightarrow H^1(T_{\overline{X}_s}(-\log D_s))$.
\end{lemma}
We will conclude in \autoref{cor:base-case-Hodge-extends} from this (using degeneration of the Hodge-de Rham spectral sequence) that the Hodge filtration on $\mathscr{E}$ extends Griffiths-transversally to the isomonodromic deformation of $(\mathscr{E},\nabla)$ over the first-order neighborhood of $s$. Strictly speaking this step is contained in the induction step of \autoref{subsection:induction-step}, but it is simpler in a number of ways, so we include it as a warm-up, since \autoref{subsection:induction-step} is mathematically (and notationally) much more involved.

Let $\mathfrak{p}$ be a closed point of $\on{Spec}(R)$ of residue characteristic $p>0$, and $s_{\mf p}$ the reduction of $s$ modulo $\mf{p}$. We use the subscript $s_{\mf{p}}$ to denote objects obtained by basechanging along $s_{\mf{p}}$, so for example $\overline{X}_{s_{\mf{p}}}$ is obtained by basechanging $\overline{X}_{s}\rightarrow \Spec(R)$ along $s_{\mf{p}}$. 
\begin{proposition}\label{prop:filtrationcomplexes}
	We have a natural map of de Rham complexes, induced by the $p$-curvature map $\psi_p$
 \[
 F^*_{\text{abs}}T_{\overline{X}_{s_{\mf p}}}(-\log D_{s_{\mf p}})_{dR} \rightarrow \End(\mscr{E})_{dR}.
 \]
 where here $ F^*_{\text{abs}}T_{\overline{X}_{s_{\mf p}}}(-\log D_{s_{\mf p}})$ is given the canonical (Frobenius-pullback) connection.
 
 Moreover, this map is compatible with  filtrations $\tau^{i}$ and $F^{i}_{conj}$ on these complexes, defined  as follows:
 {\footnotesize \begin{align*}
 &\ \ \ \ \ \ \ \ \ \ \ \ \ \ \ \ \  \cdots &  \ \ \ \ \ \ \ \ \ \ \ &\ \ \ \ \ \ \ \ \ \ \ \ \ \ \ \ \ \cdots \\
&\tau^{-1} := [F^*T_{\Xbar_{s_{\mf p}}} \rightarrow F^*T_{\Xbar_{s_{\mf p}}} \otimes \Omega_{\Xbar_{s_{\mf p}}}^1(\log D_{s_{\mf p}}) \rightarrow \ker(d)\rightarrow \cdots  ] &\longrightarrow \ \ \ \ \ \ \ \ \ \ \   &F^{-1}_{conj}\End(\mscr{E})_{dR}  \\
&\tau^{0} := [F^*T_{\Xbar_{s_{\mf p}}} \rightarrow \ker(d) \rightarrow 0] &\longrightarrow \ \ \ \ \ \ \ \ \ \ \ &F^{0}_{conj}\End(\mscr{E})_{dR} \\
&\tau^{1}=[\ker(d)\rightarrow 0] &\longrightarrow \ \ \ \ \ \ \ \ \ \ \ &F^{1}_{conj}\End(\mscr{E})_{dR}\\ 
&\ \ \ \ \ \ \ \ \ \ \ \ \ \ \ \ \ \cdots & \ \ \ \ \ \ \ \ \ \ \ &\ \ \ \ \ \ \ \ \ \ \ \ \ \ \ \ \ \cdots 
\end{align*}}
\end{proposition}
\begin{proof}
	That the map is a map of complexes follows from the fact that the $p$-curvature map $$\psi_p: F^*_{\text{abs}}T_{\overline{X}_{s_{\mf p}}}(-\log D_{s_{\mf p}}) \rightarrow \End(\mscr{E})$$ is flat \cite[5.2.3]{katz1970nilpotent}. That it is compatible with the given filtrations is immediate from the fact that the associated graded of $\on{gr}_{F_{\text{conj}}}\mathscr{E}$ has vanishing $p$-curvature.
\end{proof}

\begin{definition}
    For the filtered complex $(\End(\mscr{E})_{dR}, F_{conj}^{\bullet})$, we denote the terms in the two spectral sequences from \autoref{section:recallss} simply by $\mscr{E}'^{p,q}_1$ and $E_2'^{p,q}$, respectively. Similarly, for the filtered complex $(F^*T_{\Xbar_{s_{\mf p}}}(-\log D_{s_{\mf p}})_{dR}, \tau^{\bullet})$ (with the filtration as in \autoref{prop:filtrationcomplexes}), we obtain the terms of their respective spectral sequences $\mscr{E}^{p,q}_1, E_{2}^{p,q}$.  
\end{definition}
We now describe these spectral sequences explicitly. 

\begin{proposition}\label{prop:sspages} 
%\item 
   % The map of complexes induces a map on spectral sequences $E_1^{pq}\rightarrow E_1'^{pq}$, where
    %\[E_1'^{pq} =  H^{p+q}_{dR}(\on{gr}_{conj}^p\End(\mscr{E}) \otimes \Omega^{\bullet}).\]

    %Moreover, the complex $gr_{conj}^p\mscr{E} \otimes \Omega^{\bullet}$ is formal; more precisely, it is quasi-isomorphic (non-canonically)  to 
    %\[
    %\on{gr}_{conj}^p\End(\mscr{E})^{\nabla} \oplus  \on{gr}_{conj}^p\End(\mscr{E})^{\nabla} \otimes \Omega^1[-1] \oplus \cdots .
    %\]
    Let $\mathscr{E}'$ be the pullback of $\mathscr{E}$ to the Frobenius twist $\overline{X}_{s_{\mf p}}'$ of $\overline{X}_{s_{\mf p}}$.
    \begin{enumerate}
  \item   We have  
   
    \begin{align*}
    \mscr{E}_{1}'^{p,q}&= \on{gr}^{-p}_{F_{conj}}\End(\mscr{E})_{dR}\\
    E_2'^{p,q} &= \mb{H}^{p+q}(\overline{X}_{s_{\mf p}}', \on{gr}^{q}_{F_{Hodge}}(\End(\mscr{E}')_{dR}))\\
    \end{align*}
Note that, when $(p,q)=(2,-1)$, we have 
\[E_{2}'^{p,q} = \mb{H}^{1}(\overline{X}_{s_{\mf p}}', \on{gr}^{-1}_{F_{Hodge}}(\End(\mscr{E}')_{dR})). \] 

If the Hodge-de Rham spectral sequence for $\End{\mscr{E}}$ (where the de Rham complex of the latter is equipped with the Hodge filtration) degenerates at the $E_1$-page, then the conjugate spectral sequence
\[E_{2}'^{p,q} \Rightarrow \mb{H}^{p+q}(\End(\mscr{E})_{dR}) \]
degenerates at the $E_2$-page. 
    \item We have 
    \begin{equation}\mscr{E}_1^{p,q} = \begin{cases}
    T_{\Xbar'_{s_{\mf p}}}(-\log D_{s_{\mf p}}')\otimes \Omega^{-p+1}_{\Xbar'_{s_{\mf p}}}(\log D_{s_{\mf p}}') \ \ \text{if} \ \ p+q= -p+1, \\
    0 \ \ \ \ \text{otherwise}.
    \end{cases}
    \end{equation}

    We have 
    \begin{align*}E_{2}^{p,q} &= \mb{H}^{p+q}(T_{\Xbar'_{s_{\mf p}}}(-\log D_{s_{\mf p}}')\otimes \Omega_{\Xbar'_{s_{\mf p}}}^{-(p+k)+1}(\log D_{s_{\mf p}}')[-(p+k+q)])\\
    &= H^{-k}(T_{\Xbar'_{s_{\mf p}}}(-\log D_{s_{\mf p}}')\otimes \Omega_{\Xbar'_{s_{\mf p}}}^{-(p+k)+1}(\log D_{s_{\mf p}}'),
    \end{align*}
    where $k$ is defined by the equality $2(p+k)=1-q$, and our convention is that this group vanishes when $q$ is even. Explicitly, we have 
    \[E_{2}^{p,q} = H^{-p+(1-q)/2}(T_{\Xbar'_{s_{\mf p}}}(-\log D_{s_{\mf p}}')\otimes \Omega_{\Xbar'_{s_{\mf p}}}^{-(p+k)+1}(\log D_{s_{\mf p}}').\]
    Note that, when $(p,q)=(2,-1)$, we have  $E_{2}^{p,q} = H^1(T_{\Xbar'_{s_{\mf p}}}(-\log D_{s_{\mf p}}))$.

    \item the map  
    \[T_{\Xbar'_{s_{\mf p}}}(-\log D_{s_{\mf p}})=\mscr{E}^{1,-1}_1\rightarrow \mscr{E}_{1}'^{1,-1} = \on{gr}^{-1}_{F_{Hodge}}\End(\mscr{E}')\]
    is given by the Higgs field $\theta'$ of $\on{gr}^{\bullet}_{F_{Hodge}}\on{End}(\mathscr{E}')$.
    \end{enumerate}
\end{proposition}

\begin{proof}
    This is essentially a routine computation using the definitions in \autoref{section:recallss}, for which we now provide the details. 

    For (1), we have 
    \begin{align*}
    \mscr{E}_{1}'^{p,q} &= \mc{H}^{p+q}(\on{gr}^p_{conj}(\End(\mscr{E})_{dR})) \\
    &= (\on{gr}^p_{conj}(\End(\mscr{E})))^{\nabla} \otimes \Omega^{p+q}_{\overline{X}_{s_{\mf p}}'}(\log D_{s_{\mf p}}') \\
    &= \on{gr}_{F_{Hodge}}^{-p} \End(\mscr{E}') \otimes \Omega^{p+q}_{\overline{X}_{s_{\mf p}}'}(\log D_{s_{\mf p}}')
    \end{align*}
    by Cartier theory, using that the associated graded of the conjugate filtration has zero $p$-curvature, and its Frobenius descent is identified with the associated graded of the Hodge filtration; see \autoref{subsection:conjugate-preliminaries}. 
    
Following \autoref{defn:conjss}, $E_{2}'^{p,q}$    is given by the hypercohomology $\mb{H}^{p+q}$ of the complex
\[[\cdots \rightarrow \on{gr}^{-p+1}\End(\mscr{E}')\otimes \Omega^{p+q-1} \rightarrow \on{gr}^{-p}\End(\mscr{E}')\otimes \Omega^{p+q} \rightarrow \on{gr}^{-p-1}\End(\mscr{E}')\otimes \Omega^{p+q+1} \rightarrow \cdots ] \]
with the middle term placed in degree $p+q$ (here we suppress subscripts and $(\log D)$ for space reasons; the associated graded is with respect to the Hodge filtration). By \autoref{remark:conj-ss-paranjape}, the differential is given by the Higgs field $\theta'$. This is precisely the complex $\on{gr}^{q}_{F_{Hodge}}(\End(\mscr{E}')_{dR})$, as claimed, i.e.~the $q$-th graded component of the Higgs complex of $\on{gr}_{F_{Hodge}}\on{End}(\mathscr{E}')$, appropriately shifted.

The statement about degeneration of spectral sequences follows from the same statement in \autoref{theorem:katz-big-diagram}.
%note that the Hodge-de Rham spectral sequence has as $E_1$-page
%\[E_{1}^{p,q}\defeq \mb{H}^{p+q}(\gr^{p}_F(\End((\mscr{E})\otimes \Omega^{\bullet}))).\]
%If this degenerates at the $E_1$-page, then for each $n$,
%\begin{align*}\dim E_{2}^{p,q}&= \bigoplus_{p+q=n} \dim \mb{H}^{p+q}(\gr^{q}_F(\End(\mscr{E})\otimes \Omega^{\bullet}))\\
%&= \bigoplus_{p+q=n} \dim \mb{H}^{p+q}(\gr^{p}_F(\End(\mscr{E})\otimes \Omega^{\bullet}))\\
%&= \dim \mb{H}^1(\End(\mscr{E})\otimes \Omega^{\bullet}),
%\end{align*}
%and hence the conjugate spectral sequence degenerates at the $E_2$-page. 

For (2), by definition, we have $\mscr{E}_{1}^{p,q} = \mc{H}^{p+q}(\on{gr}^p_{\tau}( F^*_{\text{abs}}T_{\overline{X}_{s_{\mf p}}}(-\log D_{s_{\mf p}})_{dR}))$. From the definition of $\tau$, we see that  $\on{gr}^p_{\tau}( F^*_{\text{abs}}T_{\overline{X}_{s_{\mf p}}}(-\log D_{s_{\mf p}})_{dR})$ has cohomology only in degree $-p+1$, and hence $\mscr{E}_{1}^{p,q}$ vanishes unless $p+q=-p+1$. When $p+q=-p+1$, we see that 
\begin{align*}
\mc{H}^{p+q}(\on{gr}^p_{\tau}( F^*_{\text{abs}}T_{\overline{X}_{s_{\mf p}}}(-\log D_{s_{\mf p}})_{dR})) &= \mc{H}^{-p+1}( F^*_{\text{abs}}T_{\overline{X}_{s_{\mf p}}}(-\log D_{s_{\mf p}})_{dR})\\
&=  T'_{\overline{X}_{s_{\mf p}}}(-\log D'_{s_{\mf p}}) \otimes\Omega^{-p+1}_{\overline{X}_{s_{\mf p}}}(-\log D'_{s_{\mf p}}),
\end{align*}
as required. Finally, we compute $E_{2}^{p,q}$; this is the hypercohomology of the complex 
\[[\cdots \rightarrow \mscr{E}_{1}^{p-1,q} \rightarrow  \mscr{E}_{1}^{p,q} \rightarrow \mscr{E}_{1}^{p+1,q} \rightarrow \cdots ]\]
with the middle term placed in degree $p+q$. There is at most one non-zero term in this complex, which is the one in degree $p+k+q$, where $k$ satisfies
\[p+k+q=-(p+k)+1.\]
In other words, 
\begin{align*}
    E_{2}^{p,q} &= \mb{H}^{p+q}(T_{\Xbar'_{s_{\mf p}}}(-\log D_{s_{\mf p}}')\otimes \Omega_{\Xbar'_{s_{\mf p}}}^{-(p+k)+1}(\log D_{s_{\mf p}}'[-(p+k+q)]) \\
    &= H^{-k}(T_{\Xbar'_{s_{\mf p}}}(-\log D_{s_{\mf p}}')\otimes \Omega_{\Xbar'_{s_{\mf p}}}^{-(p+k)+1}(\log D_{s_{\mf p}}'),
\end{align*}
as claimed. Part (3) follows from Katz's theorem (\autoref{theorem:katz-big-diagram}) comparing the associated graded of the $p$-curvature to the Higgs field $\theta'$. 
\end{proof}

The following is a straightforward consequence of the Griffiths transversality of $F_{Hodge}$ on $\mscr{E}$.
\begin{proposition}\label{prop:ksfactorize}
    We have the following natural map of complexes, induced by $T_{\overline{X}_s}(-\log D_s)\rightarrow \gr^{-1}_{F_{Hodge}}(\End(\mscr{E}))$:
    \[[T_{\overline{X}_s}(-\log D_s)\rightarrow 0] \longrightarrow \gr^{-1}_{F_{Hodge}}(\End(\mscr{E})_{dR}).\]
   The map $KS$ of \autoref{lemma:ksvanishing} factors through the induced map $H^1(T_{\overline{X}_s}(-\log D_s)) \rightarrow \mb{H}^1(\on{gr}_{F_{Hodge}}^{-1}(\End(\mscr{E})_{dR}))$.
\end{proposition}
\begin{proof}[Proof of \autoref{lemma:ksvanishing}]
    Let $\kappa \in H^1(T_{\overline{X}_s}(-\log D_s))$ be a class in the image of $T_{S,s}$. Let $\kappa_p\in H^1(T_{\overline{X}_{s_{\mf p}}}(-\log D_{s_{\mf p}}))$ be the reduction of $\kappa$ modulo $\mf p$. We will show that $KS(\kappa_p)=0$ for $p$ sufficiently large, which will prove the lemma. More precisely, we take $p$ large enough so that the Hodge-de Rham spectral sequence for $\End(\mscr{E})$ degenerates at the first page, and hence the same is true for the conjugate spectral sequence at the second page. Finally, we may replace $R$ with a further finitely-generated localization so that the cohomological condition of \autoref{cor:unique-integral-isomonodromic} is satisfied, i.e.~$H^1_{dR}(\on{End}(\mathscr{E}_0, \nabla_0))$ is locally free and  its formation  commutes with arbitrary base change. 
    
    Let $Z_{\mf p}=F_{\text{abs}}^{-1}(s_{\mf p})\subset S_{\mf p}$.  We also assume that the given ($\omega(p)$-)integral isomonodromic deformation of $(\mathscr{E}, \nabla)$ is defined over $Z_{\mf p}$, which we may do for $p$ sufficiently large by the definition of ($\omega(p)$-)integrality. Let $(\mathscr{G}, \nabla)$ be the corresponding flat bundle on $(\overline{X}, D)_{Z_{\mf p}}/Z_{\mf p}$. By \autoref{cor:vanishing-p-curvature-for-p-power-leaves-isomonodromy}, the $p$-curvature of $\mathscr{M}_{dR}(X/S)^\natural$ vanishes when restricted to $Z_{\mf p}$ (along the $Z_{\mf p}$-point of $\mathscr{M}_{dR}(X/S)^\natural$ corresponding to $[(\mathscr{G}, \nabla)]$). Note that the hypotheses of \autoref{cor:vanishing-p-curvature-for-p-power-leaves-isomonodromy} are again satisfied by the assumption of ($\omega(p)$-)integrality.\footnote{Recall that we choose coordinates $t_1, \cdots, t_{\dim S}$ on $\widehat S$. In applying \autoref{cor:vanishing-p-curvature-for-p-power-leaves-isomonodromy} we may take $Z=\widehat{R}[[t_1, \cdots, t_{\dim S}]]/(t_1^p, \cdots, t_{\dim S}^p)$ and $m=p\dim S$.}

    Consider the map on the conjugate spectral sequences (as defined in \autoref{defn:conjss}) associated with the map of filtered complexes
    \[(F^*_{\text{abs}}T_{\overline{X}_{s_{\mf p}}}(-\log D_{s_{\mf p}})_{dR}, \tau^{\bullet})  \rightarrow (\End(\mscr{E})_{dR}, F^{\bullet}_{conj}).\]

    The terms of this spectral sequence are computed in \autoref{prop:sspages}; specifically, for $(p,q)=(2,-1)$, we have a map
    \begin{equation}\label{eqn:maponrelevantterms}
    E_{2}^{2,-1} = H^1(T_{\Xbar'_{s_{\mf p}}}(-\log D'_{s_{\mf p}})) \rightarrow \mb{H}^{1}(\overline{X}_{s_{\mf p}}', \on{gr}^{-1}_{F_{Hodge}}(\End(\mscr{E}')_{dR})) = E_{2}'^{2,-1}. \end{equation}

    Note that, by construction,  this map is obtained by taking $\mb{H}^1$ of the map of complexes
\[[\cdots \rightarrow \mscr{E}_{1}^{1,-1} \rightarrow \mscr{E}_{1}^{2,-1} \rightarrow \cdots] \longrightarrow  [\cdots \rightarrow \mscr{E}_{1}'^{1,-1} \rightarrow \mscr{E}_{1}'^{2,-1} \rightarrow \cdots] \]
    which, by \autoref{prop:sspages}, is precisely the map of complexes
\[[T_{\overline{X}'_{s_{\mf p}}}(-\log D'_{s_{\mf p}})\rightarrow 0 \rightarrow \cdots] \longrightarrow \on{gr}^{-1}_{F_{Hodge}}(\End(\mscr{E}')_{dR})  \]
induced by the Higgs field $T_{\overline{X}'_{s_{\mf p}}}(-\log D'_{s_{\mf p}})\rightarrow \on{gr}^{-1}_{F_{Hodge}}\End(\mscr{E}').$

By our assumption, the spectral sequence $(E_{2}'^{p,q})$ degenerates at the $E_2$-page, and hence $E_{2}'^{p,q}=E_{\infty}'^{p,q}$ for all $p,q$. We claim that the map  
\begin{equation}\label{eqn:map-of-es}
	E_{\infty}^{2,-1} \rightarrow E_{\infty}'^{2,-1}
\end{equation}
vanishes on $\kappa_p$. Since the spectral sequence $(E_{2}^{p,q})$  also degenerates at the second page, we conclude that $\kappa_p$ is sent to zero under the map \eqref{eqn:maponrelevantterms}.
%\[
% H^1(T_{X'}) \rightarrow \mb{H}^1(\on{gr}_F^{-1}(\End(\mscr{E}')\otimes \Omega^{\bullet}_X)).
%\]
By \autoref{prop:ksfactorize}, $KS$ factors  through this map, and hence $KS(\kappa_p)=0$ as required. 

It remains to show that the map  \eqref{eqn:map-of-es}
vanishes on $\kappa_p$. Observe that by \autoref{lemma:cocycles-agree}, this map is closely related to the $p$-curvature of the isomonodromy foliation, as discussed in the Appendices, as we now explain. The $p$-curvature of the isomonodromy foliation restricted to $s_{\mf p}$, $\Psi_p|_{s_{\mf p}}$, is by \autoref{lemma:cocycles-agree}  precisely 
\begin{equation}\label{eqn:p-curvature-at-a-point}
 	\Psi_p|_{s_\mf p}: T_{S_{\mf p}, s_{\mf p}}\to H^1(T_{\overline{X}'_{s_{\mf p}}}(-\log D'_{s_{\mf p}})) \rightarrow H^1(F_{\text{abs}}^*T_{\overline{X}'_{s_{\mf p}}}(-\log D'_{s_{\mf p}})_{dR})\overset{H^1(\psi_p)}{\longrightarrow} \mb{H}^1(\End(\mscr{E})_{dR})
 \end{equation}
where the first map is the Kodaira-Spencer map, the second is induced by the natural inclusion and the third is the one studied in \autoref{prop:filtrationcomplexes}. By that proposition and degeneration of the conjugate spectral sequence at $E_2$, the map \eqref{eqn:p-curvature-at-a-point} factors through $F^1_{conj}\mb{H}^1(\End(\mscr{E})_{dR})$, and the map \eqref{eqn:map-of-es} is the induced map 
\begin{align*} H^1(T_{\Xbar'_{s_{\mf p}}}(-\log D'_{s_{\mf p}}))&\overset{	\Psi_p|_{s_\mf p}}{\longrightarrow} F^1_{conj}\mb{H}^1(\End(\mscr{E})_{dR})\\
&\twoheadrightarrow \on{gr}^1_{F_{conj}}\mb{H}^1(\End(\mscr{E})_{dR})\\
&\simeq \mb{H}^1(\on{gr}^1_{F_{conj}}(\End(\mscr{E})_{dR}))\\
&\simeq \mb{H}^1(\on{gr}^{-1}_{F_{Hodge}}\End(\mscr{E}'))
\end{align*}
It thus suffices to show that the map \eqref{eqn:p-curvature-at-a-point} vanishes, or in other words that the $p$-curvature of the isomonodromy foliation vanishes at $[(\mathscr{E}, \nabla)]_{s_{\mf p}}$. But this is immediate from what we arranged at the start of the proof, i.e.~that the $p$-curvature of the isomonodromy foliation vanishes when restricted to $Z_{\mf p}$.
\end{proof}

\begin{corollary}\label{cor:base-case-Hodge-extends}
We keep the notation and assumptions of \autoref{setup: induction-base-case}.	Then the Hodge filtration on $(\mathscr{E}, \nabla)$ extends to a Griffiths-transverse filtration on the isomonodromic deformation of $(\mathscr{E}, \nabla)$ to $(\overline{X}_{2,\mathbb{Q}}, D_{2, \mathbb{Q}})/S_{2, \mathbb{Q}}$, where, as in \autoref{setup: induction-base-case}, $S_{n, \mb{Q}}$ denotes the $n$-order neighborhood of $s_{\mb{Q}}$, and $\overline{X}_{2, \mb{Q}}, D_{2, \mb{Q}}$  denote the schemes obtained by basechanging $\overline{X}, D$ to $S_{2, \mb{Q}}$.
\end{corollary}
\begin{proof}
This is immediate from \autoref{prop: gt-transverse-obstruction} and the degeneration of the Hodge-de Rham spectral sequence, as we now explain. Indeed, we wish to show that the map $$T_{S,s}\to H^1(T_{\overline{X}_{s_\mb{Q}}}(-\log D_{s_\mb Q}))\to \mb{H}^1(\on{End}(\mscr{E})_{dR}/F_{Hodge}^0)$$ of \autoref{prop: gt-transverse-obstruction} vanishes. This map factors through $F^{-1}_{\text{Hodge}}$, so by degeneration of the Hodge-de Rham spectral sequence it suffices to show that the map $$T_{S,s}\to H^1(T_{\overline{X}_{s_\mb{Q}}}(-\log D_{s_\mb Q}))\to \mb{H}^1(\on{gr}^{-1}_{F_{\text{Hodge}}}\on{End}(\mscr{E})_{dR})$$ induced by the Higgs field vanishes. But this is precisely the content of \autoref{lemma:ksvanishing}.
\end{proof}

\subsection{Induction step}\label{subsection:induction-step}
We now proceed with the induction step, which is similar to the induction base case, although there are several new  complications. Our goal is to extend the Hodge filtration on $(\mathscr{E}, \nabla)$ to a Griffiths-transverse filtration on its isomonodromic deformation to $(\widehat{\overline{X}}_{\mathbb{Q}}, \widehat{D}_{\mathbb{Q}})/\widehat{S}_{\mathbb{Q}}$. In this section we will assume we have found such an extension of the Hodge filtration over $S_{n, \mathbb{Q}}$, and handle deformations over split extensions of $S_{n, \mathbb{Q}}$, i.e.~extensions of the form $S_{n, \mathbb{Q}}[\epsilon]/\epsilon^2$. We will deduce the analogous extension for all deformations in \autoref{section:smoothness} using some version of the $T^1$-lifting theorem. 

As in \autoref{subsection:base-case}, we will need to verify that the obstruction of \autoref{prop: gt-transverse-obstruction} vanishes. In \autoref{subsection:base-case} we did so by comparing it to an associated graded component of the $p$-curvature of the isomonodromy filtration, with respect to the conjugate filtration. However, for $n>1$ we are not given a conjugate filtration, since we are no longer in the setting of a Picard-Fuchs equation. We will thus have to construct one ourselves. We will do so using Ogus-Vologodsky's non-abelian Hodge theory, but this will cause an additional complication: it will not be clear anymore that the $p$-curvature of the isomonodromy foliation vanishes when restricted to the corresponding section of $\mathscr{M}_{dR}(X/S)^\natural$. We will prove this via an argument in mixed characteristic relying on the Higgs-de Rham flow.

Our setup is as follows. 
\subsubsection{Setup}\label{induction_setup}
    Recall that we have  a point $s\in S(R)$, and a flat bundle of Picard-Fuchs type $(\mscr{E}, \nabla)$ on $(\overline{X}_s, D_s)$ with nilpotent residues along $D_s$. We assume in addition that $(\mscr{E}, \nabla)$ has $\omega(p)$-integral isomonodromic deformation over $\widehat S$. Note that (possibly after replacing $R$ with some localization) we may assume $(\mscr{E}, \nabla)$ is semisimple, as it is of Picard-Fuchs type.

   The goal of this section is to prove the \enquote{induction step}, which is recorded in \autoref{lemma: ks-vanish-induction-step}.  %Namely, this means that,  given a deformation as above of $(\mscr{E}, \nabla, F^{\bullet})$ to $\X_{\fracr}$, there is no obstruction to further deforming along  a split thickening $\Spec(A_{\fracr}) \xhookrightarrow{}\Spec(A_{\fracr}\oplus M)$. 
   
   Our induction hypothesis is as follows. We have an affine closed subscheme $\on{Spec}(A)$ of $S$ with $\on{Spec}(A)_{\text{red}}=s$, such that
the Hodge filtration extends in a Griffiths transverse way to $(\overline{X}, D)_A$.  We denote by $(\mscr{E}_{A_{\mathbb{Q}}}, \nabla_{A_{\mathbb{Q}}})$ the flat bundle on $(\overline{X}, D)_{A_\mathbb{Q}}$ which is the isomonodromic deformation  of $(\mscr{E}, \nabla)$.  We also assume we are given  the Griffiths transverse filtration $F^{\bullet}_{A_{\mathbb{Q}}}$ on $(\mscr{E}_{A_{\mathbb{Q}}}, \nabla_{A_{\mathbb{Q}}})$  extending the given one on $(\mathscr{E},\nabla)$. The goal of \autoref{lemma: ks-vanish-induction-step} is to show that there is no obstruction to further deforming the Hodge filtration (in a Griffiths transverse way) to the isomonodromic deformation to $(\overline{X}, D)_{A_{\mathbb{Q}}[\epsilon]/\epsilon^2}$, where we are given a free $A$-module $M$ and a map  $\Spec(A[\epsilon]/\epsilon^2) \to{} S$. 

After replacing $R$ with a finitely-generated localization, we may assume $\widehat S=\on{Spf}(R[[t_1, \cdots, t_{\dim S}]])$; fix such a set of local coordinates.

The following will be convenient for us:
\begin{assumption}\label{assumption:monomial-A}
	In the local coordinates $t_1, \cdots, t_n$ chosen above, $A$ is cut out by a monomial ideal.
\end{assumption}

This assumption will be a consequence of the setup of our induction argument in \autoref{section:smoothness}. Note that as $\on{Spec}(A)_{\text{red}}=s$, we have that the ideal of $A$ is contains $(t_1, \cdots, t_{\dim S})^q$ for some $q$.

By our assumption of $\omega(p)$-integrality, we may by localizing $R$ assume that all the objects 
\[ 
\mscr{E}_{A_{\mathbb{Q}}}, \nabla_{A_{\mathbb{Q}}}, F^{\bullet}_{A_{\mathbb{Q}}}\] 
descend to $A$, i.e. there are corresponding objects
\[ 
\mscr{E}_{A}, \nabla_{A}, F^{\bullet}_{A}\] 
whose base change along $\Spec(A_{\mathbb{Q}}) \rightarrow \Spec(A)$ are the ones above. 
For each maximal ideal $\mf{p}\subset R$, with residue field $k$ of characteristic $p$, we have the $k$-algebra
\[
A_{\mf p}:= A /\mf{p}.
\]
%For sufficiently large $p$, we can spread out $A_0$ and reduce everything modulo $p$. 
%We denote the resulting Artin algebra by $A$, so that $\Spec(A) \subset S_p$ is a closed
%subscheme with closed point $s_p$. 
Base changing along $\Spec(A_{\mf p})\rightarrow \Spec(A)$ we obtain  $(\mscr{E}_{\mf p}, \nabla_{\mf p}, F^{\bullet}_{\mf p})$ on $\overline{X}_{A_\mf p}$. Let $\widehat{R}$, $A_{\widehat R}$, $S_{\widehat R}$ denote the completion of $R$, $A$, $S$ at $\mf{p}$, respectively. 

Since $(\mscr{E}, \nabla)_{\mathbb{Q}}$ is a semisimple flat bundle, we may write it as 
\[
(\mscr{E}, \nabla)_{\mathbb{Q}}=\bigoplus_i \mscr{F}_{\mathbb{Q}, i},
\]
where each $\mscr{F}_{\mathbb{Q}, i}$ is an irreducible flat bundle. As the $\mathscr{F}_{\mathbb{Q}, i}$ underly irreducible complex variations of Hodge structure they carry (unique up to reindexing) Griffiths transverse Hodge filtrations $F^\bullet_{j, \text{Hodge}}$. We choose this direct sum decomposition so the filtrations on the $\mathscr{F}_i$ are compatible with the one on $\mathscr{E}$. Explicitly (after possibly replacing $R$ with some finitely-generated extension), we may write $$(\mathscr{E}, \nabla)_{\mathbb{Q}}\simeq\bigoplus_j \mathscr{F}_{\mathbb{Q}, j}\otimes V_j,$$ where the $\mathscr{F}_{\mathbb{Q},j}$ are mutually non-isomorphic and $\overline{\mathscr{K}}$-simple. (Recall that $\mathscr{K}$ is the fraction field of $R$.) Moreover by the theorem of the fixed part, the $V_j$ carry constant Hodge filtrations making the isomorphism above one of filtered flat bundles; see e.g.~\cite[Proposition 4.1.4(2)]{landesman2024geometric}. Choosing a basis of the $V_j$ compatible with this filtration yields the desired direct sum decomposition.

Taking the isomonodromic deformations of these flat bundles, we obtain a decomposition 
\begin{equation}\label{eqn: iso-def-decomposition}
(\mscr{E}_{A_\mathbb{Q}}, \nabla_{A_\mathbb{Q}})=\bigoplus_i \mscr{F}_{A_\mathbb{Q}, i}.
\end{equation}
Replacing $R$ by a localization if necessary and applying \autoref{lemma:integral-isomonodromic-summand}, from now on we make the following:
\begin{assumption}\label{assumption: stable-special-fiber}
Each $\mscr{F}_{A_{\mathbb{Q}}, i}$ has an integral model $\mscr{F}_{A, i}$, so that we have a decomposition
\[
(\mscr{E}_{A}, \nabla_A)=\bigoplus_i \mscr{F}_{A, i},
\]
which is  compatible with \eqref{eqn: iso-def-decomposition}. Moreover, we assume that $(\mathscr{E}, \nabla)$ satisfies the assumption of \autoref{prop:special-fiber-char-0-comparison}, i.e.~that the $R$-module  $H^1_{dR}(\on{End}(\mathscr{E}, \nabla))$ is locally free, and its formation commutes with arbitrary base change. Finally, recall that $(\mathscr{E}, \nabla)$ is a Picard-Fuchs equation arising as the cohomology of some smooth projective $f: Y\to X_s$; we assume that  the modules $\mathscr{H}^i_{dR}(Y/X_s)$ are locally free and that their formation commutes with arbitrary base change. Finally we assume that all sufficiently small primes (in terms of $(\mathscr{E}, \nabla)$, $\dim S$, the dimension of $Y/X_s$, and the function $g$ appearing in the definition of $\omega(p)$-integrality) are invertible in $R$.
%\begin{itemize} 
%\item for each $i$, and each maximal ideal $\mf p$ of $R$ the restriction of   $\mscr{F}_{A, i}$ to $\overline{X}_{s_\mf p}$, which we denote by $\mscr{F}_{\mf p, i}$, is a stable flat bundle with resepct to some fixed choice of polarization on $\overline{X}$, and 
%\item for each $j$, the Higgs bundle  $(\gr_{F^{\bullet}_{j, \text{Hodge}}}\mscr{F}_j, \theta)|_{\overline{X}_{s_\mf p}}$ is stable. 
%\end{itemize}
\end{assumption}
\begin{remark}
The last assumptions of \autoref{assumption: stable-special-fiber} will allow us to apply Ogus-Vologodsky's non-abelian Hodge theory and the Fontaine-Laffaille theory of \autoref{thm:berthelot-mazur-ogus} freely.	We are always free to replace $R$ with a finitely-generated localization in the arguments here so the assumption is harmless.
\end{remark}

Finally, we assume:

\begin{assumption}\label{assumption: hodge-filtration}
	Each $\mathscr{F}_{A,i}$ is endowed with a Hodge filtration compatible with the Hodge filtration on $(\mathscr{E}_A, \nabla_A)$.
\end{assumption}

This will be automatic from the induction argument of \autoref{section:smoothness}.

%the same is true of $(\mscr{E}_{A_{\fracr}}, \nabla_{A_{\fracr}})[\frac{1}{p}]$, and therefore for $\mf{p}$ of sufficiently large characteristic that 
%\begin{center}
%     ($\star$) \ \ \ \ \ $(\mscr{E}_A, \nabla_A)$ is polystable.
%\end{center}
%We will assume this to be the case from now on.

Again after replacing $R$ by a finitely-generated localization, we may assume it is $\mathbb{Z}$-smooth. Let $\Spec(B_{\mf p})=F^{-1}_{\text{abs}}(\Spec(A_{\mf p})).$ As $S$ is smooth, we may fix a lift $F_0$ of absolute Frobenius 
in an affine open of $S_{\widehat R}$ containing $s_{\mf p}$, and let $\Spec(B_{\widehat R})=F^{-1}_0(\on{Spec}(A_{\widehat R}))$. (Note that despite the notation, $B_{\mf p}, B_{\widehat R}$ are not the base-change of some ring $R$ to $\mf p, \widehat R$, respectively.)
We let $\widetilde{A}, \widetilde{B}$ denote the mod $\mf{p}^2$ reduction of $A, B_{\widehat R}$ respectively. 
The reduction of $F_0$ mod $\mf{p}^2$ gives a map $\widetilde{F}: \Spec(\widetilde{B})\rightarrow \Spec(\widetilde{A})$ 
lifting the Frobenius map $F: \Spec(B_{\mf p})\rightarrow \Spec(A_{\mf p})$. 

The following assumption will be convenient for us to verify the hypotheses of \autoref{cor:unique-integral-isomonodromic}---in particular,  \autoref{assumption:torsion-free-ideals}. Recall \autoref{assumption:monomial-A} that $\Spec(A)\subset \hat{S}$ is cut out by a monomial ideal in the coordinates $t_1, \cdots, t_n$.
\begin{assumption}\label{assumption:monomial-frobenius}
	The map $F_0$ sends $t_i$ to $t_i^p$ for all $i$, so that $B_{\widehat R}$ is also cut out by a monomial ideal.
\end{assumption}
When we apply our arguments here in \autoref{section:smoothness}, the assumption on $A$ will always be satisfied; we choose $F_0$ now so that \autoref{assumption:monomial-frobenius} is satisfied. 

\begin{definition}\label{defn:lifting-and-conjugate-filtration}
    \begin{enumerate}\item 
    We define $\overline{X}_{\tilde{B}}'$ as the fiber product 
    \begin{equation} \begin{tikzcd}
       (\overline{X}_{\widetilde{B}}', D_{\widetilde{B}}') \arrow[r, "\widetilde{G}"] \arrow[d]
        & (\overline{X}_{\widetilde{A}}, D_{\widetilde A}) \arrow[d] \\
        \Spec(\widetilde{B}) \arrow[r, "\widetilde{F}"]
        &  \Spec(\widetilde{A}).
        \end{tikzcd}
    \end{equation}
 This is a $W_2$-lift of the left  square of the diagram
 \begin{equation} \label{eqn:reduced-situation-frobenius}\begin{tikzcd}
    (\overline{X}_{B_{\mf p}}', D_{B_{\mf p}}') \arrow[r, "G"] \arrow[d]
    &  (\overline{X}_{A_{\mf p}}, D_{A_{\mf p}}) \arrow[d] \arrow[r] & (\overline{X}_{B_{\mf p}}, D_{B_{\mf p}})  \arrow[d]\\
    \Spec(B_{\mf p}) \arrow[r, "F"]
    &  \Spec(A_{\mf p}) \arrow[r] & \Spec(B_{\mf p}),
    \end{tikzcd}
\end{equation}    
 and therefore  $ (\overline{X}_{\widetilde{B}}', D_{\widetilde{B}}') \rightarrow \Spec(\widetilde{B})$ is a $W_2$-lift of 
 $ (\overline{X}_{B_{\mf p}}', D_{B_{\mf p}}')  \rightarrow \Spec(B_{\mf p})$. Note that the bottom horizontal composite map $ \Spec(B_{\mf p})\to  \Spec(B_{\mf p})$ is the absolute Frobenius of  $\Spec(B_{\mf p})$, and the bottom right horizontal arrow is the natural closed embedding.
 \item Using the $W_2$-lift $\overline{X}_{\widetilde{B}}' \rightarrow \Spec(\widetilde{B})$ as above, 
 we define the following connection on $\overline{X}_{B_{\mf p}}$:
  \[
    (\mscr{E}_B, \nabla)\defeq C^{-1}(G^*(\on{gr}^\bullet_{F_{\text{Hodge}}}\mscr{E}_{A_{\mf p}}, \on{gr}^\bullet_{F_{\text{Hodge}}}(\nabla_{A_{\mf p}}))),
  \]
  where $C^{-1}$ is the inverse Cartier transform of \autoref{theorem:inverse-cartier-transform}. 
 The decreasing filtration on $(\on{gr}^\bullet_{F_{\text{Hodge}}}\mscr{E}_{A_{\mf p}}, \on{gr}^\bullet_{F_{\text{Hodge}}}(\nabla_{A_{\mf p}}))$ by the sub-Higgs bundles $$(\mathscr{G}_m, \theta_m)=\left(\bigoplus_{i\leq -m} \on{gr}^i_{F_{\text{Hodge}}}\mscr{E}_{A_{\mf p}}, \on{gr}^\bullet_{F_{\text{Hodge}}}(\nabla_{A_{\mf p}})\right)$$ gives rise (by applying $C^{-1}$ to $G^*(\mathscr{F}_m, \theta_m)$) to a 
 $\nabla$-stable decreasing filtration on $\mscr{E}_B$, which we denote by $F_{conj}$.

\end{enumerate}
     
\end{definition}

Our goal is to consider deformations of $\overline{X}_A$ to $A[\epsilon]/\epsilon^2$, and to compute the obstruction to deform the Hodge filtration on $\mscr{E}_A$. We will compare it to a certain graded component of the $p$-curvature of the isomonodromy foliation on the leaf associated to $(\mathscr{E}_B, \nabla)$, with respect to the conjugate filtration defined above.

\begin{definition}
\begin{enumerate}
    \item Equip   $F^*_{\text{abs}}T_{\overline{X}'_{B_{\mf p}}/B_{\mf p}}(-\log D_{B_{\mf p}})$ with its canonical (Frobenius pullback) connection. We may filter its de Rham complex as in \autoref{prop:filtrationcomplexes} with the canonical filtration, thus obtaining the filtered complex 
    \[
    (F^*_{\text{abs}}T_{\overline{X}'_{B_{\mf p}}/B_{\mf p}}(-\log D_{B_{\mf p}})_{dR}, \tau).
    \]
    \item Consider the complex $\End(\mscr{E}_B)_{dR}$. The flat bundle $\End(\mscr{E}_B)$ carries a conjugate filtration induced from the conjugate filtration on $\mscr{E}_B$ defined in \autoref{defn:lifting-and-conjugate-filtration}. 
  
    This in turn defines the filtered complex $(\End(\mscr{E}_B)_{dR}, F_{conj}^{})$.
\end{enumerate}
  
\end{definition}

As in \autoref{prop:filtrationcomplexes}, the $p$-curvature map $$\psi_p: F_{\text{abs}}^*T_{\overline{X}'_{B_{\mf p}}/B_{\mf p}}(-\log D_{B_{\mf p}})\to \on{End}(\mathscr{E}_B)$$ induces   a map of filtered complexes
\[
    (F^*_{\text{abs}}T_{\overline{X}'_{B_{\mf p}}/B_{\mf p}}(-\log D_{B_{\mf p}})_{dR}, \tau)
    \rightarrow (\End(\mscr{E}_B)_{dR}, F_{\text{conj}}). 
    \]
%which we may filter in exactly the same way as in \autoref{prop:filtrationcomplexes}.
%We therefore obtain the  filtered complexes 
%$(F_{\text{abs}}^*T_{\X_B/B} \otimes %\Omega^{\bullet}_{\X_B/B}, \tau)$ and 
%$(\End(\mscr{E}_B)\otimes \Omega^{\bullet}_{\X_B/B}, F_{conj})$. 

\begin{lemma}[Key lemma]\label{lemma:pcurv-vanish}
    We maintain  notation as in \autoref{induction_setup}.
    \begin{enumerate}
    \item 	There is an isomorphism $(\mathscr{E}_B, \nabla)|_{A_\mf p}\simeq (\mathscr{E}_{A_\mf p}, \nabla_{A_\mf p})$. That is, in the language of \cite[Definition 4.16]{ogus-vologodsky}, $(\mathscr{E}_{A_\mf p}, \nabla_{A_\mf p})$ equipped with its Hodge filtration  is a Fontaine module; in the language of \cite{lan2019semistable} it is a Higgs-de Rham fixed point.
    \item     Moreover the restriction of the $p$-curvature of the isomonodromy foliation on $\mathscr{M}_{dR}(X/S)^\natural$ to the $B_{\mf p}$-point of $\mathscr{M}_{dR}(X/S)^\natural$ corresponding to $[(\mscr{E}_B, \nabla)]$ vanishes. More precisely,  the map 
    \begin{equation}\label{eqn: h1-p-curv-Xb}
        H^1(T_{\overline{X}'_{B_{\mf p}}/B_{\mf p}}(-\log D'_{B_{\mf p}})) \to  \mb{H}^1(F^*_{\text{abs}}(T_{\overline{X}'_{B_{\mf p}}/B_{\mf p}}(-\log D_{B_{\mf p}}))_{dR}) \overset{\psi_p}{\to} 
        \mb{H}^1(\End(\mscr{E}_B)_{dR})
        \end{equation}
vanishes on the image of the composition 
\begin{equation}\label{eqn: frob-pullback-ks}
    F_{\text{abs}}^*T_{S, B_{\mf p}}\rightarrow F^*_{\text{abs}}H^1(T_{\overline{X}'_{B_{\mf p}}/B_{\mf p}}(-\log D_{B_{\mf p}})) \to   H^1(T_{\overline{X}'_{B_{\mf p}}/B_{\mf p}}(-\log D'_{B_{\mf p}})).
\end{equation}
    \end{enumerate}
    Here the first map in \eqref{eqn: h1-p-curv-Xb} is given by the natural inclusion of $T_{\overline{X}'_{B_{\mf p}}/B_{\mf p}}(-\log D_{B_{\mf p}})$ into the de Rham complex of $F^*_{\text{abs}}(T_{\overline{X}'_{B_{\mf p}}/B_{\mf p}}(-\log D_{B_{\mf p}}))$ (namely, the inclusion of flat sections), and the second is induced by the $p$-curvature of $(\mathscr{E}_B, \nabla)$.  The first map in \eqref{eqn: frob-pullback-ks} is the Frobenius pullback of the Kodaira-Spencer map, and the second is the base change morphism.

\end{lemma}
%\begin{remark}
%\begin{itemize}
%    \item The second map in \eqref{eqn: frob-pullback-ks} is the natural basechange morphism: more generally, for a ring $R$ and a complex $C^{\bullet}$ of $R$-modules, and a module $M$, there is a natural map
%    \[
%    H^i(C^{\bullet})\otimes M\rightarrow H^i(C^{\bullet}\otimes^{L}_R M).
%    \]
%    %\item Note that the vanishing in \eqref{eqn: h1-p-curv-Xb} in the statement of \autoref{lemma:pcurv-vanish} is slightly stronger than the vanishing of the $p$-curvature of the foliation at $(\mscr{E}_B, \nabla)$: indeed, this map 
%\end{itemize}
    
%\end{remark}
We postpone the proof of \autoref{lemma:pcurv-vanish} to the next section; it is the source of the vanishing we use to see that deforming the Hodge filtration is unobstructed.

%The following is an immediate consequence of the Tor 
%spectral sequence. 
%\begin{proposition}\label{lemma:basechangecomplex}
%    Let $A$ be an Artin algebra with residue field $k$, and $K$ a bounded complex
%    of free $A$-modules. If, for all $i$, $H^i(K)$ is free, then the formation of
%    the  cohomology of $K$ is compatible with all 
%base change $T\rightarrow \Spec(A)$.
%    \end{proposition} 
%    
%        \begin{remark}\label{rmk:boundedflatamplitude}
%            We will apply \autoref{lemma:basechangecomplex} to 
%            $K\defeq R\pi_*(\gr (\End \mscr{E} \otimes \Omega^{\bullet}))$.
%            Note that this complex is quasi-isomorphic to a complex of free $A$-modules by 
%            SGA 6 III 3.7.1, which is one of the hypotheses of this proposition.
%        \end{remark}
%

\begin{proposition}\label{prop:sspages_artin}
    %Suppose that the Hodge-de Rham spectral sequence for $\End(\mscr{E}_A, \nabla)$ degenerates 
%at the $E_1$-page, and that the terms of the $E_1$-page are free $A$-modules. 

\begin{enumerate}  
    \item\label{ssterms} The terms of the spectral sequences for the filtered complexes 
    $$(F^*_{\text{abs}}T_{\overline{X}'_{B_{\mf p}}/B_{\mf p}}(-\log D_{B_{\mf p}})_{dR}, \tau),
   (\End(\mscr{E}_B)_{dR}, F_{\text{conj}})$$ are 
    exactly as in \autoref{prop:sspages}. Precisely, 
    \begin{align*}
        E_{2}'^{p,q}&=
        \mb{H}^{p+q}(G^*\on{gr}^{q}_{F_{\text{Hodge}}}(\End(\mscr{E}_B)_{dR}))\\
        E_2^{p,q} &= H^{-k}(T_{\overline{X}'_{B_{\mf p}}/B_{\mf p}}(-\log D'_{B_{\mf p}}) \otimes  \Omega_{\overline{X}'_{B_{\mf p}}/B_{\mf p}}^{-(p+k)+1}(\log D'_{B_{\mf p}})), 
    \end{align*} 
    with $k$ defined by $2(p+k)=1-q$, and our convention is that the latter group is zero when  $k$ is not an integer.
%In the above, note that $\mscr{E}_B'$ is pulled back from along 
%$G: \X_B'\xrightarrow{G} \X_A$, and hence we may pullback the Hodge filtration; 
%therefore it makes sense to define
%$\on{gr}_F^p(\End(\mscr{E}_B')\otimes \Omega^{\bullet}_{\X_B'/B})$.

\item\label{ssmaps} 
We have 
\begin{align*}
    E_{2}^{2,-1} &= H^1(T_{\overline{X}'_{B_{\mf p}}/B_{\mf p}}(-\log D'_{B_{\mf p}})) \\
                &= (H^1(T_{\overline{X}_{A_{\mf p}}/A_{\mf p}}(-\log D_{A_{\mf p}})) \otimes_{A_{\mf p}, F} B_{\mf p}
\end{align*}
where we denote by $F: A_{\mf p}\rightarrow B_{\mf p}$ the ring map corresponding to the map of schemes
$F: \Spec(B_{\mf p})\rightarrow \Spec(A_{\mf p})$ of \autoref{defn:lifting-and-conjugate-filtration}. Similarly, we have 
\begin{align*}
E_{2}'^{2,-1} &= \mb{H}^1(G^*\on{gr}^{-1}_{F_{\text{Hodge}}}(\End(\mscr{E}_{A_{\mf p}})_{dR})) \\
            &= \mb{H}^1(\on{gr}^{-1}_{F_{\text{Hodge}}}(\End(\mscr{E}_{A_{\mf p}})_{dR})
            \otimes_{A_{\mf p}, F}B_{\mf p}.
\end{align*}

The induced map $E_{2}^{2,-1}\rightarrow E_2'^{2,-1}$, which is equivalently a map
\[
    H^1(T_{\overline{X}_{A_{\mf p}}/A_{\mf p}}(-\log D_{A_{\mf p}})) \otimes_{A_{\mf p}, F} B_{\mf p} \rightarrow 
    \mb{H}^1(\on{gr}^{-1}_{F_{\text{Hodge}}}(\End(\mscr{E}_{A_{\mf p}})_{dR}))
            \otimes_{A_{\mf p}, F}B_{\mf p},
    \]
is given by $\mu \otimes_{A, F} \rm{id}_B$, where $\mu$ is the map 
\begin{equation}
   \mu: H^1(T_{\overline{X}_{A_{\mf p}}/A_{\mf p}}(-\log D_{A_{\mf p}})) \rightarrow 
    \mb{H}^1(\on{gr}^{-1}_{F_{\text{Hodge}}}(\End(\mscr{E}_{A_{\mf p}})_{dR}))
\end{equation}    
obtained  by taking $\mb{H}^1$ of the map of complexes 
\[[T_{\overline{X}_{A_{\mf p}}/A_{\mf p}}(-\log D_{A_{\mf p}}) \rightarrow 0 \rightarrow \cdots ] \rightarrow 
\on{gr}^{-1}_{F_{\text{Hodge}}}(\End(\mscr{E}_{A_{\mf p}})_{dR})
\]
induced by the Higgs field on $\on{gr}_{F_{\text{Hodge}}}(\mathscr{E}_{A_\mf p})$, namely $\on{gr}_{F_{\text{Hodge}}}(\nabla_{A_{\mf p}})$.
\end{enumerate}
\end{proposition}

\begin{proof}
    Part (\ref{ssterms}) is exactly the same as \autoref{prop:sspages}, so it remains to  
    prove (\ref{ssmaps}). Note that $F: \Spec(B)\rightarrow \Spec(A)$ is flat, being a basechange of Frobenius on $S$, which is flat on a smooth scheme. The formulas for  $E_{2}^{2,-1}, E_2'^{2,-1}$ follow immediately 
    from (\ref{ssterms}) and flat basechange.  %We now prove the claim about $E_2'^{2,-1}$.
    
    %We first define the complex 
    %\[
     %   K_p\defeq R\pi_*\on{gr}_F^{p}(\End(\mscr{E})\otimes \Omega^{\bullet}_{\X_A/A}).
      %  \]
       %  Note that our  assumption that the terms in the Hodge-de Rham 
   % spectral sequence for $(\End(\mscr{E}), \nabla, F)$, namely 
   % \[H^{p+q}(K_p)= 
    %\mb{H}^{p+q}(\on{gr}^{p}(\End(\mscr{E}, \nabla)\otimes \Omega^{\bullet}_{\X_A/A})),\]
    %are free $A$-modules, implies that $K_p$
    %satisfies the hypotheses in 
    %\autoref{lemma:basechangecomplex}.     
   %Therefore the cohomology of $K_p$ is compatible with arbitrary basechange, and 
   % the claim that 
   % \[\mb{H}^1(\on{gr}^{-1}_F\End(\mscr{E}_B'\otimes \Omega^{\bullet}_{\X_{B'}/B})) \\
    %= \mb{H}^1(\on{gr}^{-1}_F\End(\mscr{E}_A'\otimes \Omega^{\bullet}_{\X_{A'}/A}))
    %\otimes_{A, F}B\] 
     %now follows by applying basechange along $\Spec(B)\xrightarrow{F} \Spec(A)$. 
     
     The final claim follows immediately from basechange and 
     the statement that the induced map 
    \[ E_{2}^{2,-1} = H^1(T_{\overline{X}'_{B_{\mf p}}/B_{\mf p}}(-\log D'_{B_{\mf p}})) \rightarrow
     \mb{H}^1(G^*\on{gr}^{-1}_{F_{\text{Hodge}}}(\End(\mscr{E}_{A_{\mf p}})_{dR})) 
     = E_2'^{2,-1}
    \]
is induced by the map of complexes 
\[ [T_{\overline{X}'_{B_{\mf p}}/B_{\mf p}}(-\log D'_{B_{\mf p}})\rightarrow 0\cdots] \rightarrow 
   G^*\on{gr}^{-1}_{F_{\text{Hodge}}}(\End(\mscr{E}_{A_{\mf p}})_{dR})\] 
    given by the Higgs field on $\on{gr}^{-1}_{F_{\text{Hodge}}}\mscr{E}_{A_{\mf p}}$. This latter statement follows from    \autoref{prop:higgs-field-vs-p-curvature}, which identifies the Higgs field on the associated graded with the flat section corresponding to the $p$-curvature defining our map of complexes. 

\end{proof}

\begin{comment}
        \begin{proposition}
            Let $T$ be a scheme over $\mb{Q}$, $h: Y\rightarrow T$  a smooth $T$-scheme;  
            suppose that   $(\mscr{E}, \nabla, F)$ is a polarizable variation of Hodge structures
            on $Y/T$. Then the Hodge-de Rham spectral sequence 
            \begin{equation}\label{eqn:hdrdegenscheme}
                E_{1}^{p,q}\defeq R^{p+q}h_{*}^{HIG}(\gr_F^p(\mscr{E}\otimes \Omega^{\bullet}_{Y/T}))
                \Rightarrow R^{p+q}h_{*}^{dR}(\mscr{E}\otimes \Omega^{\bullet}_{Y/T})
            \end{equation}
        degenerates at the $E_1$-page, and all the modules in \eqref{eqn:hdrdegenscheme} are locally free 
        over $T$. 
        \end{proposition}
\end{comment}

%   \begin{lemma}\label{claim:ssdegen} The conjugate spectral sequence 
%    $(E_{2}'^{p,q})$ degenerates at the $E_2$-page. 
%    \end{lemma}
%    \begin{proof}
%        This follows directly from  \cite[Corollary 3.25]{ogus-vologodsky} and  \autoref{remark:conj-ss-paranjape} when $D=\emptyset$; the  proof is formal from the compatibility of the inverse Cartier transform with derived pushforward and works in general. \daniel{Find reference?}
%    \end{proof}

\begin{lemma}[Vanishing of obstruction to deforming in a Griffiths-transverse manner]\label{lemma: ks-vanish-induction-step}
    %Under the same hypotheses as  \autoref{prop:sspages_artin},
    After possibly replacing $R$ by $R[\frac{1}{f}]$ for some $f\in R$, map 
    \[
        KS: H^1(T_{\overline{X}_A/A}(-\log D_A)) \rightarrow 
        {H}^1_{\text{Higgs}}(\End(\on{gr}_{F_\text{Hodge}} \mscr{E}_A, \theta)/F^0_{\text{Hodge}})
        \]
        vanishes on the image of the Kodaira-Spencer map $T_{S, A}\to H^1(T_{\overline{X}_A/A}(-\log D_A))$.
\end{lemma}

\begin{proof}
    The argument is almost identical to that of \autoref{lemma:ksvanishing}, though we provide some details here for the reader's convenience.  It suffices to show the map vanishes mod $\mf p$ for $\mf p$ in an open subset of $\on{Spec}(R)$.   
    %and the 
    %only remaining  thing to check is the following     
    %\begin{claim}\label{claim:ssdegen} The conjugate spectral sequence 
   % $(E_{2}'^{p,q})$ degenerates at the $E_2$-page. 
   % \end{claim}
    %For the reader's convenience, we first recall how to deduce the lemma from \autoref{claim:ssdegen}.
    
    Since the $p$-curvature map 
    \[
        F^*_{\text{abs}}T_{S, B_{\mf p}} \rightarrow H^1(T_{\overline{X}'_{B_{\mf p}}/B_{\mf p}}(-\log D'_{B_{\mf p}}))\rightarrow 
        \mb{H}^1(\End(\mscr{E}_B)_{dR})
        \]
    vanishes by \autoref{lemma:pcurv-vanish} 
        the same spectral sequence argument as in the proof of \autoref{lemma:ksvanishing}
    implies $E_{2}^{2,-1} \rightarrow E_2'^{2,-1}$ vanishes on the image of 
    $T_{S, B}$.

    By \autoref{prop:sspages_artin}(\ref{ssmaps}), this is precisely the map 
    $\mu\otimes_{A, F} \id_B$, where $\mu$ is the map on $H^1$ induced by 
\[[T_{\overline{X}_{A_{\mf p}}/A_{\mf p}}(-\log D_{A_{\mf p}}) \rightarrow 0 \rightarrow \cdots ] \rightarrow 
\on{gr}^{-1}_{F_{\text{Hodge}}}(\End(\mscr{E}_{A_{\mf p}})_{dR})
\]
Since the map $A_{\mf p}\xrightarrow{F} B_{\mf p}$ makes $B_{\mf p}$ into a faithfully flat $A_{\mf p}$-algebra, we deduce 
that $\mu$ vanishes on the image of $H^1(T_{S, A_{\mf p}})$; since $KS$ factors through $\mu$, 
we conclude that $KS$ also vanishes on the image of $H^1(T_{S, A_{\mf p}})$, as required.
\end{proof}
    %It therefore  remains to  prove \autoref{claim:ssdegen}.

    \begin{corollary}\label{cor:Hodge-filtration-extends}
        Under the setup of \autoref{induction_setup}, suppose the Hodge-de Rham spectral sequence for $(\mscr{F}_{A, i}, \nabla, F^\bullet)$ degenerates at $E_1$. Recall that $\mathscr{K}=\on{Frac}(R)$.
        
        Then for any $\Spec(A_\mathscr{K}[\epsilon]/\epsilon^2)\rightarrow S$ deforming the given map $\Spec(A_\mathscr{K})\rightarrow S$, the obstruction to deforming $F^\bullet$ in a Griffiths-transverse manner to the isomonodromic deformation of $(\mscr{F}_{A_{\mathscr{K}}, i}, \nabla)$ to $(\overline{X}_{A_\mathscr{K}[\epsilon]/\epsilon^2}, D_{A_\mathscr{K}[\epsilon]/\epsilon^2})/\on{Spec}(A_\mathscr{K}[\epsilon]/\epsilon^2)$ vanishes.
    \end{corollary}
    \begin{proof}
    The argument is the same as the argument for \autoref{cor:base-case-Hodge-extends}.
    
        Since the Hodge-de Rham spectral sequence for $(\mscr{F}_{A_{\mathscr{K}}, i}, \nabla, F^\bullet)$ degenerates at $E_1$, by \autoref{prop: gt-transverse-obstruction}, the obstruction to deforming $(\mscr{F}_{A_\mathscr{K}, i}, \nabla, F^\bullet)$ isomonodromically and Griffiths-transversally to $(\overline{X}_{A_\mathscr{K}[\epsilon]/\epsilon^2}, D_{A_\mathscr{K}[\epsilon]/\epsilon^2})/\on{Spec}(A_\mathscr{K}[\epsilon]/\epsilon^2)$ vanishes if and only if the image of the class in $H^1(T_{\overline{X}_{A_\mathscr{K}}/A_\mathscr{K}}(-\log D_{A_\mathscr{K}}))$ corresponding to $(\overline{X}_{A_\mathscr{K}[\epsilon]/\epsilon^2}, D_{A_\mathscr{K}[\epsilon]/\epsilon^2})$ under the map 
       \[ 
       KS_M: H^1(T_{\overline{X}_{A_\mathscr{K}}/A_\mathscr{K}}(-\log D_{A_\mathscr{K}})) \rightarrow 
        {H}^1_{\text{Higgs}}(\End(\on{gr}_{F_\text{Hodge}} \mscr{F}_{A_\mathscr{K}, i}, \theta)/F^0_{\text{Hodge}})
        \]
        vanishes. 
        But this is a component of the map in \autoref{lemma: ks-vanish-induction-step} which itself vanishes, as required.
    \end{proof}

\subsection{Proof of the key lemma}\label{subsection:proof-of-key-lemma}
We now begin preparations for the proof of \autoref{lemma:pcurv-vanish}. This will be a mixed-characteristic argument. The idea will be to lift the construction of $(\mathscr{E}_B, \nabla)$ from the previous section to characteristic zero using the Higgs-de Rham flow of Lan-Sheng-Zuo, giving rise to a bundle $(\mathscr{E}_{B_{\widehat R}}, \nabla)$ on $(\overline{X}, D)_{B_{\widehat R}}$. We will argue that this bundle is isomonodromic and then use the $1$-periodicity of $(\mathscr{E}, \nabla)$ under the Higgs-de Rham flow (or in other words, the fact that it is a Fontaine module) to conclude that it is in fact an isomonodromic deformation of $(\mathscr{E}, \nabla)$. It will then follow from \autoref{cor:unique-integral-isomonodromic} that $(\mathscr{E}_B, \nabla)$ is itself $1$-periodic, i.e.~a Fontaine module. Now the vanishing of the $p$-curvature of the isomonodromy foliation will follow from the existence of an ($\omega(p)$)-integral isomonodromic deformation extending $(\mathscr{E}_{B_{\widehat R}}, \nabla)$ by \autoref{prop:crystal-p-power-leaves}, the analogue of \autoref{prop:p-power-leaves} for a crystal of functors. 

We keep the setup as in \autoref{induction_setup}, so that we have the diagram

\begin{equation} \begin{tikzcd}
    (\overline{X}', D')_{B_{\widehat R}}\arrow[r, "G_0"] \arrow[d]
    & (\overline{X}, D)_{A_{\widehat R}} \arrow[d] \arrow[r] & (\overline{X}, D)_{B_{\widehat R}}  \arrow[d]\\
    \Spec(B_{\widehat R}) \arrow[r, "F_0"]
    &  \Spec(A_{\widehat R}) \arrow[r] & \Spec(B_{\widehat R}).
    \end{tikzcd}
\end{equation} 
Recall that, by the induction hypothesis, we have a filtered 
flat bundle $(\mscr{E}_{A_{\widehat R}}, \nabla_{A_{\widehat R}}, F_{A_{\widehat R}}^{\bullet})$ such that 
\begin{itemize}
    \item $F_{A_{\widehat R}}^{\bullet}$ is Griffiths transverse, and 
    \item $(\mscr{E}_{A_{\widehat R}}, \nabla_{A_{\widehat R}})_{\mathbb{Q}}$ is isomonodromic.
\end{itemize}

We now explain how to construct a lift of $(\mathscr{E}_B, \nabla)$ to $ (\overline{X}', D')_{B_{\widehat R}}$, following ideas of \cite{eg_revisit} and \cite{lan2019semistable}.
%\begin{proposition}\label{prop:abs_conn_griff_trans}
%    The connection $(\mscr{E}_{\Ainthat}, \nabla_{\Ainthat})$ extends to an 
%    absolute connection 
%    \[
%        \mscr{E}_{\Ainthat}\rightarrow \mscr{E}_{\Ainthat}\otimes \
%        \Omega^1_{\mc{X}_{\Ainthat}}.\]
%        Moreover, the filtration $F_0^{\bullet}$ is still 
%        Griffiths transverse for this connection.
%
%\end{proposition}

\begin{construction}
Consider the filtered connection 
$G_0^*(\mscr{E}_{A_{\widehat R}}, \nabla_{A_{\widehat R}}, F_{A_{\widehat R}}^{\bullet})$. As in 
\cite[Definition 3.26]{eg_revisit}, we can construct the 
Artin--Rees bundle $\AR$ of
 $G_0^*(\mscr{E}_{A_{\widehat R}}, \nabla_{A_{\widehat R}}, F_{A_{\widehat R}}^{\bullet})$. 

More precisely, this is the sheaf of $\mc{O}_{\overline{X}'_{B_{\widehat R}}}$-sub-modules
\[
\sum_i G_0^*F_{A_{\widehat R}}^{i}\otimes p^{-i} \subset G_0^*\mscr{E}_{A_{\widehat R}} \otimes_{\mathbb{Z}_p} \mathbb{Q}_p.
\]

Moreover, $\AR$ is equipped with its 
natural $p$-connection, coming from the connection $(\mscr{E}_{A_{\widehat R}}, \nabla_{A_{\widehat R}})$, which we denote by $\nabla_p:=p\nabla_{A_{\widehat R}}$. 
More precisely, this is a map 
\[ 
    \nabla_p: \AR \rightarrow \AR \otimes 
    \Omega^1_{\overline{X}'_{B_{\widehat R}}}(\log D_{B_{\widehat R}}), 
    \] 
    satisfying the $p$-Leibniz rule (see \cite[Definition 3.7]{eg_revisit}) and squaring to zero. 
\end{construction}

%The following is an immediate corollary of \autoref{prop:abs_conn_griff_trans}.
%\begin{corollary}
%    $\nabla_p$ extends to a $p$-connection 
%    \[ 
%    \nabla_p: \AR \rightarrow \AR \otimes 
%    \Omega^1_{\overline{X}'_{B_{\widehat R}}}(\log D_{B_{\widehat R}}).
%    \] 
%\end{corollary}
\begin{construction}\label{construction:decent-frobenii}
    We say that  a formal  affine open subscheme 
    $\iota: \mc{U} \rightarrow \overline{X}_{B_{\widehat R}} $
     has a
    \emph{decent} Frobenius lift, if it 
    is  equipped with a lift $F_{\mc{U}}$ 
    of absolute Frobenius such 
    that 
    \begin{enumerate} \item 
    the square 
    \begin{equation} \begin{tikzcd}
        \mc{U} \arrow[d] \arrow[r, "F_{\mc{U}}"]
        & \mc{U}   \arrow[d] \\
        \Spec(B_{\widehat R}) \arrow[r, "F_0"]
        &  \Spec(B_{\widehat R}) 
        \end{tikzcd}\end{equation}
commutes, and 
\item $F_{\mc{U}}^*\mc{O}_{\mc{U}}(-{D}_{\mc{U}})= \mc{O}_{\mc{U}}(-p{D}_{\mc{U}})$, where ${D}_{\mc{U}}$ is the pullback of ${D}$ to $\mc{U}$.
\end{enumerate}
For   $\mc{U}$ equipped with such  a decent Frobenius lift, 
we may form  the diagram 
\begin{equation}
\begin{tikzcd}
    \mc{U}
    \arrow[drr, bend left, "F_{\mc{U}}"]
    \arrow[ddr, bend right, ""]
    \arrow[dr, "\alpha"] & & \\
    & \mc{U}' \arrow[r, "G_0"] \arrow[d, ""]
    & \mc{U} \arrow[d, "\pi"] \\
    & \Spec(B_{\widehat R}) \arrow[r, "F_0"]
    & \Spec(B_{\widehat R}),
    \end{tikzcd}
\end{equation}
with the square being cartesian.
\begin{remark}
    Note that the data of a  decent Frobenius lift is 
    equivalent to that of a morphism of $B_{\widehat R}$-schemes
    $\alpha: \mc{U}\rightarrow 
    \mc{U}'$ lifting  the 
     relative Frobenius map; in the following, when we say a decent Frobenius lift, we will mean  either $F_{\mc{U}}$ or $\alpha$, depending on which is more convenient for the situation.
\end{remark}
We can form the connection on the sheaf $\alpha^*\AR$  by the following recipe: for a section $s$ of $\alpha^*\AR|_{\mc{U}'}$, we set 
\[ 
    \nabla_{\mc{U}, \alpha}(s):= 
    \frac{\mathrm{id}\otimes \alpha^*}{p}\nabla_p(s).
    \] 
\end{construction}
\begin{lemma}\label{lemma:absolute_conn}
    We can cover (the formal scheme) $\overline{X}_{B_{\widehat R}}$ by decent 
    affine formal open subschemes, and  
    the connections $(\alpha^*\AR, \nabla_{\mc{U}, \alpha})$ glue to 
    a connection on $(\overline{X}_{B_{\widehat R}}, D_{B_{\widehat R}})/B_{\widehat R}$. 
 \end{lemma}

\begin{proof}
    The first part follows immediately from iterating
    \cite[Corollary 9.12]{esnault-viehweg}, i.e. first 
    applying loc.cit. to construct liftings to $\Binthat/p^2$, then to 
    $\Binthat/p^3$, and so on. 

    For the gluing, for a formal affine open 
    $\mc{U}\subset \overline{X}_{B_{\widehat R}}$ and two decent Frobenius liftings
    \[
        \alpha, \beta: \mc{U} \rightarrow \mc{U}',\]
        we must provide an $\mc{O}_{\mc{U}}$-linear isomorphism 
        \begin{equation}\label{glue}
            \epsilon_{\alpha\beta}: 
            \alpha^*\AR \xrightarrow[]{\sim} \beta^*\AR.  
        \end{equation}
            Let us pick \'etale affine coordinates  (shrinking $\mc{U}$ if necessary)
            $\tau: \mc{U} \rightarrow \mb{A}^n_{B_{\widehat R}}$; and let $\alpha_j$ denote the pullback of the coordinate function $t_j$ under the composition
            \[
            \mc{U}\xrightarrow{\alpha} \mc{U}'\xrightarrow{\tau'} \mb{A}^n_{B_{\widehat R}},
            \]
            with the last map being the basechange of $\tau$ along $F_0$. We define the $\beta_j$'s similarly. 
            Following \cite[Definition 3.27]{eg_revisit}, we now define 
            $\epsilon_{\alpha \beta}$ by 
        \[
            (s)_i\otimes 1 \mapsto 
        \sum_{|\vecu{i}|\leq i}(\partial^{\vecu{i}}(s))_{i-|\vecu{i}|}
        \otimes \frac{\prod_{j=1}^n (\alpha_j-\beta_j)^{i_j}}
        {p^{|\vecu{i}|}\vecu{i}!} + 
        \sum_{|\vecu{i}|>i}p^{|\vecu{i}|-i}(\partial^{\vecu{i}}(s))_0
        \otimes \frac{\prod_{j=1}^n (\alpha_j-\beta_j)^{i_j}}
        {p^{|\vecu{i}|}\vecu{i}!}.
        \]
        
        It is straightforward to check that the 
        $\epsilon_{\alpha\beta}$'s satisfy the cocycle condition, i.e. for 
        another $\gamma$ lifting relative Frobenius, we have 
        \[ 
\epsilon_{\beta\gamma}\epsilon_{\alpha\beta} = \epsilon_{\alpha\gamma}.\]

        It remains to check that, in \eqref{glue},
         $\epsilon_{\alpha \beta}$ is flat, namely that  the following 
         \begin{equation} \begin{tikzcd}
            \alpha^*\AR \arrow[d, "\nabla_{\mc{U}, \alpha}"] \arrow[r, "\epsilon_{\alpha\beta}"]
            & \beta^*\AR   \arrow[d, "\nabla_{\mc{U}, \beta}"] \\
            \alpha^*\AR\otimes \Omega^1_{\mc{U}/B_{\widehat R}}(\log D_{\mc{U}}) 
            \arrow[r, "\epsilon_{\alpha\beta}\otimes \mathrm{id}"]
            &   \beta^*\AR\otimes \Omega^1_{\mc{U}/B_{\widehat R}}(\log D_{\mc{U}}) \\
            \end{tikzcd}
        \end{equation}
        commutes; this is a straightforward computation. 
\end{proof}
\begin{definition}\label{defn:higgs-de-rham-flow-functor}
    We denote by $\Phi(\mscr{E}_{A_{\widehat R}}, \nabla_{A_{\widehat R}}, F_{A_{\widehat R}}^{\bullet})$ the bundle with flat connection on  $(\overline{X}_{B_{\widehat R}}, D_{B_{\widehat R}})/B_{\widehat R}$
    obtained in \autoref{lemma:absolute_conn}; following \cite{eg_revisit} we refer to $\Phi$ as the \emph{flow functor}.
\end{definition}
\begin{proposition}\label{prop:lift_inverse_cartier}
    The flat bundle  $\Phi(\mscr{E}_{A_{\widehat R}}, \nabla_{A_{\widehat R}}, F_{A_{\widehat R}}^{\bullet})$ restricted to  $\overline{X}_{B_\mf p}/B_{\mf p}$ is 
    precisely   the bundle \[
    (\mscr{E}_B, \nabla)\defeq C^{-1}(G^*(\on{gr}^\bullet_{F_{\text{Hodge}}}\mscr{E}_{A_{\mf p}}, \on{gr}^\bullet_{F_{\text{Hodge}}}(\nabla_{A_{\mf p}}))),
  \]
  of \autoref{defn:lifting-and-conjugate-filtration}.
\end{proposition}

\begin{proof}
    This is \cite[Prop 3.28]{eg_revisit} when $B_{\mf p}$ is a field. Their proof relies on \cite{lan2015nonabelian}, which is only written in this generality, but all the proofs work verbatim over more general bases.
\end{proof}

The key point is now:
\begin{proposition}\label{prop:special-point-and-isomonodromic}
	The flat bundle $\Phi(\mscr{E}_{A_{\widehat R}}, \nabla_{A_{\widehat R}}, F_{A_{\widehat R}}^{\bullet})$ satisfies the following:
	\begin{enumerate}
		\item  There is an isomorphism 	$\Phi(\mscr{E}_{A_{\widehat R}}, \nabla_{A_{\widehat R}}, F_{A_{\widehat R}}^{\bullet})|_{s_{\widehat R}}\simeq (\mathscr{E}, \nabla)_{s_{\widehat R}}$.
		\item $\Phi(\mscr{E}_{A_{\widehat R}}, \nabla_{A_{\widehat R}}, F_{A_{\widehat R}}^{\bullet})_{\mathbb{Q}}$ is isomonodromic.
	\end{enumerate}
\end{proposition}
\begin{proof}
	Part (1) is precisely the part of \autoref{thm:berthelot-mazur-ogus} about Fontaine-Laffaille modules; in other words it is the statement that Picard-Fuchs equations are $1$-periodic under the Higgs-de Rham flow (here we use \autoref{assumption: stable-special-fiber}). See \cite[Theorem 6.2]{faltings1989crystalline} for a classical reference or e.g.~\cite[Theorem 1.3]{krishnamoorthy2020periodic} for a reference in the language of Higgs-de Rham flows (in a slightly more general, parabolic, setting, though only in the case of curves---which suffices for our applications).
	
	For part (2), we use \cite[Proposition 3.39]{eg_revisit}, namely, that $\Phi(\mscr{E}_{A_{\widehat R}}, \nabla_{A_{\widehat R}}, F_{A_{\widehat R}}^{\bullet})$ is an integral model of the ``Frobenius pullback" $F^*(\mathscr{E}_{A_{\widehat R}}, \nabla_{A_{\widehat R}})$ of $(\mathscr{E}_{A_{\widehat R}}, \nabla_{A_{\widehat R}})$, which is crucially a construction that does not depend on $F^\bullet_{A_{\widehat R}}$. Namely, we have an isomonodromic flat bundle $(\mathscr{E}_{A_{\widehat R}}, \nabla_{A_{\widehat R}})$ on $(\overline{X}, D)_{A_{\widehat R}}/A_{\widehat R}$, by assumption. As $F_0$ factors through $\on{Spec}(A_{\widehat R})$, for each $\mc U$ we have that the map $F_{\mc U}$ chosen in \autoref{construction:decent-frobenii} factors through $(\overline{X}, D)_{A_{\widehat R}}$, and thus it makes sense to consider $F_{\mc U}^*(\mathscr{E}_{A_{\widehat R}}, \nabla_{A_{\widehat R}})$. We claim these bundles glue to an isomonodromic bundle. They are glued via the formula of \cite[Definition 3.1]{eg_revisit} (namely, Taylor's formula---see \cite[Proposition 3.3]{eg_revisit} for the logarithmic case) giving rise to a flat bundle on $(\overline{X}_{B_{\widehat R}}, D_{B_{\widehat R}})/B_{\widehat R}$. Finally, we check that the connection is isomonodromic in characteristic zero, i.e.~that $F^*(\mathscr{E}_{A_{\widehat R}}, \nabla_{A_{\widehat R}})_{\mathbb{Q}}$ is the restriction of an isomonodromic flat bundle on $(\overline{X},D)$ on the formal completion of $S_{\widehat R_\mathbb{Q}}$ at $s_{\widehat{R}_\mathbb{Q}}$. We do this in \autoref{lemma:frobenius-pullback-isomonodromic} below, where the isomonodromy assumption is satisfied as $(\mathscr{E}_{\widehat R}, \nabla_{A_{\widehat R}})_{\mathbb{Q}}$ is isomonodromic by the induction hypothesis, and $F^*(\mathscr{E}, \nabla)|_{\overline{X}_{0,\mathbb{Q}}}\simeq (\mathscr{E}, \nabla)|_{\overline{X}_{0,\mathbb{Q}}}$ as $(\mathscr{E},\nabla)$ is a Picard-Fuchs equation, hence gives rise to an $F$-isocrystal, see \autoref{subsection:berthelot-ogus-preliminaries} (or by part (1) and \cite[Proposition 3.39]{eg_revisit}).
\end{proof}

\begin{lemma}\label{lemma:frobenius-pullback-isomonodromic}
	Let $k$ be a finite field of characteristic $p>2$ and $S=\on{Spec}(W(k)[[t_1, \cdots, t_d]]$, and let $0\in S$ be the $W(k)$-point defined by $(t_1, \cdots, t_d)$. Let $\overline{X}$ be a smooth projective formal $S$-scheme and $D\subset \overline{X}$ a simple normal crossings divisor over $S$. Let $Z\subset S$ be a closed subscheme with $Z_{\text{red}}=\on{Spec}(W(k))$. Let $(\mathscr{E},\nabla)$ be a flat bundle on $(\overline{X}, D)_Z/Z$ with $(\mathscr{E}, \nabla)_{Z_{\mathbb{Q}}}$ isomonodromic.
	
	Choose the Frobenius lift $F_S$ on $W(k)[[t_1, \cdots, t_d]]$ defined by $t_i\mapsto t_i^p$ and let $F^*$ denote the ``Frobenius pullback" used in the proof of \autoref{prop:special-point-and-isomonodromic}, defined as in \cite[Proposition 3.3]{eg_revisit}; this is a functor from flat bundles on $(\overline{X}, D)_Z/Z$ to flat bundles on $(\overline{X}, D)_{F_S^{-1}(Z)}/F_S^{-1}(Z)$. Suppose $F^*(\mathscr{E}, \nabla)|_{\overline{X}_{0,\mathbb{Q}}}\simeq (\mathscr{E}, \nabla)|_{\overline{X}_{0,\mathbb{Q}}}$. Then $F^*(\mathscr{E},\nabla)_{F_S^{-1}(Z)_{\mathbb{Q}}}$ is isomonodromic.
\end{lemma}
\begin{proof}
	This is ``formal" from the crystalline description of the Frobenius pullback via the functoriality of the crystalline site (which is made explicit via Taylor series in \cite{eg_revisit} or \cite{berthelot2000d}); this is the statement that Frobenius pullback is flat for the non-abelian Gauss-Manin connection, and is the non-abelian analogue of the flatness of the crystalline Frobenius $\Phi$ of \autoref{subsection:berthelot-ogus-preliminaries}. In fact, the condition that $F^*(\mathscr{E}, \nabla)|_{X_{0,\mathbb{Q}}}\simeq (\mathscr{E}, \nabla)|_{X_{0,\mathbb{Q}}}$ is unnecessary. However, we give an explicit proof under this assumption, as we found it somewhat illuminating.

	Let $\widehat S:=\on{Spf}(W(k)[[t_1, \cdots, t_d]])$ and let $(\widehat{\overline{X}}, \widehat{D})$ be the base change of $(\overline{X}, D)$ to $\widehat S$. As above, choose a cover of $\widehat{\overline{X}}$ by formal affines equipped with decent Frobenius liftings over $F_S$, and \'etale affine coordinates, so that we may define a Frobenius pullback functor on flat bundles, $F^*$, as in the proof of \autoref{prop:special-point-and-isomonodromic} above or as in \cite[Proposition 3.3]{eg_revisit}.
	
	Let $\widehat{S}_\mathbb{Q}=\on{Spf}(\on{Frac}(W(k))[[t_1, \cdots, t_d]])$, and let $\widehat{\overline{X}}_{\mathbb{Q}}, \widehat{D}_{\mathbb{Q}}$ be the base change of $\overline{X}, D$ to $\widehat S_\mathbb{Q}$. We will show that $F^*(\mathscr{E},\nabla)_{\mathbb{Q}}$ is the restriction of a flat bundle on $(\widehat{\overline{X}}, \widehat{D})_\mathbb{Q}$ equipped with an \emph{absolute} connection, and is hence isomonodromic in the sense of \autoref{defn:naive-isomondromy}.
	
	Let $S_n=(F_S^n)^{-1}(0)$, i.e~the subscheme cut out by $(t_1^{p^n}, \cdots, t_d^{p^n})$. Let $(\overline{X}, D)_n=(\overline{X},D)|_{S_n}$. Note that for $N$ sufficiently large, $(F_S|_{S_n})^N$ factors through $S_0=0$. Hence $(F)^{N*}(\mathscr{E}, \nabla)|_{\overline{X}_n}$ depends only on $(\mathscr{E},\nabla)|_{\overline{X}_0}$ and carries an absolute connection (as it is pulled back from $(\mathscr{E},\nabla)|_{\overline{X}_0}$, which carries an absolute connection). We claim these connections glue on tensoring with $\mathbb{Q}$; we will glue them using the formulae in \autoref{lemma:char-0-taylor-series}. Note that all of these connections are isomorphic when restricted to $(\overline{X}_0, D_0)_\mathbb{Q}$, by the assumption that $F^*(\mathscr{E}, \nabla)|_{X_{0,\mathbb{Q}}}\simeq (\mathscr{E}, \nabla)|_{X_{0,\mathbb{Q}}}$. For each $n$, choose an isomorphism $(F)^{N*}(\mathscr{E}, \nabla)|_{\overline{X}_0}\simeq (F)^{(N-1)*}(\mathscr{E}, \nabla)|_{\overline{X}_0}$.
	
	 The sheaf $\Omega^1_{\overline{X}_{n,\mathbb{Q}}/W(k)}(\log D_n)$ is not locally free, so we cannot apply these formulae directly. However, $\Omega^1_{\overline{X}_{n,\mathbb{Q}}/W(k)}(\log D_n)|_{\overline{X}_{n-1}}$ is locally free, and restricting the absolute connection $(F^N)^*(\mathscr{E}, \nabla)$ on $(\overline{X}_n, D_n)$ to $(\overline{X}_{n-1}, D_{n-1})$ yields an ``absolute connection": $$\nabla: (F^N)^*\mathscr{E}\to (F^N)^*(\mathscr{E})\otimes \Omega^1_{\overline{X}_{n,\mathbb{Q}}/W(k)}(\log D_n)|_{\overline{X}_{n-1}},$$ so, for example, it makes sense to apply the operator $\nabla_{\frac{\partial}{\partial t_j}}$ (where $\frac{\partial}{\partial t_j}$ is dual to $dt_j$) to a section to $(F^N)^*\mathscr{E}$.
	 
	 This suffices to apply the formulae of \autoref{lemma:char-0-taylor-series} after tensoring with $\mathbb{Q}$; as all the $(F)^{N*}(\mathscr{E}, \nabla)|_{\overline{X}_0}$ are identified, this yields an absolute connection $(\mathscr{G},\nabla)$ on  $(\widehat{\overline{X}},\widehat D)_\mathbb{Q}$ as desired. It remains to identify its restriction to $F_S^{-1}(Z)_{\mathbb{Q}}$ with $F^*(\mathscr{E},\nabla)_{F_S^{-1}(Z)_{\mathbb{Q}}}$. 
	 
	 But note that the restriction of $(\mathscr{G},\nabla)$ to $(\overline{X}, D)_{Z_\mathbb{Q}}/Z_{\mathbb{Q}}$ must agree with $(\mathscr{E}, \nabla)_{Z_{\mathbb{Q}}}$, as both are isomonodromic.	But for all sufficiently large $N$, $F^{N*}((\mathscr{E},\nabla)|_{\overline{X}_0})|_{(\overline{X}_Z, D_Z)/Z}$ is an integral model of $$(\mathscr{G},\nabla)|_{(\overline{X}_Z, D_Z)/Z}=(\mathscr{E}, \nabla)_{Z_{\mathbb{Q}}}.$$ It follows that $F^{(N+1)*}((\mathscr{E},\nabla)|_{\overline{X}_0})|_{(\overline{X}_{F_S^{-1}(Z)}, D_{F_S^{-1}(Z)})/{F_S^{-1}(Z)}}$ is an integral model of $F^*(\mathscr{E},\nabla)_{F_S^{-1}(Z)_{\mathbb{Q}}}$; but it is also an integral model of $(\mathscr{G}, \nabla)|_{F_S^{-1}(Z)_{\mathbb{Q}}}$, completing the proof.
	 
\end{proof}

We now proceed with the proof of  \autoref{lemma:pcurv-vanish}.

\begin{proof}[Proof of \autoref{lemma:pcurv-vanish}]
    Our setup is as in \autoref{induction_setup}.     
    
    We first check that $(\mathscr{E}_B, \nabla)|_{A_\mf p}\simeq (\mathscr{E}_{A_\mf p}, \nabla_{A_\mf p})$. It suffices by \autoref{prop:lift_inverse_cartier} to show the stronger statement that $$\Phi(\mscr{E}_{A_{\widehat R}}, \nabla_{A_{\widehat R}}, F_{A_{\widehat R}}^{\bullet})|_{A_{\widehat R}}\simeq (\mathscr{E}_{A_{\widehat R}}, \nabla_{A_{\widehat R}}).$$ By \autoref{prop:special-point-and-isomonodromic}, these two bundles are isomorphic when restricted to $s_{\widehat R}$. But both are isomonodromic in characteristic zero, whence by \autoref{assumption: stable-special-fiber} and \autoref{cor:unique-integral-isomonodromic}, they are in fact isomorphic.
    
    Finally, we check the vanishing of the $p$-curvature of the isomonodromy foliation. on restriction to our given $B_{\mf p}$-point of $\mathscr{M}_{dR}(X/S)^\natural$. By \autoref{prop:special-point-and-isomonodromic}, $\Phi(\mscr{E}_{A_{\widehat R}}, \nabla_{A_{\widehat R}}, F_{A_{\widehat R}}^{\bullet})_{\mathbb{Q}}$ is isomonodromic over $B_{\widehat R, \mathbb{Q}}$.
    
    Recall (from \autoref{assumption:monomial-A}) that we have chosen coordinates $t_1,\cdots, t_{\dim S}$ on $\widehat S$, and that the ideal cutting out $A$ is monomial with respect to these coordinates, say the ideal $(t_1^{A_1}, \cdots, t_{\dim S}^{A_{\dim S}})$ where $A_i$ are multi-indices. We moreover have some $q$ such that $$(t_1, \cdots, t_{\dim S})^q\subset (t_1^{A_1}, \cdots, t_{\dim S}^{A_{\dim S}}).$$ Let $Z$ be the subscheme of $S_{\widehat R}$ cut out by $(t_1^{pA_1}, \cdots, t_{\dim S}^{pA_{\dim S}})$; $Z$ is a lift of $\on{Spec}(B_{\mf p})$. As in the setup for \autoref{cor:vanishing-p-curvature-for-p-power-leaves-isomonodromy}, let $Z_+$ be the subscheme cut out by $$(t_1, \cdots, t_{\dim S})^{m+p+1},$$ where $m$ is minimal such that $Z\subset V((t_1, \cdots, t_{\dim S})^m)$; for example, we may take $m=pq\dim S$. By the assumption of $\omega(p)$-integrality, we have a model $(\mathscr{G}, \nabla)$ of our isomonodromic deformation of $(\mathscr{E}, \nabla)$ over $Z_+$.
    
    Again by \autoref{cor:unique-integral-isomonodromic}, $\Phi(\mscr{E}_{A_{\widehat R}}, \nabla_{A_{\widehat R}}, F_{A_{\widehat R}}^{\bullet})$ must agree with $(\mathscr{G}, \nabla)_{B_{\widehat R}}$. Thus the hypotheses of \autoref{cor:vanishing-p-curvature-for-p-power-leaves-isomonodromy} are satisfied and the $p$-curvature vanishes as desired.
\end{proof}

\begin{remark}
Instead of using the uniqueness of integral isomonodromic deformations under suitable conditions---\autoref{cor:unique-integral-isomonodromic}---we could also have argued via polystability of the special fibers of the bundles $$\Phi(\mscr{E}_{A_{\widehat R}}, \nabla_{A_{\widehat R}}, F_{A_{\widehat R}}^{\bullet})\simeq (\mathscr{E}_{A_{\widehat R}}, \nabla_{A_{\widehat R}}),$$ which will hold on an open of $\on{Spec}(R)$.
\end{remark}

\section{Smoothness of deformation space}\label{section:smoothness}
We now prepare to prove \autoref{theorem:extending-the-Hodge-filtration}.
\subsection{$T^1$-lifting theorem}\label{subsection:T1-lifting}

Let us formulate a version of the $T^1$-lifting theorem suitable for our induction argument. Let $\fracr$ be a field of characteristic zero, $\mf{R}=\fracr[[x_1, \cdots, x_n]]$ with maximal ideal $\mf{m}=(x_1, \cdots, x_n)$, $I\subset \mf{R}$ an ideal, and $\mf{T}_i=\mf{R}/I$. Let $\mf{R}_n=\mf{R}/\mf{m}^n$. For $B$ a ring and $M$ a $B$-module, we denote by $B\oplus M$ the ring whose underlying abelian group is the direct sum $B\oplus M$, with multiplication defined by $$(b_1, m_1)\cdot (b_2, m_2)=(b_1b_2, b_1m_2+b_2m_1).$$
\begin{lemma}\label{lemma: t1-lift-induct}
	Suppose that $I\subset \mf{m}^n$. Suppose that for all free  finite rank $\mf{R}_n$-modules $M$, and all maps $f: \mf{R}\to \mf{R}_n\oplus M$ such that the composition $ \mf{R}\xrightarrow[]{f} \mf{R}_n\oplus M\to \mf{R}_n$ is the natural quotient map $$q_n: \mf{R}\to \mf{R}/\mf{m}^n,$$ we have that $f(I)=0$, i.e.~$f$ descends to a (unique) map $\mf{T}_i\to \mf{R}_n\oplus M$. Then $I\subset \mf{m}^{n+1}$.
\end{lemma}
\begin{proof}
	We set $M=\Omega^1_{\mf{R}/\fracr}\otimes_{\mf{R}}\mf{R}_n$ to be the restriction of the (completed) K\"ahler differentials of $\mf{R}/\fracr$ to $\mf{R}_n$, so that $M=\oplus_i \mf{R}_ndx_i$.  Let $d_n$ be the composition $$ d_n: \mf{R}\overset{d}{\longrightarrow} \Omega^1_{\mf{R}/\fracr}\to \Omega^1_{\mf{R}/\fracr}\otimes_{\mf{R}}\mf{R}_n=M,$$ where $d$ is the exterior derivative and the second map above is the natural quotient map, and let $f: \mf{R}\to \mf{R}_n\oplus M$ be the map $(q_n, d_n).$ Note that this is a ring homomorphism, as $$f(ab)=(q_n(a)q_n(b), d_n(ab))=(q_n(a)q_n(b), q_n(a)d_n(b)+q_n(b)d_n(a)).$$
	
	By assumption, $f(I)=0$, i.e. $d_n(g)=0$ for all $g$ in $I$. Fixing $g\in I$, we thus have $$dg=\sum_i \frac{\partial g}{\partial x_i}dx_i$$ where $\frac{\partial g}{\partial x_i} \in \mf{m}^n$ for all $i$. But this implies $g$ is in $\mf{m}^{n+1}$, using that $\fracr$ has characteristic zero. As $g$ in $I$ was arbitrary, this completes the proof.
\end{proof}
It is immediate that it suffices to check the condition of \autoref{lemma: t1-lift-induct} for rank one free $\mf{R}_n$-modules:
\begin{corollary}\label{cor:t1-lifting-rank-1}
		Suppose that $I\subset \mf{m}^n$. Suppose that for all maps $f: \mf{R}\to \mf{R}_n[\epsilon]/\epsilon^2$ such that the composition $ \mf{R}\xrightarrow[]{f} \mf{R}_n[\epsilon]/\epsilon^2\to \mf{R}_n$ is the natural quotient map $$q_n: \mf{R}\to \mf{R}/\mf{m}^n,$$ we have that $f(I)=0$, i.e.~$f$ descends to a (unique) map $\mf{T}_i\to \mf{R}_n[\epsilon]/\epsilon^2$. Then $I\subset \mf{m}^{n+1}$.
\end{corollary}
\subsection{Deforming the Hodge filtration---last step}\label{subsection:deforming-last-step}
We maintain the notation and assumptions as in \autoref{subsection:main-technical-result}. Since $(\mscr{E}, \nabla)$ is semisimple, we may write it as 
\[
(\mscr{E}, \nabla)=\bigoplus_{i=1}^m (\mscr{F}_i, \nabla)
\]
with each $\mscr{F}_i$ being irreducible. Each $(\mscr{F}_i, \nabla)$ is equipped with a Hodge filtration, which we denote simply by $F^{\bullet}$ when the context is clear. We arrange the direct sum decomposition so that the Hodge filtrations on the $\mathscr{F}_i$ are compatible with that on $\mathscr{E}$, as in \autoref{induction_setup}.

\begin{definition}For each $i$, we  define the deformation problem  $\Def_{\mscr{F}_i, \nabla, F}$, where for an Artin local $\fracr$-algebra $A_{\mathscr{K}}$, $ \Def_{\mscr{F}_i, \nabla, F}(A_{\fracr})$ is the set of isomorphism classes of  tuples
    \[
    (s_{A_{\fracr}}, \mscr{F}_{i, A_{\fracr}}, \nabla_{A_{\fracr}}, F^{\bullet}_{A_{\fracr}}, \iota )
    \]
where    
\begin{itemize}
    \item $s_{A_{\fracr}}\in S(A_{\fracr})$ is  a point extending $s_{\mathscr{K}}$, 
    \item $(\overline{X}, D)_{A_{\fracr}}:=(\overline{X}, D)\times_S s_{\mathscr{K}},$
    \item $(\mscr{F}_{i, A_{\fracr}}, \nabla_{A_{\fracr}}, F^{\bullet}_{A_{\fracr}})$ is an isomonodromic   flat bundle equipped with Griffiths-transverse filtration on $(\overline{X}, D)_{A_{\fracr}}/s_{\mathscr{K}}$, and $\iota$ is an isomorphism
    \[
    \iota: (\mscr{F}_{i, A_{\fracr}}, \nabla_{A_{\fracr}}, F^{\bullet}_{A_{\fracr}})|_{\overline{X}}\simeq (\mscr{F}_i, \nabla, F^{\bullet}).
    \]
\end{itemize}
\end{definition}

\begin{theorem}[Smoothness of deformation space]\label{theorem:smoothness}
    For each $i$, the natural map $$\Def_{\mscr{F}_i, \nabla, F^{\bullet}}\to \widehat S_{\mathscr{K}}$$ defined by $$(s_{A_{\fracr}}, \mscr{F}_{i, A_{\fracr}}, \nabla_{A_{\fracr}}, F^{\bullet}_{A_{\fracr}}, \iota )\mapsto s_{A_{\fracr}}$$ is an isomorphism.
\end{theorem}
\begin{proof}
    It is straightforward to check that $\Def_{\mscr{F}_i, \nabla, F^{\bullet}}$ is representable by a noetherian $\fracr$-algebra, say $\mf{T}_i$, using that $(\mathscr{F}_i, \nabla)$ is simple. We are studying a map $$\Def_{\mscr{F}_i, \nabla, F^{\bullet}} \rightarrow \widehat{S}_{\fracr},$$ and the latter is the spectrum of a smooth  noetherian complete local $\fracr$-algebra $\mf{R}$; we write $r_i: \mf{R}\rightarrow \mf{T}_i$ for the corresponding map on algebras.\footnote{Recall from \autoref{subsection:main-technical-result} that $\widehat{S}_{\mathscr{K}}$ is the completion of $S_{\mathscr{K}}$ at $s_{\mathscr{K}}$.} Let $\mathfrak{m}_{\mf R}, \mf{m}_{\mf T_i}$ be the maximal ideals of $\mf R, \mf T_i$, respectively, and set $\mf R_n=\mf R/\mf m_{\mf R}^n$ as in \autoref{subsection:T1-lifting}. As the Hodge-de Rham spectral sequence degenerates for each $\mathscr{F}_i$, any deformation of $(\mathscr{F}_i, \nabla)$ over $\mathscr{K}[\epsilon]/\epsilon^2$ carries a unique Griffiths transverse filtration on an irreducible flat bundle if it exists, by \autoref{cor:griffiths-transverse-def-unique}; that is, the map $$\Def_{\mscr{F}_i, \nabla, F^{\bullet}}(\mathscr{K}[\epsilon]/\epsilon^2) \rightarrow \widehat{S}_{\fracr}(\mathscr{K}[\epsilon]/\epsilon^2)$$ is injective. We deduce that $r_i$ is a surjection on cotangent spaces $\mathfrak{m}_{\mf R}/\mf{m}_{\mf{R}}^2\to \mf{m}_{\mf T_i}/\mf{m}_{\mf T_i}^2$ . It is hence surjective as a map of rings; let $I_i:=\ker(r_i)$. We now prove by induction on $n$ that $I_i\subset \mf{m}_{\mf{R}}^n$. 

    We will do this simultaneously for the  $\Def_{\mscr{F}_i, \nabla, F^{\bullet}}$'s: that is, we assume that $I_i\subset \mf{m}^n$ for all $i$, and we aim to deduce $I_i\subset \mf{m}^{n+1}$ for all $i$ using \autoref{lemma: t1-lift-induct}, as we now spell out.

By the inductive hypothesis, we know that for each $i$ the natural closed embedding $q_n: \on{Spec}(\mf R_n)\to \widehat S_{\mathscr{K}}$ factors through some (unique) map $\sigma_i: \Spec(\mf{R}_n) \rightarrow \Spec(\mf{T}_i)$. By \autoref{cor:t1-lifting-rank-1}, it suffices to show that for any morphism $f: \on{Spec}(\mf R_n[\epsilon]/\epsilon^2)\to \widehat S_{\mathscr{K}}$ such that the composition $$\on{Spec}(\mf R_n)\to  \on{Spec}(\mf R_n[\epsilon]/\epsilon^2) \xrightarrow{f} \widehat S_{\mathscr{K}}$$ is $q_n$, there exists $f_{\mf T_i}: \on{Spec}(\mf R_n[\epsilon]/\epsilon^2)\to \mf T_i$ such that the diagram 
\[
    \begin{tikzcd}
        \Spec(\mf R_n[\epsilon]/\epsilon^2) \arrow[r, "f_{\mf{T}_i}"] \arrow[dr, "f"] & \Spec(\mf{T}_i) \arrow[d] \\
        & \widehat{S}_{\mathscr{K}}
    \end{tikzcd}
    \]
commutes.  In other words, let $(\mathscr{F}_{i, n}, \nabla, F^\bullet)$ be the filtered flat bundle on $(\overline{X}, D)_{\mf R_n}$ corresponding to $\sigma_i$; it remains to show that each $(\mathscr{F}_{i, n}, \nabla, F^\bullet)$ deforms to $(\overline{X}, D)_{\mf R_n[\epsilon]/\epsilon^2}/\on{Spec}(\mf R_n[\epsilon]/\epsilon^2)$ in an isomonodromic and Griffiths-transverse manner, where $(\overline{X}, D)_{\mf R_n[\epsilon]/\epsilon^2}$ denotes the fiber product $(\overline{X}, D)\times_{\widehat{S}_\mathscr{K}, f} \on{Spec}(\mf R_n[\epsilon]/\epsilon^2).$

Let 
\[
(\mscr{E}_{A_{\fracr}}, \nabla, F^{\bullet}):= \bigoplus (\mscr{F}_{i, n}, \nabla, F^{\bullet}).
\]
We now set up the situation so that the assumptions of \autoref{induction_setup} are satisfied. As in \autoref{induction_setup} replace $R$ with a finitely-generated localization so that $\widehat S\simeq R[[t_1, \cdots, t_{\dim S}]]$; fix such a set of local coordinates. After replacing $R$ with a further finitely-generated localization, we may assume that $f$ and each $(\mathscr{F}_{i,n}, \nabla, F^\bullet)$ descends to a Griffiths-transverse filtered flat bundle $(\mathscr{F}_{A,i}, \nabla, F^\bullet)$ on $A=R[[t_1, \cdots, t_{\dim S}]]/(t_1, \cdots, t_n)^n$; note that $A$ is cut out by a monomial ideal, so \autoref{assumption:monomial-A} is satisfied. The first part of \autoref{assumption: stable-special-fiber} and \autoref{assumption: hodge-filtration} is immediate from our setup---set $$(\mathscr{E}_A, \nabla, F^\bullet)=\bigoplus (\mathscr{F}_{A, i}, \nabla, F^\bullet)$$ with the induced connection and Griffiths-transverse filtration. The second part  of \autoref{assumption: stable-special-fiber} may be guaranteed by yet again replacing $R$ by a finitely-generated localization. Finally, we make a choice of $F_0$ as in \autoref{induction_setup} so that \autoref{assumption:monomial-frobenius} is satisfied.

We let  $\mf{p}$ be a maximal ideal of $R$ of sufficiently  large residue characteristic, and, as in \autoref{induction_setup}, let $\widehat R$ be the completion of $R$ at $\mf p$. Note that e.g.~by \autoref{assumption:monomial-frobenius}  we may assume that  there exists a lift of Frobenius $F_{A_{\widehat R}}: \on{Spf}(A_{\widehat R})\rightarrow \on{Spf}(A_{\widehat R})$ such that the diagram 
    \[
    \begin{tikzcd}
  \on{Spf}(A_{\widehat R}) \arrow[r, "F_{A_{\widehat R}}"] \arrow[d]
    & \on{Spf}(A_{\widehat R}) \arrow[d, ""] \\
  \hat{S} \arrow[r, "F_0"]
&  \hat{S} \end{tikzcd}
\]
commutes.  Note that this map factors through the restriction of $F_0$ to $B_{\widehat R}$.

By \autoref{lemma:pcurv-vanish}(1), $(\mscr{E}_{A_{\mf p}}, \nabla, F^\bullet)$ is a Higgs-de Rham fixed point, and hence the Hodge-de Rham spectral sequence for $(\mscr{E}_{A_{\mf p}}, \nabla, F^\bullet)$  degenerates by \autoref{lemma: degen-hdr-fixed} below. Hence the same is true for $(\mathscr{E}_A, \nabla, F^\bullet)$ and hence for each $(\mathscr{F}_{A, i}, \nabla, F^\bullet)$. 

Now the result is immediate from \autoref{cor:Hodge-filtration-extends}, which provides the required Griffiths-transverse extension of the Hodge filtration on $(\mathscr{F}_{i,n}, \nabla)$ to its isomonodromic deformation to $(\overline{X}, D)_{\mf R_n[\epsilon]/\epsilon^2}/\on{Spec}(\mf R_n[\epsilon]/\epsilon^2)$, as desired.

\end{proof}

\begin{proof}[Proof of \autoref{theorem:extending-the-Hodge-filtration}]
The proof is immediate from \autoref{theorem:smoothness}. Indeed, we've shown that for each $i$ the natural map $$\Def_{\mscr{F}_i, \nabla, F^{\bullet}} \rightarrow \widehat{S}_{\fracr},$$ is an isomorphism. The inverse of this map precisely yields a Griffiths-transverse deformation of the Hodge filtration on $\mscr{F}_i$ to its isomonodromic deformation to $\widehat{S}_{\fracr}$. Taking the direct sum of these filtered flat bundles yields the desired statement.
\end{proof}

\subsection{Hodge-de Rham degeneration for Higgs-de Rham fixed points}
In this section we record a basic fact about Hodge-de Rham degeneration for Griffiths-transverse filtered flat bundles which are Higgs-de Rham fixed points, or in more classical language, Fontaine modules, that were used in the proof of \autoref{theorem:smoothness} above. These facts are surely well-known to experts but we were unable to find a suitable reference.

Let $k$ be a perfect field of characteristic $p>0$, and $A$ an Artin local $k$-algebra with residue field $k$. Suppose  $\overline{f}: \overline{X}\rightarrow \Spec(A)$ is  a smooth projective morphism, and $D\subset \overline{X}$ is a normal crossings divisor relative to $\Spec(A)$; we write $X$ for $\overline{X}\setminus D$ and $f: X\rightarrow \Spec(A)$ for the natural map. We form the Cartesian square
\begin{equation}
    \begin{tikzcd}
  (\overline{X}', D') \arrow[r, "G"] \arrow[d]
    & (\overline{X}, D) \arrow[d, "\bar{f}"] \\
  \Spec(A) \arrow[r, "F_{\text{abs}}"]
&  \Spec(A). \end{tikzcd}
\end{equation}

The following is a version of a result of Deligne--Illusie \cite[Prop 4.1.2]{deligne-illusie} for Higgs-de Rham fixed points:
\begin{lemma}\label{lemma: degen-hdr-fixed}
    Suppose $(\mscr{E}, \nabla)$ is a flat bundle on $(\overline{X}, D)/\Spec(A)$, and $F$ a Griffiths transverse filtration and nilpotent residues along $D$; we denote by $(\gr_F \mscr{E}, \theta)$ the associated graded Higgs bundle, which we assume lies in $\on{HIG}_{p-1}(\overline{X}, D)$. 
    
    Suppose  we have  a $W_2(k)$-lift $ (\bar{X}_{W_2}', D_{W_2}') \rightarrow \Spec(\tilde{A})$ of $(\bar{X}', D')\rightarrow \Spec(A)$, which defines the inverse Cartier transform $C^{-1}$ as in \autoref{subsection:inverse-Cartier}. Suppose the length of $s$ is less than $p-\dim_A \overline{X}$.
    
    Suppose $(\mscr{E}, \nabla, F)$ is a Higgs-de Rham fixed point relative to this lift, in the sense that 
    \begin{equation}\label{eqn:hdr-fixed-point}
    C^{-1}G^*(\gr_F \mscr{E}, \theta)\simeq (\mscr{E}, \nabla).
    \end{equation}
    
    Then the Hodge-de Rham spectral sequence for $\mscr{E}_{dR}$ 
    \[
    E_1^{p,q} = \mb{H}^q(\gr^p_F(\mscr{E}_{dR})) \Longrightarrow \mb{H}^{p+q}(\mscr{E}_{dR})
    \]
    degenerates at the $E_1$-page, and the terms of the $E_1$-page are free $A$-modules, and the formation of $\mb{H}^q(\gr^p_F(\mscr{E}_{dR})), \mb{H}^{p+q}(\mscr{E}_{dR})$ commutes with arbitrary base change.
\end{lemma}

\begin{proof}

    The proof is essentially the same as the proof of \cite[Prop 4.1.2]{deligne-illusie}. We restrict ourselves to the case when $D$ is empty and therefore $\overline{X}=X$; the general case is only notationally more complicated. 

    Let $\mf{m}$ be the maximal ideal of $A$, and $N$ be the smallest integer such that $\mf{m}^N=0$. We prove the proposition by induction on $N$. The base case $N=1$, i.e. $A=k$, follows immediately from \cite[Theorem 4.17]{ogus-vologodsky} (and from \cite[Lemma 6.1]{krishnamoorthy2020deformation} in general).\footnote{In our application, the filtered flat bundle in the base case comes from a geometric situation which lifts to $W_2(k)$, in which case this is due to Deligne-Illusie \cite{deligne-illusie}.}  Let $\iota: \Spec(k)\rightarrow \Spec(A)$ be the inclusion of the closed point, and write $\iota: X_k:=X\times_{\Spec(A)}\Spec(k) \rightarrow X$. 
    Let 
    \[
    h^{p,q}:= \dim_k \mb{H}^{p+q}(\gr^p_F \iota^*(\mscr{E}_{dR})).
    \]

    We now proceed with the induction step. The Hodge-de Rham spectral sequence \eqref{eqn:hdr-fixed-point} implies that
    \begin{equation}\label{eqn: pf-hdr-fixed-1}
    \lg \mb{H}^{p}(\mscr{E}_{dR})\leq \sum_{p+q=n} 
\mb{H}^{p+q}(\gr^p_F(\mscr{E}_{dR})),
\end{equation}
where $\lg$ stands for the length; the spectral sequence degenerates if and only if this is an equality.

For each $p$,     we may write $Rf_{*}(\gr^p(\mscr{E}_{dR}))$ as a bounded complex of finitely-generated flat, and hence free, modules. It follows from basechange that 
    \begin{equation}\label{eqn: pf-hdr-fixed-2}
   \lg  \mb{H}^{p+q}(\gr^p_F(\mscr{E}_{dR})) \leq h^{p,q} \lg A,
    \end{equation}
with equality if and only if the cohomology groups are free over $A$.
On the other hand, since $(\mscr{E}, \nabla)$ is a Higgs-de Rham fixed point, by \autoref{prop:conjugate-ss-degenerates-e2},  the conjugate spectral sequence 
\[
E_2^{p,q}= \mb{H}^{p+q}(G^*\gr^q_{\text{conj}} (\mscr{E}_{dR})) \Longrightarrow 
\mb{H}^{p+q}(\mscr{E}_{dR})
\]
degenerates, where we define the conjugate filtration using the inverse Cartier transform as in \autoref{defn:lifting-and-conjugate-filtration}; therefore 
\[
\lg \mb{H}^{p+q}(\mscr{E}_{dR}) = 
\sum_{p+q=n} \mb{H}^{p+q}(G^*\gr^q_{\text{conj}} (\mscr{E}_{dR})).
\]
Now the absolute Frobenius on $A$ factors as 
\[
\Spec(A)\rightarrow T \xhookrightarrow{j} \Spec(A)
\]
where $j: T\xhookrightarrow{} \Spec(A)$ is a proper closed subscheme. It follows by our inductive hypothesis that $\mb{H}^{p+q}(j^*\gr^q_{\text{conj}}(\mscr{E}_{dR}))$ are free over $T$ of rank $h^{q,p}$. Hence by basechange $\mb{H}^{p+q}(G^*\gr^q_{\text{conj}} (\mscr{E}_{dR}))$ is free over $A$ of rank $h^{q,p}$. Therefore 
\[
\lg \mb{H}^{p+q}(\mscr{E}_{dR}) = \sum_{p+q=n} h^{p,q} \lg A.
\]
Combining this with \eqref{eqn: pf-hdr-fixed-1}, \eqref{eqn: pf-hdr-fixed-2}, we deduce 
\begin{equation}
    \lg \mb{H}^{n}(\mscr{E}_{dR})= \sum_{p+q=n} 
\mb{H}^n(\gr^p_F(\mscr{E}_{dR})),
\end{equation}
and hence the Hodge-de Rham spectral sequence degenerates. Moreover, equality holds in \eqref{eqn: pf-hdr-fixed-2}, which means that the cohomology groups $\mb{H}^{n}(\gr^p_F(\mscr{E}_{dR}))$ are free, and hence the same is true for $\mb{H}^{p+q}(\mscr{E}_{dR})$, thus completing the induction step.
\end{proof}

%\begin{corollary}
%Under the same assumptions as in \autoref{lemma: degen-hdr-fixed},   for any $A$-module $M$, the Hodge-de Rham spectral sequence for $\mscr{E}\otimes M \otimes \Omega^{\bullet}_{\bar{X}/A}(\log D)$ 
%    \[
%    E_{1, M}^{p,q} = \mb{H}^q(\gr^p(\mscr{E}\otimes M \otimes  \Omega^{\bullet}_{\bar{X}/A}(\log D))) \Longrightarrow \mb{H}^{p+q}(\mscr{E}\otimes \Omega^{\bullet}_{\bar{X}/A}(\log D))=E_{\infty, M}^{p+q}
%    \]
%    degenerates at the $E_1$-page.  Moreover, 
%    \[
%    E_{1, M}^{p,q}=E_{1, A}^{p,q}\otimes M,\ \ E_{\infty, M}^{p+q}= E_{\infty, A}^{p+q}\otimes M.
%    \]
%\end{corollary}
%\begin{proof}
%    Again, we deal only with the case when $D$ is empty and leave the general case to the reader. By \autoref{lemma: degen-hdr-fixed}, for any $p$,  the complexes $\gr^p(\mscr{E}\otimes \Omega^{\bullet}_{X/A}), \ \ \mscr{E}\otimes \Omega^{\bullet}_{X/A}$ have free cohomologies. Therefore the formation of their  cohomologies commutes with tensoring with $M$, from which we deduce both the degeneration of the spectral sequence and the formulas for $E_{1,M}^{p,q}, \ E_{\infty, M}^{p+q}$.
%\end{proof}

\section{Proof of main  theorem}\label{section:proof-of-main-theorem}
We now prove \autoref{thm:NA-main}; it will, more or less, be an immediate consequence of \autoref{theorem:extending-the-Hodge-filtration} after some simple reductions. 
\subsection{Reduction to the Deligne canonical extension}
The following lemma will be useful for us for the reduction to the case of nilpotent residues.

\begin{lemma}\label{lemma:deligne-canonical-fix}
Let $R\subset \mathbb{C}$ be a finitely-generated $\mathbb{Z}$-algebra, $S$ a smooth $R$-scheme, $s\in S(R)$, and $(\overline{X},D)$ a smooth projective $S$-scheme equipped with a simple normal crossings divisor $D$ over $S$. Let $(\mathscr{E},\nabla)$ 	be a flat bundle on $(\overline{X},D)_s$, and suppose it admits an ($\omega(p)$-integral) formal isomonodromic deformation over $\widehat S$, the completion of $S$ at $s$. Suppose the eigenvalues of the residue matrices of $(\mathscr{E},\nabla)$ along $D_s$ are rational numbers.

Then there exists a flat bundle $(\mathscr{E}',\nabla')$ on $(\overline{X},D)_s$ with $(\mathscr{E}',\nabla')|_{\overline{X}_s\setminus D_s}\simeq (\mathscr{E}',\nabla')|_{\overline{X}_s\setminus D_s}$ and whose residue matrices have eigenvalues in $[0,1)$. Moreover $(\mathscr{E}',\nabla')$ has ($\omega(p)$-integral) formal isomonodromic deformation over $\widehat S$.
\end{lemma}
\begin{proof}
The point of the proof is that one may modify $(\mathscr{E},\nabla)$ to along $D$ to obtain the desired condition on residues; as one can do this integrally and in families, the claim about ($\omega(p)$-integral) formal isomonodromic deformation follows. We briefly explain how to perform this modification. Let $D_i=\cup D_i$, with the $D_i$ irreducible.

First, we may without loss of generality assume all the eigenvalues of the residue matrices of $(\mathscr{E},\nabla)$ are negative, by tensoring with $\mathscr{O}(ND)$ for $N\gg 0$. Let $\alpha$ be the most negative eigenvalue that appears, say along some component $D_j$ of $D$, and consider the map $$\mathscr{E}\xrightarrow{\nabla} \mathscr{E}\otimes \Omega^1_{\overline{X}_s}(\log D_s)/\Omega^1_{\overline{X}_s}\simeq \bigoplus_{i} \mathscr{E}|_{D_i} \to \mathscr{E}|_{D_j}\xrightarrow{(\on{Res}_{\nabla, D_i}-\alpha)^r} \mathscr{E}|_{D_j}$$
for $r$ sufficiently large. Replacing $\mathscr{E}$ by the kernel of this map has the effect of increasing $\alpha$ by $1$. Iterating this procedure we may move all the eigenvalues of the residue matrices into $[0,1)$. As this may evidently be done over the base of our ($\omega(p)$-)integral isomonodromic deformation the proof is complete.
\end{proof}

\subsection{The proof}
We now proceed with the proof of our main theorem.
\begin{proof}[Proof of \autoref{thm:NA-main}]
We first recall the setup. We are given $R\subset \mathbb{C}$ a finitely-generated $\mathbb{Z}$-algebra, $f: \overline{X}\to S$ a smooth projective morphism of smooth $S$-schemes, $D\subset \overline{X}$ a simple normal crossings divisor over $S$, and $s\in S(R)$. We have a Picard-Fuchs equation $(\mathscr{E},\nabla)$ on $(\overline{X}_s, D_s)$ with ($\omega(p)$-)integral formal isomonodromic deformation. Let $X=\overline{X}\setminus D$, and let $\mathbb{V}$ be the complex local system on $X_{s, \mathbb{C}}^{\text{an}}$ corresponding to $\mathscr{E},\nabla$.

By \autoref{prop:finite-orbit-criterion} it suffices to show that the orbit of the isomorphism class of $\mathbb{V}$ under $\pi_1(S_{\mathbb{C}}^{\text{an}}, s)$ is finite. Choose a finite \'etale cover $Y$ of $X_{s,\mathbb{C}}$ such that $\mathbb{V}|_{Y^{\text{an}}}$ has unipotent monodromy. By spreading out, we may find a diagram 
$$\xymatrix{
(\overline{X}', D') \ar[r] \ar[d]&  (\overline{X}, D) \ar[d]\\
S' \ar[r] & S
}$$
over $R$, whose base change to $\mathbb{C}$ satisfies the hypotheses of \autoref{lemma:passing-to-cover-trick} (where we set $X'={\overline{X}}'\setminus D'$)
 and $s'\in S'$ such that $X'_{s', \mathbb{C}}=Y$. Pulling back to $X'$ (which preserves ($\omega(p)$-)integral formal isomonodromic deformation) it thus, by \autoref{lemma:passing-to-cover-trick} to check that the orbit of $\mathbb{V}$ under $\pi_1(S', s')$ is finite. We thus replace $X, S$ with $X', S'$; note that $X'$ comes equipped with an snc compactification over $S$, with which we replace $\overline{X}$. After passing to a further cover of $S$, we may assume $X\to S$ admit a section, and so is ``good" in the sense of \cite{esnault2024non} or \autoref{theorem:esnault-kerz}.

Unfortunately the pullback above may not preserve the property of being a Picard-Fuchs equation; while the pullback of $\mathbb{V}$ will still arise from the cohomology of a family of varieties, the residues of the pullback $(\mathscr{E},\nabla)$ may no longer have eigenvalues in $[0,1)$. But we modify it as in \autoref{lemma:deligne-canonical-fix}. It thus suffices to treat the case where $(\mathscr{E},\nabla)$ has nilpotent residues along $D$.

By \autoref{theorem:esnault-kerz} it suffices to show that the Hodge filtration on $(\mathscr{E},\nabla)$ extends Griffiths-transversally to its formal isomonodromic deformation over $\widehat S_\mathbb{C}$, the formal scheme obtained by completing $S_{\mathbb{C}}$ at $s_\mathbb{C}$. But this is immediate from \autoref{theorem:extending-the-Hodge-filtration}.
\end{proof}

\section{Relationship to the $p$-curvature conjecture}\label{section:isomonodromy-and-groth-katz}
We now observe that \autoref{conj:main-conjecture} implies the classical Grothendieck-Katz $p$-curvature conjecture \cite{katz-p-curvature};  in fact, to prove the $p$-curvature conjecture for all bundles of rank $r$, on all varieties, it suffices to to verify \autoref{conj:main-conjecture} for $\mathscr{M}_{dR}(\mathscr{C}_g/\mathscr{M}_g, r)$ for $g\gg 0$. Here $\mathscr{M}_g$ is the moduli space of smooth projective curves of genus $g$, and $\mathscr{C}_g$ is the universal curve.

We first show that flat bundles whose $p$-curvature vanishes for almost all $p$ admit $\omega(p)$-integral isomonodromic deformations.

\begin{lemma}\label{lemma:p-curvature-and-omega-p-integrality}
	Let $R\subset \mathbb{C}$ be a finitely-generated $\mathbb{Z}$-algebra, $\overline{X}\to S$ a smooth projective morphism, and $D\subset \overline{X}$ a simple normal crossings divisor over $S$. Let $s\in S(R)$ be a point, and $(\mathscr{E},\nabla)$ a flat bundle on $\overline{X}_s$ with regular singularities along $D$.
	
	Suppose that for almost all primes $p$, the $p$-curvature of $(\mathscr{E},\nabla)$ vanishes mod $p$. Then $(\mathscr{E},\nabla)$ admits an integral isomonodromic deformation.
\end{lemma} 
\begin{proof}
This is immediate from the proof of \autoref{prop:isomonodromy-existence-uniqueness}. We follow the recipe of that proposition, keeping track of the denominators that appear in the ``gluing data." 

Namely, let $\widehat S$ be the formal $R$-scheme obtained by completing $S$ at $s$, and let $\{U_i\}$ be a cover of $\overline{X}_{\widehat S}$ by formal affine opens. Let $D_i=D\cap U_i$. As $(\overline{X}, D)$ is a smooth simple normal crossings pair over $S$, we have that each $(U_i, D_i)\simeq (U_{i, s}, D_{i,s})\times_R \widehat S$  Fix such isomorphisms and set $(\mathscr{E}_i, \nabla_i)$ to be the flat bundle with absolute connection on $(U_i, D_i)$ obtained by pulling back $(\mathscr{E}, \nabla)$ along the projection to $(U_{i, s}, D_{i,s})$. We will try to glue these bundles together over $U_{ij}=U_i\cap U_j$ via the formula in \autoref{lemma:char-0-taylor-series}. Note that the $(\mathscr{E}_i, \nabla_i)$ have vanishing $p$-curvature for almost all $p$, as they are pullbacks of a bundle with the same property.

Choose any section $u$ to $\on{Hom}(\mathscr{E}_i|_{U_{ij}},\mathscr{E}_j|_{U_{ij}})$ restricting to the tautological isomorphism over $U_{ij,s}$. Consider the sum $$\sum_{i_1, \cdots, i_n} (-1)^{\sum_j i_j}\frac{t_1^{i_1}\cdots t_n^{i_n}}{i_1!\cdots i_n!}\left(\prod_{j=1}^n\nabla_{\frac{\partial}{\partial t_j}}^{i_j}\right)(u),$$ as in  \autoref{lemma:char-0-taylor-series}, where the $t_i$ are coordinates on $\widehat S\simeq R[[t_1, \cdots, t_d]]$. By the vanishing of $p$-curvature we have $\nabla_{\frac{\partial}{\partial t_j}^p}\equiv 0 \bmod p$. Thus for a multi-index $I$, the coefficient of $t^I$ is $p$-integral as long as each $i_j<p^2$. This gives the desired gluing data over the ring of $g$-integral elements of $R[[t_1, \cdots, t_d]]$, where $g(p)=p^2-1$ for all sufficiently large $p$, defined as in \autoref{subsection:integrality-of-leaves}. As $$\lim \frac{p^2-1}{p}=\infty,$$ this completes the proof.
\end{proof}

\begin{proposition}\label{prop:conj-for-Mg-implies-groth-katz}
	Suppose \autoref{conj:main-conjecture} holds for all $\mathscr{M}_{dR}(\mathscr{C}_g/\mathscr{M}_g, r)$, for all $g\gg 0$. Then the Grothendieck-Katz $p$-curvature conjecture holds for all flat bundles of rank $r$ on all smooth quasi-projective varieties.
\end{proposition}
\begin{proof}
	It is well known that to prove the Grothendieck-Katz conjecture for all bundles of rank $r$, it suffices to do so on smooth proper curves. We first sketch this reduction, and then prove the proposition.
	
	Let $R\subset \mathbb{C}$ be a finitely-generated $\mathbb{Z}$-algebra and $X$ be a smooth quasi-projective $R$-scheme. Let $(\mathscr{E},\nabla)$ a flat bundle on $X$ of rank $r$, almost all of whose $p$-curvatures vanish mod $p$. As flat bundles with vanishing $p$-curvatures have regular singularities at infinity \cite[Theorem 13]{katz1970nilpotent}, it suffices to prove that $(\mathscr{E},\nabla)_\mathbb{C}^{\text{an}}$ has finite monodromy. 
	By the Lefschetz hyperplane theorem, it suffices to do so after restricting to a general sufficiently ample curve in $X$, so we may assume $X$ is a curve. Again by \cite[Theorem 13]{katz1970nilpotent}, $(\mathscr{E},\nabla)_\mathbb{C}^{\text{an}}$ has finite monodromy at infinity, so by passing to a cover we may assume $(\mathscr{E},\nabla)$ extends to a smooth projective completion of $X$. Note that these various reductions have not changed the rank of $\mathscr{E}$.
	
	We now have a flat bundle with vanishing $p$-curvatures on some smooth projective curve; it suffices to show that it has finite monodromy. Choose a finite \'etale cover $f: Y\to X$, where $Y$ has genus $g\geq r^2$, and pull $(\mathscr{E},\nabla)$ back to $Y$. By \autoref{lemma:p-curvature-and-omega-p-integrality}, $f^*(\mathscr{E},\nabla)$ 
	admits an $\omega(p)$-integral isomonodromic deformation over the formal neighborhood of $[Y]$ in $\mathscr{M}_g$. Hence, assuming \autoref{conj:main-conjecture}, $f^*(\mathscr{E},\nabla)$ has algebraic isomonodromic deformation, and hence the conjugacy class of its monodromy representation has finite orbit under $\pi_1(\mathscr{M}_g)$, by \cite[Proposition 2.1.3]{landesman2022canonical}. Hence, by \cite[Theorem 1.2.1]{landesman2022canonical}, $f^*(\mathscr{E},\nabla)$ has finite monodromy, completing the proof.
\end{proof}
\begin{remark}
	Instead of relying on \cite[Theorem 1.2.1]{landesman2022canonical}, we could have instead used \cite[Theorem 1.1.1]{lawrence2022representations}, whose proof is much easier. Indeed, under the assumption that \autoref{conj:main-conjecture} holds true for $\mathscr{M}_{dR}(\mathscr{C}_g/\mathscr{M}_g, r)$, we have that the monodromy representation of $f^*(\mathscr{E},\nabla)$ has finite orbit under $\pi_1(\mathscr{M}_{g'})$ for \emph{every} finite \'etale $f: Y\to X$, with $Y$ of genus $g'$. But this is precisely the hypothesis of \cite[Theorem 1.1.1]{lawrence2022representations}.
\end{remark}

\section{A motivic variational principle}\label{section-variational-conjectures}

In this section we briefly discuss some of the philosophy behind the proof of \autoref{thm:NA-main} and its relationship with some of the basic open conjectures on local systems in algebraic geometry, as well as how it implies several predictions of those conjectures. In some sense the latter half of this paper (on isomonodromy) is a ``non-abelian" analogue of (versions of) the variational Hodge or Tate conjectures; compare to \autoref{thm:cycle-class-initial-conditions} for one ``abelian" perspective on the same story. 

\subsection{The relative Fontaine-Mazur conjecture}
We first recall the statement of the relative Fontaine-Mazur conjecture, and draw out some of its consequences. See the introduction to \cite{liu2017rigidity} for some further discussion of this conjecture, and \cite{petrov2023geometrically} for the version of it stated here. The conjecture aims to characterize local systems of geometric origin.

Let $k$ be a finitely-generated field, $\overline{k}$ a separable closure of $k$, $\ell$ a prime different from the characteristic of $k$, and $\Lambda$  a finite extension of $\mathbb{Q}_\ell$. Let $X$ be a smooth and geometrically connected $k$-scheme.  

\begin{definition}\label{defn:geometric-origin}
	 A lisse $\Lambda$-sheaf $\mathbb{V}$ on $X_{\overline{k}}$ is \emph{of geometric origin} if there exists a dense open subset $U\subset X_{\overline{k}}$, a smooth projective $\pi: Y\to U$, and an integer $i$ such that $\mathbb{V}|_U$ is a subquotient of $R^i\pi_*\Lambda$.
\end{definition}
To first approximation one may ignore the restriction to $U$; then for a local system to be of geometric origin means that it appears in the cohomology of a family of smooth projective varieties over $X_{\overline{k}}$. One could replace the word ``subquotient" in the definition above with ``summand," as the local system $R^i\pi_*\Lambda$ is necessarily semisimple

\begin{conjecture}[Relative Fontaine-Mazur conjecture]\label{conjecture:relative-fontaine-mazur}
	Let $\mathbb{V}$ be a semisimple lisse $\Lambda$-sheaf on $X_{k'}$ for some $k'/k$ finite. Then $\mathbb{V}|_{X_{\overline k}}$ is of geometric origin.
\end{conjecture}
In other words, semisimple local systems on $X_{\overline{k}}$  which descend to $X$ (or $X_{k'}$, for $k'/k$ finite) are expected to be of geometric origin. The conjecture is basically open except in the case that $X$ is a curve and $k$ is finite, where it follows from the Langlands program for function fields. The converse to \autoref{conjecture:relative-fontaine-mazur} follows from the fact that if $Y$ appears in the cohomology of some $\pi:Y\to U$, we may descend $\pi$ to some finitely-generated extension of $k$.

The condition that a semisimple local system $\mathbb{V}$ descends to $X_{k'}$ for some $k'$ has a dynamical interpretation:~it is the same as saying that, under the natural action of $\on{Gal}(\overline{k}/k)$ on local systems on $X_{\overline{k}}$, there is an open subgroup stabilizing the isomorphism class of $\mathbb{V}$ \cite[Proposition 3.1.1]{litt2021arithmetic}. We refer to such a local system $\mathbb{V}$ as \emph{arithmetic}. The following simple consequence is of primary interest to us here:
\begin{proposition}\label{proposition:variational-arithmeticity}
	Let $X\to S$ be a smooth proper morphism of smooth $k$-schemes, and let $\mathbb{V}$ be a lisse $\Lambda$-sheaf on $X_{\overline{k}}$. Suppose that for some $s\in S(k)$, $\mathbb{V}|_{X_{s, \overline{k}}}$ is absolutely irreducible and arithmetic. Then $\mathbb{V}$ is arithmetic.
\end{proposition}
\subsection{Formal and variational meta-conjectures}
The following conjecture is thus a prediction of \autoref{conjecture:relative-fontaine-mazur}:
\begin{conjecture}[A variational motivic conjecture]\label{conjecture:variational-motivic}
	Let $f: X\to S$ be a smooth proper morphism of smooth $\overline{k}$-schemes, and let $\mathbb{V}$ be a semisimple local system on $X$. Suppose there exists $s\in S(\overline{k})$ such that $\mathbb{V}|_{X_s}$ is of geometric origin. Then for all $t\in S(\overline{k})$, $\mathbb{V}|_{X_t}$ is of geometric origin.
\end{conjecture}
In the setting of the relative Fontaine-Mazur conjecture one may take ``local system" in \autoref{conjecture:variational-motivic} to mean a lisse $\Lambda$-sheaf. But it makes sense for more general categories of coefficients---in particular, for flat bundles or $\mathbb{C}$-local systems, when $\overline{k}=\mathbb{C}$. 

While this conjecture appears to be completely out of reach, it makes a number of concrete predictions, of the following form. Let $P$ be a property of a local system that is conjecturally equivalent to being of geometric origin. Then we have the following meta-conjecture:
\begin{conjecture}[Variational meta-conjecture]\label{conjecture:meta-conjecture}
	Let $f: X\to S$ be a smooth proper morphism of smooth $\overline{k}$-schemes, and let $\mathbb{V}$ be a semisimple local system on $X$. Suppose there exists $s\in S(\overline{k})$ such that $\mathbb{V}|_{X_s}$ has property $P$. Then for all $t\in S(\overline{k})$, $\mathbb{V}|_{X_t}$ has property $P$.
\end{conjecture}
For example, one might take $P$ to be one of the following:
\begin{itemize}
	\item If one takes ``local system" to mean ``lisse $\Lambda$-adic sheaf," one may take $P$ to be the property of arithmeticity, discussed above. In this case the meta-conjecture is verified by \autoref{proposition:variational-arithmeticity}.
	\item If  $k=\mathbb{C}$ and one takes local system to mean ``flat bundle" or $\mathbb{C}$-local system, one may take $P$ to be the property of ``underlying a polarizable $\mathbb{\mathscr{O}_K}$-variation of Hodge structure for some number field $K$." This is conjecturally the same as being of geometric origin by \cite[Conjecture 12.4]{simpson62hodge}. In this case \autoref{conjecture:meta-conjecture} has been studied (and resolved) in \cite{katzarkov1999non}, \cite{jost2001harmonic}, and \cite{esnault2024non}. See also \autoref{subsection:non-abelian-p,p} below.
	\item If $k$ is finitely-generated of characteristic zero, and one takes local systems to mean ``flat bundles," one may take $P$ to be the property ``admits a spreading-out over a finitely-generated $\mathbb{Z}$-algebra such that the reduction mod $p$ has nilpotent $p$-curvature for almost all $p$." This is has been conjectured by Andr\'e to be equivalent to being of geometric origin, see \cite[Appendix to Chapter V]{andre1989g}. See \autoref{subsection:variational-Andre} below for some results in this direction. 
\end{itemize}
 The main purpose of this section is to observe that \autoref{conj:main-conjecture}, in conjunction with \autoref{conjecture:meta-conjecture}, implies a \emph{formal} analogue of \autoref{conjecture:meta-conjecture}. Namely, let $P$ be a property of flat bundles arising as (summands of) the algebraic de Rham cohomology of families of smooth projective varieties.
 
 Let $f: X\to S$ be a smooth proper morphism of smooth $R$-schemes, where $R\subset \mathbb{C}$ is a finitely generated $\mathbb{Z}$-algebra, and fix $s\in S(R)$. Let $\mathscr{K}$ be the fraction field of $R$, and let $\widehat S$ be the formal scheme obtained by completing $S$ at $s$, and $\widehat{S}_\mathscr{K}$ the formal scheme obtained by completing $S_{\mathscr{K}}$ at $s_{\mathscr{K}}$.
 \begin{conjecture}[Formal meta-conjecture]\label{conjecture:formal-meta-conjecture}
	Let $(\mathscr{E},\nabla)$ be a Picard-Fuchs equation, and suppose $(\mathscr{E},\nabla)$  admits an ($\omega(p)$-)integral isomonodromic deformation over $\widehat S$. Then the formal isomonodromic deformation of $(\mathscr{E},\nabla)$ to $\widehat{S}_{\mathscr{K}}$ has property $P$.
 \end{conjecture}
 Indeed, assuming \autoref{conj:main-conjecture}, the hypotheses of  \autoref{conjecture:formal-meta-conjecture} imply that $(\mathscr{E},\nabla)$ admits an algebraic isomonodromic deformation, say over some $T\to S$. But base-changing to $T$ we expect from \autoref{conjecture:meta-conjecture} that every fiber of this isomonodromic deformation has property $P$; hence the same is true at the formal level.
\begin{remark}
The entire discussion here makes sense if $X\to S$ is smooth but not proper, and is instead equipped with a relative simple normal crossings compactification. We omit this for notational simplicity.
\end{remark}
\begin{remark}\label{remark:abelian-variational-conjectures}
	\autoref{conjecture:meta-conjecture} and \autoref{conjecture:formal-meta-conjecture} have analogues for Gauss-Manin connections. Namely, suppose property $P$ is a condition on a cohomology class conjecturally equivalent to being a cycle class (say, being a Hodge or Tate class). Let $f: X\to S$ be a smooth proper morphism, $s\in S$ and $v\in H^{2i}_{dR}(X_s)$ a class satisfying $p$.  One expects that if the solution $\widetilde{v}$ to the Gauss-Manin connection on $\mathscr{H}^{2i}_{dR}(X/S)$ through $v$ is algebraic (or $\omega(p)$-integral), then property $P$ is satisfied by $\widetilde v$ at every point of $S$ (or formally). Indeed, if $v$ is a cycle class and $P$ is the property of being a Hodge class, we verify this in the case of ($\omega(p)$-)integrality in the course of the proof of \autoref{thm:cycle-class-initial-conditions}; in the case that we already know $\widetilde{v}$ is algebraic, it is a consequence of the theorem of the fixed part.
\end{remark}

 We now remark on how our proof of \autoref{thm:NA-main} verifies some cases of \autoref{conjecture:formal-meta-conjecture}.
\subsubsection{Non-abelian $(p,p)$-classes}\label{subsection:non-abelian-p,p} 
It should be clear at this point that \autoref{theorem:extending-the-Hodge-filtration} is an instance of \autoref{conjecture:formal-meta-conjecture}. Indeed, that theorem is precisely the case where $P$ is the property ``underlies a $\mathbb{Z}$-variation of Hodge structure." 

This is already interesting if we know $(\mathscr{E},\nabla)$ has an algebraic isomonodromic deformation (as opposed to starting with the a priori weaker hypothesis of $\omega(p)$-integral isomonodromic deformation). Katzarkov-Pantev \cite{katzarkov1999non} think of this case as an analogue of the theorem of the fixed part for Hodge classes, as in \autoref{remark:abelian-variational-conjectures}. Even in this case our \autoref{theorem:extending-the-Hodge-filtration} gives a quite different proof that the Hodge filtration extends Griffiths-transversally from those in \cite{katzarkov1999non, jost2001harmonic, esnault2024non}.

\subsubsection{A variational Andr\'e conjecture}\label{subsection:variational-Andre}
In the course of the proof of \autoref{theorem:extending-the-Hodge-filtration}, we verify another instance of \autoref{conjecture:formal-meta-conjecture}, this time for $P$ the property of ``having nilpotent $p$-curvature mod $p$ for almost all $p$." Indeed, we do more, and show that the conjugate filtration on the reduction of $(\mathscr{E},\nabla)$ mod $p$ admits a canonical extension over the formal isomonodromic deformation, with respect to which the $p$-curvature has degree $1$ (and is hence nilpotent).

One may see this a posteriori from \autoref{theorem:extending-the-Hodge-filtration} as well, as we now explain. 
\begin{theorem}
	Let $R\subset\mathbb{C}$ be a finitely-generated $\mathbb{Z}$-algebra, $S$ a smooth $R$-scheme, and $X$ a smooth projective $S$-scheme. Fix $s\in S(R)$ and suppose $(\mathscr{E},\nabla)$ is a Picard-Fuchs equation on $X_s$, which admits an $\omega(p)$-integral formal isomonodromic deformation over $\widehat S$, the completion of $S$ at $s$. 
	
	Then, after possibly replacing $R$ by a finitely-generated extension, $(\mathscr{E},\nabla)$ admits an integral isomonodromic deformation $(\widetilde{\mathscr{E}}, \widetilde{\nabla})$ over $\widehat S$, and for each closed point $\mf p$ in $\on{Spec}(R)$, the conjugate filtration on $(\mathscr{E},\nabla)_{\mf p}$ extends to a (decreasing) filtration on $(\widetilde{\mathscr{E}}, \widetilde{\nabla})_{\mf p}$ with respect to which the $p$-curvature has degree $1$.
\end{theorem}
\begin{proof}[Proof sketch]
	First, that $(\mathscr{E},\nabla)$ admits an integral isomonodromic deformation $(\widetilde{\mathscr{E}}, \widetilde{\nabla})$ is immediate from \autoref{thm:NA-main}, as it admits an algebraic isomonodromic deformation. Moreover, by \autoref{theorem:extending-the-Hodge-filtration} (and spreading out), we may assume the Hodge filtration deforms to a Griffiths transverse filtration $F$ on $(\widetilde{\mathscr{E}}, \widetilde{\nabla})$. After localizing $R$ we may assume it is smooth over $\mathbb{Z}$.
	
	Now, let $\widehat R$ be the completion of $R$ at ${\mf p}$, and $\widehat{S}_{\widehat R}$ the base change (of formal schemes). Fix a lift of absolute Frobenius on $\widehat S_{\widehat R}$, so that the Frobenius twist $X'_{\mf p}$ admits a $W_2$-lifting. As in \autoref{defn:lifting-and-conjugate-filtration}, consider the bundle  $$C^{-1}(F_{\text{abs}}^*(\on{gr}_F(\widetilde{\mathscr{E}}, \widetilde{\nabla})))$$ equipped with its conjugate filtration. As in \autoref{subsection:proof-of-key-lemma} this admits a lift to characteristic zero, namely $\Phi(\widetilde{\mathscr{E}}, \widetilde{\nabla}, F)$ as defined in \autoref{defn:higgs-de-rham-flow-functor}. By an argument identical to that in \autoref{subsection:proof-of-key-lemma} it is isomonodromic and isomorphic to $(\mathscr{E},\nabla)$ when restricted to $s_{\widehat R}$. 
	
	Thus it is isomorphic to  $(\widetilde{\mathscr{E}}, \widetilde{\nabla})$ as long as $\mf p$ lies in some dense open of $\on{Spec}(R)$, by \autoref{cor:unique-integral-isomonodromic}, whence $(\widetilde{\mathscr{E}}, \widetilde{\nabla})$ carries a conjugate filtration mod $p$ with the desired property with respect to $p$-curvature. It remains to show this filtration agrees with the one we started with; but this follows from the analogous fact for the Hodge filtration, as Picard-Fuchs equations are Fontaine modules in the language of \cite[Definition 4.16]{ogus-vologodsky} (or equivalently, Higgs-de Rham fixed points, in the language of \cite{lan2019semistable}).
\end{proof}

It would be of great interest to verify other instances of \autoref{conjecture:meta-conjecture} and \autoref{conjecture:formal-meta-conjecture}, as well as, of course, \autoref{conjecture:variational-motivic}.

\appendix
\section{Crystals of functors and their $p$-curvature}\label{appendix:crystals-and-functors}
In this appendix we develop some generalities on crystals (of functors) and their $p$-curvature. In the next appendix we will make these explicit for the functor $\mathscr{M}_{dR}(X/S)^\natural$ sending an $S$-scheme $T$ to the set of isomorphism classes of flat connections on $X_T/T$ (with suitable conditions on regular singularities along the boundary of a simple normal crossings compactification).
\subsection{Tangent functors and extensions}
Given a functor $F:\on{Sch}/S\to \on{Sets}$, a quasi-coherent sheaf on $F$ is an assignment, for $T\in \on{Sch}/S$ and each $f\in F(T)$, a quasi-coherent sheaf $f^*\mathscr{F}$ on $T$, and for each $g: T'\to T,$ an isomorphism $g^*(f^*\mathscr{F})\to (f\circ g)^*\mathscr{F}$ satisfying the evident cocycle condition. 

For an arbitrary scheme (or functor) $Z$  and $M$ and a quasi-coherent sheaf on $Z$, we denote by $Z\oplus M$ the square zero thickening of $Z$ given by $M$: more precisely, $Z\oplus M$ is affine over $Z$ with structure sheaf is given by the sheaf of $\mscr{O}_Z$-algebras $\mscr{O}_Z\oplus M$, where multiplication is given by 
\[
(r_1, m_1)(r_2, m_2)=(r_1r_2, r_1m_2+r_2m_1),
\]
for local sections $r_i, m_i$ of $\mscr{O}_Z, M$ respectively. Explicitly, a $T$-point of $Z\oplus M$ is a $T$-point $f$ of $Z$ and a ring map $\mathscr{O}_T\oplus f^*M\to \mathscr{O}_T$. Set $Z[\epsilon]=Z\oplus \mathscr{O}_Z$.

\begin{definition} The tangent functor $T_{F/S}$ is the set-valued functor on $\Sch_S$ such that for $Z\rightarrow S$, 
$T_{F/S}(Z)$ is the set  $F(Z[\epsilon])$. Note that $T_{F/S}$ is equipped with a natural map to $F$, induced by the natural closed embedding $Z\hookrightarrow Z[\epsilon]$, and that there is a natural monomorphism $F\hookrightarrow T_{F/S}$, induced by the ``reduction modulo $\epsilon$ map $Z[\epsilon]\to Z$ (the ``zero section" of $T_{F/S}$).
\end{definition}

\subsection{Crystals}\label{subsection:crystals}
Let $R$ be a ring and $S$ be a smooth $R$-scheme. Let $((\widehat{S\times S})_{PD}, \bar{I}, \gamma)$ denote the formal scheme obtained as the PD-completion of the diagonal $S\xhookrightarrow{} S\times S$,  and $\pi_i: (\widehat{S\times S})_{PD} \rightarrow S$  the $i$-th projection map for $i=1, 2$. A \emph{crystal of functors} over $S$ is a functor $$F: \on{Sch}/S\to \on{Sets}$$ equipped with an isomorphism    \[
    \gamma: \pi_1^*F \xrightarrow[]{\simeq}  \pi_2^*F
    \]
    which is the identity map when restricted to the diagonal $V(\bar{I})$, and furthermore satisfies the cocycle condition when pulled back to the PD-completion of $S\times S\times S$ along the diagonal, i.e.~$\pi_{23}^*\gamma \circ\pi_{12}^*\gamma=\pi_{13}^*\gamma.$
\begin{definition}\label{definition:flat-section}
    Suppose $F$ is a crystal of functors. A $S$-point of $F$, $\sigma\in F(S)$, induces $(\widehat{S\times S})_{PD}$-points $\sigma_i\in \pi_i^*F((\widehat{S\times S})_{PD})$.  We say that $\sigma$ is flat if $\gamma 
    \circ \sigma_1=\sigma_2$.
    
%    Fix a closed subscheme $Z$ of $S$ (not necessarily smooth), with ideal sheaf $J$, and let $t$ be a $Z$-point of $F$. Let $N(Z, n)=V(\pi_1^*J, \pi_2^*J, \overline{I}^{[n]})\subset (\widehat{S\times S})_{PD}$. For each $n>0$ we obtain (from $t$) $N(Z, n)$-points of $\pi_1^*F$ and $\pi_2^*F$, denoted $\pi_1^*t$ and $\pi_2^*t$. Let $\gamma_{N(Z, n)}$ be the pullback of $\gamma$ to $N(Z, n)$. We say that $t$ is $[n]$-flat if $\gamma_{N(Z, n)}(\pi_1^*t)=\pi_2^*t$.
\end{definition}
%\begin{remark}
%Loosely speaking, $[n]$-flatness is a formalization of the notion of being a solution to a differential equation ``to order $n$." As in \autoref{prop:p-power-leaves} we will wish to check some kind of flatness over $p$-power leaves to guarantee vanishing of $p$-curvature; it will turn out (in \autoref{prop:crystal-p-power-leaves} below) that $[p+1]$-flatness will suffice to guarantee vanishing of $p$-curvature along a $p$-power section to $F$.
%\end{remark}

\begin{remark}
What we have defined in this section is a crystal on the \emph{nilpotent} crystalline site: that is, the subcategory of the crystalline site consisting of divided power thickenings whose defining ideal is pd-nilpotent. Here, we say a divided power ideal $I$ is pd-nilpotent if for every $x\in I$, there exists $n$ with $\gamma_m(x)=0$ for all $m\geq n$.
\end{remark}

\subsection{$p$-curvature of a crystal of functors}
Suppose $S$ is a  $R$-scheme, where $p=0$ in $R$. Let $F$ be a crystal of functors on $S$. We give the definition of the $p$-curvature of $F$ due to Mochizuki (see also \cite[\S~2.3]{osserman2004connections} and \cite[Proposition 1.7]{ogus-vologodsky}).  

Let $I$ be the ideal of $\mscr{O}_{\widehat{(S\times S)}_{\PD}}$ defining the usual diagonal, i.e.~the pullback of the ideal of the diagonal through the natural map $\widehat{(S\times S)}_{\PD}\to S\times S$.

We have the following result, which should be viewed as an analogue of the standard fact that $\Omega^1_S\simeq I/I^2$. Recall that for a pd-algebra $(B, J, \delta)$, and an integer $k\geq 1$, we have the ideals $J^{[k]}$ of $B$ defined as 
\[
J^{[k]} := \text{the ideal generated by } \delta_{e_1}(x_1)\cdots \delta_{e_t}(x_t) \text{ with } \sum e_i \geq k \text{ and } x_i\in J.
\]

\begin{proposition}\label{prop:frobdifferentials}
Suppose either that $S$ is smooth over $R$.
%, or that it arises as a Frobenius pullback, i.e. it is of the form $Z^{[p]}$ for some closed subscheme $Z\xhookrightarrow{} X$, with $X/\on{Spec}(R)$ smooth, as in \autoref{definition:p-power-leaves}.

Then there is a natural isomorphism 
\[F_{\text{abs}}^*\Omega^1_S\simeq \bar{I}/(\bar{I}^{[p+1]}+I),\]
given by 
\[
F_{\text{abs}}^*da \mapsto (a\otimes 1-1\otimes a)^{[p]}
\]
for local sections $a$ of $\mscr{O}_S.$
\end{proposition}
\begin{proof}
    This is \cite[Proposition~1.6]{ogus-vologodsky}.
\end{proof}

By definition of the crystal structure on $F$, we have an isomorphism 
\begin{equation}\label{eqn:crystalstructure}
     \gamma: \pi_1^*F \xrightarrow[]{\simeq}  \pi_2^*F.
    \end{equation}
  Note that the maps   
  \[
  \pi_1|_{V(I, \bar{I}^{[p+1]})}, \pi_2|_{V(I, \bar{I}^{[p+1]})}: V(I, \bar{I}^{[p+1]}) \rightarrow S
  \]
  are canonically identified with each other as $V(I, \bar{I}^{[p+1]})\subset V(I)$, and we therefore denote them both by $\pi_{eq}$. From the isomorphism  \eqref{eqn:crystalstructure}, we  therefore obtain an automorphism $\alpha$ of $\pi_{eq}^*(F)$ over $V(I, \bar{I}^{[p+1]})$.
  
  Below, given a sheaf $\mathscr{G}$ on $S$, we denote the corresponding sheaf on $F$ by $p^*\mathscr{G}$, i.e.~we think of $F$ as being equipped with a ``structure morphism" $p: F\to S$.

\begin{prop/constr}\label{prop:osserman_construction} 
The image of $\bar{I}$ in $\mscr{O}_{(S\times S)_{PD}}/(I, \bar{I}^{[p+1]})$ is square zero; we may therefore write $V(I, \bar{I}^{[p+1]})$ in the form $S\oplus \mf{m}$, where $\mf{m}$ denotes $\bar{I}/(I, \bar{I}^{[p+1]})$, viewed as an $\mscr{O}_S$-module in the natural way. 

Moreover, 
% for $S$ either a smooth scheme or a Frobenius pullback,
 $\alpha$ gives rise to a natural transformation of $(p^*F_{\text{abs}}^*\Omega^1_S)^\vee\to T_{F/S}$, where $p^*F_{\text{abs}}^*\Omega^1_S$ is the sheaf on $F$ associated to $F_{\text{abs}}^*\Omega^1_S$.
\end{prop/constr}

\begin{proof}
    The claim that the image of $\bar{I}$ is square zero is immediate from the definitions. 
From \autoref{prop:frobdifferentials}, we have that $\mf{m}\simeq F^*\Omega^1_S$.  As the natural map $V(I, \overline{I}^{[p+1]})\to S$ is split (by the diagonal map), we find that $V(I, \overline{I}^{[p+1]})=S\oplus \mathfrak{m}$ as desired. Note that $\pi_{eq}^*F\simeq F\oplus p^*\mf{m}$.

   Out of  $\alpha$,  we now construct a map  a map of functors over $F$:
   \[p^*\mf{m}^{\vee} \rightarrow T_{F/S}.\]
   We do this on $Z$-points for $Z$ an $S$-scheme equipped with an element $t\in F(Z).$ Each element in $p^*\mf{m}^{\vee}(Z)$ gives a map 
   \[Z[\epsilon] \rightarrow Z\oplus t^*p^*\mf{m}.\]

   The latter is equipped with a natural map 
   \[t\oplus \mf{m}: Z\oplus t^*p^*\mf{m}\rightarrow F\oplus p^*\mf{m}\simeq F\oplus p^*\mf m.\]
   Composing with $\alpha$, we obtain another map $\alpha \circ (t\oplus \mf{m})$ with the same source and target; finally, using the projection $F\oplus p^*\mf{m}\rightarrow F$, we obtain a  map 
   \[Z[\epsilon] \rightarrow F,\]
   i.e. an element of $T_{F/S}(Z)$, as required. It is straightforward to check that this gives a  map of functors.  
\end{proof}
\begin{definition}\label{defn:osserman_pcurv}
    We refer to the   natural transformation constructed in \autoref{prop:osserman_construction} as the $p$-curvature of the crystal $F$. For any $t\in F(Z)$,   we say that the $p$-curvature of $F$ vanishes at $t$ if the map 
    \[t^*(F_{\text{abs}}^*\Omega_S^{1})^\vee\rightarrow t^*T_{F/S}\] factors through the ``zero section" $t^*F$.
\end{definition}

Now let $Z\subset S$ be a closed subscheme and $Z^{[p]}=F^{-1}_{\text{abs}}(Z)$, as in \autoref{definition:p-power-leaves}. In this case the natural map $(F^*_\text{abs}\Omega^1_Z)^\vee\to t^*(\pi^*F_{\text{abs}}^*\Omega_S^{1})^\vee$ is an isomorphism. We immediately deduce the following analogue of \autoref{prop:p-power-leaves} for crystals of functors:

\begin{proposition}\label{prop:crystal-p-power-leaves}
	Suppose $Z$ is a closed subscheme of $S$, $Z^{[p]}=F_{\text{abs}}^{-1}(Z)$, and fix $t\in F(Z^{[p]})$. Let $W\subset V(I, \overline{I}^{[p+1]})\subset (\widehat{S\times S})_{\PD}$ be the subscheme $\pi_{eq}^{-1}(Z^{[p]})$. Let $\pi_{eq}|_W^{-1}t$ be the $W$-point of $\pi_{eq}|_W^*F$ induced by $t$. Suppose $\gamma|_W(\pi_{eq}|_W^{-1}t)=\pi_{eq}|_W^{-1}(t)$. Then the composition $$(F^*_\text{abs}\Omega^1_{Z^{[p]}})^\vee\overset{\sim}{\to} t^*(p^*F_{\text{abs}}^*\Omega_S^{1})^\vee\rightarrow t^*T_{F/S}$$ vanishes, where the first map is the Frobenius pullback of the derivative of the closed embedding $\iota: Z^{[p]}\to S$, and the second is the pullback of the $p$-curvature of $F$. In other words, the $p$-curvature of $F$ vanishes at $t$.
\end{proposition}
\begin{proof}
That we have an isomorphism $(F^*_\text{abs}\Omega^1_{Z^{[p]}})^\vee\overset{\sim}{\to} t^*(p^*F_{\text{abs}}^*\Omega_S^{1})^\vee$ is a direct computation from the fact that $d(f^p)=0$ for any regular function $f$ on $S$. The rest is immediate from the definition of $p$-curvature; we spell out why. 

By \autoref{prop:frobdifferentials} and \autoref{prop:osserman_construction}, $$W\simeq Z^{[p]}\oplus F_{\text{abs}}^*\Omega^1_S|_{Z^{[p]}}\simeq Z^{[p]}\oplus F^*_{\text{abs}}\Omega^1_{Z^{[p]}}.$$ Now unwinding \autoref{prop:osserman_construction}, restricted to $W$, the $p$-curvature measures the difference between $\gamma|_W(\pi_{eq}|_W^{-1}t)$ and $\pi_{eq}|_W^{-1}(t)$. But these two $W$-points of $\pi_{eq}|_W^*F$ are the same by assumption. 
\end{proof}

%Specializing to $\mscr{M}=\mscr{M}_{dR}(\X/S)$, we therefore obtain 
%\begin{corollary}
%On the smooth locus $\mscr{M}^{\circ}$, the $p$-curvature of $\mscr{M}_{dR}(\X/S)\rightarrow S$ is a map 
%\[
%\Psi: F^*T_S\rightarrow h^1/h^0(\tau_{\leq 1} \tilde{f}_*(\End(\mscr{E}_{univ})\otimes \Omega^{\bullet}),
%\]
%where $\tilde{f}: \X\times_S \mscr{M} \rightarrow \mscr{M}$ is the natural map, and $\mscr{E}_{univ}$ denotes the universal flat bundle on $\X\times_S \mscr{M}/\mscr{M}$. 
%
%For any $t: T\rightarrow \mscr{M}^{\circ}$, the $p$-curvature of $\mscr{M}$ vanishes at $t$, in the sense of \autoref{defn:osserman_pcurv}, if and only if $\Psi_T$ is zero. 
%\end{corollary}
%\begin{remark}
%    In \autoref{section:pcurv_formula}, we will give a formula for  the map $\Psi$ in terms of the usual $p$-curvature.
%\end{remark}
\subsection{Two notions of $p$-curvature agree} 
We will not use the following but we record it for psychological comfort.

Suppose that $S$ is a smooth scheme  over $\mb{F}_q$ and  that $\pi: Y\rightarrow S$ is a crystal of schemes, with $Y/S$ smooth. The crystal structure gives rise to a horizontal foliation on $Y/S$
\[\kappa: \pi^*T_S\rightarrow T_{Y},\] induced by a map $\pi^{-1}T_S\to T_Y$
defined as follows on $Z$-points. Given a $Z$-point of $\pi^{-1}T_S$, i.e.~a map $\delta: Z[\epsilon]\to S$ and a map $t: Z\to Y$ lifting $\delta \bmod \epsilon$, we wish to construct a canonical map $Z[\epsilon]\to Y$ extending $t$. But the existence of such a map is immediate from the fact that there is a canonical PD-structure on $Z[\epsilon]$ with respect to the ideal $\epsilon$ (as it is square zero). Indeed, let $$\iota: Z[\epsilon]\to Z\to Z[\epsilon]\overset{\delta}{\to} S$$ be the natural ``constant" map. Then $(\iota, \delta): Z[\epsilon]\to S\times S$ canonically factors through a map $$(\iota, \delta)^{\on{PD}}\to S\times^{\on{PD}}S.$$ The map $t$ gives rise to a ``constant" section $\alpha$ to $\iota^*Y/Z[\epsilon]$, i.e.~one that factors through $Z$. But $\iota^*Y={(\iota, \delta)^{\on{PD}}}^*\pi_1^*Y$; our desired map $Z[\epsilon]\to Y$ is given by applying ${(\iota, \delta)^{\on{PD}}}^*\gamma$ to this section, where $\gamma$ is the isomorphism of \autoref{subsection:crystals} in the definition of a crystal.
%\begin{proposition}
%    The $p$-th power maps on $T_Y$ and $T_R$ induce a map $\phi: F^*T_{\mc{Y}}\rightarrow T_{\mc{Y}}$. Moreover, this is map independent of the chosen presentation.  
%\end{proposition}

In \autoref{defn:p-curv-foliation}, we defined the $p$-curvature of a foliation.  We now make the same definition in the case of  a crystal of smooth schemes, and record that it agrees with the definition of $p$-curvature due to Mochizuki.
\begin{proposition}
    Denote by $\Pi$ the composition
    \[F^*\pi^*T_S\xrightarrow[]{F^*\kappa} F^*T_{{Y}} \xrightarrow[]{\phi} T_{{Y}} \xrightarrow[]{} T_{{Y}}/\kappa(\pi^*T_S)\simeq T_{Y/S},\] where $\phi$ is the map induced by taking the $p$-th power of a vector field. 
    Then $\Pi$ agrees with $\Psi$ as defined in \autoref{defn:osserman_pcurv}.
\end{proposition}
\begin{proof}
The argument is essentially the same as the one for flat bundles, which is recorded in \cite[Proposition 1.7]{ogus-vologodsky}.    
\end{proof}

\section{The non-abelian Gauss-Manin connection and its $p$-curvature}\label{appendix:NAGM}
In this appendix we explain that $\mathscr{M}_{dR}(X/S)^\natural$ is a crystal and compute its $p$-curvature. See \autoref{subsection:appendix-summary} for a less technical summary of the key material we develop here.
\begin{remark}
It would arguably be better to show that the stack $\mathscr{M}_{dR}(X/S)$ is itself a crystal of stacks over $S$, and compute the $p$-curvature of this stack. For $\pi: \mathscr{M}_{dR}(X/S)\to S$ the structure morphism, this would be a map $\pi^*F_{\text{abs}}^*T_S\to T_{\mathscr{M}_{dR}(X/S)/S}$, where here the target is the \emph{tangent stack} of $\mathscr{M}_{dR}(X/S)/S$.

The arguments required to do so would be close to those contained here. Unfortunately there is not, to our knowledge, any literature on crystals of stacks and we felt that developing this theory here would be too cumbersome to justify its inclusion, since it is unnecessary for our applications.
\end{remark}

\subsection{The tangent functor of $\mathscr{M}_{dR}^\natural$}
We fix notation as in \autoref{subsection:M_dR-and-isomonodromy}, i.e.~$R$ is a ring, $S$ is a smooth $R$-scheme, $\overline{f}: \overline{X}\to S$ is a smooth projective morphism, and $D\subset \overline{X}$ is a simple normal crossings divisor over $S$. Set $X=\overline{X}\setminus D$. Recall that $\mathscr{M}_{dR}(X/S, r)^\natural$ is the functor sending an $S$-scheme $T$ to the set of isomorphism classes of flat bundles on $\overline{X}_T/T$ with regular singularities along $D_T$.

Let $Z$ be an affine $S$-scheme and $[(\mathscr{E}, \nabla)]\in \mathscr{M}_{dR}(X/S)^\natural(Z)$. It follows immediately from \autoref{cor:de-rham-complex-deformations} that the preimage of $[(\mathscr{E}, \nabla)]$ under the natural map $$T_{\mathscr{M}_{dR}(X/S)^\natural/S}(Z)\to \mathscr{M}_{dR}(X/S)^\natural(Z)$$ is canonically identified with the quotient $$\mathbb{H}^1(\overline{X}, \on{End}(\mathscr{E})_{dR})/\on{Aut}(\mathscr{E}, \nabla).$$
In particular, if $(\mathscr{E}, \nabla)$ is simple, then this is simply $\mathbb{H}^1(\overline{X}, \on{End}(\mathscr{E})_{dR})$ itself. 
\subsection{The non-abelian Gauss-Manin connection} \label{subsection:appendix-nagm}
\begin{lemma}\label{lemma:glue_taylor}
Let $(A, I, \gamma)$ be a pd-nilpotent divided power algebra, $T=\on{Spec}(A)$, $T_0=\on{Spec}(A/I)$, and $U,V$ smooth affine $T$-schemes with $U_0:=U_{T_0}$. Let $D$ (resp.~$E$ be a relative snc divisor on $U$ (resp.~$V$) over $A$. Let $(\mathscr{E}, \nabla)$ be a flat vector bundle on $(V, E)/T$, and let $f_i, f_j: U\to V$ be two morphisms over $T$ whose restrictions to $U_0$ agree, with $f_i^{-1}(E), f_j^{-1}(E)\subset D$. Then there is a canonical isomorphism $\alpha_{ij}: f_i^*(\mathscr{E}, \nabla)\simeq f_j^*(\mathscr{E},\nabla)$ given by Taylor's formula. Explicitly, in \'etale-local coordinates $x_1, \cdots , x_d$ on $V/T,$ where $E$ is cut out by $x_1\cdots x_r$, 
\[
\alpha_{ij}: s\otimes 1 \mapsto \sum_{\underline{\ell}} \partial ^{\underline{\ell}}(s) \otimes 
\prod_{m=1}^d \gamma_{\ell_m}(f_i^*(x_{m})-f_j^*(x_{m})).
\]
Here, $s$ is a local section of $\mscr{E}$, $\underline{\ell}$ denotes a tuple of integers $(\ell_1, \cdots , \ell_d)$ with each $\ell_m\geq 0$, and $\partial^{\underline{\ell}}$ denotes the operator $\prod \nabla_{\partial_{x_m}}^{\ell_m}$. Finally, the symbol $\gamma_{\ell_m}$ denotes the p.d. structure extended to $\mscr{O}_UI,$ which is unique by  \cite[\href{https://stacks.math.columbia.edu/tag/07H1}{Lemma 07H1}]{stacks-project}.
Given $f_i, f_j,f_k$ as above, we have $\alpha_{jk}\circ\alpha_{ij}=\alpha_{ik}$.
\end{lemma}
\begin{proof}
We briefly indicate why the formula for $\alpha_{ij}$ makes sense in the logarithmic setting, where this is not obvious as $\partial^{\underline \ell}$ does not preserve $\mathscr{E}$. Suppose $\underline{\ell}$ is a multi-index, and set $E_i=V(x_i), i=1,\cdots, r$. Let $\underline{\ell}E=\sum_{i=1}^r \ell_i E_i$. Then $\partial^{\underline \ell}$ sends $\mathscr{E}$ to $\mathscr{E}(\underline{\ell}E)$. But $f_i^*, f_j^*$ send the ideal of $E$ into the ideal of $D$;  hence $$\prod \nabla_{\partial_{x_m}}^{\ell_m}$$ is contained in the ideal of $\underline{\ell}E$, and hence the tensor product above belongs to $f_j^*(\mathscr{E})$. See \cite[Proposition 3.3]{eg_revisit} for a similar argument.

It is a straightforward (classical) computation that the map, as written above, is a well-defined isomorphism of flat bundles, and moreover that the cocycle condition is satisfied; see e.g.~\cite[Lemma 2.8]{bhatt2011crystalline} for a similar computation.
\end{proof}
\begin{construction}[$\mathscr{M}_{dR}(X/S)^\natural/S$ is a crystal]\label{prop:crystal taylor}
	Let $(T, I, \gamma)$ be a divided power scheme with $I$ pd-nilpotent,  $T_0=V(I)$, $X_0/T_0$ a smooth scheme, and $(\overline{X}_i, D_i), (\overline{X}_j, D_j)$ two smooth $T$-schemes equipped with simple normal crossings divisors over $T$ whose fiber over $T_0$ is $(X_0, D_0)$. Let $X_i=\overline{X}_i\setminus D_i, X_j=\overline{X}_j\setminus D_j$. Then there is a canonical bijection $$\varphi_{ij}: \mathscr{M}_{dR}(X_i/T)^\natural(T) \overset{\sim}{\to} \mathscr{M}_{dR}(X_j/T)^\natural(T)$$ given explicitly as follows. For  $(\mscr{E}, \nabla)\in \mathscr{M}_{dR}(X_i/T)^\natural(T)$, we give a recipe to construct an object of $\mathscr{M}_{dR}(X_j/T)^\natural(T)$: 

Let $\{U_\beta\}$ be an affine cover of $\overline{X}_0$, and $U_\beta^i, U_\beta^j$ the corresponding affine covers of $\overline{X}_i, \overline{X}_j$. For each $\beta$, choose an isomorphism $\xi_\beta: U_\beta^j\overset{\sim}{\to} U_\beta^i$ sending $D_j\cap U_\beta^j$ to $D_i\cap U_\beta^i$ which is the identity mod $I$; such exists because deformations of smooth affine schemes equipped with simple normal crossings divisors are unique.
We obtain a connection $(\mscr{E}_{\beta}, \nabla):= \xi_{\beta}^*((\mscr{E}, \nabla)|_{U_{\beta}^{i}})$ on $U^j_{\beta}$. It remains to glue these together to obtain a flat bundle on $(\overline{X}_j, D_j)/T$. 

 For two indices $\beta, \gamma$,  define $U^i_{\beta\gamma}:= U^i_{\beta}\cap U^i_{\gamma}$, and similarly $U^j_{\beta\gamma}$. Then we have maps 
 \[
 \xi_{\beta}|_{U^j_{\beta\gamma}},  \xi_{\gamma}|_{U^j_{\beta\gamma}}: U^j_{\beta\gamma} \overset{\sim}{\to} U^i_{\beta\gamma}
 \]
 which agree mod $I$ by construction. For brevity, we denote these two maps by $F, G$ respectively. By \autoref{lemma:glue_taylor} applied to $F$ and $G$, we obtain an isomorphism $F^*((\mscr{E}, \nabla)|_{U^i_{\beta \gamma}}) \simeq G^*((\mscr{E}, \nabla)|_{U^i_{\beta \gamma}})$ given in coordinates by 
 \[
s\otimes 1 \mapsto \sum_{\underline{\ell}} \partial ^{\underline{\ell}}(s) \otimes 
\prod_{m=1}^d \gamma_{\ell_m}(F^*(x_{m})-G^*(x_{m})).
\]

 Therefore we have an isomorphism 
\[
f_{\beta \gamma}: \mscr{E}_{\beta}|_{U^j_{\beta \gamma}} \simeq  \mscr{E}_{\gamma}|_{U^j_{\beta \gamma}}.
\]
Again by \autoref{lemma:glue_taylor}, the maps $f_{\beta \gamma}$ satisfy the cocycle condition, and hence we obtain a flat bundle on $(\overline{X}_j, D_j)/T$. It is straightforward to check that this gives a bijection, and that it satisfies the cocycle condition: given $(\overline{X}_i, D_i), (\overline{X}_j, D_j), (\overline{X}_k, D_k)$, we have $\varphi_{jk}\circ\varphi_{ij}=\varphi_{ik}$.

We now construct descent data as in \autoref{subsection:crystals} to endow $\mathscr{M}_{dR}(X/S)^\natural$ with the structure of a crystal of functors. Namely, consider the formal scheme $(\widehat{S\times S})_{\PD}$, and let $\pi_i: (\widehat{S\times S})_{\PD}\to S$ the two projections. We wish to construct an isomorphism $\pi_1^*\mathscr{M}_{dR}(X/S)^\natural\xrightarrow{\sim} \pi_2^*\mathscr{M}_{dR}(X/S)^\natural$ satisfying the cocycle condition. On $T$-points, for $(T, I, \gamma)$ a divided power scheme with $I$ pd-nilpotent, this boils down to the following: given two maps $f_1, f_2: T\to S$ which agree modulo $I$, set $(\overline{X}_i, D_i)=(\overline{X}, D)\times_{S, f_i}  T$ for $i=1, 2$. We need to construct a natural bijection between isomorphism classes of flat bundles on $(\overline{X}_1, D_1)/T$ and isomorphism classes of flat bundles on $(\overline{X}_2, D_2)/T$ satisfying the cocycle condition. But this is precisely what we've done above.
\end{construction}
\begin{remark}
In fact all of these constructions give rise to equivalences of categories, endowing the stack $\mathscr{M}_{dR}(X/S)$ with the 	structure of a crystal of stacks. But we will not use this fact here.
\end{remark}

Let $(T, I, \gamma)$ be a divided power (formal) scheme with $I$ (topologically) pd-nilpotent. If the natural inclusion $T_0\hookrightarrow T$ is equipped with a splitting $s: T\to T_0$ (for example, if $T=(\widehat{S\times S})_{PD}$ and $I$ is the ideal sheaf of the diagonal), and $(\overline{X}, D)/T$ is an snc pair, with $(\overline{X}_0, D_0)=(\overline{X}, D)|_{T_0}$, \autoref{prop:crystal taylor} gives rise to a map from isomorphism classes of flat bundles on $(\overline{X}_0, D_0)/T_0$ to isomorphism classes of flat bundles on $(\overline{X},D)/T$, as follows. Given a flat bundle $(\mathscr{E}_0, \nabla_0)$ on $(\overline{X}_0, D_0)/T_0$, we may pull it back to a flat bundle on $(\overline{X}_0, D_0)\times_{T_0} T/T$, and then apply the construction of \autoref{prop:crystal taylor} to the pair of deformations of $(\overline{X}_0, D_0)$ given by $(\overline{X}_0, D_0)\times_{T_0} T, (\overline{X},D)$. A flat bundle on $X/T$ obtained this way is \emph{isomonodromic}. 

Note that in the case that $T=(\widehat{S\times S})_{PD}$, this agrees with the notion of  flatness from \autoref{definition:flat-section}. To compare this notion to the characteristic zero definition of \autoref{remark:general-isomonodromic-in-char-0} note that the flat bundles obtained via \autoref{prop:crystal taylor} are locally pulled back from $(\overline{X}_0, D_0)$.

\begin{remark}
In \cite[\S8]{simpson1994moduli}, Simpson gives an essentially identical construction of the crystal structure on $\mathscr{M}_{dR}(X/S)$ (in the case where $D=\emptyset$), working over $\mathbb{C}$. He compares it with the usual analytic description of isomonodromy (i.e.~that the conjugacy class of the monodromy of a family of connections is locally constant) \cite[Theorem 8.6]{simpson1994moduli}. This shows, for example, that our description here also agrees with that of \autoref{defn:naive-isomondromy}, though one could also see this directly, by e.g.~comparing the formulas involved in the explicit constructions.
\end{remark}

\subsection{The Kodaira--Spencer cocycle}\label{section:ks cocycle}
\begin{construction}\label{constr:ks cocycle}
    Suppose $\overline{f}: \overline{X}\rightarrow S$ is a smooth projective family equipped with a simple normal crossing divisor $D/S$. Then the Kodaira--Spencer class, viewed as a section $\KS$ of $\Omega^1_S\otimes R^1\fbar_*T_{\overline{X}/S}(-\log D)$, may be represented by a cocycle as follows. Let $S^{(2)}\subset S\times S$ denote the first order thickening of the diagonal inside $S\times S$. The natural thickening   $S \xhookrightarrow{} S^{(2)}$  may be split by the first projection map $\pi_1: S^{(2)}\rightarrow S$.  

    Define $(\overline{X}, D)_{S^{(2)}}:= (\overline{X}, D) \times_{S, \pi_2} S^{(2)}$, so that we have a Cartesian diagram 
    \begin{equation} \begin{tikzcd}
        (\overline{X}, D) \arrow[r] \arrow[d, "\overline{f}"]
        & (\overline{X}, D)_{S^{(2)}}\arrow[d] \\
        S \arrow[r]
        &  S^{(2)}.
        \end{tikzcd}
    \end{equation}
    In other words, $(\overline{X}, D)_{S^{(2)}}$ is a thickening of $(\overline{X}, D)$, over the split extension $S\xhookrightarrow{} S^{(2)}$. Note that this is, in general, not the trivial thickening with respect to the splitting $\pi_1: S^{(2)}\rightarrow S$. Pick open affines $U_i\subset \overline{X}$ covering $\overline{X}$, so that we obtain $\widehat{U}_{i, 1}:= U_i\times_{S, \pi_1}S^{(2)}$, and similarly $\widehat{U}_{i,2}$; note that the latter is an open affine inside $\overline{X}_{S^{(2)}}$. Let $D_i=U_i\cap D$, and define $\widehat D_{i,1}, \widehat D_{i,2}$ to be the pullback of $D_i$ to $S^{(2)}$.

    Since both $(\widehat{U}_{i,1}, \widehat D_{i,1})$ and $(\widehat{U}_{i,2}, D_{i,2})$ are deformations of $(U_i, D_i)$ (a smooth snc pair) over $S^{(2)}$, there exists an isomorphism of deformations
    \[
    \widetilde{f}_i: (\widehat{U}_{i,2}, \widehat D_{i,2}) \simeq (\widehat{U}_{i,1}, \widehat D_{i,1}). 
    \]
    Denote by $f_i$ the composition 
    \[\widehat{U}_{i,2}\xrightarrow[]{\widetilde{f}_i} \widehat{U}_{i,1} \rightarrow U_i,\]
    with the last map being the natural retraction. Note that $f_i$ is the identity modulo $I$.

Let $U_{ij}$ denote the intersection $U_i\cap U_j$.
    Let $c_{ij}\in \fbar^*\Omega^1_S \otimes T_{\overline{X}/S}(-\log D)(U_{ij})$ be the section given by 
    \[
    dx_k \mapsto f_i(x_k)-f_j(x_k)
    \]
    in local coordinates $x_1, \cdots, x_d$ on $\overline{X}/S$, where $D$ is cut out by $x_1\cdots x_r$. 
    Note that we view $f_i(x_k)-f_j(x_k)$ as a section of $\fbar^*\Omega^1_S$, since it is a function on $X_{S^{(2)}}$ vanishing on $X\subset X_{S^{(2)}}$, and that this map vanishes on $x_1, \cdots, x_r$ and hence defines a map out of $\Omega^1_{\overline{X}/S}(\log D)$.
    
    It is straightforward to check that $c_{ij}$ is well-defined, and the following result is a matter of unwinding definitions:
\end{construction}
\begin{proposition}\label{prop:ks cocycle}
    The sections $c_{ij}$ constructed in \autoref{constr:ks cocycle} form a cocycle representing the Kodaira--Spencer class.
\end{proposition}

\subsection{The $p$-curvature of $\mscr{M}_{dR}$: cocycles}
Write $\mscr{M}=\mscr{M}_{dR}(X/S)^\natural$; given a sheaf $\mathscr{G}$ on $S$, let $p^*\mathscr{G}$ denote the corresponding sheaf on $\mathscr{M}$. Recall that we defined the non-abelian $p$-curvature to be a map 
\[
\Psi: p^*F_{\text{abs}}^*T_S\rightarrow T_{\mscr{M}/S}.
\]
%or equivalently a section $\Psi$ of $\pi^*F^*\Omega^1_S\otimes T_{\M/S}$.
Recall that we have described the fiber of $T_{\mathscr{M}/S}$ over a point $[(\mathscr{E}, \nabla)]$ of $\mathscr{M}$ as $$\mathbb{H}^1(\overline{X}, \on{End}(\mathscr{E})_{dR})/\on{Aut}(\mathscr{E}, \nabla).$$
%T_{\M/S}\simeq R^1\tilde{f}_{*, dR}\End(\mscr{E}_{univ})$, with $\tilde{f}: \X\times_S \M\rightarrow \M$ is the projection map, and $\mscr{E}_{univ}$ is the universal flat bundle on $\X\times \M/\M$.

Let $h: T\rightarrow S$ be a $S$-scheme and $t: T \rightarrow \mscr{M}$ be a $T$-point of $\mscr{M}$, which we may think of as being represented by a flat bundle $(\mathscr{E}_T,\nabla)$ on $(\overline{X}, D)_T/T$. To give  a description of  $\Psi$, it suffices to give sections $\Psi_T\in h^*F_{\text{abs}}^*\Omega^1_S\otimes R^1\overline{f}_{T*}\End(\mscr{E}_T)_{dR}$ compatibly as $T$ varies;  here 
$(\overline{X}, D)_T:= (\overline{X}, D)\times_S T$, and $\overline{f_T}: \overline{X}_T\rightarrow T$ is the basechange of $\overline{f}: \overline{X}\rightarrow S$ fits into the diagram
    \begin{equation}
        \begin{tikzcd}
            (\overline{X}, D)_T \arrow[r, "\tilde{h}"] \arrow[d, "\overline{f_T}"] & (\overline{X}, D) \arrow[d, "\overline{f}"] \\
            T \arrow[r, "h"] &S.
        \end{tikzcd}
    \end{equation}
    
   \begin{construction}
    For brevity, denote by $\widehat{S}$ the PD-completion of the diagonal inside $S\times S$, i.e.~$(\widehat{S\times S})_{\PD}$ and by $\bar{I}$ its  natural  p.d. ideal; let $\pi_i: \widehat{S} \rightarrow S$ be the $i$-th projection, for $i=1,2$. 
    
    Pick open affines $U_i\subset \overline{X}$ covering $\overline{X}$, and define $\widehat{U}_{i, 1}:= U_i\times_{S, \pi_1}\widehat{S}$, and similarly $\widehat{U}_{i,2}$. Similarly define $D_i=U_i\cap D$, and $\widehat{D}_{i,1}, \widehat{D}_{i,2}$ the pullback to $\widehat{U}_{i,1}, \widehat{U}_{i,2}$.

    Since both $(\widehat{U}_{i,1}, \widehat D_{i,1})$ and $(\widehat{U}_{i,2}, \widehat D_{i,2})$ are deformations of $(U_i, D_i)$ over $\hat{S}$---a smooth affine snc pair---there exists an isomorphism  over $\widehat{S}$
    \[
    \widetilde{f}_i: (\widehat{U}_{i,2}, \widehat D_{i,2}) \simeq (\widehat{U}_{i,1}, \widehat D_{i,1}). 
    \]
    Denote by $f_i$ the composition 
    \[\hat{U}_{i,2}\xrightarrow[]{\tilde{f}_i} \hat{U}_{i,1} \rightarrow U_i,\]
    with the last map being the natural retraction. 

    In local coordinates $x_1, \cdots, x_d$ on $U_i$, consider  elements $c_{ij, \mscr{E}_T}\in f_T^*h^*\bar{I}\otimes \End(\mscr{E}_T)$ defined as 
    \[
s\otimes 1 \mapsto \sum_{\underline{\ell}} \partial ^{\underline{\ell}}(s) \otimes 
\prod_{m=1}^d \tilde{h}^*\gamma_{\ell_m}(f_i^*(x_{m})-f_j^*(x_{m})),
\]
for $s$ a local section of $\mscr{E}_T$, with notation as in \autoref{lemma:glue_taylor}. 
\end{construction}

Recall that, for a flat bundle $(\mscr{N}, \nabla)$ on a scheme $Y$, a hypercocycle representing a class in $\mb H^1(\mscr{N}_{dR})$ is given by:
\begin{itemize}
    \item a covering of $Y$ by open subsets $U_i$, and 
    \item after setting $U_{ij}:= U_i\cap U_j$,  $U_{ijk}:= U_i\cap U_j\cap U_k$, a collection of sections $(f_{ij}, \omega_i)$, with $f_{ij}\in \mscr{N}(U_{ij})$, $\omega_i\in \mscr{N}(U_i)$ such that  $f_{ij}|_{U_{ijk}}-f_{ik}|_{U_{ijk}}+f_{jk}|_{U_{ijk}}$, and 
    \[
d\omega_i=0, \ \ \    \omega_i|_{U_{ij}}-\omega_j|_{U_{ij}} = df_{ij}.
    \]
    In particular, if we have flat sections $(f_{ij})$ satisfying the cocycle condition, then the collection $(f_{ij}, 0)$ defines  a hypercocycle.
\end{itemize}
\begin{proposition}\label{prop:p-curv cocycle}
\begin{enumerate}
    \item The collection of tuples  $(c_{ij, \mscr{E}_T}, 0)$ form a hypercocycle representing a global section of  $(h\times h)^*\bar{I} \otimes R^1\overline{f_T}_{*}\End(\mscr{E}_{T})_{dR}$. 
    \item Moreover,  after restricting to $V(\bar{I}^{[p+1]}+I)$, this hypercocycle represents the class $\Psi$,
    \item and the restriction of $c_{ij, \mscr{E}_T}$ is given by 
     \[
    c_{ij, \mscr{E}_T, V}: s\otimes 1\mapsto \sum_{k=1}^d \nabla_{\partial_{x_k}}^{p}(s) \otimes \gamma_p(f_i^*(x_k)-f_j^*(x_k)),
    \]
    for $s$ a local section of $\mscr{E}_T$.
\end{enumerate}
\end{proposition}

\begin{proof}
    The first follows from the fact that the $(c_{ij, \mscr{E}_T})$ themselves form a cocycle, as well as the fact that each map $c_{ij, \mscr{E}_T}$ is flat. The  second claim follows from the construction of $\Psi$ and the explicit recipe given in \autoref{prop:crystal taylor}. Finally, the formula for $c_{ij, \mscr{E}_T, V}$ is immediate from the definition of $V(\bar{I}^{[p+1]}+I)$.
\end{proof}

Recall that in \autoref{section:ks cocycle}, we defined the Kodaira--Spencer class $\KS$, which is a global section of  $\Omega^1_S\otimes R^1\fbar_*T_{\overline{X}/S}(-\log D).$ Therefore for any $h: T\rightarrow S$, we obtain the pulled-back class $\KS_T$, which is a global section of  $h^*\Omega^1_S\otimes R^1\overline{f_T}_*T_{\overline{X}_T/T}(-\log D_T)$. 
\begin{lemma}\label{lemma:cocycles-agree}
Suppose $p>2$. Let $\psi_{\mscr{E}_T}: F_{\text{abs}}^*T_{\overline{X}_T}(-\log D)\rightarrow \End(\mscr{E}_T)$ be the $p$-curvature map of the flat bundle $\mscr{E}_T$.  Then the section
\[
\Psi_T = \psi_{\mscr{E}_T}(F_{\text{abs}}^*\KS_T)
\]
of  $h^*F_{\text{abs}}^*\Omega^1_S\otimes R^1\overline{f_T}_*\End(\mscr{E}_T)_{dR}$ represents the $p$-curvature $t^*\Psi$ of $\mathscr{M}$ restricted to $T$.

Note that, in the above equality, we view $F_{\text{abs}}^*\KS_T$ as a flat section of $$h^*F_{\text{abs}}^*\Omega^1_S \otimes R^1\overline{f_T}F_{\text{abs}}^*T_{\overline{X}_T/T}(-\log D_T)_{dR},$$ with the canonical (Frobenius pullback) connection on $F_{\text{abs}}^*T_{\overline{X}_T/T}(-\log D_T)$.

In other words, $t^*\Psi$ is represented by the composite map $$h^*F_{\text{abs}}^*T_S\to h^*F_{\text{abs}}^*R^1\overline{f_T}_{*}T_{\overline{X}_T/T}(-\log D_T)\xrightarrow{\psi_{\mathscr{E}_T}} R^1\overline{f_T}_*\on{End}(\mathscr{E}_T)_{dR}$$
where the first map above is the pullback of the Kodaira-Spencer map, and the second is the $p$-curvature of $(\mathscr{E}_T, \nabla)$.
\end{lemma}
\begin{proof}
    We choose an affine covering by $U_i\subset \overline{X}$, and denote by $\widehat{U}_{i,1}:=U_i\times_{S, \pi_1} \widehat{S}$ and similarly $\widehat{U}_{i,2}$ as before, along with maps 
    \[
    f_i: \hat{U}_{i,2}\simeq \hat{U}_{i,1}\rightarrow U_i.
    \]
    Since we have a canonical identification $\mscr{O}_{\widehat{S}}/\bar{I}^{[2]}\simeq \mscr{O}_S/I^2$, the maps $f_i \bmod \bar{I}^{[2]}$ give maps as in \autoref{constr:ks cocycle}. The claim now follows straightforwardly by inspecting the explicit cocycles representing these classes, as given in  \autoref{prop:p-curv cocycle} and \autoref{prop:ks cocycle}. Indeed, by \autoref{prop:ks cocycle}, the class $\KS_T$ has a cocycle representative 
    \[
    c_{ij} = \sum  (f_i(x_k)-f_j(x_k)) \otimes \partial_{x_k} \in h^*\Omega^1_S\otimes \overline{f}_{T*}T_{\overline{X}_T/T}(-\log D_T)(U_{ij}).
    \]
     Therefore $F^*_{abs}(\KS_T)$ is represented by 
\begin{equation}\label{eqn: frob-pullback-ks-class}
    F^*_{abs} c_{ij} = \sum  F^*_{abs}(f_i(x_k)-f_j(x_k)) \otimes F^*_{abs}\partial_{x_k} \in h^*F^*_{abs}\Omega^1_S\otimes F^*_{abs}\overline{f}_{T*}T_{\overline{X}_T/T}(-\log D_T)(U_{ij}).
    \end{equation}
   Under the  explicit isomorphism $F^*_{abs}\Omega^1_S\simeq \bar{I}/(\bar{I}^{[p+1]}+I)$ of \autoref{prop:frobdifferentials}, the section \eqref{eqn: frob-pullback-ks-class} becomes 
   \[
   \sum  \gamma_p(f_i(x_k)-f_j(x_k)) \otimes F^*_{abs}\partial_{x_k},
   \]
     and therefore $\psi_{\mscr{E}_T}(F^*_{abs}\KS_T)$ is represented by the cocycle 
     \[
    \psi_{\mscr{E}_T}(F^*_{abs}c_{ij}): s\otimes 1\mapsto \sum_{k=1}^d \nabla_{\partial_{x_k}}^{p}(s) \otimes \gamma_p(f_i^*(x_k)-f_j^*(x_k)),
    \]
    which gives precisely the global section $\Psi_T$ by \autoref{prop:p-curv cocycle}.

    The last claim is a matter of unwinding definitions.
\end{proof}

\subsection{Vanishing of $p$-curvature}
In this section we give a criterion for vanishing of the $p$-curvature of $\mathscr{M}_{dR}(X/S)^\natural$ in terms of integral (or $\omega(p)$-integral) isomonodromic deformations. This will boil down to checking the hypotheses of \autoref{prop:crystal-p-power-leaves}. 

Recall that ($\omega(p)$-)integral isomonodromic deformations are simply families of flat bundles (over certain special bases) that are isomonodromic when restricted to characteristic zero.  Loosely speaking, the idea of this section will be that in good situations we can check the hypotheses of \autoref{prop:crystal-p-power-leaves} by lifting to characteristic zero, using \autoref{cor:unique-integral-isomonodromic}. This is the analogue for isomonodromy of the fact that if a power series $f(z)\in \mathbb{Z}_p[[t]]$ satisfies a differential equation to some order, the same is true of its mod $p$ reduction.

Let $R$ be a Noetherian domain of characteristic zero and $S=\on{Spf}R[[t_1, \cdots, t_d]]$. Let $\overline{X}$ be an adic smooth projective formal $R$-scheme equipped with a simple normal crossings divisor $D$ over $S$. For a multi-index $A=(a_1, \cdots, a_d),$ let $t^A=t_1^{a_1}\cdots t_d^{a_d}$ and $|A|=\sum_i a_i$. Let $Z\subset S$ be a zero-dimensional monomial subscheme cut out by some ideal $J=(t^{A_1}, \cdots, t^{A_d})$, and fix $m$  such that $Z$ is contained in $V((t_1, \cdots, t_d)^m)$. Let $Z_{+}$ be the (monomial) subscheme cut out by $(t_1, \cdots, t_d)^{m+p+1}$; note that $Z\subset Z_{+}$. Let $0\subset S$ be the closed subscheme corresponding to the ideal $(t_1, \cdots, t_d)$.

The following lemma checks the hypotheses of \autoref{prop:crystal-p-power-leaves}  by lifting to characteristic zero.
\begin{lemma}\label{cor:vanishing-p-curvature-for-p-power-leaves-isomonodromy}
Suppose each $A_i$ is a multiple of $p$, i.e.~the reduction $Z_p$ of $Z$ mod $p$ is a $p$-power subscheme of $S$. Let $(\mathscr{E}, \nabla)$ 
be a flat vector bundle on $(\overline{X}, D)_{Z_+}/Z_+$. 
Suppose there exists an isomonodromic flat bundle $(\mathscr{F}, \nabla)$ on $(\overline{X}, D)_{\mathbb{Q}}/S_{\mathbb{Q}}$ and an isomorphism $$(\mathscr{E}, \nabla)_{\mathbb{Q}}\simeq (\mathscr{F}, \nabla)|_{Z_{+, \mathbb{Q}}}.$$ Finally, suppose $(\mathscr{E}, \nabla)|_{\overline{X}_0}$ satisfies the cohomological conditions of \autoref{cor:unique-integral-isomonodromic}, i.e.~that $H^1(\on{End}(\mathscr{E})_{dR}|_{\overline{X}_0})$ is locally free and its formation commutes with arbitrary base change. Then the $p$-curvature of $\mathscr{M}_{dR}(X/S)^{\natural}$ vanishes at the $Z_p$-point of $\mathscr{M}_{dR}(X/S)^{\natural}$ corresponding to $(\mathscr{E}, \nabla)|_{Z_p}$.
\end{lemma}
\begin{proof}
	We check the hypotheses of \autoref{prop:crystal-p-power-leaves}. That is, set $V(I, \overline{I}^{[p+1]})\subset (\widehat{S\times S})_{\PD}$ to be as in \autoref{prop:osserman_construction}, $\pi_{eq}: V(I, \overline{I}^{[p+1]})\to S$ the projection map, and $W=\pi_{eq}^{-1}(Z_p)$. We must show that the $W$-point of $\mathscr{M}_{dR}^\natural(X/S)$ induced by $(\mathscr{E}, \nabla)|_{Z_p}$ is fixed by the map \eqref{eqn:crystalstructure} in the definition of a crystal. We will check this on a larger mixed-characteristic subscheme $W_+$ of $(\widehat{S\times S})_{\PD}$ containing $W$.
	
	We first let $W_+'$ be the subscheme of $(\widehat{S\times S})_{\PD}$ given by  the preimage of $Z_+\times Z_+$ under the natural map $(\widehat{S\times S})_{\PD}\rightarrow  S\times S$. We let $W_+$ be the subscheme of $W_+'$ cut out by the ideal consisting of $R$-torsion. One may give $W_+$ the following explicit description. Letting $J_+$ be the ideal of $Z_+$, the functions on $W_+$ are given by the image of the natural map from
	 $$(R[[t_1, \cdots t_d]]\otimes_R R[[t_1, \cdots t_d]])_{\PD}/(J_+\otimes 1, 1\otimes J_+)$$
	 	 to 
	  $$(R_\mathbb{Q}[[t_1, \cdots t_d]]\otimes_{R_\mathbb{Q}} R_\mathbb{Q}[[t_1, \cdots t_d]])/(J_+\otimes 1, 1\otimes J_+).$$
	  Note that this map is injective in degrees at most $m+p$.
	  
	  We first check that $W_+$ contains $W$. By \autoref{prop:frobdifferentials} and \autoref{prop:osserman_construction}, $$W\simeq Z_p\oplus F^*_{\text{abs}}\Omega^1_{S}|_{Z_p} \simeq Z_p\oplus F^*_{\text{abs}}\Omega^1_{Z_p};$$ by \autoref{prop:frobdifferentials} the ring of functions on $W$ are free as an $\mathscr{O}_{Z_p}$-module with basis $\{1, (t_1\otimes 1-1\otimes t_1)^{[p]}, \cdots, (t_d\otimes 1-1\otimes t_d)^{[p]}\}.$ As $Z_p$ is contained in $V((t_1, \cdots, t_d)^m)$, the description of $W_+$ from the previous paragraph yields that $W\subset W_+$ as desired.
	
	Let $\pi_1, \pi_2: W_+\to Z_+$ denote the two projections. It suffices by \autoref{prop:crystal-p-power-leaves} to show that the map $\gamma$ of \eqref{eqn:crystalstructure} in the definition of the crystal structure on $\mathscr{M}_{dR}^\natural(X/S)$ sends $\pi_1^*(\mathscr{E},\nabla)$ to $\pi_2^*(\mathscr{E},\nabla)$. By assumption we have $$\gamma(\pi_1^*(\mathscr{E},\nabla)_{\mathbb{Q}})=\pi_2^*(\mathscr{E},\nabla)_{\mb{Q}};$$ we would like to conclude by \autoref{cor:unique-integral-isomonodromic}.
		 
	 Now consider the filtration of $\mathscr{O}_{W_+}$ by the ideals ${I}_n$ induced by the filtration on $$(R_\mathbb{Q}[[t_1, \cdots t_d]]\otimes_{R_\mathbb{Q}} R_\mathbb{Q}[[t_1, \cdots t_d]])$$ by multi-degree, i.e.~$I_n$ consists of the $R_{\mathbb{Q}}$-span of all monomials in $1\otimes t_i, t_i\otimes 1$ of degree at least $n$. This satisfies the hypotheses of \autoref{cor:unique-integral-isomonodromic} (for example, by the explicit description of $\mathscr{O}_{W_+}$ above) and hence, letting $\pi_i$ be the two projections $W_+\to Z_+$, we have $\gamma(\pi_1^*(\mathscr{E}, \nabla))\simeq \pi_2^*(\mathscr{E}, \nabla)$. Restricting this isomorphism to $W$ gives the claim by \autoref{prop:crystal-p-power-leaves} as desired.
\end{proof}

%\begin{corollary}\label{cor:vanishing-p-curvature-for-p-power-leaves-isomonodromy}
%	Let $p$ be a prime, and $\mf p$ a maximal ideal of $R$ of residue characteristic $p$. Suppose the hypotheses of \autoref{lemma:checking-n-flatness} are satisfied, and furthermore that each $A_i$ is divisible by $p$, so $Z_{\mf p}$ is a $p$-power subscheme of $S$ in the sense of \autoref{definition:p-power-leaves}. Then the $p$-curvature of the isomonodromy foliation vanishes along the $Z_{\mf p}$-point of $\mathscr{M}_{dR}(X/S)^\natural$ corresponding to $(\mathscr{E}, \nabla)|_{Z_{\mf p}}$.
%\end{corollary}
%\begin{proof}
%This is immediate by combining \autoref{lemma:checking-n-flatness} and \autoref{prop:crystal-p-power-leaves}.
%\end{proof}

%\section{Basic properties of isomonodromic deformation}

\bibliographystyle{alpha}
\bibliography{bibliography-p-painleve}

\end{document}